\documentclass[a4paper,12pt,twoside]{article}
\usepackage[latin1]{inputenc}                 
\usepackage{amsfonts,amssymb,amsmath,exscale} 
\usepackage[dvips]{graphicx,color,psfrag}     
\usepackage{fancyhdr}                         
\usepackage{caption}                          
\usepackage{rotating}
\usepackage{defcml}                           
\Floatboxname{Box}
\usepackage{verbatim}                         

\usepackage{amsmath,amsthm,amsfonts,amssymb,amscd}
\usepackage{subfig}
\usepackage{graphicx,bm,color}
\usepackage{algorithm}
        \usepackage{setspace}
\usepackage{algpseudocode}
\usepackage{stmaryrd}
\usepackage[framed,numbered,autolinebreaks,useliterate]{mcode}
\usepackage{blindtext}
\usepackage{multicol}

\usepackage{xargs}                      
\usepackage[color=gray!20, backgroundcolor=yellow, textwidth=1.85cm]{todonotes}
\newcommandx{\ju}[2][1=]{\todo[linecolor=blue,backgroundcolor=blue!25,bordercolor=blue,#1]{#2}}
\usepackage[numbers, sort, comma, square]{natbib}
\usepackage{mathtools}

\usepackage{booktabs,ctable,multirow,longtable} 
\usepackage[latin1]{inputenc}                 
\usepackage{amsfonts,amssymb,amsmath,exscale} 
\usepackage[dvips]{graphicx,color,psfrag}     
\usepackage{fancyhdr}                         
\usepackage{caption}                          
\usepackage{rotating}
\usepackage{defcml}                            
\usepackage{url}
\usepackage{lscape}
\Floatboxname{Box}
\newcolumntype{\xx}{\bm{x}}

\newcommand{\xx}{\boldsymbol{x}}
\newcommand{\R}{\mathbb{R}}

\definecolor{gray}{gray}{0.6}

\newtheorem{Remark}{Remark}[section]

\newtheorem{form}{Formulation}[section]

\def\R{\mathbb R}
\newcommand{\noii}[1]{{\textcolor{black}{#1}}}
\newcommand{\grm}[1]{{\textcolor{black}{#1}}}

\fancypagestyle{plain}{%
	\fancyhead{}%
	\fancyfoot[c]{\sffamily\thepage}%
}
\makeatletter
 
\usepackage{scalerel,stackengine}
\stackMath
\newcommand\reallywidecheck[1]{%
	\savestack{\tmpbox}{\stretchto{%
			\scaleto{%
				\scalerel*[\widthof{\ensuremath{#1}}]{\kern-.6pt\bigwedge\kern-.6pt}%
				{\rule[-\textheight/2]{1ex}{\textheight}}
			}{\textheight}%
		}{0.5ex}}%
	\stackon[1pt]{#1}{\scalebox{-1}{\tmpbox}}%
}

\begin{document}
	\Titel{
		Global-local techniques for hydraulic fracture
	}
	\Autor{N. Noii, H. A. Jahangiry, H. Waisman}
	\Report{02--I--17}
	\Journal{
		
	}
	%
	
	
	
\thispagestyle{empty}

\vspace*{1cm}
\ce{\bf \Large Level-Set Topology Optimization for Ductile and Brittle}
\vspace*{0.6cm}
\ce{\bf\Large Fracture Resistance Using the Phase-Field Method}

\vskip .35in

\ce{Nima Noii\(^{a,}\)\footnote{Corresponding author.\\[1.3mm]
		E-mail addresses:
		noii@ikm.uni-hannover.de (N. Noii); hassan.jahangiry@gmail.com (H. A. Jahangiry); waisman@civil.columbia.edu (H. Waisman).
	}, Hassan Ali Jahangiry \(^{b}\), Haim Waisman \(^{c}\)} \vskip .2in

\ce{\(^a\) Institute of Continuum Mechanics} \ce{Leibniz Universit\"at Hannover, Appelstrasse 11, 30167 Hannover, Germany}\vskip .25in

\ce{\(^b\) Faculty of Civil Engineering, Semnan University, Semnan, Iran} \vskip .2in

\ce{\(^c\) Computational Mechanics Group} \ce{Department of Civil Engineering \& Engineering Mechanics,}\vskip -.08in{\hspace*{4cm}Columbia University, New York, USA} 

\vspace*{0.5cm}
\ce{\bf \large \textcolor{blue}{\underline{Accepted 15.02.2023}}}
\vspace*{0.5cm}
\ce{\bf \large \textcolor{blue}{\underline{Computer Methods in Applied Mechanics and Engineering}}}
\vspace*{0.5cm}
\ce{\bf\large \textcolor{blue}{ \underline{ISSN 0045-7825}}}

\vspace*{0.45cm}
\vskip .2in

\begin{Abstract}
	This work presents a rigorous mathematical formulation for \textit{topology optimization of a macro structure} undergoing ductile failure. The prediction of ductile solid materials which exhibit dominant plastic deformation is an intriguingly challenging task and plays an extremely important role in various engineering applications. Here,  we rely on the phase-field approach to fracture which is a widely adopted framework for modeling and computing the fracture failure phenomena in solids. The first objective is to optimize the topology of the structure in order to minimize its mass, while accounting for structural damage. To do so, the topological phase transition function (between solid and void phases) is introduced, thus resulting in an extension of all the governing equations. Our second objective is to additionally enhance the fracture resistance of the structure.   Accordingly, two different formulations are proposed. One requires only the residual force vector of the deformation field as a constraint, while in the second formulation, the residual force vector of the deformation and phase-field fracture simultaneously have been imposed. An incremental minimization principles for a class of gradient-type dissipative materials are used to derive the governing equations. Thereafter, to obtain optimal topology to enhance the structural resistance due to fracture, the level-set-based formulation is formulated. The level-set-based topology optimization is employed to seek an optimal layout with smooth and clear boundaries. Sensitivities are derived using the analytical gradient-based adjoint method to update the level-set surface for both formulations. Here, the evolution of the level-set surface is realized by the reaction-diffusion equation to maximize the strain energy of the structure while a certain volume of design domain is prescribed. Several three-dimensional numerical examples are presented to substantiate our algorithmic developments.

\end{Abstract} 
\newpage
\vspace{-0.5cm}
{\small\tableofcontents}

	
\sectpa[Section1]{Introduction} 

 Topology optimization has been of great interest in both academia \cite{sigmund2013topology,deaton2014survey}, and industry \cite{zhu2016topology,aage2017giga}, since the seminal pioneering research work \grm{in} \cite{bendsoe1988generating}, and the recent significant developments in additive manufacturing. It has been used as a powerful tool in various engineering fields to achieve smart, robust, and lightweight structures. In fact,  topology optimization is a mathematical method that seeks to find the optimal material distribution that satisfies the equilibrium, objective, and constraint functions. By definition, topology optimization lies in the conceptual design phase, which can effectively determine the number, connection pattern and presence of holes in the design domain and evolve design members to improve the expected performance.

In the scientific literature, there are several different numerical methods to perform topology optimization \cite{sigmond2,sigmund2013topology,noii2017new}. In the earliest research study \cite{bendsoe1}, the maximum stiffness of structures was obtained by the microstructure or homogenization method, so to determine the optimum layout of the structure. Since then, topology optimization has been widely applied to a large variety of scientific fields. 

The most popular topology optimization methods are the Solid Isotropic Material with Penalization (SIMP) technique \cite{bendsoe2}, the Evolutionary Structural Optimization (ESO) technique \cite{Xie} and the Level Set Method (LSM) \cite{Osher1, Allaire1, Wang}. The major \grm{advantages for}  using the level-set method \grm{is} due to \grm{the direct description of} the optimal geometry, clear and smooth interface representation, interface merging or splitting, suppressed checkerboard pattern and \grm{circumvention of } the \textit{islanding} phenomenon \cite{Yamada}. Also, since the obtained layouts can be used directly on 3D printers without any post-processing techniques, this speeds up the industrial production process. The level-set method proposed in \cite{Osher1}, has been a versatile method for the implicit representation of evolutionary interfaces in an Eulerian coordinate system. Many efforts have been made to develop and improve level-set-based topology optimization,  see for examples \cite{Allaire2, Allaire1, Burger, Osher2, Sethian, Wang, JAHANGIRY1}, for more details. The key idea employing the level-set method topology optimization is to represent the interfaces by a discretized implicit dynamic hypersurface, which then evolves under a velocity field toward optimality. 

Nevertheless, the conventional level-set method is unable to create holes in the design domain. Furthermore, at every iteration, one must also, reshape the level-set function \grm{to satisfy the signed distance} characteristic which may not be convenient. \cite{Jahan_IGA_LSM_RDE1}. Accordingly, the Allen-Cahn equation \cite{allen1979microscopic} together with  phase-field approach \cite{bourdin2006phase} have been applied to optimal design problems. An advantage to the phase-field type models is that  the nucleation of voids in the design domain can be readily achieved. This type of regularized formulation is based on a Lyapunov energy functional, it is formulated based on a smooth double-well potential which takes a global minimum value at every phase  \cite{takezawa2010shape,bourdin2006phase}. In contrast, the major limitation of using the phase-field approach is two-fold. Since, it is a surface tracking method, it does not allow the number of voids in structure to be increased \cite{takezawa2010shape}. Additionally, double-well potential term results in more non-linearity in the formulation, and thus requires a sufficient number of iterations until the convergence of model \grm{is} reached \cite{choi2011topology,takezawa2010shape}.
Thereafter, the Reaction-diffusion (R-D) equation \cite{Yamada} which is inspired by the phase-field concept, was introduced to update the surface of the level-set method, i.e., the design field. Since then, many efforts have been devoted to use R-D instead of the conventional Hamilton-Jacobi (H-J) equation \cite{choi2011topology, otomori2015matlab, Jahan_IGA_LSM_RDE1}. As an advantage, the R-D equation compared to the phase-field approach no longer \grm{needs} the double-well potential, so it supports flexible changes in structural design \cite{choi2011topology}. In order to update the design variable,  the pointwise gradient-based sensitivity analysis is developed to derive the normal velocity of the R-D equation.

In the literature, there are still very few research works that take material non-linearity, e.g. elastic-plastic constitutive response, and fracture into account. This is mainly due to the challenges in deriving path-dependent sensitivities \grm{when} considering plasticity, and fracture \noii{response}, \grm{as well as} the expensive computational simulation it requires. \noii{Alberdi et al. \cite{alberdi2018unified} developed an adjoint sensitivity analysis framework to evaluate  the path-dependent design sensitivities for problems related to inelastic materials 
and dynamic responses affecting topology optimization, see also \cite{michaleris1994tangent}.}

Evidently, fracture arises in the form of evolving crack surfaces and ductile solid materials one can also expect dominant plastic deformation  \cite{noii2021bayesian}. In these types of materials, the crack evolves at a slow rate and is accompanied by significant plastic distortion, thus resulting in high material non-linearity. The prediction of such failure mechanisms due to crack initiation and growth coupled with elastic-plastic deformations is an intriguingly challenging task and plays an extremely important role in various engineering applications. 

Recently, the variational approach to fracture by Francfort and Marigo and the related regularized formulation \cite{bourdin2008variational}, which is also commonly referred to as a phase-field model of fracture, see e.g. the review paper \cite{ambati2015review, WU20201}, is a widely accepted framework for modeling and computing the fracture failure phenomena in elastic solids. 
Variational approaches are introduced based on energy minimization principles \cite{francfort+marigo98, miehe+welschinger+hofacker10a, kuhn2015}, and their regularization is obtained by
$\Gamma$-convergence, which is fundamentally inspired by the work of image segmentation conducted by Mumford and Shah \cite{mumford+shah89}. Then the model is improved by formulating a Ginzburg-Landau-type evolution equation of the fracture phase-field \cite{hakim+karma09}. Such models incorporate non-local effects based on length scales, which reflect properties of the material \textit{micro-structure size} with respect to the \textit{macro-structure size}.
A variety of research studies have recently extended the phase-field approach to fracture toward the ductile case. The essential idea is to couple the evolution of the crack phase-field to an elasto-plasticity model. Initial works on this topic include~\citep{alessi2014,ambati2015,alessi2017,borden2016,kuhn2016,ulloa2016} (see~\cite{noii2021bayesian,alessi2018comparison}  for an overview). Phase-field models for ductile fracture were subsequently developed in the context of cohesive-frictional materials~\citep{kienle2019}
including thermal effects~\citep{dittmann2020}, fiber pullout behavior \cite{storm2021comparative}, hydraulic fracture~\citep{heider2020phase,ulloa2022variational,noii2019phase}, \grm{stochastic} analysis \cite{noii2022probabilistic, noii2021bayesian,noii2022bayesian1}, degradation of the fracture toughness \cite{yin2020ductile}, multi-surface plasticity~\citep{fang2019}, and multi-scale approach \cite{liu2022phase,noii2020adaptive,gerasimov2018non,aldakheel2021multilevel} among others. 

Despite many investigations on topology optimization, failure behaviors due to the fracture, and more precisely ductile failure have seldom been considered. The present contribution is aimed at investigating \textit{topology optimization of structures enhancing the fracture resistance while modeling fracture using the phase-field framework}. Thus, we concentrate our efforts on structures consisting primarily of metals or other materials with relatively stiff elastic regions such that small deformation plasticity. A successful extension of the topology optimization approach to this setting would pave the way for the wide adoption of considering the fracture response in optimizing a design domain of industrial additive manufacturing technologies. 

In the literature, one of the earliest investigations on topology optimization considering fracture-induced softening behavior is by Challis et al. \cite{Chalis}, who used the virtual crack extension model to maximize fracture resistance and the level-set method to optimize the layout. Kang et al. \cite{Kang} utilized the J-integral to predict crack nucleation and propagation at predefined locations in the design domain. Failure mitigation in optimal topology design using a coupled nonlinear continuum damage model under fixed and variable loading are subjected to many studies \cite{james2015topology, james2014failure}. 
\noii{Interested readers for topology optimization due to fracture resistance of structures are further referred to
\cite{amir2013reinforcement,li2017topology,noel2017level,li2017design,wu2022topology}.}
Recently, the phase-field fracture method has been used to increase the fracture resistance of quasi-brittle composite materials using BESO for work maximization \cite{Xia_Fracture, Da_Fracture, russ2020novel, russ2021novel}. Russ and Waisman \cite{Russ} employed the phase-field fracture method to increase the brittle fracture resistance using SIMP to minimize the weight in the presence of compliance and fracture surface energy functional constraints. Wu et al. \cite{Wu} used reaction-diffusion  level-set method-based topology optimization to maximize the fracture resistance of brittle composite materials subject to the volume constraint using the phase-field fracture model. Hu et al. \cite{Hu} employed the Extended Finite Element Method (XFEM) to increase the fracture resistance using the BESO for topology optimization, in which the J-integral criterion and the mean compliance are considered as a multi-objective function while a given volume of Design area is prescribed.  

In the present study, following \cite{noii2021bayesian}, we present a variational formulation for brittle and ductile phase-field fracture based on variational principles, rooted in incremental energy minimization, for gradient-extended dissipative solids \cite{miehe+hofacker+schaenzel+aldakheel15}. The coupling of plasticity to the crack phase-field is achieved by a constitutive work density function, which is characterized by a degraded stored elastic energy and the accumulated dissipated energy due to plasticity and damage. Accordingly, first, we focus on the development of governing equations based on the topological field. Thus, we define the transition rule between the solid region (material counterpart), and the void region (non-material counterpart) for constitutive equations,  using the quadratic function which depends on the Heaviside function of the topological field. \noii{Here, an exact Heaviside step function is used, which helps to eliminate the negative effects arising from intermediate densities and checkerboard patterns in the final results and only cut elements  will be created by the interfaces (i.e., zero level-set).}

Desai et al. \cite{desai_Allaire} employed the phase-field model of fracture for brittle materials, in which the level-set method was used to maximize the elastic energy functional while a certain volume of design domain was prescribed. \noii{In this paper allocating/initializing the domain with multiple voids has been done for topology optimization due to fracture resistance. Because, initializing the domain with multiple voids, leads to multiple singular areas (due to new topology layout with additional weak elements as voids), thus multiple crack initiation appears in the inappropriate region. Additionally, fracture response and topology analysis have inverse effects on optimization problems. Since, fracture response makes structure \textit{weaker}, while topology analysis makes structure \textit{stiffer} with imposed constraints. As such, the idea of the reaction-diffusion-based level-set method is to enable void nucleation within the design domain. Thus, this helps to eliminate the   destructive effects of the fracture caused by the initialization of the Hamilton-Jacobi-based level-set method.}

The second objective is to introduce a mathematical formulation of the topological optimization problems dealing with brittle and ductile phase-field fracture. \noii{This is realized by means of the reaction-diffusion-based level-set method.} We mainly concentrate on the three-dimensional setting. The level-set based optimization problem is defined by maximizing the total mechanical work under a certain prescribed structural volume fraction (as a prescribed quantity). Since, we are dealing with a rate-dependent nonlinear boundary value problem (i.e., ductile phase-field fracture), an incremental topology optimization approach is used to reach the target volume. Thereafter, two different types of formulation are proposed. One equipped only with the residual force vector of the displacement field as a constraint. Additionally in the next formulation, it is required that \textit{simultaneously} the residual form of the displacement and the phase-field fracture be imposed. The first type enables computations performed with less implementation effort and indeed it mimics the \textit{classical} implementation of topology optimization for elastic-plastic problems \cite{fritzen2016topology, huang2008topology,gangwar2022thermodynamically}. Since the second formulation requires both equilibrium equations to be imposed in the optimization problem, one may obtain more accuracy in optimizing a design domain. To the best of our knowledge, none of the aforementioned works have included these formulations for ductile fracture. Even in the case of topological optimization for brittle fracture, only the first type of formulation is applied, see \cite{desai_Allaire, Xia_Fracture, Wu}. The accuracy of these formulations is further tackled in this contribution, and the superiority of the second formulation is highlighted.

In literature, most of the nonlinear problems are resolved through the gradient-based topology optimization e.g. \cite{schwarz2001topology, yoon2007topology, huang2007topology, huang2008topology, xia2017evolutionary, maury2018elasto, zhao2019material, zhao2020topology, Jahan_IGA_LSM_RDE_EP, Jahan_IGA_LSM_RDE_GN} among many others. To this end, herein we employ a displacement-controlled adjoint sensitivity that is used for gradient-based topology optimization. The key requirement for realizing this coupled adjoint-based sensitivity analysis (due to coupled plasticity and phase-field fracture), is derivatives of the objective and constraint functions with respect to the topological field (so-called design variable) which are consistently formulated to evolve the structural topology.

Additionally, we examine the efficiency of the optimum layout due to brittle/ ductile fracture as compared to the optimum layout obtained by assuming linear elasticity, and also with the original domain. Hence, we present quantitative and qualitative indicators to highlight the role of considering fracture equilibrium in the topology optimization framework. These indicators correspond to: $(i)$ The load-displacement curve to examine the maximum load capacity before crack initiation. $(ii)$ Qualitative response of crack phase-field pattern to highlight the effects of damage response within every new topology. So, we expect less crack area exists in the final topological configuration. $(iii)$ Lastly, to examine the convergence of the optimization process by the objective function and volume constraint, which implies the stiff response of material through the total mechanical energy of structure under a specific volume. 

The main objective of this paper is summarized as follows:

\begin{itemize}
\item {A novel} theoretical topology optimization formulation to be derived for the phase-field brittle/ductile fracture;

\item Level-set-based method for fracture resistance topology optimization of 3D structures;

\item The reaction-diffusion equation to capture the evolution of the level-set function;

\item An analytical adjoint method (through the shape derivatives) to derive the point-wise sensitivities of the design variables;

\item An exact (binary) Heaviside step function to accurately resolve the material phases (solid/void).
\end{itemize}

The paper is structured as follows.  In Section 2, we provide a brief overview of topology optimization coupled with phase-field fracture method within elastic-plastic materials. We further outline the theoretical framework for the ductile phase-field fracture models in a variational setting, making use of incremental energy minimization. Next in Section 3, we first introduce the reaction-diffusion equation based level-set method, and then adjoint-based sensitivity analysis is presented. In Section 4, four numerical simulations (associated with brittle and ductile fractures) are performed to demonstrate the correctness of our algorithmic developments. Finally, the conclusion with some remarks and suggestions for future research are provided.


\sectpa[Section2]{Phase-Field Modeling of Ductile Fracture in Elastic-Plastic Materials}
This section outlines a mathematical framework for topology optimization of fracture-resistance of  structure undergoing ductile failure. We first elaborate on the primary fields and function spaces. Here, the fracture process is modeled by employing the well-developed phase-field formulation to resolve the sharp crack surface topology in the regularized concept. To this end, by introducing a topological field, we modify the governing equations for phase-field thus suitable for topology optimization framework. Thus, governing equations associated with topology optimization of ductile phase-field fracture \grm{are} derived.

\begin{figure}[!t]
	\centering
	\caption*{\hspace*{0cm}\underline{\texttt{time}=0}\hspace*{5.5cm}\underline{\texttt{time}=$t_m$}\hspace*{0.5cm}}
	{\includegraphics[clip,trim=3cm 29cm 1cm 0cm, width=17cm]{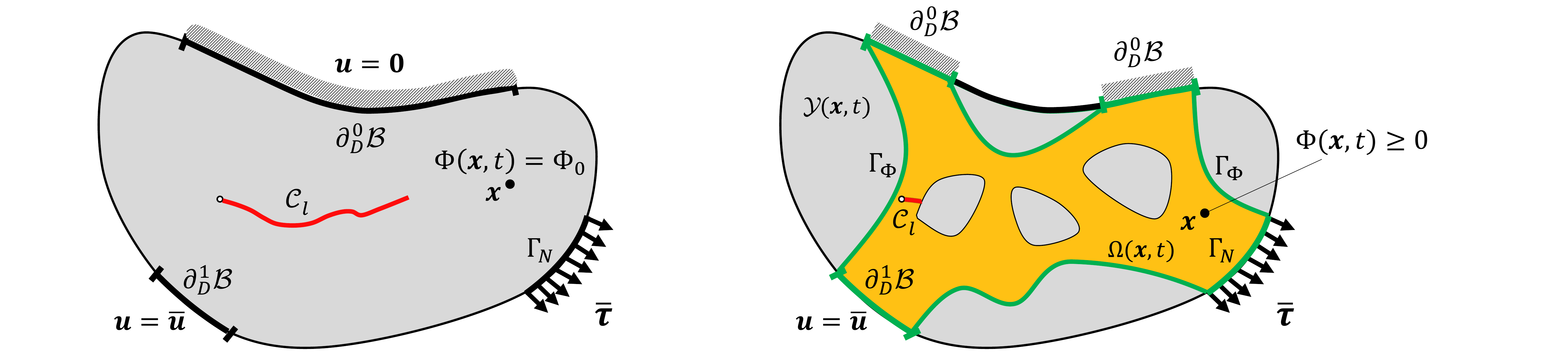}}  
	\vspace*{-0.35cm}
	\caption*{\footnotesize{\hspace*{1cm}$\calB=\Omega(\Bx,0)$ and $\calY(\Bx,0)=\varnothing$\hspace*{2.7cm}$\calB=\Omega(\Bx,t)\cup\calY(\Bx,t)\cup\Gamma_\Phi(\Bx,t)$\hspace*{1.5cm}}}
	\caption{Setup for design domain $\calB$ which consists of an admissible solid design domain $\Omega(\Bx,t)$ (material domain), and void domain $\calY(\Bx,t)$ (non-material domain).}
	\label{F1}
\end{figure}

\sectpb[Section2]{Primary fields and function spaces}
Let $\calB\subset{\mathbb{R}}^{\delta}$ be an arbitrary solid domain, $\delta=\{2,3\}$ with a smooth  boundary $\partial\calB$ as depicted in Figure \ref{F1}. Here, we denote $\calB$ as a design domain. We assume Dirichlet boundary conditions on $\partial_D\calB $ and Neumann boundary conditions on $\partial_N \calB := \Gamma_N \cup \mathcal{C}$, where $\Gamma_N$  denotes the outer domain boundary and {$\calC\in \mathbb{R}^{\delta-1}$} is the crack boundary, as illustrated in Figure \ref{F1}. Next, we introduce following primary fields to state our variational formulation.

\sectpd[Section27]{Displacement field and crack phase-field}

The response of the fracturing solid at material points $\Bx\in\calB$ and time $t\in \calT = [0,T]$ is described by the displacement field $\Bu(\Bx,t)$ and the crack phase-field $d(\Bx,t)$ as
\begin{equation}
	\Bu: 
	\left\{
	\begin{array}{ll}
		\calB \times \calT \rightarrow \mathbb{R}^\delta \\[2mm]
		(\Bx, t)  \mapsto \Bu(\Bx,t)
	\end{array}
	\right.
	\AND
	d: 
	\left\{
	\begin{array}{ll}
		\calB \times \calT \rightarrow [0,1] \\[2mm]
		(\Bx, t)  \mapsto d(\Bx,t)
	\end{array}
	\right.
	\WITH
	\dot{d} \ge 0 .
	\label{s2-fields}
\end{equation}
Intact and fully fractured states of the material are characterized by $d(\Bx,t)=0$ and $d(\Bx,t)=1$, respectively. In order to derive the variational formulation, the following space is first defined. For an arbitrary $A\subset\mathbb{R}^\delta$, we set
\begin{align}
	\mathrm{H}^1(\calB,A):=\{v:\calB\times\calT\rightarrow A\quad:\quad v\in \mathrm{H}^1(A)\}.
\end{align}
We also denote the Sobolev vector valued space $\mathbf{H}^1(\calB,A):=\left[\mathrm{H}^1(\calB,A)\right]^\delta$ and define
\begin{equation}
	\calW_{\overline{\Bu}}^{\Bu}:=\{\Bu\in\mathbf{H}^1(\calB,\mathbb{R}^\delta)\;\colon\;
	\Bu=\bm{0} \ \text{on} \ \partial^0_D\calB \AND \Bu=\overline{\Bu} \ \text{on} \ \partial^1_D\calB\}.
	\label{eq:spaces_u}
\end{equation}
Concerning the crack phase-field, we set
\begin{equation}
	\calW^{d}:=\mathrm{H}^1(\calB) \AND  \calW^{d}_{d^n}:=\{d\in \mathrm{H}^1(\calB,{\color{black}[0,1]}) \; \colon\; \ d \geq d^n \},
	\label{eq:spaces_d}
\end{equation}
where $d^n$ is the damage value in a previous time instant. Note that $\calW^{d}_{d^n}$ is a non-empty, closed, and convex subset of $\calW^{d}$, and introduces the evolutionary character of the phase-field, incorporating an irreversibility condition in incremental form. 

\sectpd[Section27]{Topological field}

\begin{figure}[!t]
	\centering
	{\includegraphics[clip,trim=0cm 29cm 0cm 0cm, width=16cm]{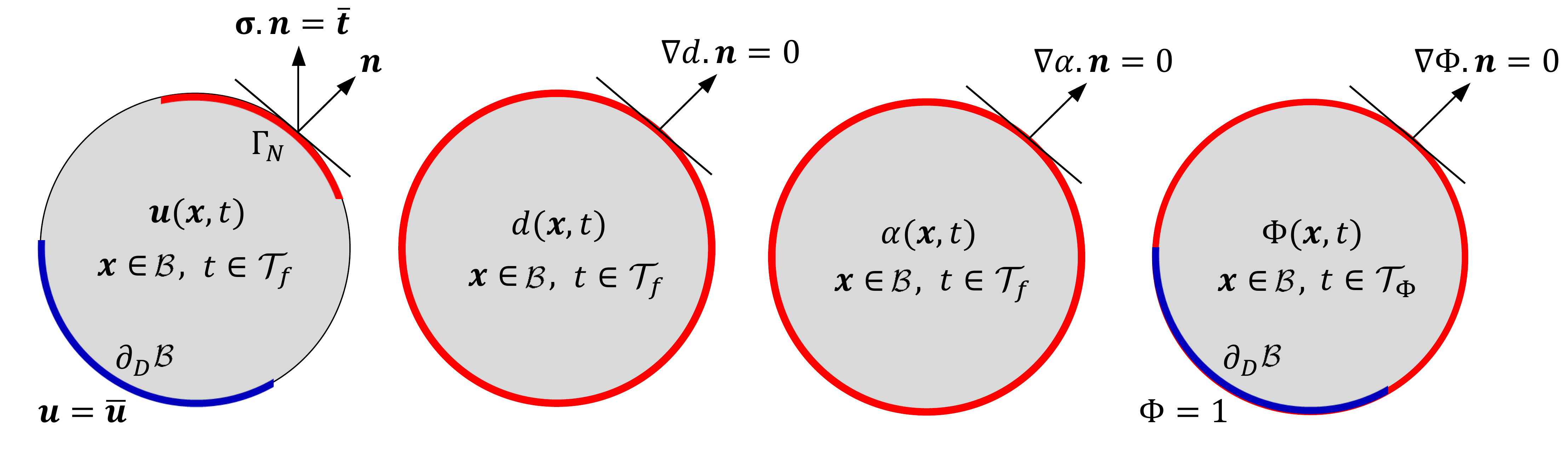}}  
	\caption{Primary variable fields in topology optimization ductile phase-field fracture, for a solid
		body $\calB \subset {\mathbb{R}}^{\delta}$ with dimension $\delta \in [2,3]$. $(a)$
		The displacement field $\Bu$ defined on $\calB$ and Neumann-type
		boundary condition for traction $\bar{\Bt}=\Bsigma \cdot
		\Bn$. $(b)$ The crack phase-field is determined by
		Dirichlet-type boundary condition $d=1$ on $\mathcal{C}$ and Neumann-type
		boundary condition $\nabla d \cdot \Bn = 0$ on the full surface $\partial{\calB}$. $(c)$ 
		The hardening variable $\alpha$ in $\calB$ that is
		continuous in entire domain with Neumann-type
		boundary condition $\nabla \alpha \cdot \Bn = 0$. $(d)$ The topological field $\Phi$ is determined by
		Dirichlet-type boundary condition $\Phi=1$ on $ \partial^1_D\calB$ (non-zero prescribed load surface), and Neumann-type
		boundary condition $\nabla \Phi \cdot \Bn = 0$ on the full surface $\partial{\calB}$.}
	\label{F2}
\end{figure}

To elaborate on the topological field (also known as design variable), let the design domain $\calB$ (reference domain) be divided into  \textit{admissible solid design domain} (also refers to as a material domain) denoted as $\Omega\subset\calB \in \mathbb{R}^\delta$, and the remaining region as \textit{void domain} (also referred as a non-material domain) as $\calY:=\calB \backslash \Omega\subset\calB \in \mathbb{R}^\delta$, thus $\calB=\Omega\cup\calY$, see Figure \ref{F2}. Next, we can assume admissible solid design domain $\Omega$ to be bounded by the given two volume fraction $0\le\Theta_1\le\Theta_2\le1$, such that
\begin{equation}\label{addmis_domain1}
	\Theta_1 \text{V}(\calB)\le \text{V}(\Omega)\le\Theta_2 \text{V}(\calB)\WITH 
	\text{V}(\bullet)=\int_{\bullet} \text{d}\Bx
\end{equation}
In topological optimization, one aims to find an optimal solid design domain $\Omega$ which will minimize an objective function for a given (theoretical/ practical) set of constraints. Any set $\Omega$ with finite perimeter could be represented by a continuous indicator field, so-called \textit{topological field}. Thus, the minimization of an objective function  will be reduced from the determination of an optimal solid design domain $\Omega$ to finding a topological field denoted as $\Phi(\Bx,t)$. Thus, we define a topological field $\Phi(\Bx,t): \calB \times \calT\rightarrow [-1,1]$, which following \req{addmis_domain1} is reduced to:
\begin{equation}\label{addmis_domain2}
	\Theta_1 \text{V}(\calB)\le \widehat{\text{V}}(\Omega)\le\Theta_2 \text{V}(\calB)\WITH 
	 \widehat{\text{V}}(\Omega)=\int_{\calB} \text{H}(\Phi(\Bx,t))\text{d}\Bx
\end{equation}
We define a Heaviside function for a design variable as our topological $\text{H}(\Phi(\Bx,t)): \calB \times \mathbb{T}\rightarrow [0,1]$. In which, a direct mapping between two phases, i.e., material and non-material phases, in $\calB$, is approximated through \textit{ an exact Heaviside function} for characteristic point $\Bx$ defined as a $\text{H}(\Phi(\Bx,t))=1$ if $\Bx\in\Omega$, and so $\text{H}(\Phi(\Bx,t))=0$ if $\Bx\in\calY$, respectively. We note that $\TH^n(\Phi)=\TH(\Phi)$ for all $n\ge1$. Remarkably, to determine the volume integral $\text{V}(\Omega)$, one can use  $\widehat{\text{V}}(\Omega)$ from \req{addmis_domain2} to get the volume integral over the entire domain of the function $\Phi(\Bx,t)$, hence, as opposed to the SIMP method, one need not explicitly determine an optimal solid domain $\Omega$, like in \req{addmis_domain1}, which results in $	\text{V}(\Omega)=\widehat{\text{V}}(\Omega)$.

It should be noted that the characteristic material point for $\Bx\in\Omega$ is represented by  $\Phi(\Bx,t)>0$, and a characteristic void (non-material point) for $\Bx\in\calY$ is represented by  $\Phi(\Bx,t)<0$. Additionally, $\Phi(\Bx,t)=0$ (i.e., zero level-set) represents a material interface (i.e, material boundary) so-called $\Gamma_{\Phi}$ that is a surface to distinguish between a set of material and non-material points in given domain $\calB$, see Figure \ref{F2}. Finally, the set of the topological field is defined as the following space:
\begin{equation}
	\calW_{\overline{\Phi}}^{\Phi}:=\{\Phi\in\mathrm{H}^1(\calB,[-1,1])\;\colon\;
	\Phi=1 \ \text{on} \ \partial^1_D\calB\}.
	\label{eq:spaces_phi}
\end{equation}
%
\sectpd[Section27]{Hardening field and plastic strain}
Focusing on the isochoric setting of von Mises plasticity theory, we define the plastic strain tensor $\Bve^p(\Bx,t)$ and the hardening variable $\alpha(\Bx,t)$ as
\begin{equation}
	\Bve^p: 
	\left\{
	\begin{array}{ll}
		\calB \times \calT \rightarrow \mathbb{R}^{\delta\times\delta}_\mathrm{dev} \\[2mm]
		(\Bx, t)  \mapsto \Bve^p(\Bx,t)
	\end{array}
	\right.
	\AND
	\alpha: 
	\left\{
	\begin{array}{ll}
		\calB \times \calT \rightarrow \mathbb{R}_+ \\[2mm]
		(\Bx, t)  \mapsto \alpha(\Bx,t)
	\end{array}
	\right.
	\WITH
	\dot{\alpha} \ge 0 ,
	\label{s2-fields}
\end{equation}
where $\mathbb{R}^{\delta\times\delta}_\mathrm{dev}:=\{\Be\in\mathbb{R}^{\delta\times\delta} \ \colon \ \Be^T=\Be,\ \tr{[\Be]}=0\}$ is the set of symmetric second-order tensors with vanishing trace. The plastic strain tensor is considered as a local internal variable, while the hardening variable can be considered as a non-local internal variable. In particular, $\alpha$ may be introduced to incorporate phenomenological hardening responses and/or non-local effects, for which the evolution equation reads
\begin{equation}
	\dot\alpha = \sqrt{\frac{2}{3}}\,\vert\dot{\Bve}^p\vert,
	\label{s2-evol-alpha}
\end{equation}
is considered. As such, $\alpha$ can be viewed as the equivalent plastic strain, which starts to evolve from the initial condition $\alpha(\Bx,0) = \text{0}$, for which we assume $\Bve^{p}\in\mathbf{Q}:=\mathrm{L}^2(\calB;\mathbb{R}^{\delta\times\delta}_\mathrm{dev})$. Moreover, in view of~\eqref{s2-evol-alpha}, it follows that $\alpha$ is irreversible. Assuming in this section the setting of gradient-extended plasticity, we define the function spaces
%
\begin{equation}
	\calW^\alpha_{\alpha_n,\,\Bq}:=\{\alpha\in\calW^\alpha \quad \colon \quad \alpha = \alpha_n + \sqrt{2/3}\,\vert \Bq \vert, \ \Bq \in \mathbf{Q} \},
	\label{Walpha}
\end{equation}
{\color{black} where $\calW^\alpha=\mathrm{L}^2(\calB)$ for local plasticity, while $\calW^\alpha=\mathrm{H}^1(\calB)$ for gradient plasticity.} The hardening law~\eqref{s2-evol-alpha} is thus enforced in incremental form by setting $\alpha\in \calW^\alpha_{\alpha_n,\,\Bve^p-\Bve^p_n}$.

The gradient of the displacement field defines the symmetric strain tensor of the geometrically linear theory as
\begin{equation}
	\Bve = \nabla_s \Bu = \sym[ \nabla \Bu ] := \frac{1}{2} [\nabla\Bu + \nabla\Bu^T]
	.
	\label{s2-disp-grad}
\end{equation}
In view of the small strain hypothesis and the isochoric nature of the plastic strains, the strain tensor is additively decomposed into an elastic part $\Bve^e$ and a plastic part $\Bve^p$ as
\begin{equation}
	\Bve = \Bve^e + \Bve^p 
	\WITH
	\tr{[\Bve]} = \tr{[\Bve^e]}
	.
	\label{s2-strain-e-p}
\end{equation}

The solid $\calB$ is loaded by prescribed deformations and external traction on the boundary, defined by time-dependent Dirichlet conditions and Neumann conditions
\begin{equation}
	\Bu = \overline{\Bu} \ \textrm{on}\ \partial_D\calB
	\AND
	\Bsigma \cdot \Bn 
	= \overline{\Btau} \ \textrm{on}\ \partial_N\calB ,
	\label{s2-bcs}
\end{equation}
where $\Bn$ is the outward unit normal vector on the surface $\partial \calB$. The stress tensor $\Bsigma$ is the thermodynamic dual to $\Bve$ and $\bar{\Btau}$ is the prescribed traction vector.

\sectpb[Section27]{Geometry Projection}

An important part of the level-set-based topology optimization is the geometry mapping to a structural model that affects the predicted behavior of a converged design. The geometry mapping projects the geometry defined by the discretized level-set function onto the structural model that provides the structural response. The most common approach is the Eulerian approach, referred to as the density-distribution method, which leads to the evolution of materials in the fixed discretized domain. In the level-set-based topology optimization, the Heaviside function and Dirac-$\delta$ function are two important ingredients that are correspondingly used for volume and boundary integration, respectively. The exact  Heaviside function based on the topological field reads:
\begin{equation}\label{eq:heavydide_def}
	\begin{aligned}
		\operatorname{H}\big( \Phi (\Bx,t) \big)=\left\{ 	\begin{aligned}
			& 1 \qquad \text{if} \;\;\Bx\in\Omega\rightarrow\Phi \left( \Bx \right)\ge 0 \\ 
			& 0 \qquad \text{if} \;\;\Bx\in\calY\rightarrow\Phi \left( \Bx \right)<0, \ \\ 
		\end{aligned}
		\right.
	\end{aligned}
\end{equation} 
By definition \cite{osher2004level}, Dirac-$\delta$ function is the directional derivative of  $\text{H}(\Phi(\Bx,t))$ in the normal direction $\widehat{\Bn}_{\Phi}$ which is defined as
\begin{equation}\label{eq:dirac_def}
	\begin{aligned}
		\delta(\Bx)=\nabla(\Phi(\Bx,t))\cdot\widehat{\Bn}_{\Phi} \WITH
		\widehat{\Bn}_{\Phi} = \frac{\nabla \Phi }{\noii{|\nabla \Phi|} }
	\end{aligned}
\end{equation}
Thus, following \req{eq:dirac_def} a material interface is related to the Dirac-$\delta$ function by 
\begin{equation}
	\Gamma_{\Phi}(\Bx)=\{ \Bx\;\colon\; \delta (\Phi (\Bx))> 0 \}\;,
\end{equation}
Subsequently, by having a definition of the exact Heaviside function and Dirac-$\delta$ function, we are now able to define the integration rule for some arbitrary function $F(\Bx)$ in $\Omega$ and $\calY$ as:
\begin{equation}
	\begin{aligned}
		\int_{\Omega}{F}(\Bx)\mathrm{d}{\bm{x}}\text{=}\int_{\calB}{F}(\Bx)\text{H}(\Phi) \mathrm{d}{\bm{x}}, \AND
		 \int_{\calY}{F}(\Bx)\mathrm{d}{\bm{x}}\text{=}\int_{\calB}{F}(\Bx)\big(1-\text{H}(\Phi)\big) \mathrm{d}{\bm{x}}.
	\end{aligned}
\end{equation}
Additionally, the integral in the full domain $\calB$ is approximated toward its outer surface $\partial \calB$ through Dirac-$\delta$ function as:
\begin{equation}\label{eq:geo_transfer}
	\begin{aligned}
		\int_{\partial\calB}{F}(\Bx)\mathrm{d}{\bm{a}}\text{=}\int_{\calB}{F}(\Bx)\delta(\Phi)\left| \nabla\Phi\right|\mathrm{d}{\bm{x}}. \\
	\end{aligned}
\end{equation}

\sectpb[Section22]{Energy quantities for ductile phase-field fracture}
\label{sec:energy_functions}
%
Here, the ductile fracture model is formulated based on a geometrically conceived approach to the phase-field modeling, following the local plasticity theory described in~\cite{miehe+hofacker+schaenzel+aldakheel15,noii2021bayesian}. The original model is constructed within a variationally consistent framework, in agreement with the incremental energy minimization principle. In this case, to formulate the topology-based ductile fracture model, we define
\begin{equation}
	\mbox{Global Primary Fields}: \ 
	\BfrakU := \{\Bu, d,\Phi\}\;,
	\label{primary_m2}
\end{equation}
and the set of constitutive state variables as
\begin{equation}
	\mbox{Constitutive State Variables}: \
	\BfrakC := \{ \Bve, \Bve^p, \alpha, d, \nabla d,\nabla \Phi \}
	.
	\label{state_m2}
\end{equation}
(see Figure \ref{F2}). In order to derive the  variational formulation, we define an energy density function per unit volume $ W(\BfrakC)$. Thus, at a fixed point,  the  total energy function is additively decomposed into elastic contribution $W_{elas} (\cdot\,)$,  a plastic contribution $W_{plas} (\cdot\,)$, and a fracture contribution~$W_{frac} (\cdot\,)$, resulting
\begin{equation}
	{W(\BfrakC):= W_{elas}(\Bve,\Bve^p ,d ,\alpha,\Phi) + W_{plas}(\alpha, d, \nabla \alpha,\Phi) +
		W_{frac}(d, \nabla d,\Phi )}. 
	\label{psuedo-energy}
\end{equation}
Next, by having pseudo-energy density per unit volume functions given at hand, a rate-dependent pseudo potential energy functional can be written as
\begin{equation}
	\begin{aligned}
		{\calE}(\Bu,\Bve^p,\alpha,d,\Phi) &:= 
		\int_{\calB} {W}(\BfrakC) \,\text{d}\Bx \;, 
		\label{potential-energy-functional}
	\end{aligned}
\end{equation}
Accordingly, to derive the variational formulation one has to define the  constitutive energy density functions, namely $W_{elas}$, $W_{plas}$, and $W_{frac}$, which is discussed below. 

\sectpc{Elastic contribution}
The elastic energy density $W_{elas}$ in \eqref{psuedo-energy}
is expressed in terms of the effective strain energy density $\psi_e(\Bve^e)$. In our formulation to preclude fracture in compression, a decomposition of the effective strain energy density into \textit{damageable} and \textit{undamageable} parts \grm{is} employed. Thus, we perform
additive decomposition of the strain tensor into \textit{volume-changing}
(volumetric) and \textit{volume-preserving} (deviatoric) counterparts
\[
\bm\varepsilon^e(\Bu,\Bx)=\bm\varepsilon^{e,vol}(\Bu,\Bx)+\bm\varepsilon^{e,dev}(\Bu,\Bx),
\]
where
\begin{equation}
	\begin{aligned}
	&\bm\varepsilon^{e,vol}(\Bu,\Bx):=\mathbb{P}^{vol}:\bm \varepsilon^e,\; \WITH \mathbb{P}^{vol}_{ijkl}:=\frac{1}{3}\delta_{ij}\delta_{kl} ,\\
	&\bm\varepsilon^{e,dev}(\Bu,\Bx):=\mathbb{P}^{dev}:\bm \varepsilon^e,
	\WITH \mathbb{P}^{dev}:=\mathbb{I}-\frac{1}{3}\text{\BI}\otimes\text{\BI} \AND 
	\mathbb{I}_{ijkl}:=\frac{1}{2}\big(\delta_{ik}\delta_{jl}+\delta_{il}\delta_{jk}\big)\;,
	\end{aligned}
\end{equation}
\noii{along with the first and second invariants of the elastic strain denoted as $I_1(\Bve^e):=\text{tr}[\Bve^e]$  and $I_2(\Bve^e):=\tr[(\Bve^e)^2]$.} 
The effective  strain energy function $\psi_e(\Bve^e)$ is additively decomposed into  damageable and undamageable contributions:
\begin{equation}
	\psi_e(I_1(\Bve^e),I_2(\Bve^e)(\Bx))=\psi_e^{+}(I_1,I_2)+\psi_e^{-}(I_1),
\end{equation}
such that
\begin{equation}
	{\psi_e^{+}}={H{^+}[I_1]}\psi_e^{vol}\big(I_1\big)
	+\psi_e^{dev}\big(I_1,I_2\big)~\AND 
	{\psi_e^{-}}=\big(1-{H{^+}[I_1]}\big)\psi_e^{vol}\big(I_1\big).~	
	\label{eq232_5}
\end{equation}
{Therein, $H{^+}[I_1(\Bve^e)]$ is a \textit{positive Heaviside function} which returns one and zero for $I_1(\Bve^e)>0$ and $I_1(\Bve^e)\leq0$, respectively.} We note that in this paper the volumetric and deviatoric counterparts of the energy take the following forms:
\noii{\begin{equation}
	\psi_e^{vol}\big(I_1\big)
	=\frac{{K}}{2}I^2_1=\frac{K}{2}\Big(\Bve^{e,vol}:\bm I\Big)^2
     \AND
	\psi_e^{dev}\big(I_1,I_2\big)
	={\mu}\Big(\frac{I_1^2}{3}-I_2\Big)
	=\mu {\Bve}^{e,dev}:{\bm \varepsilon}^{e,dev},
	\label{eq:psi_iso11}
\end{equation}}
Here, $K= \lambda+\frac{2}{3}\mu>0$ is the bulk modulus which includes the elastic Lam\'e's first constant denoted by $\lambda$  and the shear modulus $\mu$. The total elastic contribution to the pseudo-energy \eqref{psuedo-energy} finally reads
\begin{equation}
	W_{elas}(\Bve,\Bve^p,d,\alpha):=
	f(\Phi)\Big[g(d)\psi_e^{+}(I_1,I_2)+ \psi_e^{-}(I_1)\Big],
	\label{elas-part0}
\end{equation}
where $g(d)$ is the \textit{elastic degradation function}, which takes in a simple quadratic form as: 
\begin{equation}
	\begin{aligned}
		g(d)=(1-\kappa)(1-d)^2+\kappa,
	\end{aligned}
\end{equation}
where $\kappa$ is so-called ersatz material parameter (set as a very small quantity) which is used  to avoid numerical instabilities, and mathematically it is also dependent on the discretization space, see
\cite{khodadadian2020bayesian}.
Additionally, to define the transition rule of point $\Bx\in\calB$ between the solid region in $\Omega$, and the void region $\calY$ for the constitutive equations,  we define the quadratic function $f(\Phi)$ so-called \textit{topological phase transition function} as:
\begin{equation}\label{fun:topological_phase}
	\begin{aligned}
		f(\Phi)=(1-\kappa)\text{H}(\Phi)^2+\kappa,
	\end{aligned}
\end{equation}
such that 
\begin{itemize}
	\item In the solid part of the domain ($\Phi\ge0$), and the intact region ($d\approx0$) yields ${W}_{elas}(\Bve, d)\approx{\psi}_{e}(\Bve)$,
    \item In the void part of the domain ($\Phi<0$) yields ${W}_{elas}(\Bve, d)\approx0$.
\end{itemize}
The first equation corresponds to the balance of linear momentum within the quasi-static response defined as
\begin{equation}
		\fterm{
	\div\,[\Bsigma(\Bve,d,\Phi)] + f(\Phi) \overline\Bf = \bm{0}\ ,
}
	\label{s2-equil:defo}
\end{equation}
where dynamic effects are neglected, and $\overline\Bf$ is the given body force. Following the Coleman-Noll procedure, the stress tensor is obtained from the potential ${W}_{elas}$ in \req{elas-part0} by
\begin{equation}
	\Bsigma := \partial_{\Bve^e} {W}_{elas}
	=  f(\Phi) \widehat{\Bsigma}(\Bve_e)
	=f(\Phi)\left[g(d) \widetilde\Bsigma_{+} + \widetilde\Bsigma_{-}\right],
	\label{vari-bf-23}
\end{equation}
therein,
\begin{equation}\label{eq232_8}
	{\widetilde{\Bsigma}}^{+}({\Bve_e})
	=K\text{H}(I_1)(\bm \varepsilon^e:\textbf{I})\textbf{I}
	+2\mu{\bm \varepsilon}^{e,dev},\quad \text{and} \quad {\widetilde{\Bsigma}^{-}}({\Bve_e})=K\big(1-\text{H}(I_1)\big)(\bm \varepsilon^e:\textbf{I})\textbf{I}.
\end{equation}
where $\widehat{\Bsigma}(\Bve_e)$ is the solid material stress tensor, and $\widetilde\Bsigma_{\pm}$ is the effective positive and negative stress tensor, respectively. The decoupled representation of the fourth-order elasticity tensor (to relate the work into conjugate pairs of stress and strain tensors) is obtained through the additive decomposition of the stress tensor, which reads as follows:
\begin{equation}
	\mathbb{C}:=\frac{\partial {\bm \sigma}( \Bve_e)}{\partial {\bm \varepsilon}}
	=f(\Phi)\left[g(d)\frac{\partial {\bm \sigma}^{+}( \Bve_e)}{\partial {\bm \varepsilon}}
	+\frac{\partial {\bm \sigma}^{-}(\Bve_e)}{\partial {\bm \varepsilon}}\right]
	=:f(\Phi)\left[g(d)\widetilde{\mathbb{C}}^++\widetilde{\mathbb{C}}^-\right],
	\label{eq:elast_tensor1}
\end{equation}
where $\widetilde{\mathbb{C}}^+$, \grm{and} $\widetilde{\mathbb{C}}^-$ are the fourth-order elasticity tensor corresponds to damageable and undamageable counterparts, see Section 2.4.1.

\sectpc{Plastic contribution}
The plastic contribution $W_{plas}$ is expressed in terms of an effective plastic energy density $\psi_p$, whose form will depend on the adopted phenomenological model. In line with previous works~\citep{noii2021bayesian,miehe2017phase}, we consider a potential function in the context of local von Mises plasticity:   
\begin{equation}
	{\psi}_{p}(\alpha,\nabla\alpha) :=\frac{1}{2}h \alpha^2(\Bx,t) ,
	\label{psi_p}
\end{equation}
with the isotropic hardening modulus $h\ge 0$.
Thus, the plastic contribution to the pseudo-energy density~\eqref{psuedo-energy} then reads
\begin{equation}
	\begin{aligned}
		W_{plas}(\alpha,d ,\nabla \alpha):=
		f(\Phi)g(d) {\psi}_{p}(\alpha),
	\end{aligned}
	\label{plas-part0}
\end{equation}
where $g_p(d)=g(d)$ is the \textit{ plastic degradation function}, along with topological phase transition function $f(\Phi)$. In case of different \grm{degradation} function for plastic, and elastic contributions, see \cite{noii2021bayesian}.

\sectpc{Fracture contribution} 
The phase-field contribution $W_{frac}$ is expressed in terms of the crack surface energy density $\gamma_l$ and the fracture length-scale parameter $l_f$ that governs the regularization. In particular, the sharp-crack surface topology $\calC$ is regularized by a functional $\calC_l$, as outlined in ~\cite{Wi20_book,heister2015}. The regularized functional reads
\begin{equation}
	\calC_l(d) = \int_{\calB} \gamma_l(d, \nabla d) \,\text{d}\Bx.
	\label{s2-gamma_l}
\end{equation}
with positiveness for crack dissipation as:
\begin{equation}
	\frac{d}{dt} \calC_l(d)  \ge 0\;.
	\label{gamma-evol}
\end{equation}
In line with standard phase-field models \cite{miehe+welschinger+hofacker10a,ambati+kruse+lorenzis16}, a general surface density function for the $\gamma_l(d, \nabla d)$ is defined as 
\begin{equation}
	\gamma_l(d, \nabla d):=\frac{1}{2}\, \bigg(\frac{d^2}{l_f} + l_f \nabla d \cdot \nabla d \bigg).
\end{equation}
Finally, the fracture contribution of pseudo-energy density given in \req{psuedo-energy} is modified for our topology optimization problem by introducing the topological phase transition function $f(\Phi)$ through the following explicit form
\begin{equation}
	{W_{frac}}(d, \nabla d,\Phi) =  
	f(\Phi)[1 - {g}(d)]\; \psi_c +
	2 f(\Phi){\psi_c}\; l_f \;{\gamma}_l(d, \nabla d)\;,
	\label{frac-part}
\end{equation}
where ${\psi}_c > 0$ is the so-called critical fracture density energy. 
The critical elasticity density function depends on the critical effective stress $\sigma_c$ or Griffith's energy release rate $G_c$, as outlined in \cite{aldakheel+blaz+wriggers18}
\begin{equation}
	\psi_c = \frac{\sigma_c^2}{2E} = \frac{3}{8 l_f \sqrt{2}} G_c \; .
\end{equation}
By taking the variational derivative $\delta_d W$ of \req{psuedo-energy}, and employing some additional algebraic manipulation, the  second PDE is derived for topology optimization of  fracture-resistance  in the rate-dependent setting as
%
\begin{equation}
		\fterm{
	\eta_f f(\Phi)\dot{d}=(1-d)f(\Phi) \calH-\Big[f(\Phi)d - l_f^2 \div[f(\Phi).\nabla d]\Big]\;, 
}
	\label{frac-eqs}
\end{equation}
along with $f(\Phi)\nabla d\cdot \Bn=0$ on $\partial \calB$. \noii{Here, $\eta_f$ is a material parameter (also known as mobility parameter) that characterizes the viscous response of the fractured state, see \cite{miehe+welschinger+hofacker10a,miehe+schaenzel+ulmer15}.}
Additionally, we note that the following identity has been used:
\begin{equation}
	\delta_d\left[ f(\Phi)\gamma_l \right] =f(\Phi)\delta_d\gamma_l-l_f\nabla d\cdot f(\Phi) \WITH \delta_d\gamma_l=\frac{1}{l_f}(d-l^2_f\Delta d)
	\label{frac-eqs1}
\end{equation}
Here, the crack driving force function is shown as  $\calH(\Bx,t)$ reads:
\begin{equation}
	\calH(\Bx,t):=\max_{s\in[0,t]}\widetilde{D}\big(\BfrakC(\Bx,s)\big) \WITH  \widetilde{D}:=\zeta{\Big\langle}\frac{\psi_e^{+} + \psi_p}{{\psi}_c} - 1 {\Big\rangle} ,
	\label{histfield_m3}
\end{equation}
where, the Macaulay brackets denotes the ramp function $\langle x \rangle:=(x+|x|)/2$. \noii{We note that the crack driving force $\widetilde{D}$ is computed based on effective quantities (i.e., no degradation function is included), see \cite{miehe+schaenzel+ulmer15} Sections 3.1.2-3.1.3.} Additionally, it is worth noting that if we set $f(\Phi)=1$ in \req{frac-eqs}, we recover the standard evolution equation for phase-field fracture as given in \cite{miehe+schaenzel+ulmer15}.
Also, note that $\calH(\Bx,t)$ ensures 
that the \textit{local irreversibility} condition (positivity of the fracture dissipation) {for avoiding crack healing} is satisfied, i.e., $\dot{d} \geq 0$ \cite{noii2021quasi}. Additionally, $\zeta\geq0$ is a scaling parameter that introduces further flexibility in the formulation, allowing to tune the post-critical range to better match with experimental results (cf.~\cite{noii2021bayesian}).

\sectpb{Dissipation for the rate of minimization principle}
%
Gradient theories for standard dissipative solids in elastic-plastic materials undergoing fracture are governed by two scalar constitutive functions consisting of the energy storage and the dissipation functions. In the previous section, we described the internal energy storage, and in this section, the dissipative response of the material is further elaborated.  Similar to  the energy density functional, the dissipative material response is additively decomposed into plastic deformations, and viscosity part due to fracture evolution through:
\begin{equation}\label{eq:dissipation}
	\mathcal{D}(\dot{\BfrakC})=\mathcal{D}_{plas}(\dot{\BfrakC})+\mathcal{D}_{vis}(\dot{\BfrakC}).
\end{equation}
Next, we explain individual contributions of the dissipation potential and their ingredients.

\sectpc[Section232]{Plastic dissipation}

The plastic dissipation-potential density function provides a major restriction on constitutive equations for elastic-plastic and dissipative materials based on the principle of maximum dissipation. This restriction is due to  the second law of thermodynamics (Clausius-Planck inequality) within in \textit{a reversible (elastic) domain} in the space of the dissipative forces.
Defining a dual driving force $\{\Bs^p , -\Bfrakh^p\}$ with respect to the primary fields $\{\Bve^p,\alpha\}$, and applying the Coleman-Noll procedure to the free energy density function \req{psuedo-energy} yields the following thermodynamic conjugate variables:
\begin{align}
	\Bs^p := -\partial_{\Bve^p} {{W}}
	= \Bsigma(\Bve_e,\Phi) =f(\Phi) \widehat{\Bsigma}(\Bve_e)
	\AND \Bfrakh^p:=\delta_\alpha{{W}}
	=f(\Phi) g(d)\,h\alpha \;.
	\label{dual}
\end{align}
We note that the constitutive equations in \req{dual} are modified through topological phase transition function $f(\Phi)$. In agreement with the classical setting of elasto-plasticity, the yield function is defined as
\begin{equation}
	\beta^p(\Bs^p,\Bfrakh^p,d,\Phi):= \hbox{$\sqrt{3/2}$}\;\vert \BF^p \vert - \Bfrakh^p- f(\Phi) g(d)\sigma_Y 
	\WITH
	\BF^p:= \dev[\Bs^p] = \Bs^p - \frac{1}{3} \mbox{tr} [\Bs^p] \text{\BI} .
\end{equation}
With the yield function at hand, the dissipation-potential  density function for plastic response reads
\begin{equation}
	\widehat{{\Phi}}_{plas}(\dot\Bve^p,\dot\alpha,d,\Phi) = \sup_{\{\Bs^p,\Bfrakh^p\}} \{ \Bs^p:\dot{\Bve}^p - \Bfrakh^p\dot\alpha  \mid \ \beta^p(\Bs^p,\Bfrakh^p;d)\leq 0 \} ,
	\label{phiP}
\end{equation}
Following \cite{noii2021bayesian} after some algebraic manipulation, we can thus write~\eqref{phiP} as:
\begin{equation}
	\widehat{\Phi}_{plas}(\dot{\bm{\varepsilon}}^p,\dot\alpha,d,\Phi)=
	\begin{dcases}
		\sqrt{2/3} f(\Phi)\,  g(d)\sigma_Y  \vert \dot{\bm{\varepsilon}}^p \vert \quad &\text{if} \ \  \vert \dot{\bm{\varepsilon}}^p \vert = \sqrt{3/2}\,\dot\alpha,  \\
		+\infty \quad &\text{otherwise}.
	\end{dcases}
	\label{eq:pdiss_relation}
\end{equation}
which follows from the principle of maximum plastic dissipation, and is modified through $f(\Phi)$. Taking the supremum of inequality  function \req{eq:pdiss_relation}  yields, as a necessary condition, the primal representation of the plasticity evolution problem in the form of a Biot-type equation:
\begin{equation}
	\{\Bs^p , -\Bfrakh^p\}\in\partial_{\{\dot\Bve^p,\dot\alpha\}}\,\widehat{{\Phi}}_{plas}(\dot\Bve^p,\dot\alpha,d,\Phi) .
	\label{eq:pbiot}
\end{equation}
The Euler equations of the maximization principle~\eqref{phiP} yield the corresponding the flow rule and hardening law
\begin{equation}
	\dot\Bve^p = \lambda^p \hat\Bn
	\WITH
	\hat\Bn=\frac{\partial \beta}{\partial \Bs^p}  ,
	\AND
	\dot\alpha = -\lambda^p\frac{\partial \beta}{\partial \Bfrakh^p} ,
	\label{plast-evol-eqs}
\end{equation}
together with the KKT conditions
\begin{equation}
	\beta^p \le 0, \quad\quad  \quad\quad \lambda^p \ge 0, 
	\quad\quad \mbox{and} \quad\quad
	\beta^p\; \lambda^p = 0 .
	\label{palst-kkt}
\end{equation}
Equations~\eqref{plast-evol-eqs} and~\eqref{palst-kkt} constitute the so-called \textit{dual form} of the elasto-plastic problem in strong form. Finally,  the dissipation potential functional for plastic flow is defined as
\begin{equation}
		\begin{aligned}
	{\mathcal{D}}_{plas}(\dot\Bve^p,\dot\alpha,d,\Phi)
	=  \int_\calB{	\widehat{{\Phi}}_{plas}}~\text{d}\Bx.
 	\end{aligned}
\end{equation}
 \noii{For sake of completeness, by means of \req{plast-evol-eqs}, the fourth-order elasticity tensor given in \req{eq:elast_tensor1}  reads:
	\begin{equation}\label{eq_C_vol_dev}
		\begin{aligned}
			&\widetilde{\mathbb{C}}^+=KH{^+}(\Bve)\bm I\otimes\bm I+2\mu(1-3\mu\delta_1)\mathbb{P}^{dev}
			+6\mu^2(\delta_1-\delta_2){\hat\Bn}\otimes{\hat\Bn},
			\\ 
			&\mathbb{C}^-=K\big(1-H{^+}(\Bve)\big)\bm I\otimes\bm I~,
		\end{aligned}
	\end{equation}
	with scaling factors
	\begin{equation}\label{eq_C_vol_dev1}
		\delta_1=\frac{\lambda^p}{\sqrt{3/2}||{\widetilde{\Bsigma}}^{dev}||+3\mu\lambda^p}
		\AND
		\delta_2=\frac{1}{3\mu+H}\;,
	\end{equation}
	see \cite{de2011computational} for further details on computational methods for plasticity.}
\begin{Remark}
	\label{rem:redu_yield_fun}
	Note that for the implementation purpose, we let  any $\Bx\in\calY$  be in elastic region, thus we simply set $\beta^p(\Bs^p,\Bfrakh^p;d)<0$ if $f(\Phi)=\kappa$ or if simply $\text{H}(\Phi)<0$ holds.
\end{Remark}
%
\sectpc[Section2320]{Dissipative power density due to fracture viscous}

In line with \citep{noii2021bayesian,miehe+welschinger+hofacker10a} ${\Phi}_{vis}$ denotes the dissipative power density due to viscous resistance forces for fracture part denoted as ${\Phi}^{frac}_{vis} (\dot{d})$ which reads:
\begin{equation}
			\begin{aligned}
	\widehat{\Phi}_{vis} (\dot{d},\Phi) =f(\Phi)\left[ \frac{\eta_f}{2}\dot{d}^{\,2} (\Bx)+I_+(\dot d)\right],
	 	\end{aligned}
	\label{Dvis}
\end{equation}
\noii{which is already given in \req{frac-eqs}. Following \cite{miehe+welschinger+hofacker10a} the indicator function $I_+\colon \mathbb{R} \to \mathbb{R} \cup \{+\infty\}$ has been introduced to impose the irreversibility condition embedded in $d\in\calW^d_{d_n}$.} As before, a global rate potential of the dissipative power density due to viscous reads
\begin{equation}
	\mathcal{D}_{vis}(\dot{d},\Phi) := 
\int_\calB \widehat{\Phi}_{vis} (\dot{d},\Phi) \,\text{d}\Bx,   
	\label{potential-functional-visco}
\end{equation}
%

\sectpb[Section23ext]{Potential energy of the external loading}

The macroscopic continuum domain is assumed to be loaded by macroscopic external field actions. Thus, we assumed the macroscopic  body force \grm{is} denoted by $\overline\Bf$, and  traction field on the surfaces $\partial_N\calB$ is denoted by $ \overline{\Btau}$ for the mechanical contribution by means of topological phase transition function. Thus, the work of external loads due to mechanical contribution reads:
\begin{equation}\label{eq:pot_external}
	\mathcal{E}_{ext} (\dot{\bm u},\Phi) := 
	\int_{\calB} f(\Phi) \overline \Bf \cdot \dot{\bm u}\,\text{d}\Bx +
	\int_{\partial_N\calB} f(\Phi) \overline{\Btau} \cdot \dot{\bm u}~\mathrm{d}a,
\end{equation}
%

\sectpb[Section232]{Potential functional for the fracturing elasto-plastic media}
%
Here, the governing equations of the fracturing elasto-plastic solid can be derived from the basis of the energy functional~\eqref{potential-energy-functional}, dissipation functional \req{eq:dissipation}, and external potential \req{eq:pot_external} by invoking rate-type variational principles~\citep{miehe2011,miehe2016bphase}.
Thus, the energy functional for the fracturing elastic-plastic solid material is required for the following potential
\begin{align}
	\Pi ( \dot{\Bve},\dot{\Bve^p}, d,\dot\alpha,\Phi)
	:&=\frac{d}{dt}{{\mathcal{E}}} ( {\Bve},{\Bve^p},d,\alpha)
	+\calD( \dot{\Bve},\dot{\Bve^p}, d,\dot\alpha)
	-\mathcal{E}_{ext}({\dot{\bm u}}),
	\label{eq:energy_rate}
\end{align}
The minimization principle of the  rate-dependent gradient-extended energy functional \textit{at fixed topological configuration} (thus $\Phi(\Bx,t)$ is given a-priory, \noii{so that one can distinguish $\Phi(\Bx,t)$ with other primary fields in minimization problem}) can be written in the following compact form
\begin{equation}
	\boxed{ (\dot{\Bve},\dot{\Bve^p}, d,\dot\alpha) = \arg\big\{ 
		\inf_{\Bu\in\calW_{\overline{\Bu}}^{\Bu}} \ 
		\inf_{ d\in \calW^{d}} \	
		\inf_{\{\Bve^p,\alpha\}\in\mathbf{Q}\times \calW^\alpha_{\alpha_n,\,\Bq}} \,  \Pi(\dot{\Bve},\dot{\Bve^p}, d,\dot\alpha;\Phi) \big\} , }
	\label{minimzation_rate_energy1}
\end{equation}
%
\sectpb[Section23]{Potential energy for the incremental
	minimization principle}
\label{sec:Inc_specifi energy}
%
In this section, the incremental minimization principle for a class of gradient-type dissipative materials in \req{minimzation_rate_energy1} is employed to derive the governing equations. Let, the loading
interval {$\calT_f := (t_0,T_f)$}
be discretized using the discrete-time (loading) points
\begin{equation}\label{eq:time_int_u}
	0 = t_0 < t_1 < \ldots < t_n < t_{n+1} < \ldots < t_N = T_f,
\end{equation}
{with the end time value $T_f>0$.} {The parameter $t\in \calT_f$ denotes the time. We note that for rate-independent problems the time corresponds to an incremental loading parameter.} In order to advance the solution within a specific time step, we focus on the finite time increment $[t_n,t_{n+1}]$, with
\begin{equation}
	\tau_f=t_{n+1}-t_{n}>0,
	\label{eq:time_interval}
\end{equation}
denoting the time step. To formulate the incremental variational principles at the current time $t_{n+1}$ associated with known fields at $t_n$, we need to elaborate  the incremental energy storage, dissipation potential, and external work done at a finite time step $[t_{n},t_{n+1}]$.
The global rate potential form $\Pi(\Bu,\Bve^p,\alpha,d;\Phi)$ given in \req{minimzation_rate_energy1}, inline with our recent study in \cite{noii2021bayesian}, it can be written in an  incremental energy form $\Pi^\tau(\Bu,\Bve^p,\alpha,d;\Phi)$ for a time increment $[t_{n}, t_{n+1}]$, as	
\begin{equation}\label{final_inc_pot0}
	\begin{aligned}
	\Pi^\tau( {\Bve},{\Bve^p}, d,\alpha,\Phi)
	:&=\int^{t_{n+1}}_{t_{n}} \Pi ( \dot{\Bve},\dot{\Bve^p}, d,\dot\alpha,\dot\Phi) ~\mathrm{d}t\\
	&={\mathcal{E}}^\tau( {\Bve},{\Bve^p},d,\alpha,\Phi)
	+\calD^\tau( \dot{\Bve},\dot{\Bve^p}, d,\dot\alpha,\Phi)
	+\mathcal{E}^\tau_{ext}(\bm u,\Phi),\\
\end{aligned}
\end{equation}	
which can be further expanded as
	\begin{equation}\label{final_inc_pot}
		\begin{aligned}
	\Pi^\tau( {\Bve},{\Bve^p}, d,\alpha,\Phi)	&=\int_{\calB} \Big({W}(\BfrakC)-{W}(\BfrakC_n)\Big)~\text{d}\Bx\\
	&+\int_\calB
	\Big({\tau_f\widehat{{\Phi}}^\tau_{vis}}
	+\tau_f\widehat{{\Phi}}^{{\tau}}_{plast}
	+ I_+( d-d_n)\Big)~\text{d}\Bx\\
	&+\int_{\calB} f(\Phi)\overline\Bf \cdot ({\Bu-\Bu_n})~\text{d}\Bx
	+\int_{\partial_N\calB} f(\Phi)\overline{\Btau} \cdot ({\Bu-\Bu_n})~\mathrm{d}a,\\
\end{aligned}
\end{equation}
where the incremental dissipative power density due to viscous effect is given by
\begin{equation}
   \Phi^{\tau}_{vis}= f(\Phi)\left[\frac{\eta_f}{2\tau_f^2}(d-d_n)^{\,2} +I_+(d-d_n)\right],
	\label{phiP_tau123}
\end{equation}
and the incremental plastic dissipation potential is given as
\begin{equation}
	\widehat{{\Phi}}^\tau_{plas}(\Bve^p,\alpha,d;\Bve_n^p,\alpha_n,\Phi) = 
	\frac{1}{\tau_f}\sup_{\{\Bs^p,\Bfrakh^p\}} \{ \Bs^p:\big({\Bve}^p-{\Bve}_n^p\big) - \Bfrakh^p\big(\alpha-\alpha_n\big)  \ \mid \ \beta^p(\Bs^p,\Bfrakh^p,d)\leq 0 \} .
	\label{eq:pdiss_incr1}
\end{equation}
An incremental form of  Euler equations written in terms of the flow rule and hardening law is given by:
\begin{equation}\label{eq:flow_rule_incr}
	\Bve^p = \Bve_n^p+\tau_f\lambda^p \frac{\partial \beta^p}{\partial \Bs^p}  
	\AND
	\alpha = \alpha_n+\tau_f\lambda^p\frac{\partial \beta^p}{\partial \Bfrakh^p} ,
\end{equation}
at time $t_n$. Additionally, the direction of the plastic flow is further defined through
\begin{equation}
	{\hat\Bn}:=\frac{\partial \beta^p}{\partial \Bs^p}=\hat\Bn^{trial}=\frac{\BF^{p,trial}}{\vert \BF^{p,trial} \vert}
	=\frac{\Bve^p-\Bve_n^{p}}{\vert \Bve^p-\Bve^p_n \vert}
\end{equation}
where, from standard arguments of von Mises plasticity, $\BF^{p,trial}:=\BF^{p}(\Bve,\Bve^p_n,d)$. Thus, incremental plastic dissipation potential given \req{eq:pdiss_incr1} becomes
\begin{equation}
	\widehat{{\Phi}}^\tau_{plas}(\dot{\bm{\varepsilon}}^p,\dot\alpha;d,\Phi)=
	\begin{dcases}
		\sqrt{2/3} \,  g_p(d)f(\Phi)\sigma_Y  \vert {\bm{\varepsilon}}^p-{\bm{\varepsilon}}_n^p \vert \quad 
		&\text{if} \ \  \vert \dot{\bm{\varepsilon}}^p-{\bm{\varepsilon}}_n^p \vert = \sqrt{3/2}\,(\alpha-\alpha_n),  \\
		+\infty \quad &\text{otherwise}.
	\end{dcases}
	\label{eq:pdiss_incr2}
\end{equation}
for a detail derivation, see \cite{noii2021bayesian} Section 2.3. Following the incremental potential given in \req{final_inc_pot}, the time-discrete counterpart of the canonical rate-dependent variational principle in \req{minimzation_rate_energy1} takes the following compact form:
\begin{equation}
	\boxed{ ({\Bve},{\Bve^p}, d,\alpha) = \arg\big\{ 
		\inf_{\Bu\in\calW_{\overline{\Bu}}^{\Bu}} \ 
		\inf_{ d\in \calW^{d}} \	
		\inf_{\{\Bve^p,\alpha\}\in\mathbf{Q}\times \calW^\alpha_{\alpha_n,\,\Bq}} \,  \Pi^\tau({\Bve},{\Bve^p}, d,\alpha;\Phi)  \big\} . }
	\label{minimiaz2}
\end{equation}
Finally, we seek to find the stationarity conditions for the minimization problem~\eqref{minimiaz2}, with following formulation in abstract from.

\begin{form}[Weak form of the Euler-Lagrange equations for ductile phase-field fracture]
	\label{form_3}
	Let constants $(K, \mu,H,$ $\sigma_Y, \psi_c,l_f,\zeta,\eta_f) \geq 0$ be given with the initial conditions $\Bu^0=\Bu(\bm{x},0)$, $\Bve^p_0=\Bve^p({x},0)$, $\alpha_0=\alpha({x},0)$, and $d^0=d(\bm{x},0)$ for a fixed topological field with $\Phi$. For the time increments $n\in\{1,2,\ldots, N\}$, we solve a displacement equation where we seek $\Bu:= \Bu^n\in \calW_{\overline{\Bu}}^{\Bu}$ for a fixed $\Phi\in\calW_{\overline{\Phi}}^{\Phi}$ such that
	\begin{equation}\label{eq:couple_1}
		\begin{aligned}
			&\int_\calB\big[\Bsigma(\Bve,\Bve^p,d,\Phi):\Bve(\delta\Bu)-f(\Phi)\overline{\Bf}\cdot\delta\Bu\big]\,\text{d}\Bx - \int_{\partial\calB_N}f(\Phi)\overline{\Btau}\cdot\delta\Bu\,da =0 \quad  \forall \,\delta\Bu\in\calW^{\Bu}_0\  . 
		\end{aligned}
	\end{equation}
	The phase-field variational formulation is needed to find $d:= d^n: \calB \in \mathcal{W}^d_{d_{n-1}}$ for a fixed $\Phi\in\calW_{\overline{\Phi}}^{\Phi}$ such that
	\begin{equation}\label{eq:couple_2}
		\int_\calB \bigg[ \bigg(f(\Phi)(1-d)\calH -  {d} + \frac{\eta_f}{\Delta t}(d-d^n)\bigg)\delta d - l_f^2f(\Phi)  \nabla d \cdot \nabla (\delta d) \bigg]\,\text{d}\Bx = 0   \;\;\forall \,{\color{black}\delta d\in\calW^{d}}  . 
	\end{equation}
	along with $\nabla d \cdot \Bn=0 \; \text{on}\; \partial \calB$. This leads to a coupled multi-field problem is  defined for equations \req{eq:couple_1}, and \req{eq:couple_2} together with \req{eq:flow_rule_incr} and solved for $(\Bu,d)$  for a fixed $\Phi$.
\end{form}

\begin{Remark}
	\label{rem:redu_brit_frac}
	Up to this point, we have considered a ductile phase-field fracture for a solid exhibiting an elastic-plastic response followed by damage (hereafter $\texttt{E-P-D}$); or elastic, followed by elastic-plastic, and then plastic-damage (hereafter $\texttt{E-P-DP}$). To reduce the model into a brittle fracture response, it is sufficient that one imposes $\sigma_Y\rightarrow\infty$ (hereafter $\texttt{E-D}$) for which the yield surface becomes negative.
\end{Remark}
\sectpa[Section3]{Topology Optimization of Fracturing Elasto-Plastic Material}
%
The mathematical formulation for topology optimization of a solid subject to the equality/inequality \grm{constraints} reads:
%
\begin{equation}\label{compat_argmin}
		\Phi = 
		\substackrel{\Phi \in 	\calW_{\overline{\Phi}}^{\Phi} }{\mbox{argmin}}\,
		\Big\{[\;\BfrakJ(\BfrakU)] \quad |\quad 
	\BfrakR(\BfrakU)=0 \AND \BfrakG(\BfrakU)\leq \bar{G} \Big\}
\end{equation}
for unknown parameters $\Phi$ with
\begin{equation} 
	\begin{aligned}
		&\BfrakJ(\BfrakU)
		=\int_{\calB}{\widehat{F}}(\Bu,d,\Phi)  \mathrm{d}{\bm{x}} 
		\WITH \widehat{F}(\Bu,d,\Phi)=F(\Bu,d,\Phi)\text{H}^{c_1}(\Phi),\\
		&\BfrakR(\BfrakU)
		=\int_{\calB}{\widehat{R}}(\Bu,d,\Phi)  \mathrm{d}{\bm{x}} 
		\WITH \widehat{R}(\Bu,d,\Phi)=R(\Bu,d,\Phi)\text{H}^{c_2}(\Phi),\\
		&\BfrakG(\BfrakU)
		=\int_{\calB}{\widehat{G}}(\Bu,d,\Phi)  \mathrm{d}{\bm{x}} 
		\WITH \widehat{G}(\Bu,d,\Phi)=G(\Bu,d,\Phi)\text{H}^{c_3}(\Phi).
	\end{aligned}
\end{equation}		

Here, $ \BfrakJ(\BfrakU)$ is a functional (so-called objective functional) with its density function $ \widehat{F} $ for a global primary field $ \BfrakU$ within the structural design domain $ \calB $.  The minimization problem is subjected to equality  constraint functional as $\BfrakR(\BfrakU)$ with its density $ \widehat{R} $, along with  inequality  constraint functional denoted as $\BfrakG(\BfrakU)$ with its density $  \widehat{G}$. 
\noii{Here, the set of penalty parameters $(c_1,c_2,c_3)$ for the Heaviside step function is used to for $(i)$ resolving the structural analysis by further penalizing the material point to void transition rule, and $(ii)$ to avoid the singularity in sensitivity analysis. Nevertheless, we note it is required to highlight that since an exact Heaviside step function is employed in this manuscript, the only advantage of using penalty parameters here is to avoid singularity and not its effect on the structural analysis. We note that in our case following topological phase transition function in \req{fun:topological_phase}, we set $c_1=c_2=c_3=2$}.
Here, to  resolve \req{compat_argmin} the level-set method is employed together with a gradient-based approach, so-called the adjoint  method.

\sectpb[Section30]{Reaction-diffusion based level-set method}

The reaction-diffusion based topology optimization is a well-developed technique for  modeling the evolution of the interfaces. In this approach, as opposed to the typical SIMP formulation \cite{bendsoe2, bendsoe1} the material distribution is characterized by the scalar level-set function $\Phi(\Bx,t) \in \calW_{\overline{\Phi}}^{\Phi}$, see Section 2.1. Considering a closed moving front for the  boundary material $ \Gamma_\Phi(\Bx,t)=\Bx\in \R^n\ $ which is represented by zero level-set of the hyper surface $ \Phi (\Bx,t) $ i.e. $ \Gamma_\Phi(\Bx,t)=\left\{ \Bx\;\colon\; \Phi (\Bx, t)=0 \right\}\ $. To determine the motion/velocity of the interface, one can take the chain rule of  $\Gamma_\Phi(\Bx,t)$, which results in the following Hamilton-Jacobi equation:
\begin{equation}
	\begin{aligned}
		\frac{\partial \Phi (\Bx,t)}{\partial t}=\frac{\partial \Phi (\Bx,t)}{\partial \Bx}\cdot\frac{\partial \Bx}{\partial t}=\Bv\cdot\nabla \Phi, \
	\end{aligned}
\end{equation}
where $\Bv$ is the velocity field of the closed moving front.
Since, outward normal direction is given by $ \widehat{\Bn}=\frac{\nabla \Phi }{\left| \nabla \Phi  \right|}\ $the Hamilton-Jacobi equation can be replaced by:
\begin{equation}\label{eq33}
	\begin{aligned}
		\frac{\partial \Phi }{\partial t}
		=\Bv\cdot\nabla \Phi =\Bv\cdot\nabla \Phi \cdot\frac{\left| \nabla \Phi  \right|}{\left| \nabla \Phi  \right|}
		=\Bv\cdot\left| \nabla \Phi  \right|\cdot\frac{\nabla \Phi }{\left| \nabla \Phi  \right|}
		=\Bv\cdot\widehat{\Bn}_{\Phi}\cdot\left| \nabla \Phi  \right|
		={\widehat{v}_{\Phi}}\left| \nabla \Phi  \right|,
	\end{aligned}
\end{equation}
with $ {\widehat{v}_{\Phi}}=\Bv.\widehat{\Bn}_{\Phi} $ being the normal velocity field of the interface. Equation \req{eq33} is the standard level-set evolution equation of the closed front by changing the zero level-set. The conventional level-set method with the Hamilton-Jacobi equation has many limitations, such as the need for frequent initialization as the starting point that affects the final topology, the need for reshaping as a signed-distance function, time step dependence with Courant- Friedrichs-Lewy (CFL) Conditions and potential instabilities when a hole is nucleated. These limitations affect the accuracy of the final topology, and to this end, following  \cite{wang2003level, Yamada} the reaction-diffusion evolution equation has been suggested as:
\begin{equation}\label{eq:RD_phi}
	\fterm{
		\underbrace{\eta_{\Phi} \dot{\Phi}}_{\text{{topological field update}}} = \underbrace{{\widehat{v}_{\Phi}}(\Phi)}_{\text{topological sensitivity }} + \underbrace{l^2_{\Phi} {{\nabla }^{2}}\Phi}_{\text{\grm{regularization}}}\quad in \;\calB.
	}
\end{equation}
The topological reaction-diffusion equation in \req{eq:RD_phi} is defined in the domain $\calB$ that is augmented with an imposed Neumann homogeneous boundary condition as
\begin{equation}\label{eq30}
	{\nabla \Phi \cdot \Bn} = 0 \quad on \; \partial\calB.
\end{equation}
Here, $ l_{\Phi} $ is the regularization parameter, so-called topological length-scale. Also, $ \eta_{\Phi} $ is denoted as a topological viscosity. Adding the diffusion term $ l_\Phi {{\nabla }^{2}}\Phi \ $ as a Tikhonov regularization technique, provides more stabilization, smoothness, and well-defined boundaries \cite{Tikhonov, Yamada}. 
\begin{Remark}
	\label{rem:regulariaze_topology}
	As a matter of fact, since topological length-scale $l_\Phi$ acts as a regularization indicator for the solution of reaction-diffusion equation, choosing the smaller value for $l_\Phi$ (which depends on discretization size) results in smoother topology, but some narrow members will be removed. Thus, the final topology contains more global voids and so less local voids. For a detailed discussion, see \cite{Yamada}.
\end{Remark}

Next temporal discretization of \req{eq:RD_phi} is performed to derive an incremental topology optimization, which is needed due to the non-linearity of ductile phase-field fracture. 
Here, in the topology optimization process, we deal with pseudo-time for temporal discretization in \req{eq:RD_phi} compared with real-time which needs for path-dependent ductile fracture given in \req{eq:time_int_u}. Let, the pseudo-time interval {$\calT_\Phi := (t_0,T_\Phi)$} be discretized using the discrete pseudo-time through:
\begin{equation}\label{eq:time_int_phi}
0 = t_0 < t_1 < \ldots < t_m < t_{m+1} < \ldots < t_M = T_\Phi,
\end{equation}
with the end pseudo-time value $T_\Phi>0$. In order to advance the solution within a specific time step, we focus on the finite time increment $[t_m,t_{m+1}]$, where
\begin{equation}
	\tau_\Phi=t_{m+1}-t_{m}>0,
	\label{eq:time_interval_phi}
\end{equation}
denotes the pseudo-time step within the topology optimization iterations. To formulate the incremental topology optimization, we aim to determine topology field $\Phi(\Bx,t)$ at the current pseudo-time $t_{m+1}$ associated with known fields at $t_m$. Thus, the following  variational formulation for the reaction-diffusion equation in \req{eq:RD_phi} reads:
\begin{form}[Weak form of the Euler-Lagrange equation for reaction-diffusion equation]
	\label{form_4}
	Let constants $(l_\Phi,\eta_\Phi) \geq 0$ be given with the initial conditions $\Phi^0=\Phi(\bm{x},0)$ for a fixed $(\Bu,d,\alpha)$. For the time increments $m\in\{1,2,\ldots, M\}$, we solve a reaction-diffusion equation where we seek $\Phi:= \Phi^m\in \calW_{\overline{\Phi}}^{\Phi}$ such that
	\begin{equation}\label{eq:couple_3}
		\int_\calB \bigg[ \bigg({\widehat{v}_{\Phi}}(\Phi)- \frac{\eta_\Phi}{\tau_\Phi }(\Phi-\Phi^m)\bigg)\delta \Phi - l_\Phi^2 \nabla \Phi \cdot \nabla (\delta \Phi) \bigg]\,\text{d}\Bx = 0   \;\;\forall \,{\color{black}\delta \Phi\in\calW^{\Phi}_0}  . 
	\end{equation}
	along with $\nabla \Phi \cdot \Bn=0 \; \text{on}\; \partial \calB$.
\end{form}
An implicit geometry parameterization based on Formulation \ref{form_4} for topological field $-1\leq\Phi(\Bx,t)\leq1$ results in material point (i.e., $\Bx\in\Omega$) if $\Phi(\Bx,t)>0$, non-material point (i.e., $\Bx\in\calY$) if $\Phi(\Bx,t)<0$, and represented as a material boundary interface (i.e., $\Bx \in\Gamma_{\Phi}$) if $\Phi(\Bx,t)=0$, which implies zero level-set.
\begin{figure}[!t]
	\centering
	{\includegraphics[clip,trim=1cm 20cm 0cm 0cm, width=16.5cm]{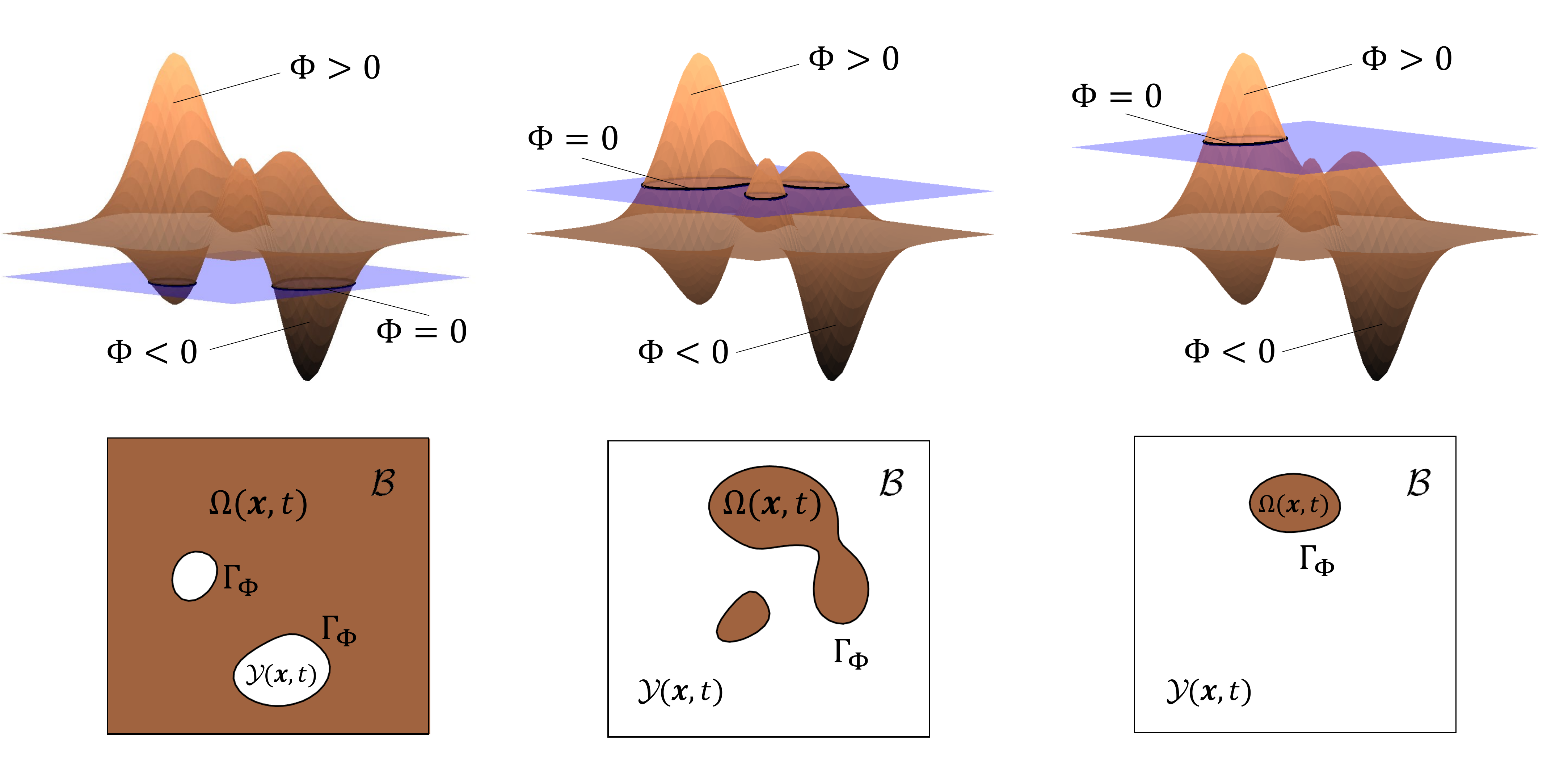}}  
	\caption{Level-set based geometric parameterization based on topological field $\Phi(\Bx,t)$.
	}
	\label{Figure5}
\end{figure}

\sectpb[Section29]{Space-time finite element discretization}
%
In this section, we provide spatial discretizations of the variational forms for Formulation (2.1), and \req{form_4}. This result in a discretized multi-field problem to be solved for three-field unknowns represented by $(\Bu, d,\Phi)$. Here, we use a Galerkin finite element method to discretize the equations with first-order isoparametric elements. In particular, \noii{to demonstrate the versatility} of the proposed model, two types of elements have been used. 
Examples 1-4 are discretized using first-order linear three-dimensional tetrahedral elements $P_1$. Examples 2-3 are discretized using $Q^1$-conforming trilinear hexahedral elements for the three-dimensional problems, i.e., the ansatz and test space uses $Q_1^c$--finite elements, e.g., for details, we refer readers to the \cite{wriggers2008nonlinear}. Let the continuous domains $\calB$ (solid domain), be approximated by $\calB_h$ such that $\calB\approx\calB_h$. Approximated domain $\calB_h$ is then decomposed into non-overlapping finite numbers of elements $\calB_e\subset\calB_h$ such that 
\begin{equation}
	\calB\approx\calB_h = \bigcup_{e}^{n_{e}}\calB_e\;.
\end{equation}
The finite element discretized solutions are approximated by
\begin{equation}
	\Bu^h=\BN_\Bu\hat{\bm u}, \quad
	d^h=\bm{N}_d\hat{\bm d}, \quad
	\Phi^h=\bm N_\Phi\hat{\Phi}.
	\label{Discr1}
\end{equation}
with set of nodal solution values $(\hat{\bm u},\hat{\bm d},\hat{\BA})$, and linear basis functions $\BN_\Bu$, $\BN_d$, and $\BN_\Phi$ correspond to the three primary fields. Also, their derivative follow
\begin{equation}
	\bm\varepsilon(\bm u_h)=\bm B_u\hat{\bm u}, \quad
	\nabla d_h=\bm B_d\hat{\bm d}, \quad
	\nabla \Phi_h=\bm B_d\hat{\Phi}.
	\label{Discr2}
\end{equation}
which represent constitutive state variables by $(\bm\varepsilon(\bm u_h), \nabla d_h, \nabla \Phi_h)$. Here, $\BB_{\Bu}$, $\BB_{d}$, and $\BB_{\Phi}$ are the matrix representation of the shape function's derivative, corresponding to the global deformation, crack phase-field, and topological field, respectively. Next, for each primary field, the set of discretized equilibrium equations based on the \textit{residual} force vector denoted by $\mathbf{R}^{\bullet}$ has to be determined. First, for the displacement deformation field, we have
\begin{equation}
    {{\bf{R}}_\Bu} (\Bu,d,\Phi)=
    \sum_{e} \int_{\calB_e}\big[{\BB^T_\Bu} \Bsigma_h(\Bve,\Bve^p,d;\Phi)- f(\Phi) \overline{\Bf}\cdot\BN_\Bu\big]\,\text{d}\Bx - \int_{\partial\calB_N} f(\Phi) \overline{\Btau}\cdot\BN_\Bu\,d\Ba ={\bm 0}, 
	\label{phf_fem_u}
\end{equation}
next, for the crack phase-field part we have,
\begin{equation}
	\begin{aligned}
		{{\bf{R}}_d}(\Bu,d,\Phi)=
		\sum_{e}\int_{\calB_e} \bigg[ \bigg((1-\kappa)(d_h-1)\calH+  {d} + \frac{\eta_d}{\tau_f}(d_h-d_h^{n})\bigg)\cdot\BN_d
		+   l_d^2\BB_d^T\cdot\nabla d_h \bigg]\,\text{d}\Bx ={\bm 0}. 
		\label{phf_fem_d}
	\end{aligned}
\end{equation}
and lastly for the topological field is given by,
\begin{equation}
	{{\bf{R}}_\Phi}(\Bu,d,\Phi)=
	\sum_{e}\int_{\calB_e} \bigg[ \bigg( {\widehat{v}_{\Phi}} - \frac{\eta_{\Phi}}{\tau_\Phi}(\Phi_h-\Phi^m_h)\bigg)\cdot\BN_\Phi- l_\Phi^2  \BB_\Phi^T\cdot\nabla \Phi_h \bigg]\,\text{d}\Bx ={\bm 0}. 
	\label{phf_fem_phi}
\end{equation}
%

\sectpb[Section29]{Fracture-resistance topology optimization}
%
In the following section, a mathematical formulation for the optimization problem given in \req{compat_argmin} is specified. More specifically, the objective functional $\BfrakJ(\BfrakU)$ along with  equality, and inequality constraint functions are determined. Furthermore, a gradient-based adjoint optimization approach is employed. Hereinafter, to derive the adjoint equations we need to derive sensitivities at time $n\in\calT_N$, see \req{eq:time_int_u}, at fixed topology field $\Phi_{m-1}$ at time ${m-1}\in\calT_M$, see \req{eq:time_int_phi}. Thus, we aim to obtain $\Phi_{m}$ at time ${m}\in\calT_M$. To this end, an objective functional  is defined as a \textit{structural stiffness} of material which has to be maximized while retaining a certain \textit{volume ratio}. To compute the structural stiffness of material specifically for the non-linear response, we determine the accumulation of energy from scratch until the current time by looking at the area under the load-displacement curve \cite{huang2008topology,xia2017evolutionary}. To this end, we consider  the rate of total work functional for $\BfrakJ(\BfrakU)$ which takes the following form:
\begin{equation}
	\begin{aligned}
		\BfrakJ
		=-\int^{t_N}_{0}\int_{\calB}{\bm{\sigma}(\bm u,d,\Phi )}:\dot{\Bve}(\bm u,\Phi)\;\text{d}{\bm{x}}  \;\text{d}{t}=-\int^{t_N}_{0}\int_{\calB}{\Bf(\Phi )}.\dot{\Bu}\;\text{d}{\bm{x}} \;\text{d}{t}
	\end{aligned}
\end{equation}
Since, we are dealing with a rate-dependent nonlinear boundary value problems (i.e., ductile phase-field fracture), we consider incremental topology optimization, see Section \ref{Section3: LG_vol}. Thus, $\BfrakJ(\BfrakU)$ is approximated by the trapezoidal rule up to the prescribed displacement as:
\begin{equation}\label{eq:fun_obj}
	\BfrakJ=\sum_{\text{n}=1}^{N} {\BfrakJ^i} \WITH \BfrakJ^n=-\frac{1}{2}(\mathbf{P}^{n}_\Bu+\mathbf{P}^{n-1}_\Bu).(\Bu^n-\Bu^{n-1})
\end{equation}
which follows an inequality constraint due to the volume ratio by:
\begin{equation}\label{eq:vol_const}
	\BfrakG(\Phi)=\int_{\calB}{\widehat{G}}(\Phi)  \mathrm{d}{\bm{x}} 
	=\int_{\calB}\operatorname{H}(\Phi)\text{d}{\bm{x}}=\widehat{\text{V}}(\Phi)=V(\Omega)\le \bar{\Omega}
\end{equation}
Here, $ \widehat{\text{V}}(\Phi)=V(\Omega)\le\bar{\Omega}\ll V(\calB) $ and $\bar{\Omega}$ are the existing volume and prescribed volume constraint, respectively. Following, \req{eq:time_int_u}, we note that, $ 1\le n\le N $ in \req{eq:fun_obj} indicates the prescribed load increment onto the solid body. 

So, it is now only equality constraints that are left to complete the optimization problem given in \req{compat_argmin}. Here, two
different types of formulations are proposed. The first one requires only the deformation force vector as a constraint, thus optimization problem given in \req{compat_argmin} is reduced to:
\begin{equation*}\label{topology_BVP1}
	\texttt{Formulation 1}:\quad \Phi = 
	\substackrel{\Phi \in \calW_{\overline{\Phi}}^{\Phi}  }{\mbox{argmin}}\,
	\Big\{[\;\BfrakJ(\BfrakU)] \quad |\quad 
	\mathbf{P}_\Bu-{\mathbf{R}_\Bu}={\bm 0}, \AND \BfrakG(\Phi)\leq \bar{\Omega} \Big\}\;.\hspace*{1.5cm}
\end{equation*}
Since, two distinguished fields are required to describe ductile phase-field fracture ($\Bu,d$), thus the second formulation is constrained with residual force vector of both
deformation and phase-field fracture. Thus, we have:
\begin{equation*}\label{topology_BVP2}
	\texttt{Formulation 2}:\quad \Phi = 
	\substackrel{\Phi \in \calW_{\overline{\Phi}}^{\Phi}  }{\mbox{argmin}}\,
	\Big\{[\;\BfrakJ(\BfrakU)] \quad |\quad 
	(\mathbf{P}_\Bu-{\mathbf{R}_\Bu},\;{\mathbf{R}_d})=({\bm 0},{\bm 0}), \AND \BfrakG(\Phi)\leq \bar{\Omega} \Big\}\;.
\end{equation*}
It is worth noting that Formulation 2 is a generalization of Formulation 1, as it includes additional constraints of the fracture in the optimization process. Accordingly, in the following, we will derive the sensitivity analysis due to Formulation 2, and then we will reduce it to Formulation 1. In the first numerical example in Section 4, we investigate the efficiency and accuracy of Formulations 1 and 2.
%
\sectpb[Section29]{Sensitivity analysis}
%
In order to derive the gradient-based fully coupled adjoint-based sensitivity analysis, derivatives of the objective and constraint functions with respect to the design variables need to be provided. First, we consider here Formulation 2, which is more general, and next, we will reduce it to Formulation 1. Let \noii{$ \bm{\mu }\ $ and $ \bm{\lambda }\ $} describe the adjoint vectors for time $n-1$, and $n$, respectively. The Lagrangian functional of the optimization problem corresponds to Formulation 2 reads:
\begin{equation}\label{LGforopt}
	\begin{aligned}
	\BfrakL(\Phi,{\bm\lambda_{\Bu}},{\bm\lambda_{d}},{\lambda_{V}},\noii{{\bm\mu_{\Bu}},{\bm\mu_{d}}};{\Bu},{d},{\alpha})
	=\sum_{\text{n}=1}^{N}\Big({\BfrakJ}^{\text{n}}
	&+({\bm\mu^{n}_{\Bu}})^T(\mathbf{P}^{n-1}_\Bu-\mathbf{R}^{n-1}_\Bu)
	+({\bm\lambda^{n}_{\Bu}})^T(\mathbf{P}^{n}_\Bu-\mathbf{R}^{n}_\Bu)\\
	&+({\bm\mu^{n}_{d}})^T \mathbf{R}^{n-1}_\text{d}
	+({\bm\lambda^{n}_{d}})^T \mathbf{R}^{n}_\text{d}\Big)
	+ {\lambda^{T}_{V}} \BfrakG(\Phi) . \ \\
	\end{aligned}\	
\end{equation}
The minimization problem for the given Lagrangian functional of the ductile phase-field fracture in \req{LGforopt} takes the following compact form:
\begin{equation}\label{compat_argmin2}
	\fterm{ 
		\{ \Phi,{\bm\lambda_{\Bu}},{\bm\lambda_{d}},{\bm\lambda_{V}},\noii{{\bm\mu_{\Bu}},{\bm\mu_{d}}}\} = 
		\mbox{arg} \{\; 
		\substackrel{\bm{u} \in \calW_{\overline{\Phi}}^{\Phi}}{\mbox{min}}
		\;\substackrel{{\bm\lambda_{\Bu}},{\bm\lambda_{d}},{\bm\lambda_{V}}\in \BL^2 }{\mbox{max}} \,
		[\; \BfrakL(\Phi,{\bm\lambda_{\Bu}},{\bm\lambda_{d}},{\lambda_{V}};{\Bu},{d},{\alpha}) \; ]
		\; \}.
	}
\end{equation}
To perform the adjoint sensitivity analysis, the derivatives of the Lagrangian function with respect to the design variables can be expressed as:
\begin{equation}\label{eq:sens_G}
	\begin{aligned}
		\mathcal{G}(\Phi,{\bm\lambda_{\Bu}},{\bm\lambda_{d}},{\lambda_{V}},\noii{{\bm\mu_{\Bu}},{\bm\mu_{d}}};{\Bu},{d},{\alpha})
		:=D_{\Phi} \BfrakL
		&=\sum_{\text{n}=1}^{N} \Big(
		D_{\Phi} \BfrakJ^i
		+({\bm\mu^{n}_{\Bu}})^TD_{\Phi} \mathbf{R}_\Bu^{n-1}
		+({\bm\lambda^{n}_{\Bu}})^TD_{\Phi} \mathbf{R}_\Bu^{n}\\[2 mm]
		&+({\bm\mu^{n-1}_{d}})^TD_{\Phi} \mathbf{R}_d^{n-1}
		+({\bm\lambda^{n}_{d}})^TD_{\Phi} \mathbf{R}_d^{n}
		\Big)\\[2 mm]
		&+{\lambda^{T}_{V}} \frac{d \widehat{G}(\Phi)}{d \Phi}
	\end{aligned}
\end{equation}
To formulate \req{eq:sens_G}, directional derivative of the objective functional $\BfrakJ$ reads:
\begin{equation}\label{eq:sens_J}
	\begin{aligned}
		D_{\Phi} \BfrakJ^i=
	    -\frac{1}{2}\Big(   D_{\Phi} \mathbf{P}_\Bu^{n-1} +D_{\Phi} \mathbf{P}_\Bu^{n} \Big)\Delta \Bu^i
	    -\frac{1}{2}\Big( \mathbf{P}_\Bu^{n-1}+\mathbf{P}_\Bu^{n}\Big) \frac{d \Delta \Bu^i}{d\Phi} 
	\end{aligned}\
\end{equation}
Since, on the Dirichlet boundary $\partial_D\calB $, the displacement field is independent to the topological field  (i.e., the topological field is fixed, so  $\Phi=1$ on $\partial_D\calB $), and additionally, for any internal nodes within domain $\mathbf{P}(\Bx)=\bm{0}$ for $\Bx\in\calB$, so we have the following identity:
\begin{equation}\label{eq:sens_J1}
	\begin{aligned}
		\mathbf{P}_\Bu^{n-1}\frac{d \Delta \Bu^n}{d\Phi}=0  \AND \mathbf{P}_\Bu^{n}\frac{d \Delta \Bu^n}{d\Phi}=0 
	\end{aligned}\
\end{equation}
thus, \req{eq:sens_J} is reduced to
\begin{equation}\label{eq:sens_J2}
	\begin{aligned}
		D_{\Phi} \BfrakJ^n=
		-\frac{1}{2}\Big(   D_{\Phi} \mathbf{P}_\Bu^{n-1} +D_{\Phi} \mathbf{P}_\Bu^{n} \Big)\Delta \Bu^n
	\end{aligned}\
\end{equation}
We next continue by performing the directional derivative of the residual force vector for displacement as:
\begin{equation}\label{eq:sens_R_u_i}
	\begin{aligned}
		D_{\Phi} \mathbf{R}_\Bu^{\bullet}
		&=\frac{d \mathbf{R}_\Bu^{\bullet}}{d \Phi}
		+\frac{\partial \mathbf{R}_\Bu^{\bullet}}{\noii{\partial} \Delta \Bu^{\bullet}}\cdot\frac{d\Delta \Bu^{n}}{d\Phi}
		+\frac{\partial \mathbf{R}_\Bu^{\bullet}}{\noii{\partial}  \Delta d^{\bullet}}\cdot\frac{d \Delta d^{\bullet}}{d\Phi}\\[2 mm]
		&=\frac{d \mathbf{R}_\Bu^{\bullet}}{d \Phi}
		+\textbf{K}_{\Bu,\Bu}^{\bullet} \cdot\frac{d\Delta \Bu^{\bullet}}{d\Phi}
		+\textbf{K}_{\Bu,d}^{\bullet} \cdot\frac{d \Delta d^{\bullet}}{d\Phi}\;,
		\quad\text{for}\quad\bullet\in\{n-1,n\}
	\end{aligned}
\end{equation}
which follows by directional derivative of residual force vector for crack phase-field as: 
\begin{equation}\label{eq:sens_R_d_i}
	\begin{aligned}
		D_{\Phi} \mathbf{R}_d^{\bullet}
		&=\frac{d \mathbf{R}_d^{\bullet}}{d \Phi}
		+\frac{\partial \mathbf{R}_d^{\bullet}}{\noii{\partial} \Delta \Bu^{\bullet}}\cdot\frac{d\Delta \Bu^{\bullet}}{d\Phi}
		+\frac{\partial \mathbf{R}_d^{\bullet}}{\noii{\partial}  \Delta d^\bullet}\cdot\frac{d \Delta d^{\bullet}}{d\Phi}\\[2 mm]
		&=\frac{d \mathbf{R}_d^{\bullet}}{d \Phi}
		+\textbf{K}_{d,\Bu}^{\bullet} \cdot\frac{d\Delta \Bu^{\bullet}}{d\Phi}
		+\textbf{K}_{d,d}^{\bullet} \cdot\frac{d \Delta d^{\bullet}}{d\Phi}\;.%
				\quad\text{for}\quad\bullet\in\{n-1,n\}
	\end{aligned}
\end{equation}
Additionally, the sensitivity analysis of the volume constraint in \req{eq:vol_const} reads:
\begin{equation}\label{eq:sens_vol}
	\begin{aligned}	
		\frac{d }{d \Phi }\widehat{G}(\Phi)=\frac{d }{d \Phi }\int_{\calB}\text{H}(\Phi)=\int_{\calB}{\delta \left( \Phi \right)\mathrm{d}{\bm{x}} }. 
	\end{aligned}\
\end{equation}
So, the sensitivity analysis given in \req{eq:sens_G} using (\req{eq:sens_J2}-\req{eq:sens_vol}),  results in
\begin{equation}\label{eq45}
	\begin{aligned}
		\mathcal{G}^{\;n}(\Phi,{\bm\lambda_{\Bu}},{\bm\lambda_{d}},{\lambda_{V}},{\bm\mu_{\Bu}},{\bm\mu_{d}};{\Bu},{d},{\alpha})=
		\sum_{\text{n}=1}^{N} \Big(
		&({\bm\mu^{n}_{\Bu}}-\frac{1}{2}\Delta \Bu^{n})^T\frac{d \mathbf{P}_\Bu^{n-1}}{d \Phi}
		+({\bm\lambda^{n}_{\Bu}}-\frac{1}{2}\Delta \Bu^{n})^T\frac{d \mathbf{P}_\Bu^{n}}{d \Phi}\\[2 mm]
		&+({\bm\mu^{n}_{\Bu}})^T\frac{d \mathbf{R}_\Bu^{n-1}}{d \Phi}
		+({\bm\lambda^{n}_{\Bu}})^T\frac{d \mathbf{R}_\Bu^{n}}{d \Phi}\\[2 mm]
		&+({\bm\mu^{n}_{d}})^T\frac{d \mathbf{R}_d^{n-1}}{d \Phi}
		+({\bm\lambda^{n}_{d}})^T\frac{d \mathbf{R}_d^{n}}{d \Phi}\\[2 mm]
		&+\Big(
		({\bm\mu^{n}_{\Bu}})^T \textbf{K}_{\Bu,\Bu}^{n-1}
		+({\bm\mu^{n}_{d}})^T \textbf{K}_{d,\Bu}^{n-1}
		\Big)\cdot \frac{d\Delta \Bu^{n-1}}{d\Phi}\\[2 mm]
		&+\Big(
		({\bm\lambda^{n}_{\Bu}})^T \textbf{K}_{\Bu,\Bu}^{n}
		+({\bm\lambda^{n}_{d}})^T \textbf{K}_{d,\Bu}^{n}
		\Big)\cdot \frac{d\Delta \Bu^{n}}{d\Phi}\\[2 mm]
		&+\Big(
		({\bm\mu^{n}_{\Bu}})^T \textbf{K}_{\Bu,d}^{n-1}
		+({\bm\mu^{n}_{d}})^T \textbf{K}_{d,d}^{n-1}
		\Big)\cdot \frac{d\Delta d^{n-1}}{d\Phi}\\[2 mm]
		&+\Big(
		({\bm\lambda^{n}_{\Bu}})^T \textbf{K}_{\Bu,d}^{n}
		+({\bm\lambda^{n}_{d}})^T \textbf{K}_{d,d}^{n}
		\Big)\cdot \frac{d\Delta d^{n}}{d\Phi} \Big)\\[2 mm]
	\end{aligned}\
\end{equation}
To simplify \req{eq45} which helps for further condensation, we split all degrees of freedom  into prescribed nodes (corresponds to non-zero Dirichlet boundary), and interface nodes denoted as $\{a,b\}$. So, following \req{eq45} the first two terms can be omitted, if for the nodal points within non-zero Dirichlet boundary, the following constraints hold:
\begin{equation}
	\texttt{C}^n_1=[\bm\mu^{n}_{\Bu}]^a-\frac{1}{2}[\Delta \Bu^{n}]^a={\bm 0} \AND
	\texttt{C}^n_2=[\bm\lambda^{n}_{\Bu}]^a-\frac{1}{2}[\Delta \Bu^{n}]^a={\bm 0}\;,
\end{equation}
thus yields
\begin{equation}
	[\bm\mu^{n}_{\Bu}]^a=[\bm\lambda^{n}_{\Bu}]^a=\frac{1}{2}[\Delta \Bu^{n}]^a
\end{equation}
This is a classical way of formulating an evolutionary topology optimization of elastic-plastic
structures, see for detailed discussion \cite{xia2017evolutionary}.

\begin{figure}[t!]
	\centering
	\vspace{-0.1cm}
	\subfloat{\includegraphics[clip,trim=2cm 8.7cm 2cm 9cm, width=8.2cm]{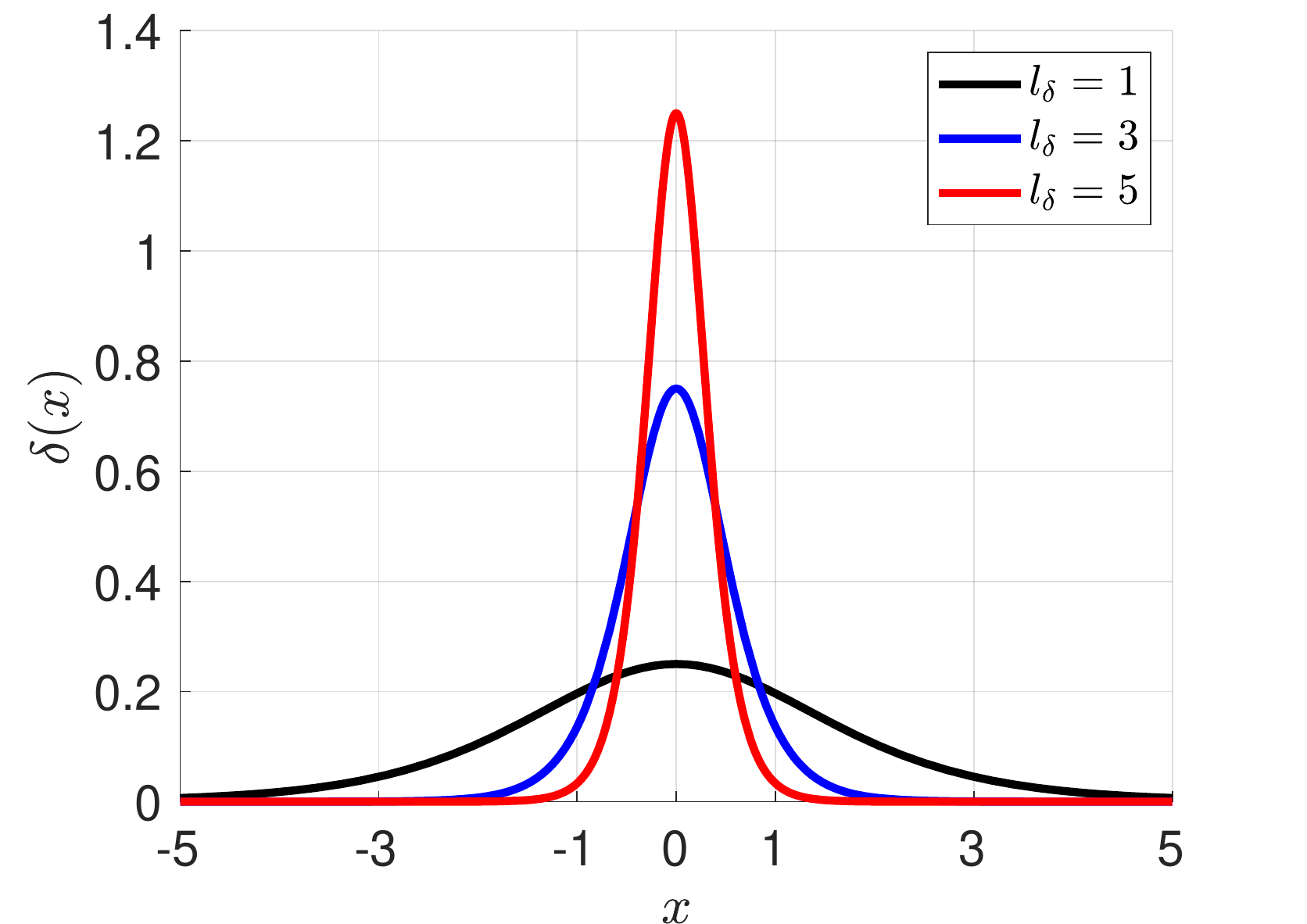}}   	\subfloat{\includegraphics[clip,trim=2cm 8.7cm 2cm 9cm, width=8.2cm]{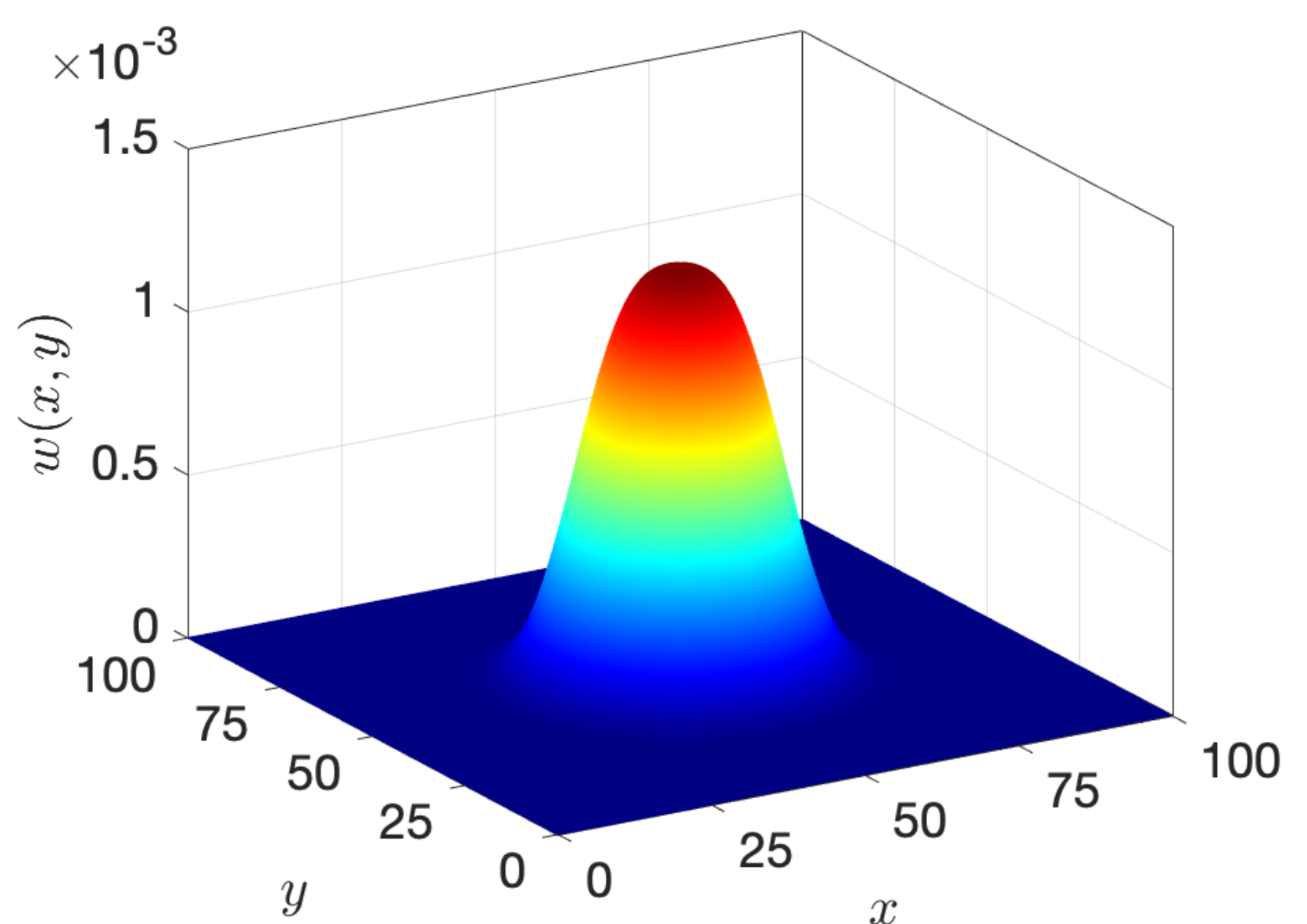}}		
		\vspace{-0.1cm}
	\caption*{\hspace*{1.3cm}(a)\hspace*{8cm}(b)}
	\caption{Schematic description of the (a) approximation of Dirac delta function given in \req{eq:dirac_approx}, and (b) the weight function used in sensitivity filtering given in \req{eq59}.}
	\label{Figure_dirac_filter}
\end{figure}

In order to eliminate the implicit derivatives for both deformation and crack phase-field $ \frac{\partial \Delta \Bu}{\partial \Phi }\ $and $ \frac{\partial \Delta d}{\partial \Phi}\ $ given in \req{eq45} (the last four lines), the terms between square brackets must be set as zero, to derive adjoint sensitivity equations. This leads to solving the following coupled system of equations for the adjoint vectors $ \bm{\mu}\ $and $ \bm{\lambda }\ $ by:
\begin{equation}\label{eq:sol_adjoint}
	[\textbf{K}^{n-1}]^T{\bm\mu}^n
	={\bm 0} \WITH \texttt{C}^n_1={\bm 0}
	\AND
	[\textbf{K}^{n}]^T{\bm\lambda}^n ={\bm 0} \WITH \texttt{C}^n_2={\bm 0}
\end{equation}
such that the coupled tangent stiffness matrix due to ($\Bu,d$) reads:
\begin{equation}\label{eq11}
	{{\bf{K}}^{n}}=\begin{bmatrix}
		{{\bf{K}}^{n}_{\Bu \Bu}} & {{\bf{K}}^{n}_{\Bu,d}} \\[5pt]
		{{\bf{K}}^{n}_{d, \Bu}} & {{\bf{K}}^{n}_{d,d}}
	\end{bmatrix},
	\AND
	{{\bf{K}}^{n-1}}=\begin{bmatrix}
		{{\bf{K}}^{n-1}_{\Bu \Bu}} & {{\bf{K}}^{n-1}_{\Bu,d}} \\[5pt]
		{{\bf{K}}^{n-1}_{d, \Bu}} & {{\bf{K}}^{n-1}_{d,d}}
	\end{bmatrix},
\end{equation}
with solution fields for Lagrange multipliers
\begin{equation}
	{\bm\mu}^n=\begin{bmatrix}
		{\bm\mu^{n}_{\Bu}}  \\[5pt]
		{\bm\mu^{n}_{d}}
	\end{bmatrix} 
   \AND
   {\bm\lambda}^n=
    \begin{bmatrix}
		{\bm\lambda^{n}_{\Bu}}  \\[5pt]
		{\bm\lambda^{n}_{d}}
	\end{bmatrix} 
\end{equation}
It is worth noting that since we are using monothonic loading thus we \grm{have}:
\begin{equation}
	[\Delta \Bu^{n-1}]^a=[\Delta \Bu^{n}]^a \qquad \rightarrow \qquad {\bm\mu}^n={\bm\lambda}^{n-1}
\end{equation}
Finally, the sensitivity analysis given in \req{eq:sens_G} by using adjoint variables derived in \req{eq:sol_adjoint} at time $n$ for Formulation 2 reads:
\begin{equation}\label{eq:final_sestivity}
	\begin{aligned}
		\mathcal{G}^n(\Phi,{\bm\lambda_{\Bu}},&{\bm\lambda_{d}},{\lambda_{V}},{\bm\mu_{\Bu}},{\bm\mu_{d}};{\Bu},{d},{\alpha})=\\[5pt]
		&\underbrace{	\displaystyle\sum_{\text{n}=1}^{N} \Big(({\bm\mu^{n}_{\Bu}})^T\frac{d \mathbf{R}_\Bu^{n-1}}{d \Phi}
		+({\bm\lambda^{n}_{\Bu}})^T\frac{d \mathbf{R}_\Bu^{n}}{d \Phi}
		+({\bm\mu^{n}_{d}})^T\frac{d \mathbf{R}_d^{n-1}}{d \Phi}
		+({\bm\lambda^{n}_{d}})^T\frac{d \mathbf{R}_d^{n}}{d \Phi}\Big)}_{\mathcal{G}_{S}^n}
		+\underbrace{{\lambda^{T}_{V}}\int_{\calB}{\delta \left( \Phi \right)\mathrm{d}{\bm{x}}}}_{\mathcal{G}_{V}^n}\;,
	\end{aligned}
\end{equation}
whereas the sensitivity analysis $\mathcal{G}$ given in \req{eq:final_sestivity} is additively decomposed into the solid and volume counterpart shown as $\mathcal{G}_S$, and $\mathcal{G}_V$, respectively. We further require this decomposition in Section \ref{Section3:sec_filitering_G_s} for the filtering method.
Additionally, Dirac delta function given in  \req{eq:final_sestivity}  is approximated as a regularized form through:
\begin{equation}\label{eq:dirac_approx}
	\delta\Big(\Phi(\Bx)\Big)\approx \frac{l_{\delta} e^{-l_{\delta}\Phi(\Bx)}}{(1+e^{-l_{\delta}\Phi(\Bx)})^2} 
	\WITH 
	\delta\Big(\Phi(\Bx)\Big)\approx\frac{d \text{H}\Big(\Phi(\Bx)\Big)}{d \Phi(\Bx)}\;.
\end{equation} 
Here, the regularization length-scale parameter for the Dirac delta function is set as $l_{\delta}=5$, see Figure \ref{Figure_dirac_filter}(a).

Let us now reduce the sensitivity analysis of \req{eq:final_sestivity} for Formulation 2 toward Formulation 1. With the same procedure, following \req{eq:sol_adjoint} the adjoint sensitivity equation for Formulation 2 reads:
\begin{equation}\label{eq:sol_adjoint2}
	[\textbf{K}_{\Bu,\Bu}^{n-1}]^T{\bm\mu}^n
	={\bm 0} \WITH \texttt{C}^n_1={\bm 0}
	\AND
	[\textbf{K}_{\Bu,\Bu}^{n}]^T{\bm\lambda}^n ={\bm 0} \WITH \texttt{C}^n_2={\bm 0}
\end{equation}
So, only deformation contribution exists in adjoint equations. Thus, the total sensitivity analysis for Formulation 1 at time $n$ reads:
\begin{equation}\label{eq:final_sestivity2}
	\begin{aligned}
		\mathcal{G}^n(\Phi,{\bm\lambda_{\Bu}},{\bm\lambda_{d}},{\lambda_{V}},{\bm\mu_{\Bu}},{\bm\mu_{d}};{\Bu},{d},{\alpha})=
		\underbrace{	\displaystyle\sum_{\text{n}=1}^{N} \Big(({\bm\mu^{n}_{\Bu}})^T\frac{d \mathbf{R}_\Bu^{n-1}}{d \Phi}
			+({\bm\lambda^{n}_{\Bu}})^T\frac{d \mathbf{R}_\Bu^{n}}{d \Phi} \Big)}_{\mathcal{G}_{S}^n}
		+\underbrace{{\lambda^{T}_{V}}\int_{\calB}{\delta \left( \Phi \right)\mathrm{d}{\bm{x}}}}_{\mathcal{G}_{V}^n}\;,
	\end{aligned}
\end{equation}
Thus, there exists no contribution for crack-phase-field in the total adjoint sensitivity analysis given in \req{eq:final_sestivity2}. Evidently, we required fewer adjoint equations for \req{eq:final_sestivity2} \grm{which} have to be resolved compared to the \req{eq:final_sestivity}.
It is now required that by using the sensitivity analysis derived in \req{eq:final_sestivity} to find a relationship with normal velocity field ${\widehat{v}_{\Phi}}$, and thus solve Formulation \req{form_4}. At this point, it is necessary to remark that the level-set velocity field which is required for closed front propagation is obtained from the steepest descent method \cite{allaire2004structural}. For this reason, the derivative of the Lagrangian function $\BfrakL$ in \req{LGforopt} with respect to pseudo time gives:
\begin{equation}\label{eq:steep_decent}
	\begin{aligned}
		\frac{\partial \mathcal{L}}{\partial t}=\frac{\partial \mathcal{L}}{\partial \Phi }.\frac{\partial \Phi }{\partial t}=\mathcal{G}{\widehat{v}_{\Phi}}\left| \nabla \Phi  \right|
		\WITH \frac{\partial \mathcal{L}}{\partial t}<0,\\
	\end{aligned}
\end{equation}
here \req{eq33} is used. An important observation is that if the negative quantity for the time derivative of the Lagrangian function $\mathcal{L}$ to be held, it is sufficient that the normal velocity field of the \req{eq:steep_decent}, defines through:
\begin{equation}\label{eq:veloc_phi}
	\begin{aligned}
			\fterm{ 
		{\widehat{v}_{\Phi}}\left| _{\Phi = 0} \right.=-\mathcal{G}\;,
	}
	\end{aligned}
\end{equation}
such that following properties holds:
\begin{equation}
	\begin{aligned}
		\frac{\partial \mathcal{L}}{\partial t}=\int_{\calB}{\mathcal{G}\delta \left( \Phi  \right){{\widehat{v}_{\Phi}}}\left| \nabla \Phi  \right|\mathrm{d}{\bm{x}} }=\int_{\partial\calB}{\mathcal{G}{{\widehat{v}_{\Phi}}}\mathrm{d}{\bm{a}}}\;, \\
	\end{aligned}
\end{equation}
here, the identity \req{eq:geo_transfer} is used, see for a detailed discussion \cite{allaire2004structural}.

By extracting the normal velocity field ${\widehat{v}_{\Phi}}$ through the sensitivity analysis in \req{eq:veloc_phi}, we are now able to solve Formulation \req{form_4}, and thus update the topology of structure. 

\begin{Remark}
	\label{rem:augumentedLG}
	We note that in this work, \grm{Lagrangian} method is used, while one can simply extend the method toward augumented \grm{Lagrangian} method. To do so, modified Lagrange variable is introduced as ${\widehat{\bm\lambda}^{T}_{\Bu}}={\bm\lambda^{T}_{\Bu}}+\frac{p_u}{2}(\mathbf{R}^{n-1}_\Bu
	+\mathbf{R}^{n}_\Bu)$ for mechanical part, ${\widehat{\bm\lambda}^{T}_{d}}={\bm\lambda^{T}_{d}}+\frac{p_d}{2}( \mathbf{R}^{n}_\text{d})$ for phase-field part, and ${\widehat{\lambda}^{T}_{V}}={\lambda^{T}_{V}}+\frac{p_V}{2}( G(V))$. So, \req{LGforopt} will be changed toward $\BfrakL(\Phi,{\widehat{\bm\lambda}_{\Bu}},{\widehat{\bm\lambda}_{d}},{\widehat{\lambda}_{V}};\bar{\Bu},\bar{d},\bar{\alpha})$ for set of unknowns $(\Phi,{\bm\lambda_{\Bu}},{\bm\lambda_{d}},{\lambda_{V}})$.
\end{Remark}

\noii{\begin{Remark}
	It is of great importance to mention that the sensitivity procedure is performed in Section 3.4 is based on displacement increments $\Delta \Bu$ which aligns with \cite{fritzen2016topology,huang2008topology,xia2017evolutionary,gangwar2022thermodynamically}. But, alternative approach  to derive adjoint equations is by performing sensitivity analysis based on total $\Bu$ for a path-dependent problem which discussed in \cite{michaleris1994tangent,alberdi2018unified,Russ,russ2021novel}. 
\end{Remark}}

\sectpb[Section3: LG_vol]{Lagrange multiplier due to volume constraint}
A gradually convergence bi-sectioning algorithm is used to satisfy the volume constraint in \req{eq:vol_const}. To do so, Lagrange multiplier  $\lambda_{V}$ has to be used to enforce the constraint inequality into the Lagrangian function in \req{LGforopt}. Since, we are dealing with a rate-dependent nonlinear boundary value problem (i.e., ductile phase-field fracture), thus, an incremental topology optimization approach is used to reach the target volume $ \bar{\Omega}$, as the following set:
\begin{equation}\label{eq:incr_volume}
	\calB = \bar{\Omega}^0<  \ldots < \bar{\Omega}^m < \bar{\Omega}^{m+1} < \ldots < \bar{\Omega}^m= \bar{\Omega},
\end{equation}
in which it is so-called \textit{expected volume} $\bar{\Omega}^m$ at time $m$, and it is defined through:
\begin{equation}\label{eq:specifed_vol}
	\bar{\Omega}^m=\chi^{m-1}_v(\Phi)-\Theta_v( \chi^{m-1}_v(\Phi)-\bar{\Omega})
	\WITH
	\chi^m_v(\Phi)=\frac{{\text{V}}(\Omega^m)}{V(\calB)}=\frac{\widehat{\text{V}}(\Phi^m)}{V(\calB)}
\end{equation}
such that $\chi^m_v(\Phi)$ is volume ratio at time $m$. Here, we set $\Theta_v=0.05$
. Then by means of \req{eq:incr_volume}, where we expect  the topology iteration $m$ to approach the volume $\bar{\Omega}^m$. By doing so, the non-linearity of the problem is reduced, and one avoids approaching the target value instantly in one step. We now elaborate on the bi-sectioning method to determine $\lambda_{V}$. First, we set a lower- $ {{\lambda }_{l}} $ and an upper-bound $ {{\lambda }_{u}} $ for the Lagrange multiplier \cite{Sigmond1}. Then, we update $\lambda_{V}$ as follows:
\begin{equation}\label{eq56}
	\begin{aligned}
		\lambda_V =\sqrt{{{\lambda }_{l}}{{\lambda }_{u}}}. \\
	\end{aligned}
\end{equation}
As a result, by having $\lambda_V$ we could determine the volume sensitivity and together with the solid sensitivity $\widetilde{\mathcal{G}}_{S}$ it results in $\mathcal{G}$. Accordingly, we can compute normal velocity field ${\widehat{v}_{\Phi}}$ from \req{eq:veloc_phi}, and thus solving the reaction-diffusion equation for $\Phi(\Bx,t)$. By having the volume ratio $\chi_v(\Phi)$ for the new topology, one can assess if $\chi_v(\Phi)$ is greater than the expected volume $\bar{\Omega}$, then $\lambda_l$ will be updated based on $\lambda_V$, otherwise  $\lambda_u$ will be replaced by $ \lambda_V$. The detailed bi-sectioning algorithm  to determine the Lagrange multiplier  $\lambda_{V}$ due to volume constraint in the optimization process is depicted in Algorithm \ref{Table_bisection}.

\newpage	

\begin{algorithm}
	\small
	\caption{\em Bi-sectioning algorithm for Lagrange multiplier due to volume constraint.}
	\label{Table_bisection}
	{
		\begin{tabular}{l}
			{\bf Input:} topology optimization data $(\Phi^{m-1}, \lambda^{m-1}_{V},\bar{\Omega})$  from step $m-1$, and $\widehat{\mathcal{G}}^{\;m}_{S}$ from step $m$. \\[2mm]
			
			\underline{Initialization}:\\[0.3cm]
			
			 \hspace*{0.75cm}\textbullet\;we set: $\bar{\Omega}^m=\chi^{m-1}_v(\Phi)-\Theta_v(\chi^{m-1}_v(\Phi)-\bar{\Omega})$\;, \\[0.3cm]
			 
			 \hspace*{0.75cm}\textbullet\;we set for topological field as $\Phi^{m-1,1}=\Phi^{m-1}$\;,\\[0.3cm]
			
			 \hspace*{0.75cm}\textbullet\;choose lower-, and upper-bound for $\lambda_{V}$, we set $\lambda_{l}=10^{-8}$, and $\lambda_{u}=\max(10^{8},10^{4}\lambda^{m-1}_{V})$\;,\\[0.3cm]
			
			{{\bf Bi-sectioning iteration} $k\geq1$}:\\[0.3cm]
			
		     \hspace*{0.75cm}\textbullet\;we set $\lambda_V =\sqrt{{{\lambda }_{l}}{{\lambda }_{u}}}$\\[0.3cm]
		     
		     \hspace*{0.75cm}\textbullet\;given $(\lambda_V,\Phi^{m-1,k})$ compute ${\mathcal{G}}^{m}_{V}$ from \req{eq:final_sestivity}-\req{eq:dirac_approx} together with the given $\widehat{\mathcal{G}}^{\;m}_{S}$ results in $\mathcal{G}^m$,\\[0.3cm]		
			
			 \hspace*{0.75cm}\textbullet\;given ${\mathcal{G}}^{m}$, find the normal velocity field :${\widehat{v}_{\Phi}}\left| _{\Phi = 0} \right.=-\mathcal{G}^m(\Phi,{\bm\lambda_{\Bu}},{\bm\lambda_{d}},{\lambda_{V}},{\bm\mu_{\Bu}},{\bm\mu_{d}};{\Bu},{d},{\alpha})$,\\[0.3cm]			
			
			 \hspace*{0.75cm}\textbullet\;${\text{normalizing part:}}$ given ${\widehat{v}_{\Phi}}| _{\Phi = 0}$, find normalize velocity field ${{v}_{\Phi}}| _{\Phi = 0}=\frac{{\widehat{v}_{\Phi}}| _{\Phi = 0}}{\texttt{mean}({{\widehat{v}_{\Phi}}| _{\Phi = 0}})}$, \\[0.3cm]
			
			 \hspace*{0.75cm}\textbullet\;${\text{topological part:}}$ given ${{v}_{\Phi}}| _{\Phi = 0}$, solve ${{\bf{R}}_\Phi}(\Bu,d,\Phi)=\bm{0}$ for $\Phi$, set $\Phi=:\Phi^{m-1,k}$, \\[0.3cm]
			
			 \hspace*{0.6cm} \textbullet\;given $\Phi^{m-1,k}$, find volume of material domain $\Omega^k(\Bx,t)$, with ${\text{V}}(\Omega^k)=\widehat{\text{V}}(\Omega^k)=\displaystyle\int_{\calB} \text{H}(\Phi^k)\text{d}\Bx$\;,\\[0.3cm]
						
			 \hspace*{0.75cm}\textbullet\; given $\widehat{\text{V}}(\Omega^k)$, find volume ratio by $\chi^k_v=\frac{\widehat{\text{V}}(\Omega^k)}{V(\calB)}$ with ${\text{V}}(\calB)=\displaystyle\int_{\calB} 1\text{d}\Bx$\;,\\[0.3cm]
			
			 \hspace*{0.75cm}\textbullet\; if $\chi^k_v\geq\bar{\Omega}^m$ fulfilled, set  $\lambda_l=\lambda^k_V$ otherwise  $\lambda_u=\lambda^k_V$\;,\\[0.3cm]
			
			 \hspace*{0.75cm}\textbullet\; check convergence of bi-sectioning algorithm:\\[0.3cm]
			
			\hspace{2cm}$\mathrm{Res}_\mathrm{V}^k:=\frac{|\lambda^k_V-\lambda^{k-1}_V|}{|\lambda^k_V+\lambda^{k-1}_V|}$ \\[0.3cm]
			
			 \hspace*{0.75cm}\textbullet\; if fulfilled, set $(\Phi^{m-1,k}, \lambda^{k}_{V})=:(\Phi^{m}, \lambda^{m}_{V})$ then stop; \\[0.2cm]
			
			\hspace{1.9cm}\; else $k+1\rightarrow k$. \\[0.5cm]
						
			{\bf Output:} solution $(\Phi^{m}, \lambda^{m}_{V})$.
		\end{tabular}
	}
\end{algorithm}

We note that in Algorithm \ref{Table_bisection}, $\texttt{mean}(\bullet)$ denotes a standard mean value function. Let us write
Algorithm \ref{Table_bisection} in the following abstract form:
\begin{equation}
	\Bs^m=\texttt{LSM}(\Bs^{m-1})
\end{equation}
with $\Bs=(\Phi,\lambda _{V})$  at time $t^{m}$. We will use this notation in our final proposed topology optimization algorithm.
%
\sectpb[Section3:sec_filitering_G_s]{Filtering method due to numerical instabilities}
Since the point-wise exact Heaviside step function is employed for geometry mapping when a hole is nucleated in an appropriate region, its sensitivity is always zero and no material can be added to that domain which is referred to as evolutionary locking, see \cite{huang2010evolutionary} Section 3.3.3. Therefore, a sensitivity filter scheme is used to provide the potential for material creation in the void regions. To the best of the authors' knowledge, the filter schemes are used to overcome the numerical instabilities of topology optimization such as intermediate densities, mesh dependency, checkerboard patterns, and local minima that occur in element-wise constant density distribution techniques such as SIMP and ESO \cite{diaz, sigmond2, sigmund4, sigmund3}.

To apply the filtering method, recall the sensitivity analysis $\mathcal{G}$ given in \req{eq:final_sestivity} is additively decomposed through:
\begin{equation}\label{eq:addtive_split_solid_sens}
	\mathcal{G}^{\;m}=\mathcal{G}^{\;m}_{S}+\mathcal{G}^{\;m}_{V}
\end{equation}
We aim here to modify the solid counterpart of the sensitivity analysis $\mathcal{G}^m_S$. To do so, following \cite{Sigmund5} the symmetric sensitivity filter formulation without density weighting scheme is written as:
\begin{equation}\label{eq:modify_G_S0}
	\begin{aligned}
		{\mathcal{\widetilde{G}}^{\;m}_{S}}=\frac{\displaystyle\sum\limits_{j\in {|node|}}{w\left( {\Bx_{mj}} \right){\mathcal{G}^m_{S,j}}}}{\displaystyle\sum\limits_{j\in {|node|}}{w\left( {\Bx_{mj}} \right)}}, \\
	\end{aligned}
\end{equation}
where $ |node|$ and $ w\left( {\Bx_{mj}} \right)\ $ denote the number of nodes and the weight factor, respectively.
The smooth continuous weight factor is defined as:
\begin{equation}\label{eq59}
	\begin{aligned}
		w\left( {\Bx_{mj}} \right)={{e}^{-3{{(\frac{\left\| {\Bx_{mj}}-{\Bx_{nm}} \right\|}{{r_{\min }}})}^{3}}}}, \\
	\end{aligned}
\end{equation}
where $ {\Bx_{nm}}\ $and $ {r_{\min }}\ $ denote the spatial coordinates of $m_{th}$ nodal point and filter radius, respectively. Figure \ref{Figure_dirac_filter}(b) illustrates the nonlinear filter factor in a 2D mesh using \req{eq59}, schematically.
The result obtained from the sensitivity filter is still a chaotic phenomenon and the optimization process may fail to converge. To overcome this obstacle, sensitivities are averaged with their historical information, which is written as:
\begin{equation}\label{eq:modify_G_S}
	\begin{aligned}
		\widehat{\mathcal{G}}^{\;m}_{S}=\frac{\widetilde{\mathcal{G}}^{\;m}_{S}+\widehat{\mathcal{G}}_{S}^{\;m-1}+\widehat{\mathcal{G}}_{S,n}^{\;m-2}}{3}\,\,\,\,\,with\,\,\,\,m>2.\
	\end{aligned}
\end{equation}
Here, $1\leq m\leq M$ indicates the iteration of the optimization process, see \req{eq:time_int_phi}.
Thus, by means of  \req{eq:modify_G_S0}, and \req{eq:modify_G_S}, the sensitivity analysis given in \req{eq:addtive_split_solid_sens} is finally modified through:
\begin{equation}
	\mathcal{G}^m=\widehat{\mathcal{G}}_{S}^m+\mathcal{G}_{V}^m.
\end{equation}
Thus,  $ \widehat{\mathcal{G}}^m_{S}$ is employed as the normal velocity field to update the topological field  through Formulation \ref{form_4}. The detailed topology optimization for fracture resistance of ductile material to determine topological field $\Phi(\Bx,t)$  is depicted in Algorithm \ref{TableGL}.

\begin{Remark}
	\label{rem:solwrate_top}
    It is worth noting that for the highly nonlinear problem, like in our case, the slow rate of convergence for topology optimization problems needs to be taken into account. To do so, different ways can be highlighted. ($i$) If we set $\Theta_V$ in \req{eq:specifed_vol} as a small quantity thereafter the speed rate of volume of $\Omega$ also will be reduced. ($ii$) If we set $\tau_\Phi$ in (3.1) as small value to maintain the stability of the solution thereafter the speed rate $\Phi$ will be reduced (due to a necessary CFL condition, see e.g. \cite{wang2003level})
\end{Remark}

\begin{algorithm}
	\small
	\caption{\em Topology optimization for fracture resistance of ductile material.}
	\label{TableGL}
	{
		\begin{tabular}{l}
		{{\bf Topology optimization iteration} $1\leq m\leq T_M=T_\Phi$}.\\[0.2cm]	
		
		{\bf Input:} solution data $(\Phi^{m-1},\lambda^{m-1}_{V},\bar{\Omega})$ from step $m-1$. \\[0.4cm]
			
			{{\bf Ductile phase-field fracture iteration } $1\leq n\leq T_N=T_f$}\\[0.2cm]
			
			 {\bf Input:} loading data $(\bar{\bm u}_{n},\bar{\bm \tau}_{n})$ on $\partial_D \calB$, respectively; \\[2mm]
			\hspace{1.3cm} solution $(\Bu^{n-1},\Bve^{n-1}_{p} ,d^{n-1} ,\alpha^{n-1})$ from time step $n-1$ with fixed $\Phi^{m-1}$. \\[2mm]
			
			\quad\quad \underline{Staggered iterative solution for ductile fracture}:\\[0.5cm]
			
			\quad\quad\quad\;\; \textbullet\;${\text{phase-field part:}}$ given $(\bm u^k; \Phi^{m-1})$, solve ${{\bf{R}}_d}(\Bu,d;\Phi)=\bm{0}$ for $d$, set $d=:d^k$,\\[0.2cm]			
			\quad\quad\quad\;\; \textbullet\;${\text{mechanical part:}}$ given $(d^k; \Phi^{m-1})$, solve ${{\bf{R}}_\Bu}(\Bu,d;\Phi)=\bm{0}$ for $\Bu$, set $\Bu=:\Bu^k$, \\[0.2cm]

			\quad\quad\quad\;\; \textbullet\; for the obtained pair $(\bm u^k,d^k;\Phi^{m-1})$, check {staggered residual} by \\[0.2cm]
			\hspace{2cm}$\mathrm{Res}_\mathrm{Stag}^k:=|\mathcal{E}_{\Bve}(\bm u^k,d^k,\Phi^{m-1};\delta \Bu)|
			+|\mathcal{E}_{d}(\bm u^k,d^k,\Phi^{m-1};\delta d)|
			\leq\texttt{TOL}_\mathrm{Stag}, \; \forall \; (\delta \Bu,\delta d)\in(\calW_{\bm{0}}^{\Bu},\calW_{d})$\\[0.2cm]
			
			\quad\quad \quad \;\; \textbullet\; if fulfilled, set $(\bm u^k,d^k;\Phi^{m-1})=:(\bm u^n,d^n;\Phi^{m-1})$ then stop\;,\\
			
			\hspace{1.9cm}\; else $k+1\rightarrow k$. \\[2mm]
			
			\quad\quad\quad\;\; \textbullet\;
			{\bf Output:} solution $(\Bu^{n},\Bve^{n}_{p} ,d^{n} ,\alpha^{n})$ at $n^{\text{th}}$ time-step. \\[0.5cm]
			
			\quad\quad \underline{Adjoint sensitivity method}:\\[0.5cm]	
			
			\quad\quad\quad\;\; \textbullet\; given $(\Bu^{n},\Bve^{n}_{p} ,d^{n} ,\alpha^{n})$; solve$^\star$ \\[0.3cm]
			
			\hspace{2cm}$[{\bf{{\widehat{K}}}}^{n}]^T{\bm\lambda}^n ={\bm 0} \WITH \texttt{C}^n_2={\bm 0} \AND {\bm\mu}^n={\bm\lambda}^{n-1}\;,$\\[0.3cm]
		
			\quad\quad\quad\;\; {\color{white}\textbullet}\; set $({\bm\mu},{\bm\lambda})=:({\bm\mu}^n,{\bm\lambda}^n)$, \\
			
			\quad\quad\quad\;\; \textbullet\; given $({\bm\mu}^n,{\bm\lambda}^n)$ update solid sensitivity analysis $\mathcal{G}_{S}^m=\displaystyle\sum_{\text{n}=1}^{N}\mathcal {G}_{S}^n$ $^{\star\star}$ ,\\[0.3cm]
			
			{\bf Output:} solid sensitivity analysis with $\mathcal{G}_{S}^m$ \\[0.3cm]
			
			\quad\quad \underline{Filtering method}:\\[0.5cm]
						
	     	\quad\quad\quad\;\; \textbullet\; predictor: compute symmetric sensitivity filter for adjoint solid counterpart $\mathcal{G}_{S}^m$ by\\[0.2cm]

				\hspace{3cm} $		{\mathcal{\widetilde{G}}_{S}^m}=\frac{\displaystyle\sum\limits_{j\in {|node|}}{w\left( {\Bx_{mj}} \right){\mathcal{G}_{S,j}}}}{\displaystyle\sum\limits_{j\in {|node|}}{w\left( {\Bx_{mj}} \right)}}$  \\[0.2cm]
			
            \quad\quad\quad\;\; \textbullet\; corrector: add filtering into structural sensitivity functional: ${\mathcal{\widetilde{G}}_{S}^{\;m}}$ through\\[0.2cm]
            
            		\hspace{3cm} $\widehat{\mathcal{G}}^{\;m}_{S}=\frac{\widetilde{\mathcal{G}}^{\;m}_{S}+\widehat{\mathcal{G}}_{S}^{\;m-1}+\widehat{\mathcal{G}}_{S}^{\;m-2}}{3}\WITH m>2.$\\[0.5cm]
			
			\quad\quad \underline{Update topological field: Level-set method}\\[0.5cm]	
			
			\quad\quad\qquad \textbullet\; given ${\mathcal{\widehat{G}}^{\;m}_{S}}$, solve 
	        $\Bs^m=\texttt{LSM}(\Bs^{m-1})$, for $\Bs=(\Phi,\lambda _{V})$ set $\Bs=:\Bs^m$ \\[0.5cm]
			
			{\bf Output:} solution $(\Phi^{m},\lambda^{m}_{V})$.\\[0.4cm]
		 
		 $^{\star}$ compute jacobian  $\bf{{\widehat{K}}}$ for \texttt{Formulation 1} by ${\bf{{\widehat{K}}}}=\textbf{K}_{\Bu,\Bu}$, and \texttt{Formulation 2} by ${\bf{{\widehat{K}}}}=\textbf{K}$ from \req{eq11}\\
		 $^{\star\star}$ compute sensitivity ${\mathcal{\widehat{G}}^{\;n}_{S}}$ for \texttt{Formulation 1} by \req{eq:final_sestivity2}, and for \texttt{Formulation 2} by \req{eq:final_sestivity}. \\  	
		\end{tabular}
	}
\end{algorithm}

\sectpa[Section5]{Numerical Examples}
\label{Section5}
This section demonstrates the performance of the proposed topology optimization due to brittle and ductile phase-field fracture.  Four boundary value problems are investigated, in which the first two are related to brittle fracture, and the last two deal with ductile phase-field fracture. More specifically, in all the numerical examples, we examine the efficiency of the optimum layout when taking into account brittle/ ductile fracture and compare it with the optimum layout obtained due to pure linear elasticity (excluding fracture), and also with the original domain. To do so, we present the quantitative and qualitative indicators to highlight the role of fracture in the topology optimization framework. These indicators are as follows: $(i)$ The load-displacement curve to examine the maximum load capacity before crack initiation. $(ii)$ Qualitative response of crack phase-field pattern to highlight the effects of damage response within every new topology. $(iii)$ Numerical indicators such as the objective function and the volume constraint which imply a stiff response of the material. 

\begin{table}[!ht]
	\small
	\caption{Material/model parameters used in the numerical examples based on \cite{ambati2015review,noii2021bayesian}}	\vspace{2mm}
	\centering
	\begin{tabular}{ccllllll}
		No.  &Parameter & Name                   & Exm. 1  & Exm. 2  & Exm. 3  & Exm. 4& Unit            \\\hline 
		1.   &$K$        & Bulk modulus      & $17.3$    & $13.46$ & $175$    & $175$ & $\mathrm{MPa}$ \\
		2.   &$\mu$      & Shear modulus      & $8$  & $10.95$  & $80.76$  & $80.76$  & $\mathrm{MPa}$                 \\
		3.   &$h$      & Hardening modulus     & --  & -- & $200$  & $200$  & $\mathrm{MPa}$                 \\
		4.   &$\sigma_Y$      & Yield stress    & --  & -- & $543$  & $243$ & $\mathrm{MPa}$                 \\
		5.   &${\psi}_c$      & Specific fracture energy    & $5\times 10^{-4}$  & $1\times 10^{-5}$   & $13$  & $2$  &  $\mathrm{MPa}$                 \\
		6.   &$\zeta$      & Scaling factor     & $1$  & $1$  & $10$  & $5$  &  --                \\
		7.   &$\eta_f$      & Fracture viscosity    & $10^{-6}$  & $10^{-6}$ & $10^{-6}$  & $10^{-6}$  &  $\mathrm{N/m^{2}s}$                \\
		8.   &$\kappa$      & Stabilization parameter   & $10^{-8}$  & $10^{-8}$ & $10^{-8}$  & $10^{-8}$  &  $\mathrm{MPa}$                 \\
		9.   &$l_f$      & Fracture length-scale   & $0.3$  & $(3.6,2.8,2.4)$  & $0.18$  & $0.04$  &    $\mathrm{mm}$               \\
		11.   &$\eta_\Phi$      & Topological viscosity    & $1$  & $1$ & $1$  & $1$  &  --                \\
		12.   &$l_\Phi$      & Topological length-scale   & $10^{-2}$  & $10^{-2}$& $10^{-2}$  & $10^{-2}$  &   $\mathrm{mm}$               \\
		13.   &$r_{min}$      & Filtering radius   & $3l_f$  & $3l_f$& $3l_f$  & $3l_f$  &   $\mathrm{mm}$               \\
		14.   &$\tau_f$      & Fracture time increment   & $10^{-4}$  & $10^{-4}$ & $10^{-4}$  & $10^{-4}$  &   $\mathrm{mm}$               \\
		15.   &$\tau_\Phi$      & Topology time increment   & $10^{-4}$  & $10^{-4}$ & $10^{-4}$  & $10^{-4}$  &   $\mathrm{mm}$               \\
		\hline
		\label{material-parameters}
	\end{tabular}
\end{table}

It is noteworthy that, \textit{non-optimized} result refers to the initial structure which has not been optimized, or alternatively, refers to the first iteration of the topology optimization scheme. Additionally, by \textit{elasticity} result, we mean that the topology optimization is performed on a linear elastic model of the given boundary value problem, i.e. no plasticity is assumed and the failure mechanism is purely brittle.

In the last example, all the optimal design options that arise from linear elasticity, brittle, and ductile fracture are examined. Thus, here we mainly aim to illustrate the efficiency of the proposed topology optimization formulation which basically demonstrates various designs to prevent crack propagation.  The material parameters used in the numerical examples are listed in Table~\ref{material-parameters}. 

At this point, it is necessary to remark that the optimization process is terminated when the advancement of the objective function for three consecutive iterations becomes less than a specified threshold while the volume constraint is satisfied during the optimization process.
Here, we set optimization threshold as $10^{-4}$, see for more detail \cite{Jahan_IGA_LSM_RDE_GN}.

\noii{For the implementation, we have used both \textsc{MATLAB R2018}b \cite{MATLAB18b} and \textsc{Fortran 90}. More precisely, user elements including the constitutive modeling at each Gaussian quadrature point are written in \textsc{Fortran 90}. The general framework for the ductile fracture approach {is} then implemented in \textsc{MATLAB} as a parent/main program such that all subprograms in \textsc{Fortran 90} {are} called as a \textsc{Mex-file}.} 

\begin{figure}[!t]
	\centering
	{\includegraphics[clip,trim=1cm 28cm 0cm 1cm, width=16cm]{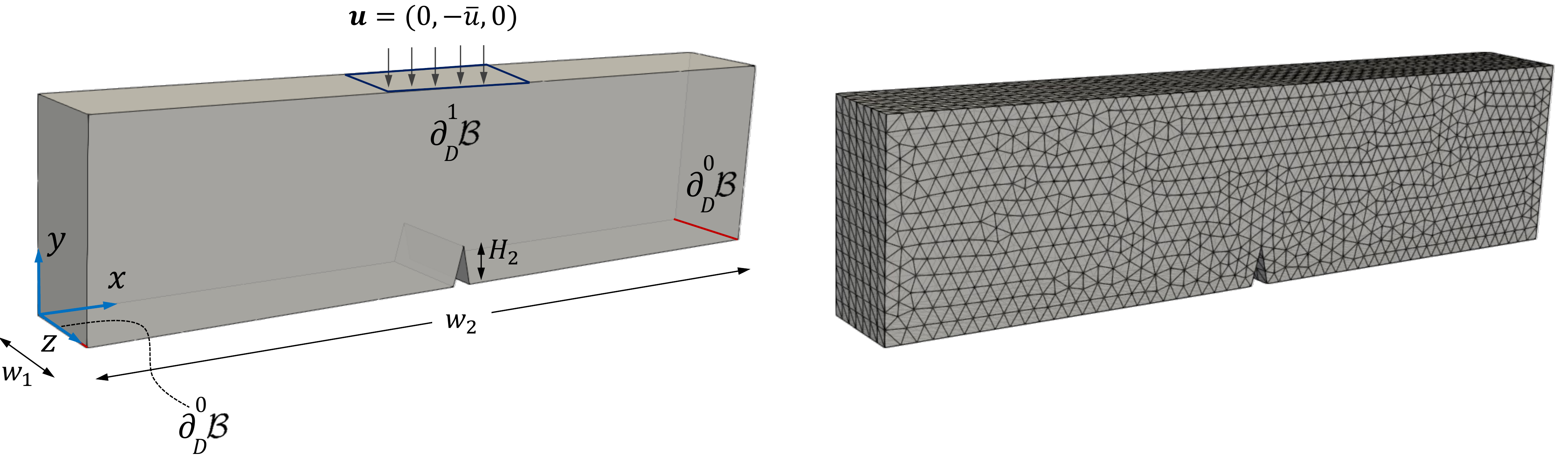}}  
	\vspace*{-0.7cm}
	\caption*{\hspace*{4.3cm}(a)\hspace*{8cm}(b)\hspace*{2cm}}
	\caption{Example 1. The representation of the (a) geometry and boundary conditions, and (b) finite element discretization.}
	\label{Exm1_bvp}
\end{figure}
%

\sectpb[Section51]{Example 1: Three-point bending test under compression loading}

To gain insight into the performance of the proposed topology optimization scheme toward failure mechanics, a three-point bending test of a solid beam under compression is concerned. In this example, we consider brittle fracture only.
The two aforementioned topology optimization procedures denoted as Formulation 1, and Formulation 2 are examined with respect to the total energy functional (stiffer structural response). For this example, the configuration is shown in Figure \ref{Exm1_bvp}(a). The left edge is constrained in all directions while the right edge is fixed for displacement in $y-z$ directions. The geometrical dimensions for Figure \ref{Exm1_bvp}(a) are set as $w_1=1\;mm$, $w_2=8\;mm$, $H_1=2\;mm$, and $H_2=0.40\;mm$.

A monotonic displacement increment  ${\Delta \bar{u}}_y=-1\times10^{-3}\;mm$ is applied in a vertical direction in a part of the top page $area=(3,2,2)\times (5,2,0)$ of the specimen for 180-time steps. Thus, the final prescribed displacement load for this optimization problem is set as ${\bar{u}}_y=-1.8\;mm$. The minimum finite element size in the solid domains is $h_{min}=0.15\;mm$, which, in turn, the heuristic requirement $h<l/2$ inside the localization zone is fulfilled. The design domain $\calB$ is discretized into 33140 linear tetrahedral elements. The material and numerical parameters are those given in Table  \ref{material-parameters}, respectively. 

\begin{figure}[t!]
	\caption*{\hspace*{2cm}\underline{$\chi_v=0.88$}\hspace*{4cm}\underline{$\chi_v=0.69$}\hspace*{3.5cm}\underline{$\chi_v=0.40$}\hspace*{1.5cm}}
	\vspace{-0.1cm}
	\subfloat{\includegraphics[clip,trim=5cm 11.5cm 6cm 15cm, width=5.5cm]{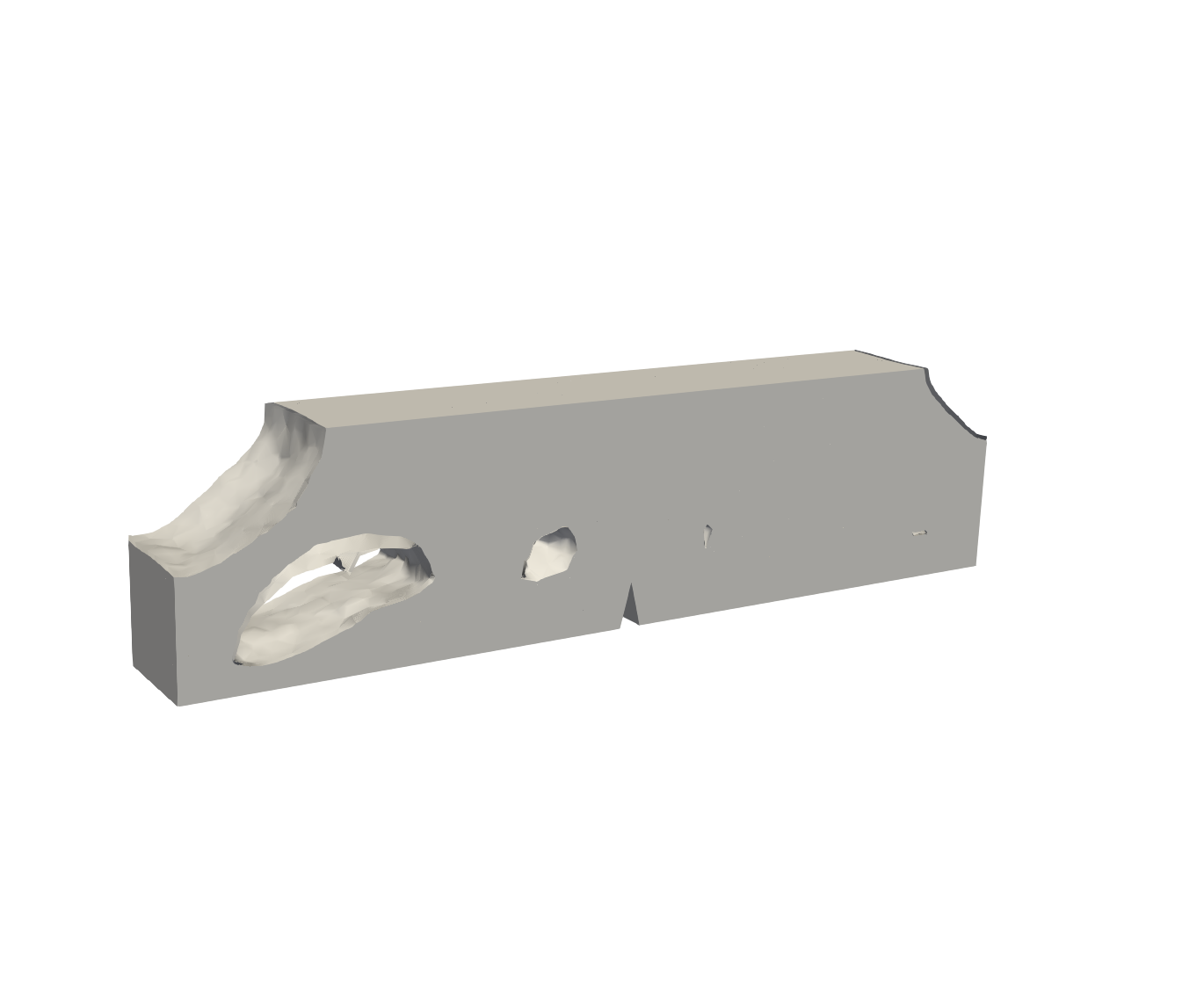}}   		\subfloat{\includegraphics[clip,trim=5cm 11.5cm 6cm 15cm, width=5.5cm]{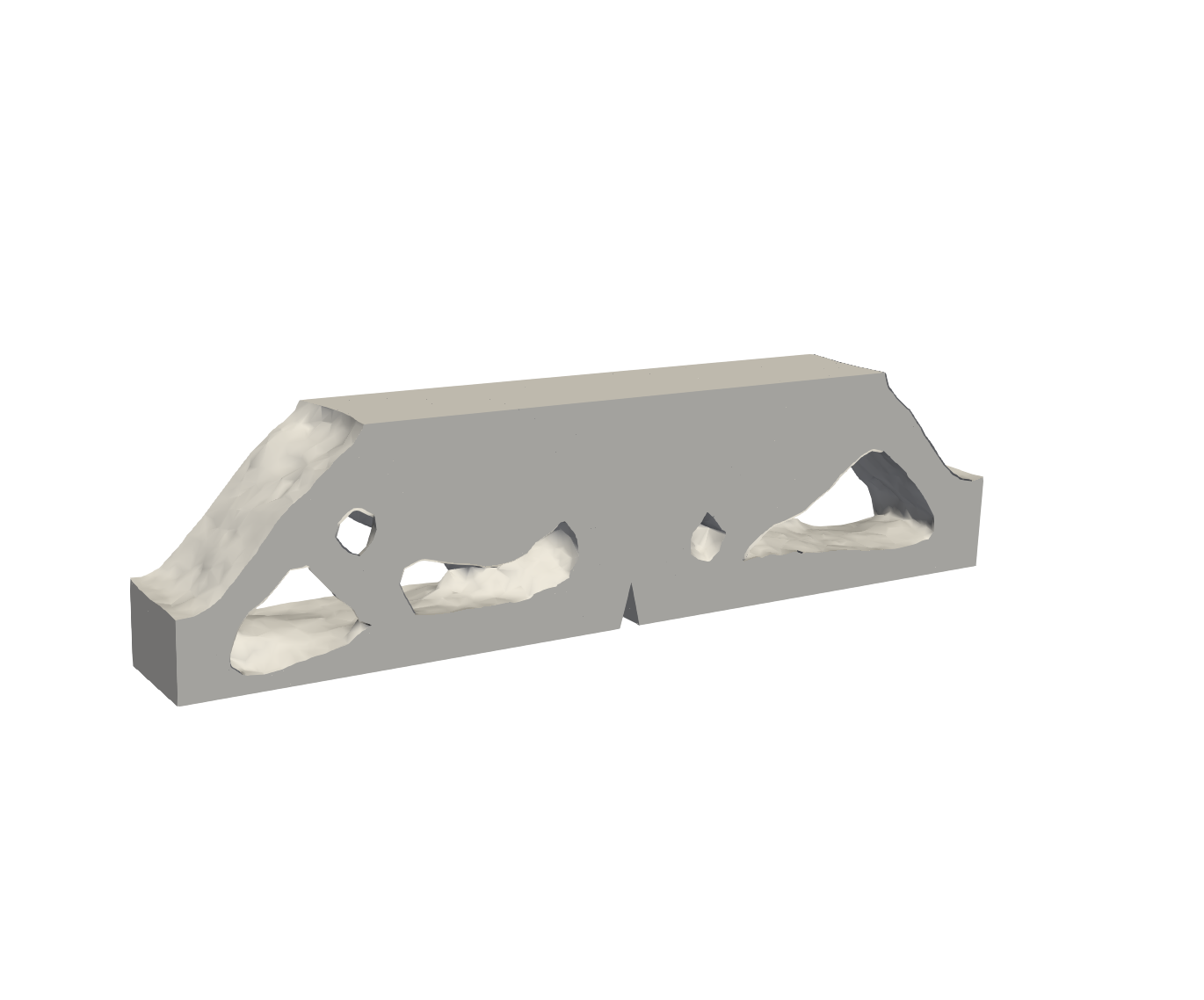}}   			\subfloat{\includegraphics[clip,trim=5cm 11.5cm 6cm 15cm, width=5.5cm]{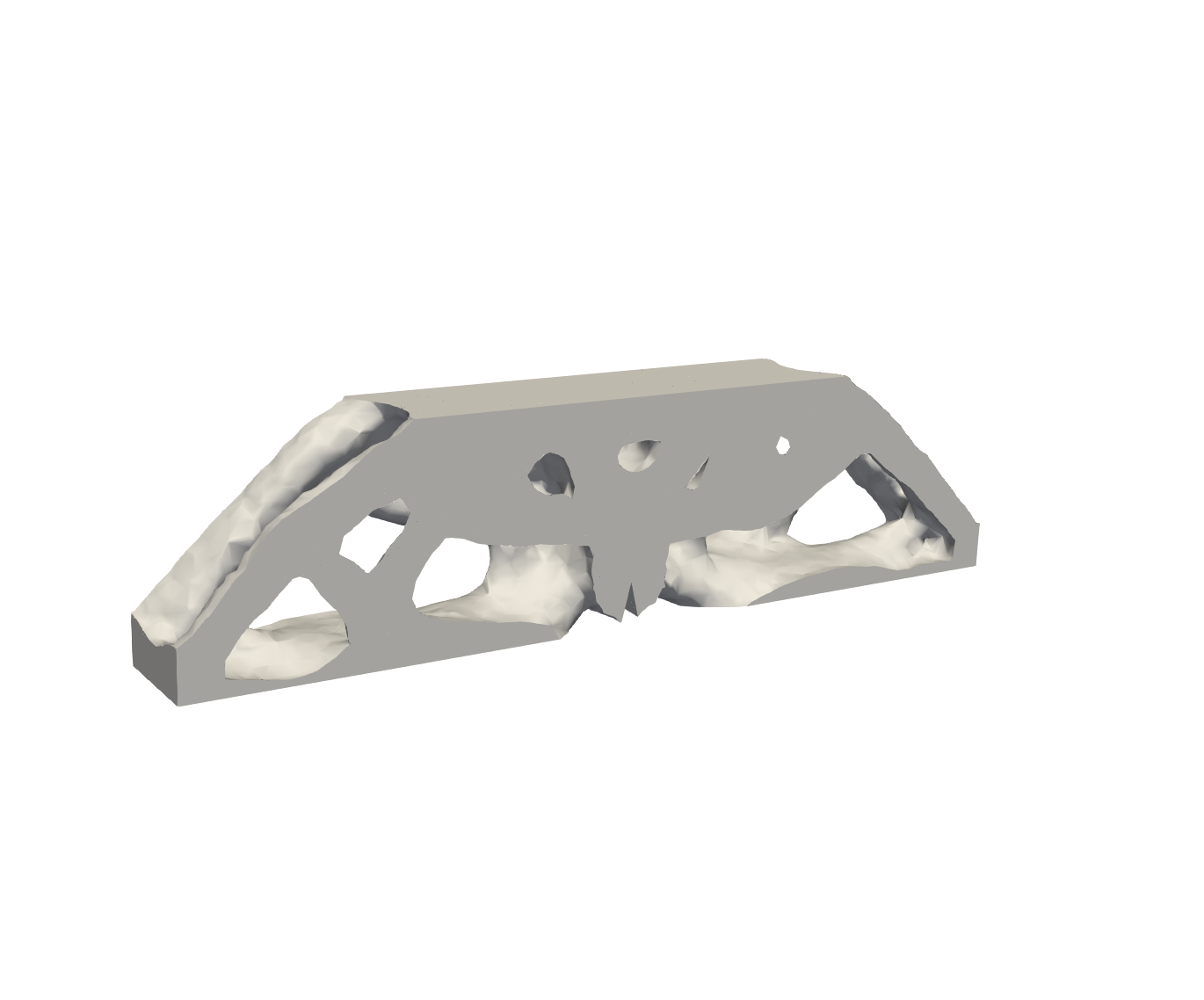}}
	\caption{Example 1. Evolution history of the optimal layouts  for different volume ratio of  the three-point bending test based on Formulation 1.}
	\label{Exm1_C2}
\end{figure}

\begin{figure}[t!]
	\caption*{\hspace*{2cm}\underline{$\chi_v=0.63$}\hspace*{4cm}\underline{$\chi_v=0.61$}\hspace*{3.5cm}\underline{$\chi_v=0.58$}\hspace*{1.5cm}}
	\vspace{-0.1cm}
	\subfloat{\includegraphics[clip,trim=5cm 11.5cm 6cm 15cm, width=5.5cm]{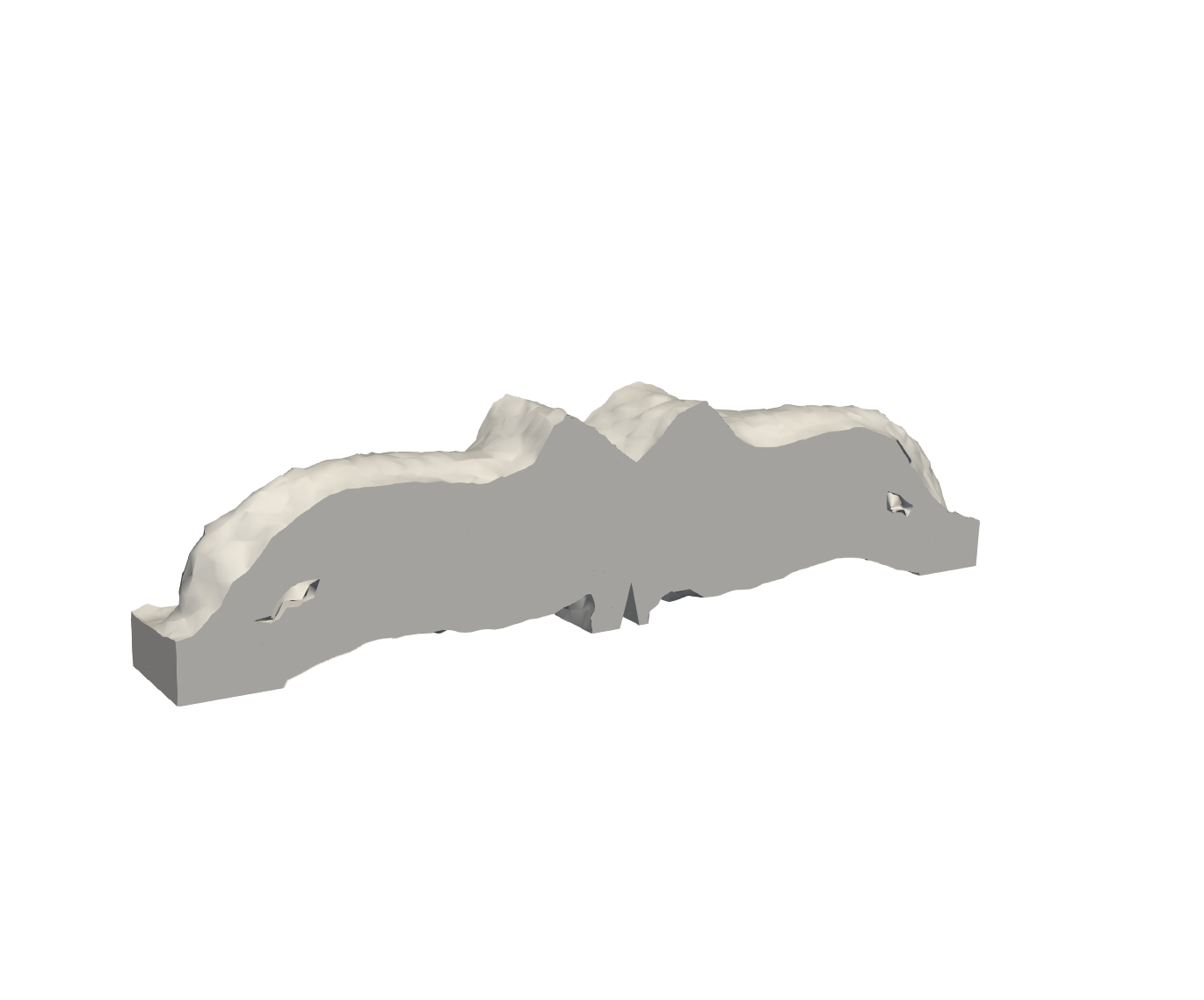}}   		\subfloat{\includegraphics[clip,trim=5cm 11.5cm 6cm 15cm, width=5.5cm]{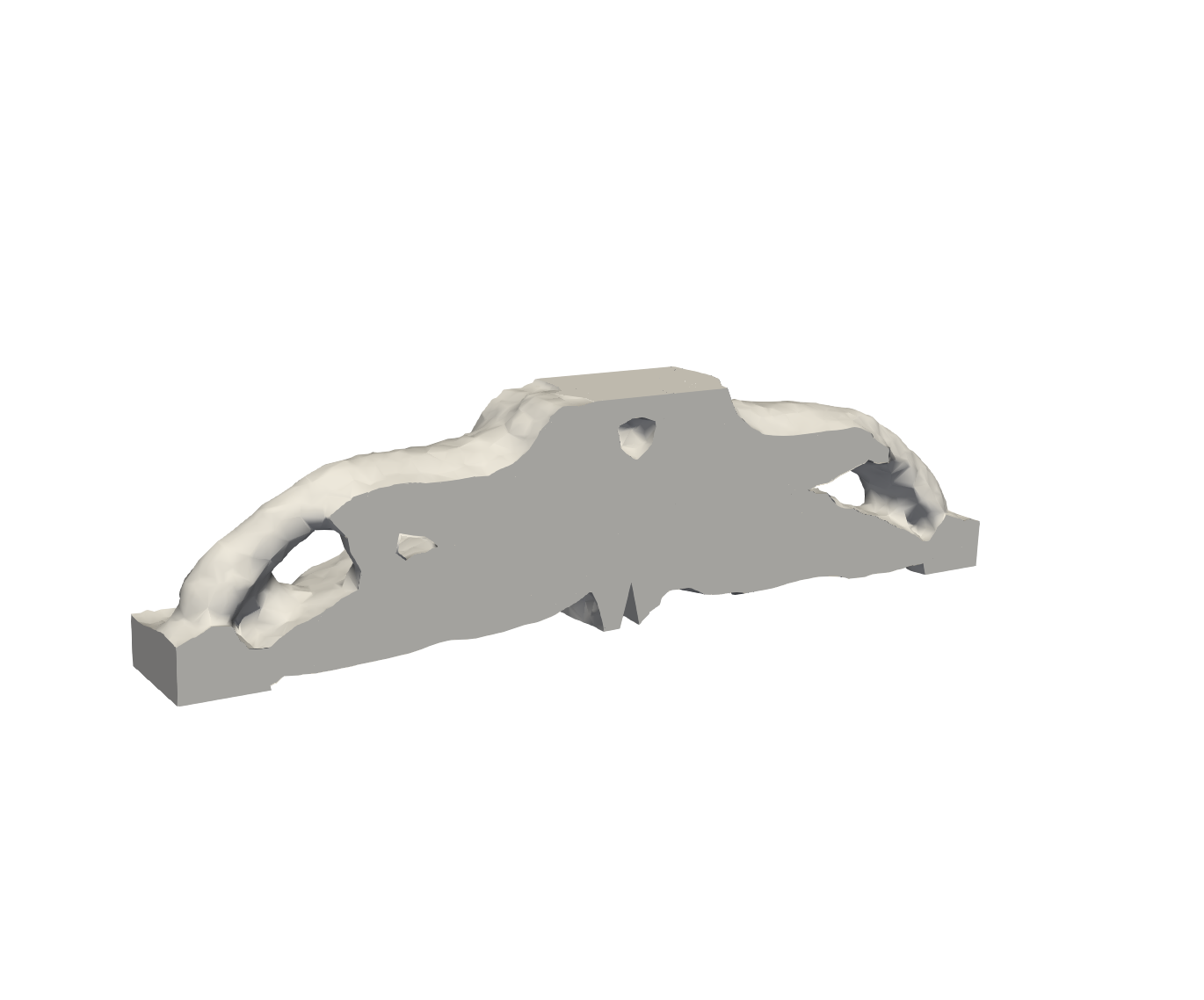}}   			\subfloat{\includegraphics[clip,trim=5cm 11.5cm 6cm 15cm, width=5.5cm]{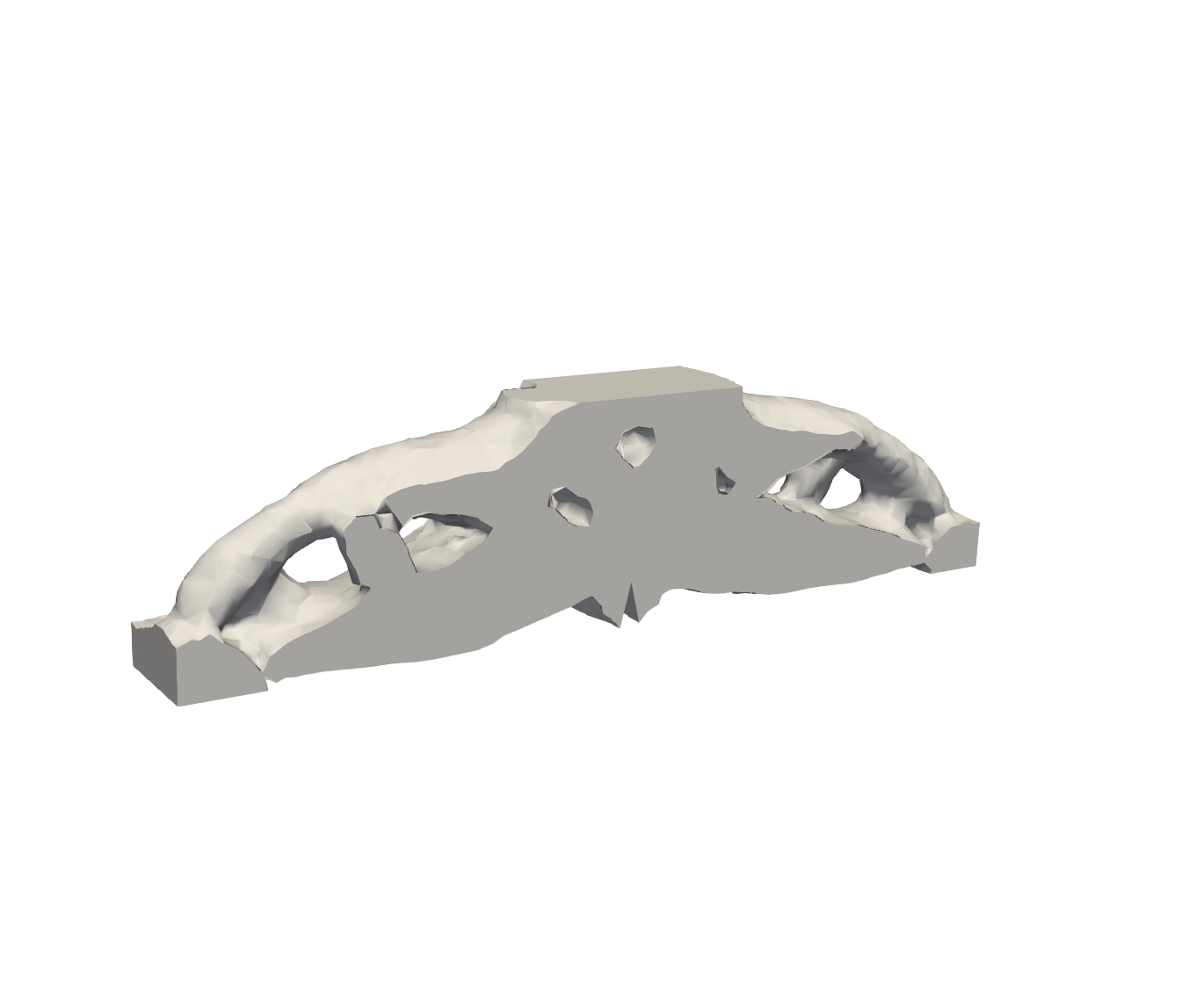}}
	\caption*{\hspace*{2cm}\underline{$\chi_v=0.55$}\hspace*{4cm}\underline{$\chi_v=0.50$}\hspace*{3.5cm}\underline{$\chi_v=0.40$}\hspace*{1.5cm}}
	\vspace{-0.1cm}
	\subfloat{\includegraphics[clip,trim=5cm 11.5cm 6cm 15cm, width=5.5cm]{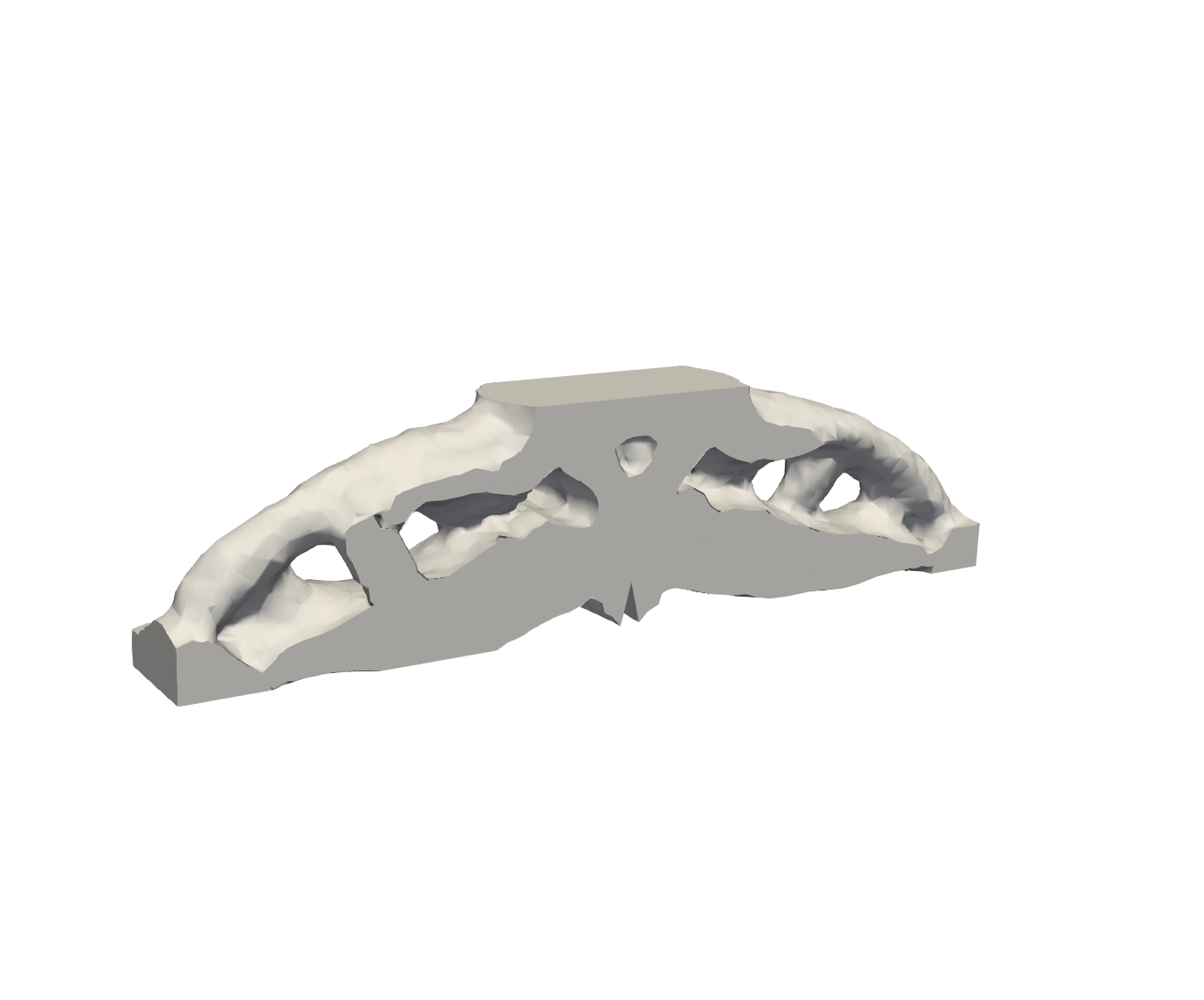}}   		\subfloat{\includegraphics[clip,trim=5cm 11.5cm 6cm 15cm, width=5.5cm]{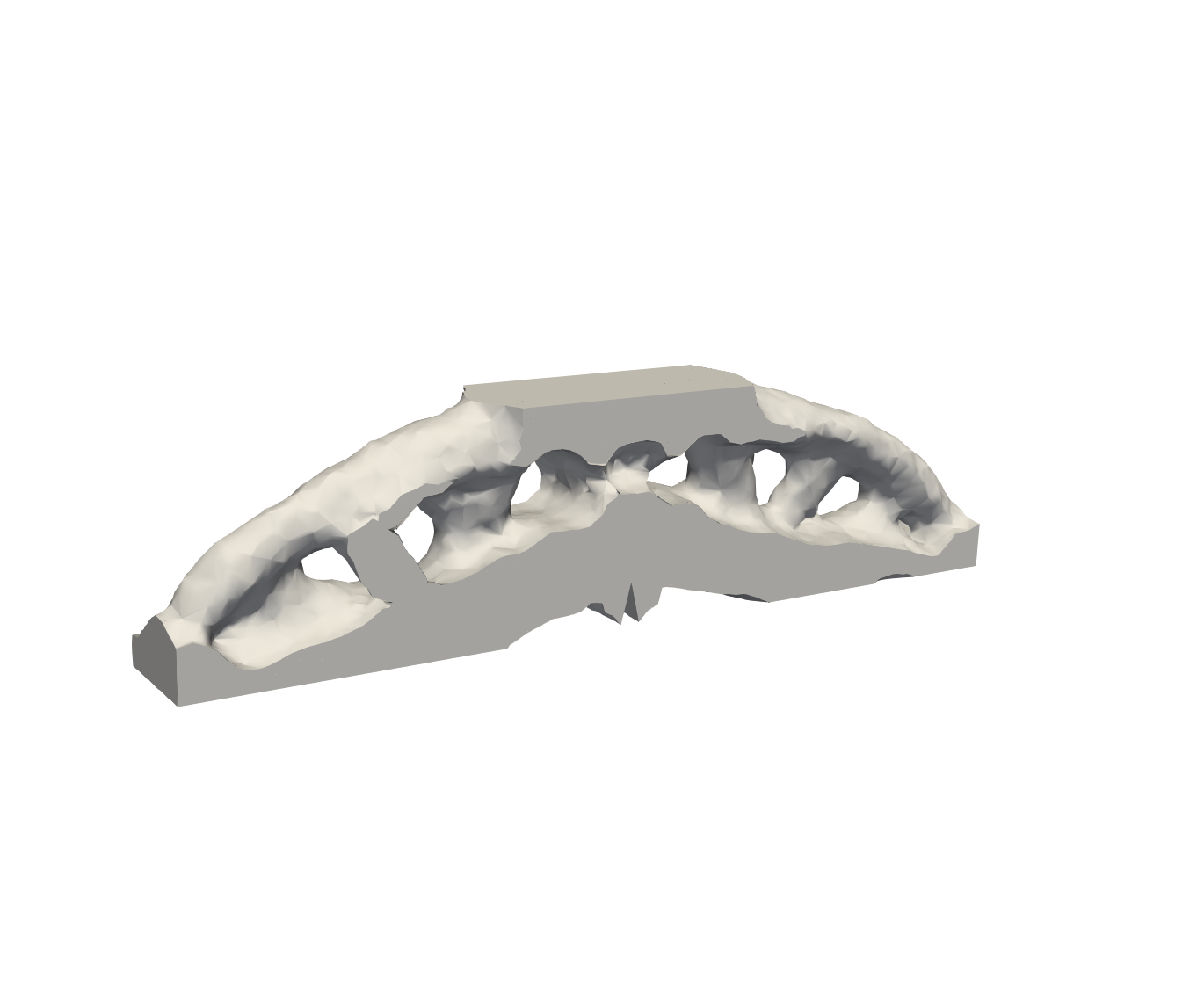}}   			\subfloat{\includegraphics[clip,trim=5cm 11.5cm 6cm 15cm, width=5.5cm]{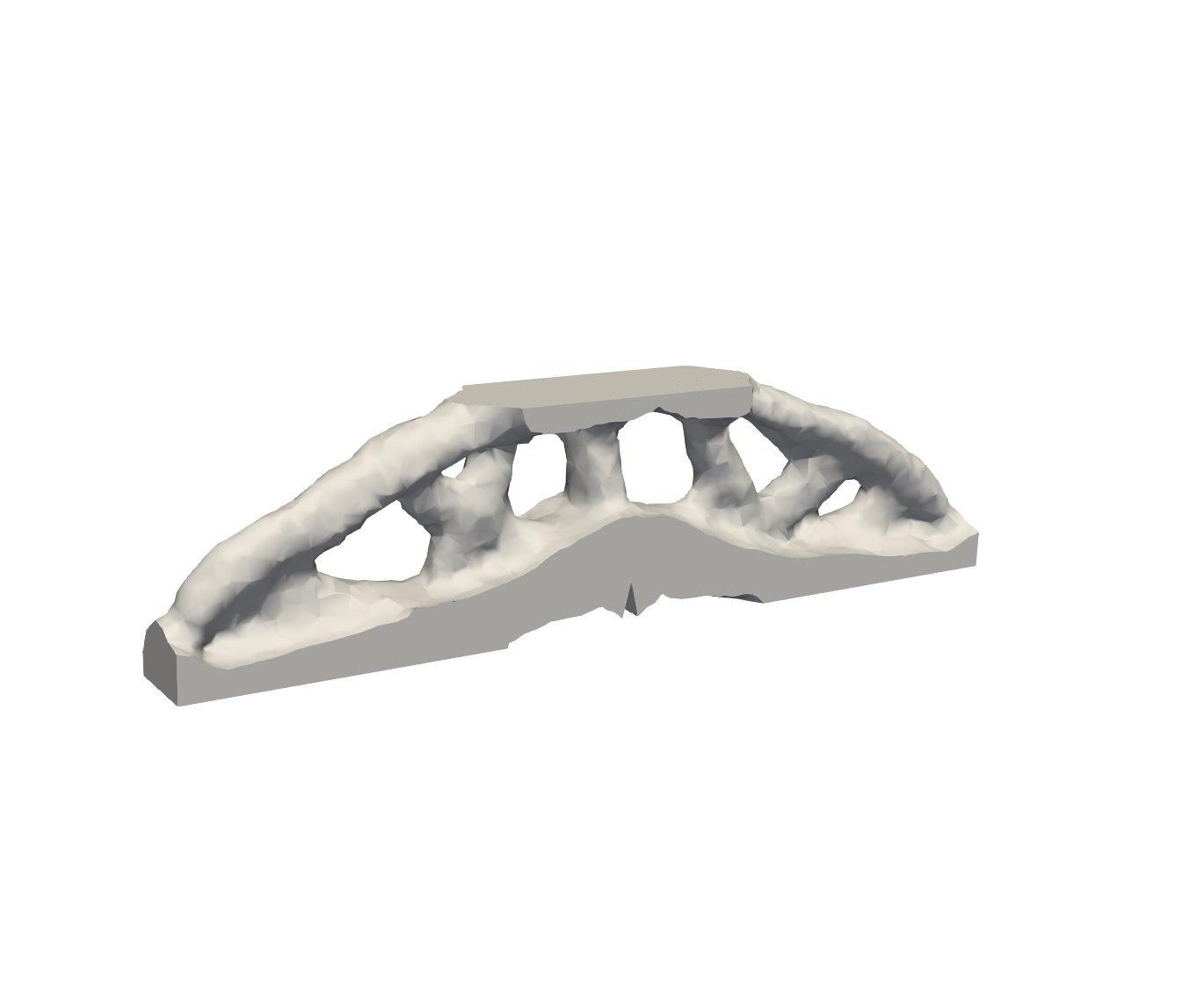}}
	\caption{Example 1. Evolution history of the optimal layouts for different volume ratio of  the three-point bending test based on Formulation 2.}
	
	\label{Exm1_C1}
\end{figure}

\begin{figure}[t!]
	\caption*{\hspace*{2cm}\underline{$\chi_v=0.89$}\hspace*{4cm}\underline{$\chi_v=0.69$}\hspace*{3.5cm}\underline{$\chi_v=0.40$}\hspace*{1.5cm}}
	\vspace{-0.1cm}
	\subfloat{\includegraphics[clip,trim=5cm 11.5cm 6cm 15cm, width=5.5cm]{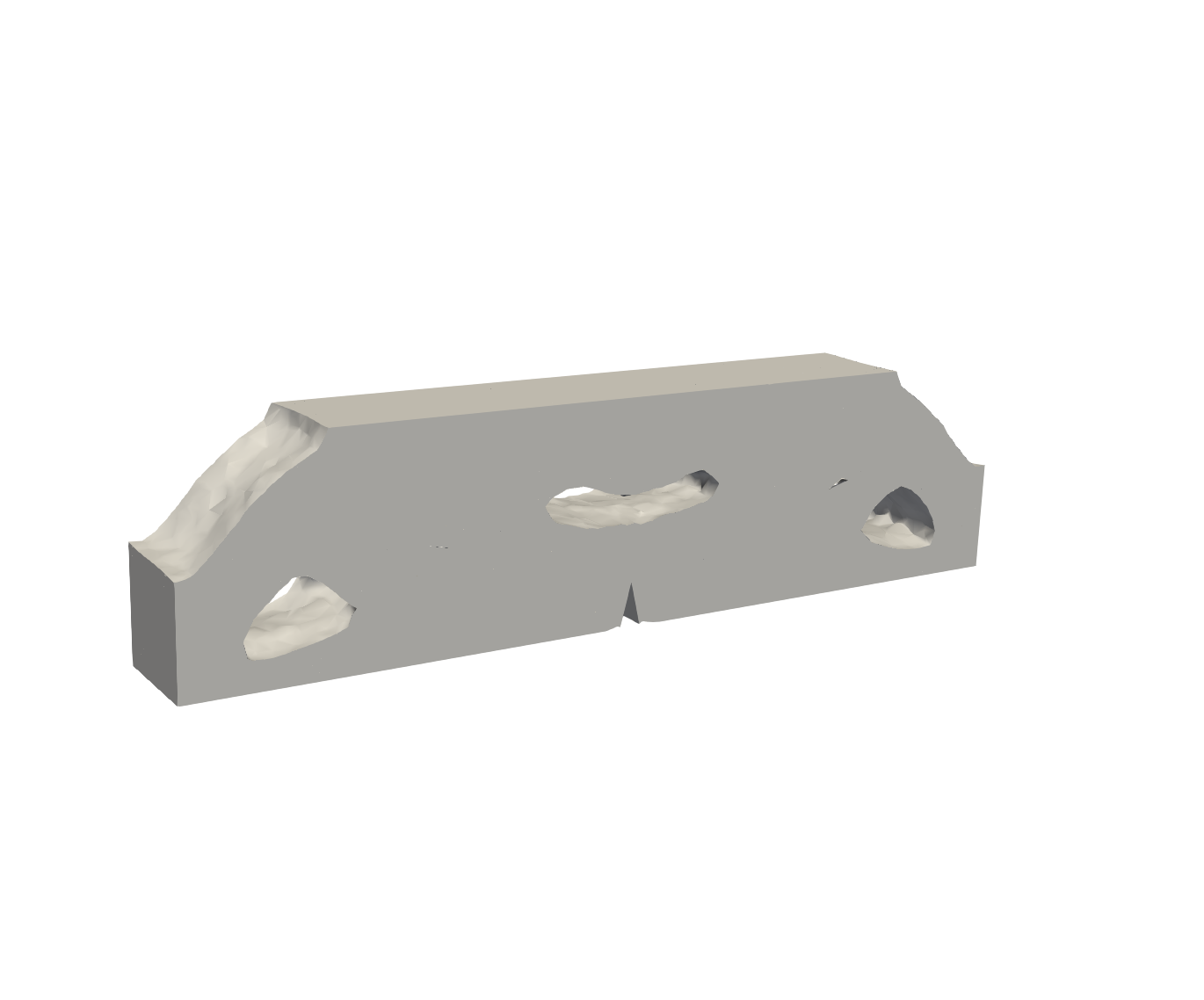}}   		\subfloat{\includegraphics[clip,trim=5cm 11.5cm 6cm 15cm, width=5.5cm]{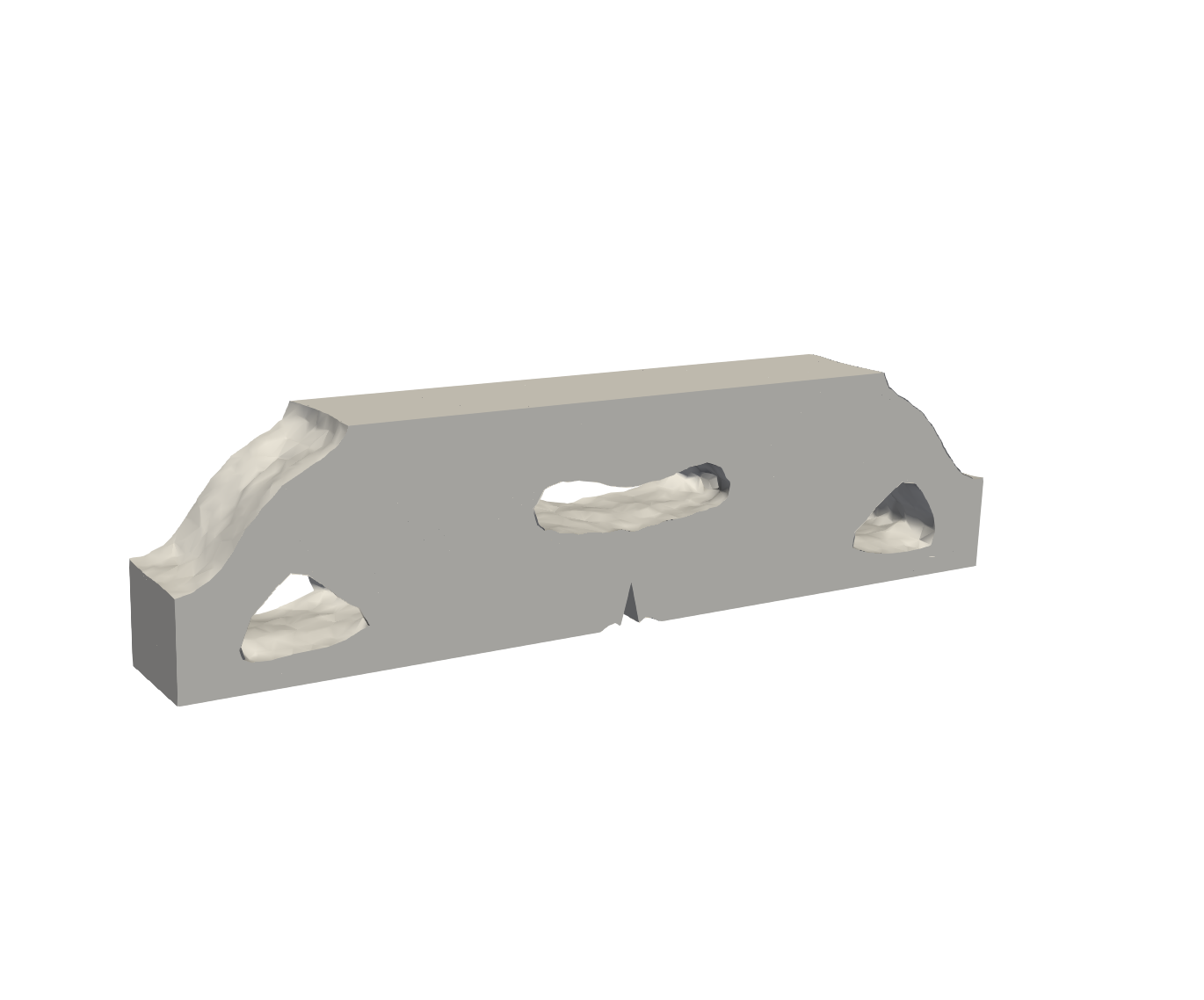}}   			\subfloat{\includegraphics[clip,trim=5cm 11.5cm 6cm 15cm, width=5.5cm]{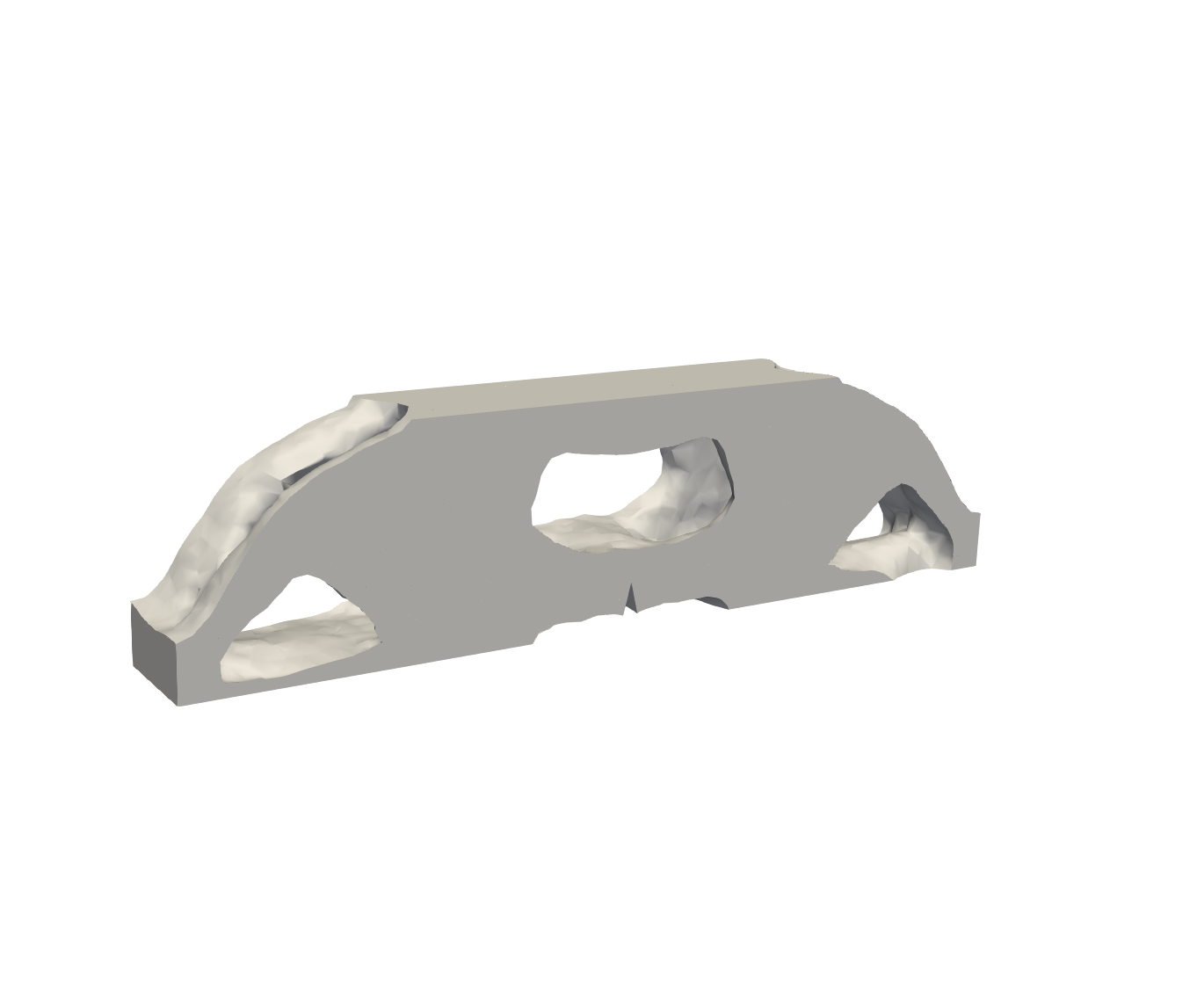}}
	\caption{Example 1. Evolution history of the optimal layouts for different volume ratio of  the three-point bending test based on linear elasticity results (excluding fracture phase-field).}
	\label{Exm1_C1_E}
\end{figure}

We start with the presentation of the quantitative and qualitative results from Formulation 1. The evolution history of the optimal layouts based on the topological field of the Formulation 1 and  2 are shown in Figure \ref{Exm1_C2}, and Figure \ref{Exm1_C1}, respectively. Additionally, for the sake of comparison, the final topologies for different volume ratios in the case of the pure linear elasticity are depicted in Figure \ref{Exm1_C1_E}. At this point, it is necessary to remark that a clear difference in the final layouts between Formulations 1, 2 and pure linear elasticity can be observed.  

\begin{figure}[!t]
	\centering
	{\includegraphics[clip,trim=1cm 30cm 0cm 0cm, width=18cm]{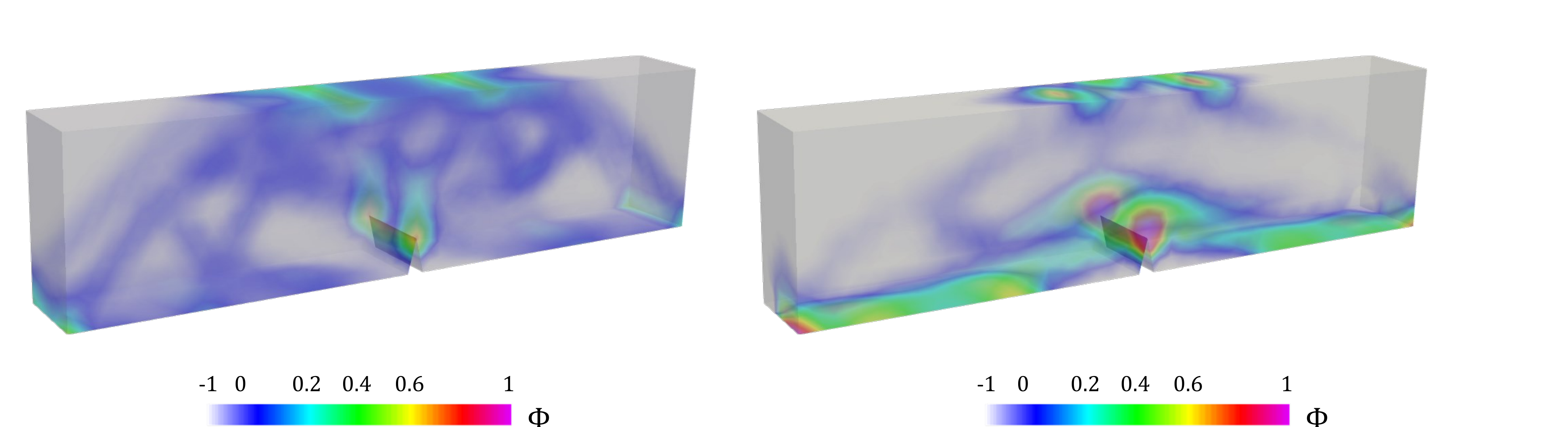}}  
	\vspace*{-0.4cm}
	\caption*{\hspace*{4.3cm}(a)\hspace*{8cm}(b)\hspace*{2cm}}
	\caption{Example 1. Computed topological field $\Phi$ at the final topology optimization iteration at $\chi_v=0.4$ based on (a) Formulation 1, and (b) Formulation 2.}
	\label{Exm1_LSM}
\end{figure}

The computed regularized topological field $\Phi(\Bx,t) \in [-1,1]$ by means of reaction-diffusion evolution through Formulation \req{form_4} for both Formulations 1 and 2 are illustrated in Figures \ref{Exm1_LSM}(a-b), respectively. The gray color shows non-material points such that  $\Phi(\Bx,t)<0$, while the colorful area represents material points with $\Phi(\Bx,t)>0$.  Additionally, the purple area denoted as $\Gamma_{\Phi}$ distinguishes between the material and the non-material sets as a zero-level set contour given by $\Phi(\Bx,t)=0$.

\begin{figure}[!b]
	\centering
	\vspace{-0.1cm}
	\includegraphics[clip,trim=1cm 33.5cm 1.5cm 1cm, width=16cm]{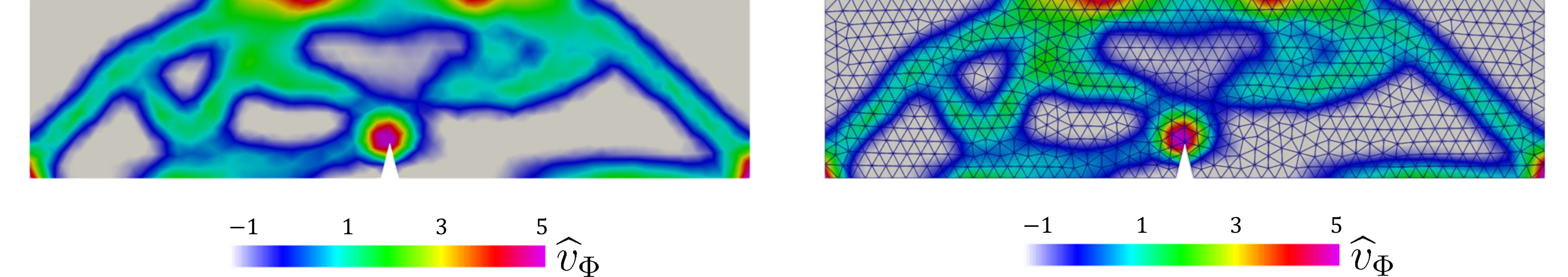}		
	\vspace*{-0.7cm}
	\caption*{\hspace*{2.3cm}(a)\hspace*{7.7cm}(b)\hspace*{2cm}}
	\caption{\noii{Example 1. Velocity distribution when $\chi_v=0.4$ on the (a) continuum space, and (b) asymmetric finite element discretization.}}
	\label{Exm1_velocity}
\end{figure}

\noii{The positive normal velocity indicates that the material distribution could be increased, while the negative normal velocity leads to nucleation of more voids in those regions. The normal velocity of the final optimum layout ${\widehat{v}_{\Phi}}\left| _{\Phi = 0} \right.$ for Example 1 (Formulation 1) is investigated). By means of Figure \ref{Exm1_velocity}(a), it is clear that in the vicinity of supports, loading areas, as well as the crack initiation zones (sharp corners regions) resulting the large positive quantity of the normal velocity field which prevents material removal. Also, Figure \ref{Exm1_velocity}(b) shows the normal velocity distribution overlaid on the finite element discretization such that the split elements can be clearly observed.}

\begin{figure}[!t]
	\centering
	{\includegraphics[clip,trim=0cm 28cm 0cm 0cm, width=17cm]{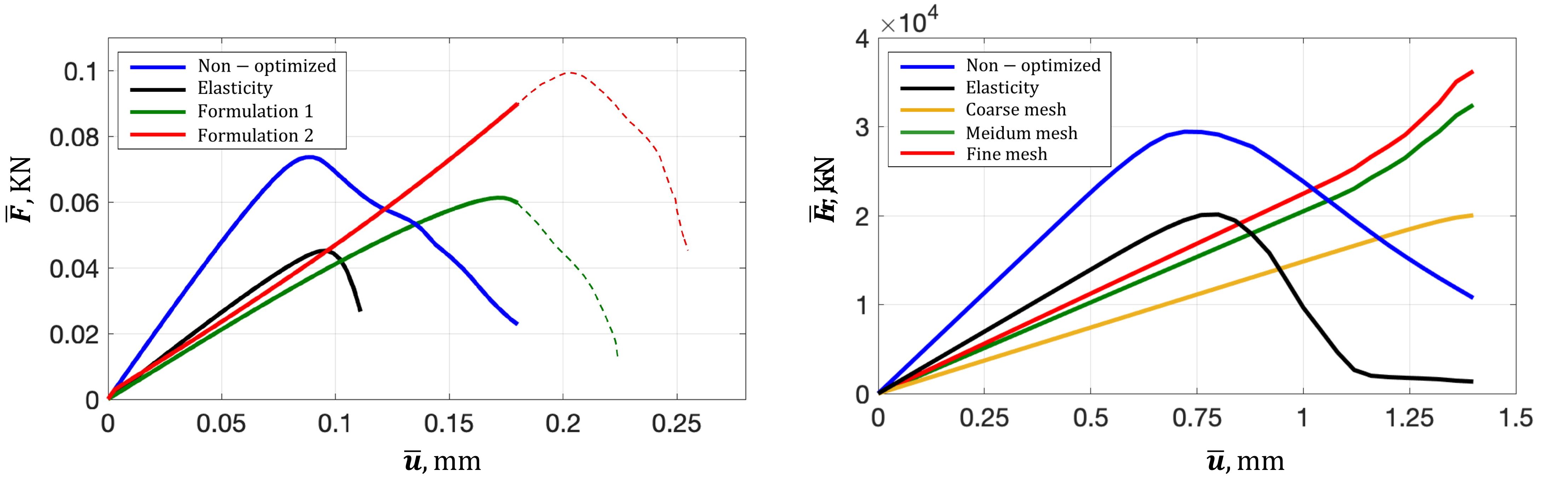}}  
	\vspace*{-0.7cm}
	\caption*{\hspace*{4.3cm}(a)\hspace*{8cm}(b)\hspace*{2cm}}
	\caption{ Comparison load-displacement curves for the, (a) three-point bending test under compression loading in Example 1, and (b) L-shaped panel test under mixed-mode fracture in Example 2.	}
	\label{Exm1_LD}
\end{figure}

For further comparison, the load-displacement curve for the non-optimized, elasticity, Formulation 1 and 2 results are depicted in Figure \ref{Exm1_LD}(a). A slight stiffness degradation can be observed in the final layout obtained from Formulation 1, which indicates the propagation of local cracks in the final optimal layout, during the time increment. It is worth noting that by means of Formulation 2, the maximum load-carrying capacity of the final layout has significantly increased compared with the other formulations. Thus, following Figure \ref{Exm1_LD}(a) the maximum load-carrying capacity employing Formulation 2 is $52.7 \%$ greater than the pure linear elasticity case and $35.7 \%$ greater than Formulation 1. It can be also grasped that the failure displacement (when the crack initiation is observed) for Formulation 1 is $51.4 \%$ greater than the non-optimized case and $49.7 \%$ greater than the pure linear elasticity case. \noii{Additionally, the failure softening region for both formulations is depicted in Figure \ref{Exm1_LD}(a) as a dashed line, to illustrate fractured state.}

\begin{figure}[!t]
	\centering
	{\includegraphics[clip,trim=0cm 28cm 0cm 0cm, width=17cm]{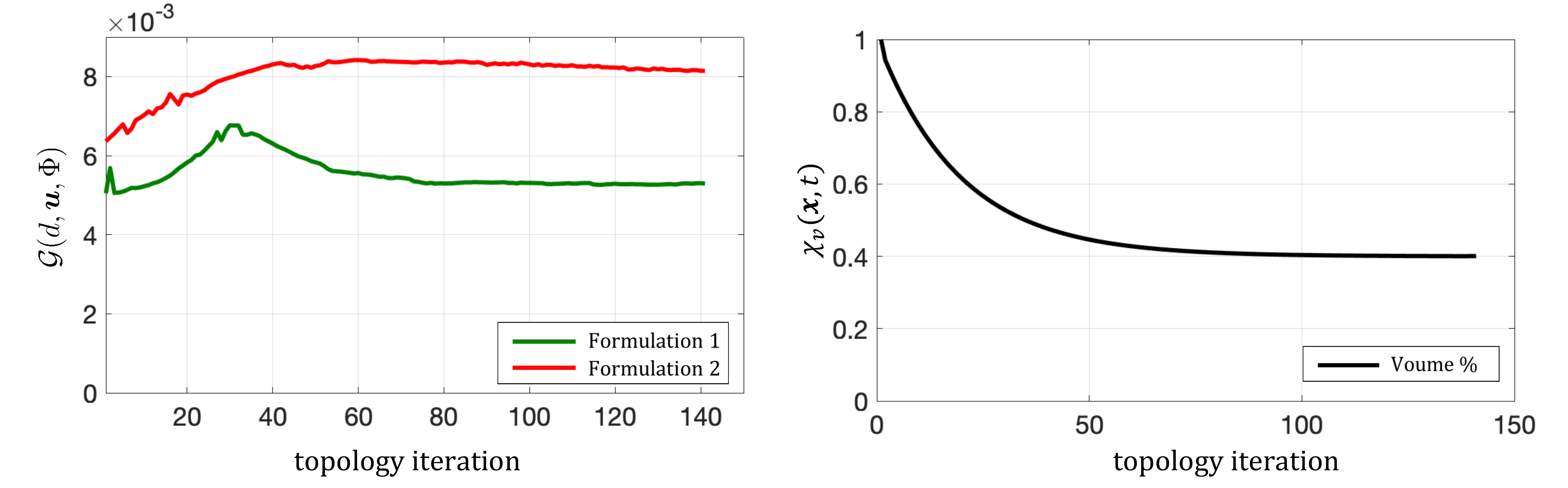}}  
	\vspace*{-0.7cm}
	\caption*{\hspace*{4.3cm}(a)\hspace*{8cm}(b)\hspace*{2cm}}
	\caption{Example 1. Convergence history for three-point bending of the (a) objective function for the Formulation 1-2, and (b) volume constraint function.}
	\label{Exm1_conv}
\end{figure}

Another impacting factor that should be noted is the convergence history of objective (i.e., \req{eq:fun_obj}) and volume constraint (i.e., \req{eq:vol_const}) functions which are shown in Figures \ref{Exm1_conv}(a-b), respectively. We note that, by considering the objective function in Figure \ref{Exm1_conv}(a), evidently Formulation 2 shows higher objectivity, i.e., a stiffer structural response due to larger stored mechanical energy in structure, \grm{compared} to Formulation 1 with the same volume ratio. Thus, Formulation 2 underlines its efficiency by adding extra constraints to the minimization problem. At this point, it is necessary to remark that in Formulation 1, one requires only one adjoint sensitivity equation to be solved through  \req{eq:sol_adjoint2} compared with Formulation 2 in \req{eq:sol_adjoint}, thus keeping the computational cost reasonably low. We note that the volume constraint in Figures \ref{Exm1_conv}(b) for all the methods indicates similar response and asymptotically approaches to its desired value, i.e., $\chi_v=0.40$.

\begin{figure}[!t]
	\centering
	{\includegraphics[clip,trim=1cm 14.3cm 0cm 3cm, width=17cm]{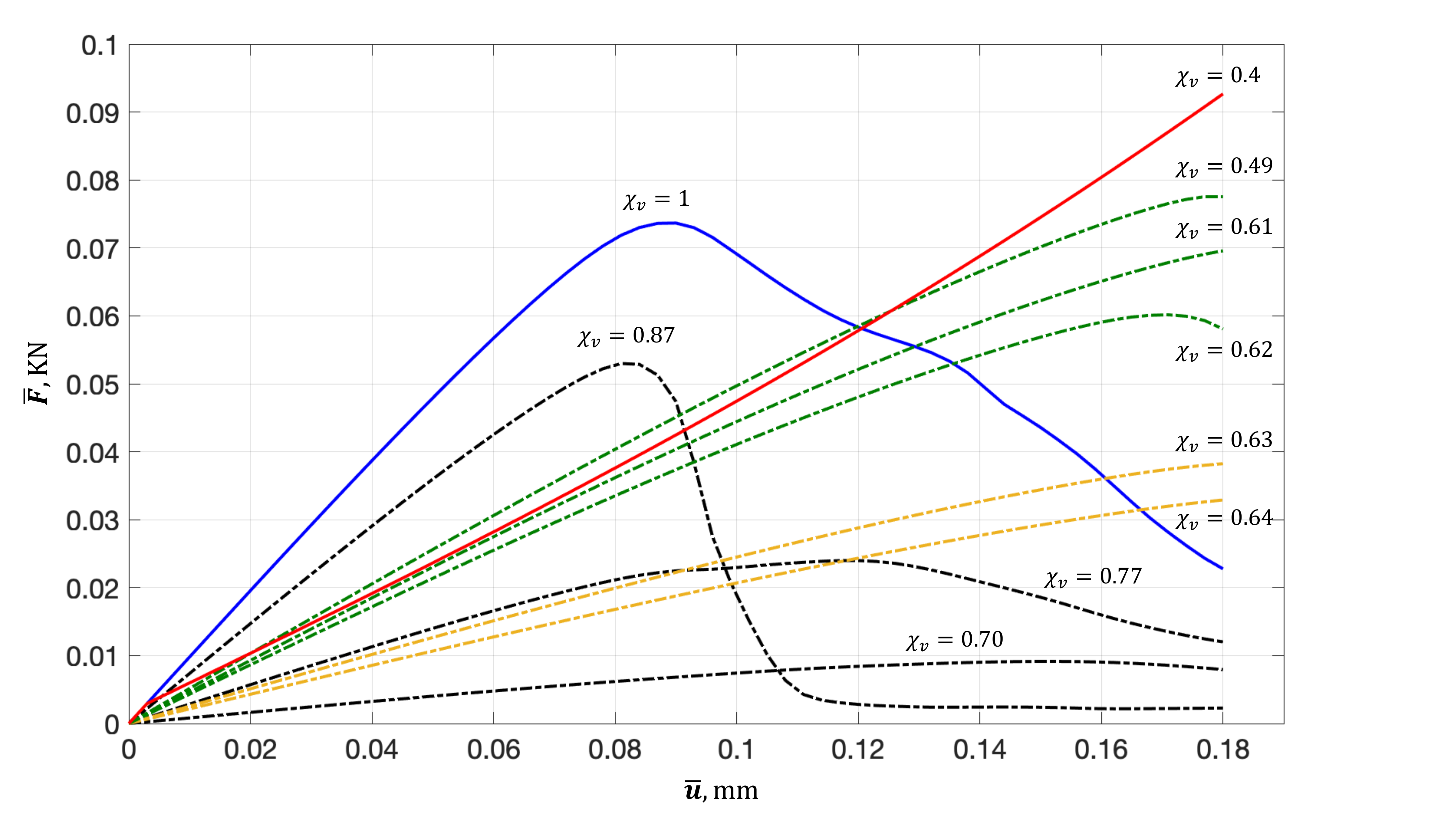}}  
	\vspace*{-0.7cm}
	\caption{Example 1. Comparison load-displacement curves for different volume ratio of topology optimization iteration up to final optimum layout through Formulation 2.}
	\label{Exm1_inv}
\end{figure}

Additionally, the load-displacement curves for different volume ratios, during the evolution of the topology optimization process based on the Formulation 2 are depicted in Figure \ref{Exm1_inv}. Indeed, this highlights the effects of the damage response within every new topology, and thus by approaching $\chi_v=0.40$ the topology is improved significantly, and less fracture is observed. Finally, at $\chi_v=0.40$ there will be no softening due to fracture (no crack initiation). Thus, the proposed model showed its proficiency. 

\begin{figure}[!t]
	\vspace{-0.1cm}
	\subfloat{\includegraphics[clip,trim=5cm 11.5cm 6.7cm 14cm, width=8.4cm]{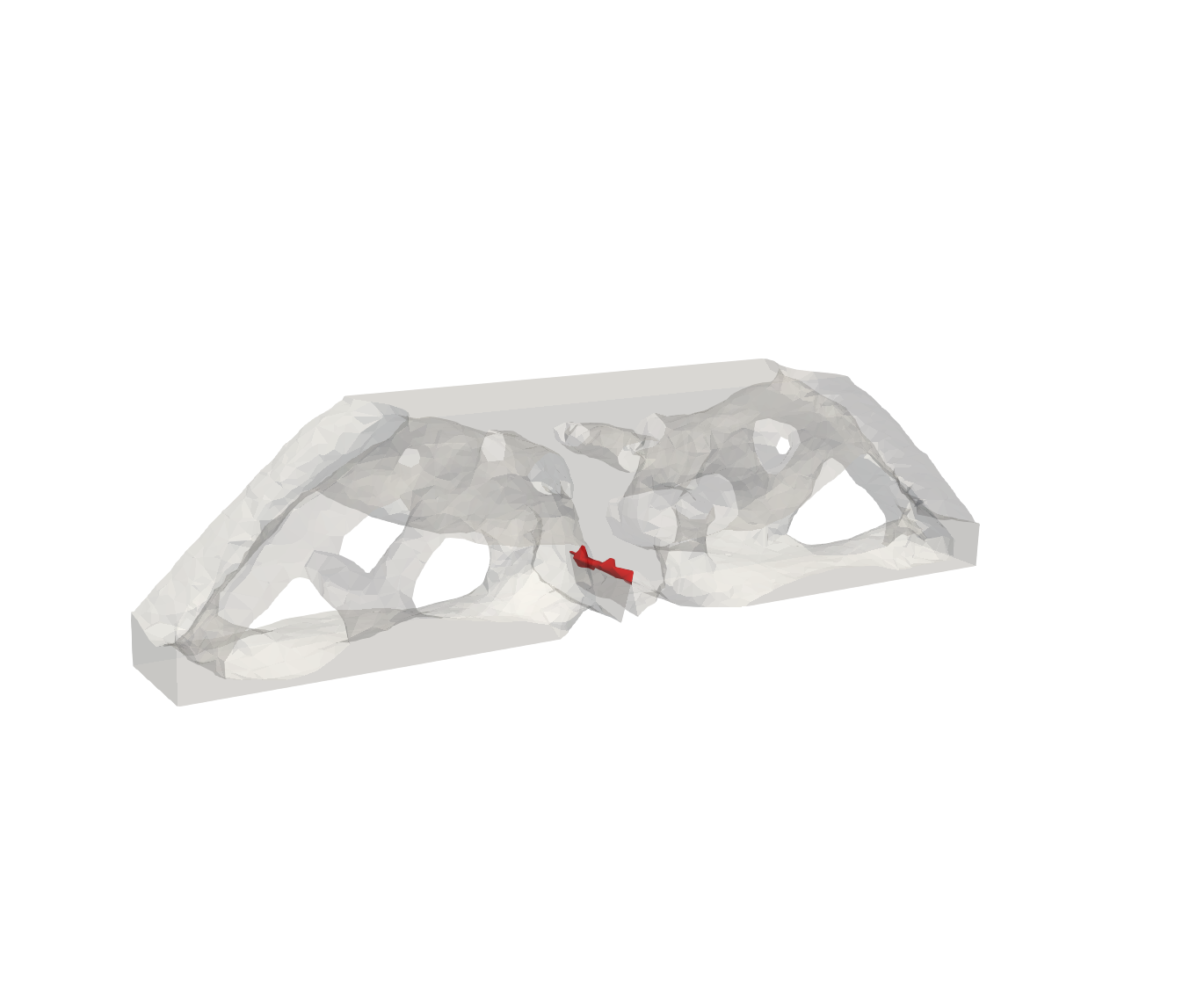}}   		
	\subfloat{\includegraphics[clip,trim=5cm 11.5cm 6.7cm 14cm, width=8.4cm]{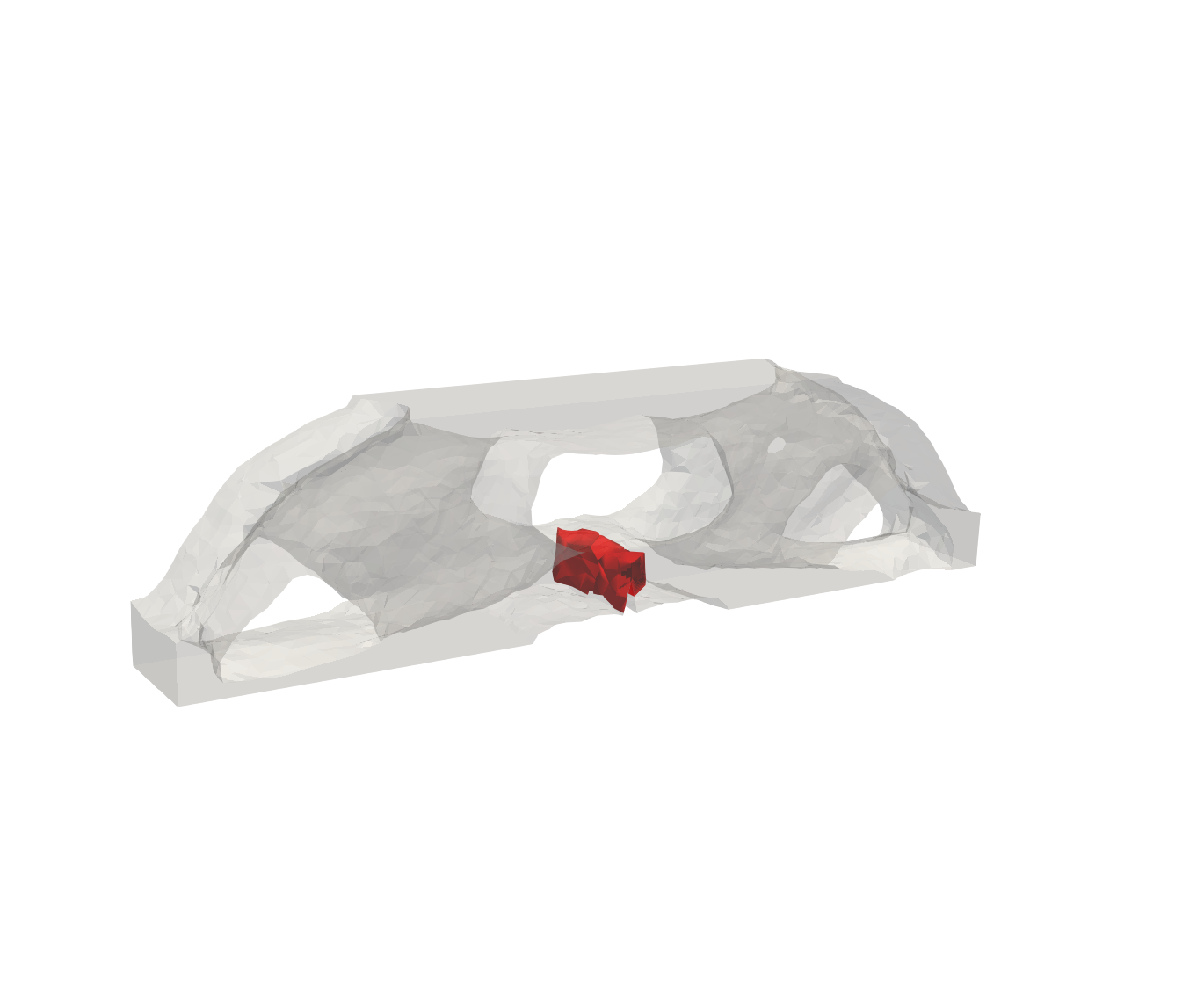}}   		
	\caption*{\hspace*{1.8cm}(a)\hspace*{8cm}(b)}
	\caption{Example 1. The crack phase-field response for optimum layout based on (a) Formulation 1, and (b) linear elasticity results undergoing brittle fracture.}
	\label{Exm1_C1_elast_crack}
\end{figure}

\begin{figure}[!t]
	\vspace{-0.1cm}
	\subfloat{\includegraphics[clip,trim=5cm 11.5cm 6cm 12cm, width=5.5cm]{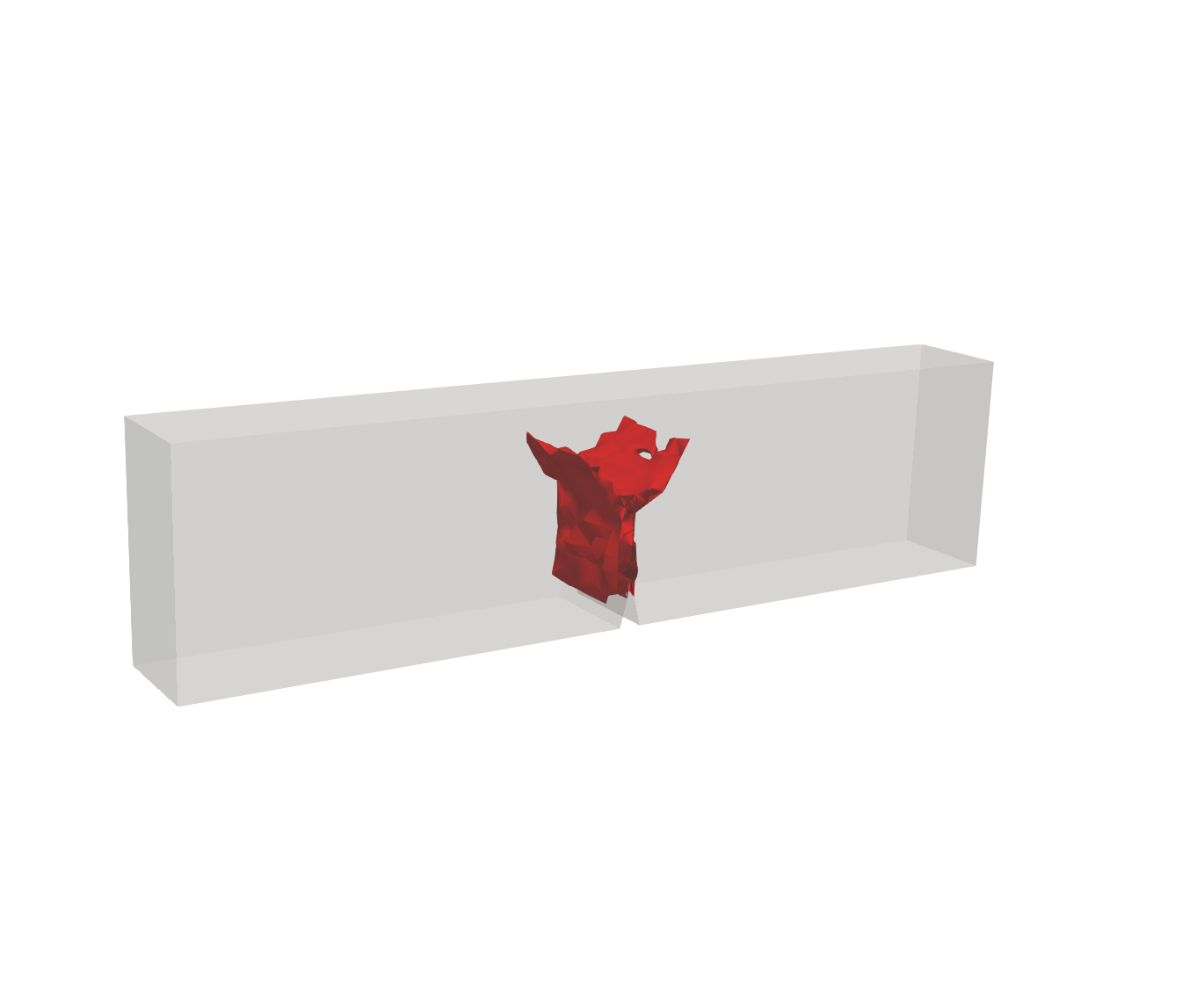}}   		\subfloat{\includegraphics[clip,trim=5cm 11.5cm 6cm 12cm, width=5.5cm]{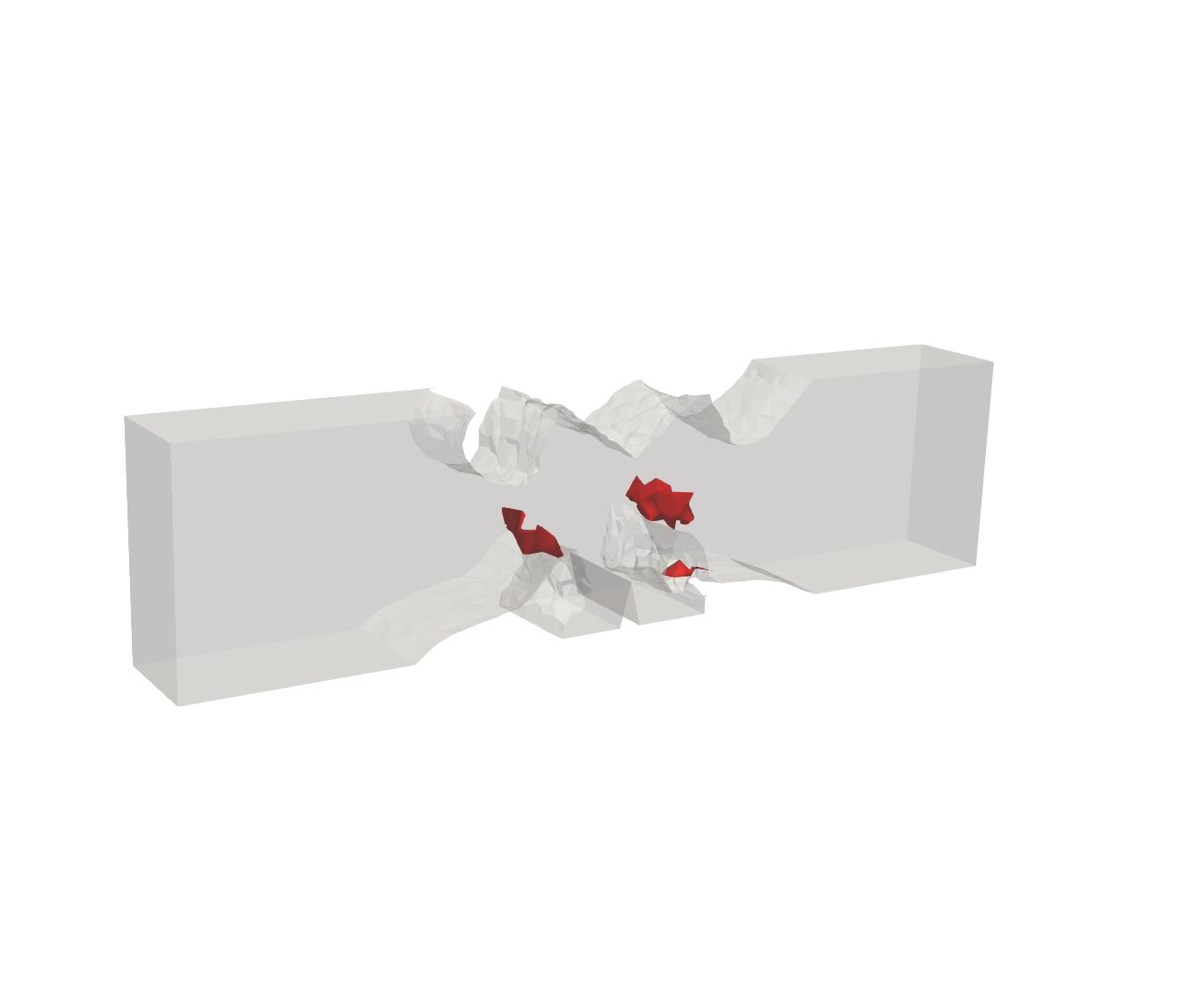}}   			\subfloat{\includegraphics[clip,trim=5cm 11.5cm 6cm 12cm, width=5.5cm]{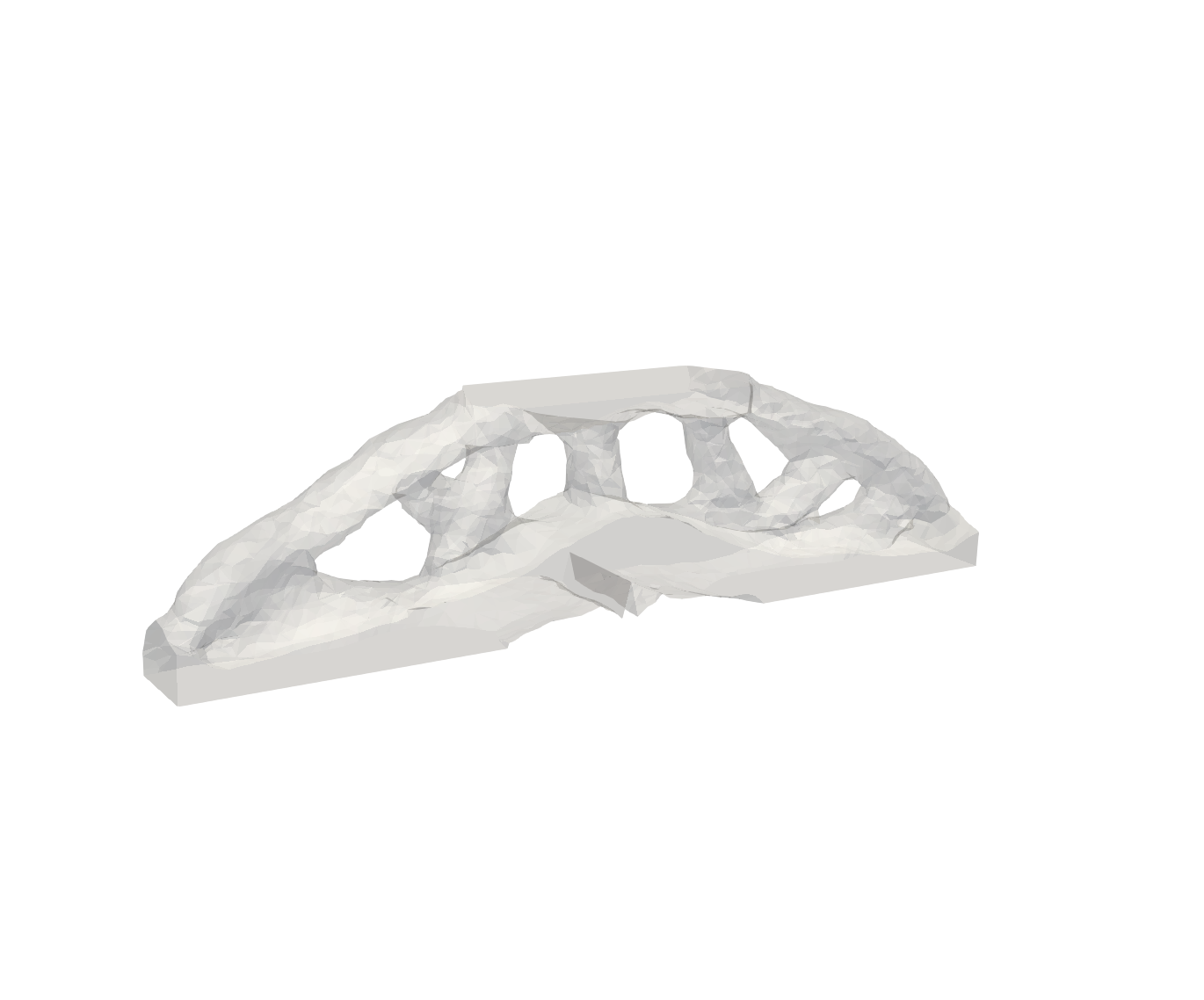}}
	    \caption*{\hspace*{1.3cm}(a)\hspace*{5cm}{(b)}\hspace*{5cm}{(c)}}
	\caption{Example 1. The crack phase-field response for different optimal topology layout for (a) non-optimized result $\chi_v=1$, (b) $\chi_v=0.71$, and (c) final optimal topology layout $\chi_v=0.40$ through Formulation 2.}
	\label{Exm1_C1_crack}
\end{figure}

The crack phase-field profiles due to optimal layout for Formulation 1, and elasticity results are depicted in Figure \ref{Exm1_C1_elast_crack}, while, for  Formulation 2 similar crack profiles are shown for different time stages in Figure \ref{Exm1_C1_crack}. The first important observation is that the final layouts obtained from Formulations 1 and 2 have overcome the fracture limitations caused by stiffness degradation (compare with Figure \ref{Exm1_C1_crack}(a) for the non-optimized result). It is of great importance that all approaches reach failure due to the softening, while the optimal layout due to Formulation 2 has not reached failure, nor any damage is visible. 

\begin{figure}[t!]
	\subfloat{\includegraphics[clip,trim=5cm 11.5cm 6.7cm 14cm, width=8.4cm]{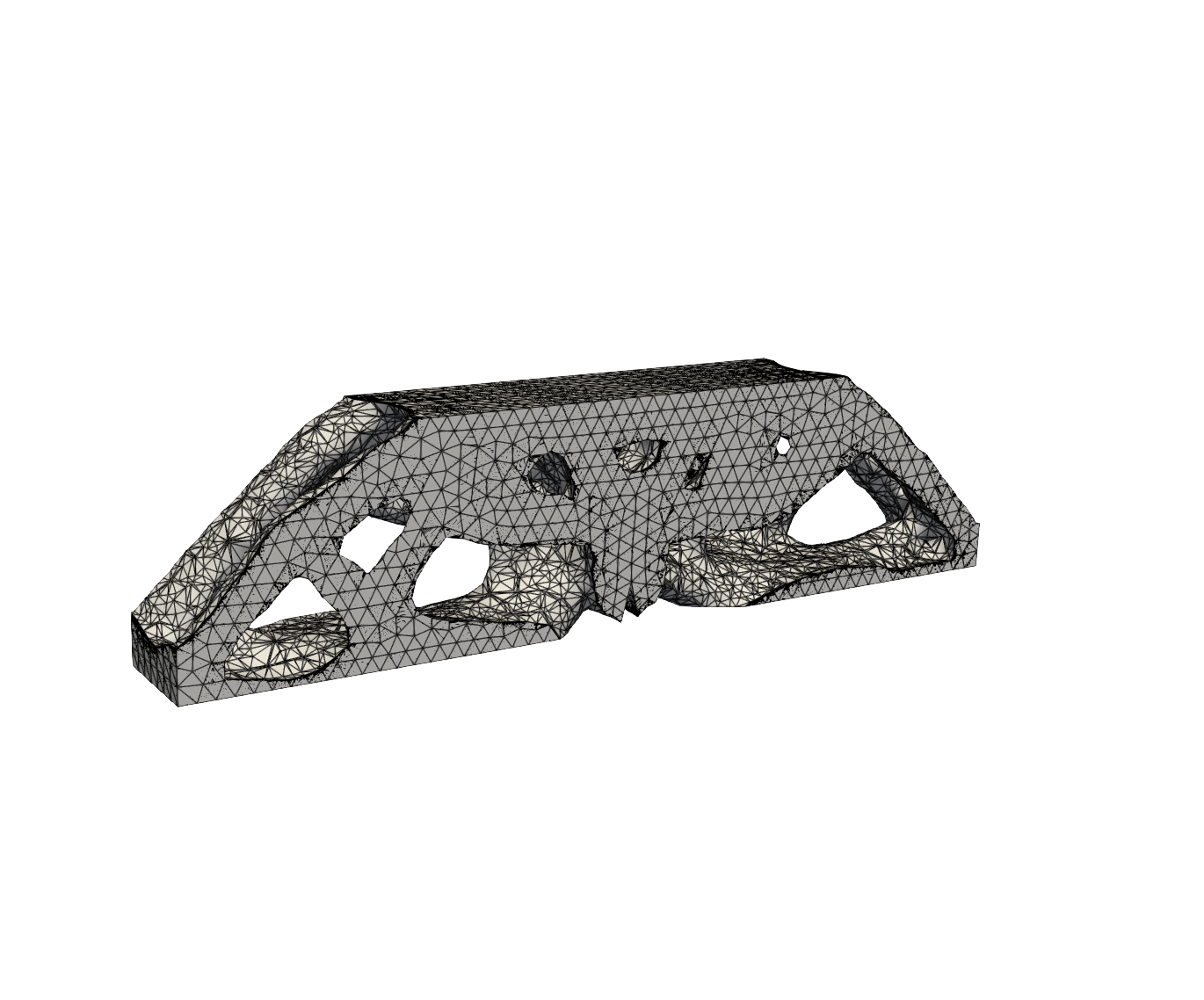}}   		
	\subfloat{\includegraphics[clip,trim=5cm 11.5cm 6.7cm 14cm, width=8.4cm]{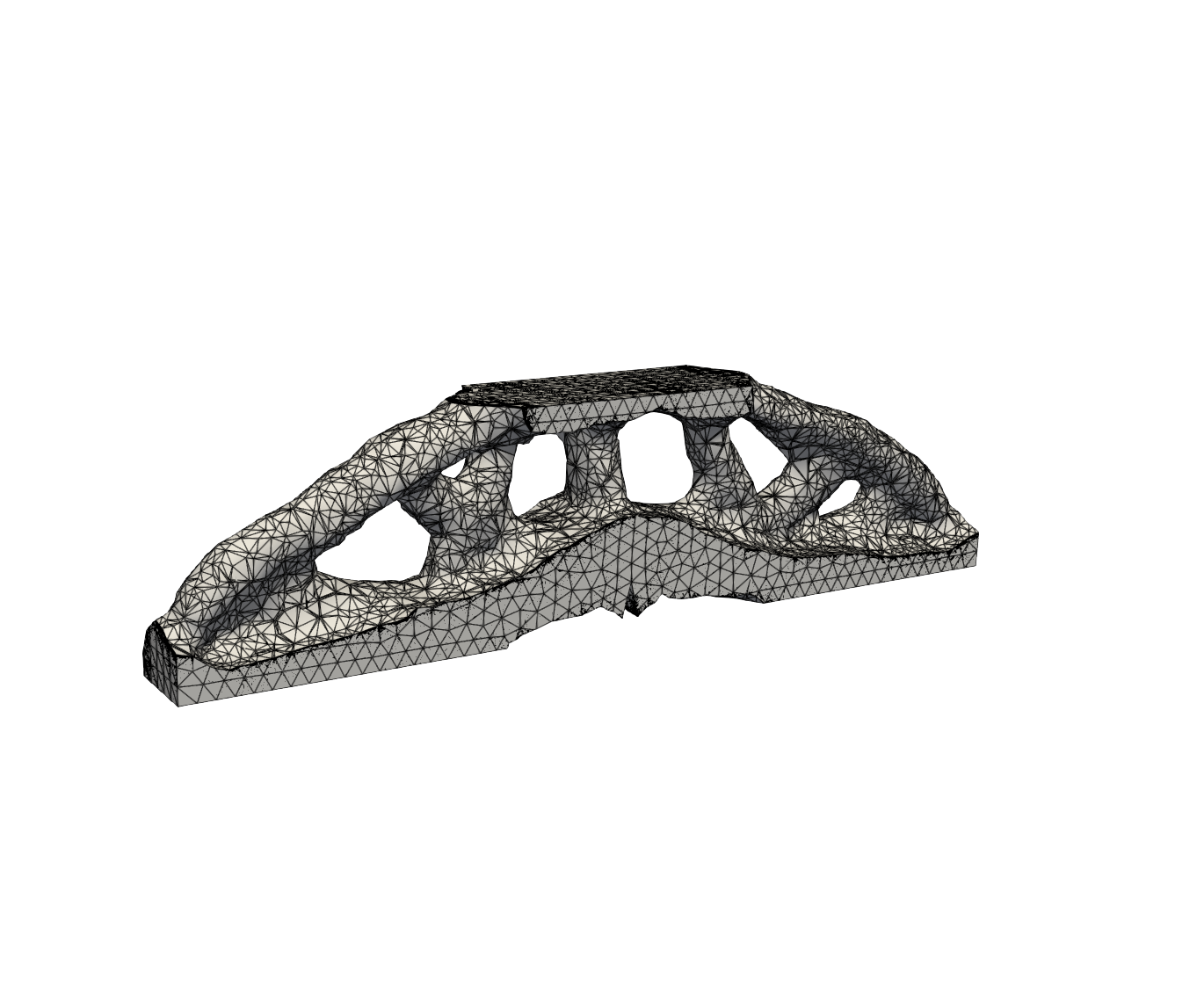}}   		
	\caption*{\hspace*{1.8cm}(a)\hspace*{8cm}(b)}
	\caption{Example 1. Finite element discretization of the optimal layouts when $\chi_v=0.4$ based on (a) Formulation 1, and (b) \noii{Formulation 2}.}
	\label{Exm1_mesh}
\end{figure}

Lastly, the discretized final optimum layouts for Formulations 1 and 2 are depicted in Figures \ref{Exm1_mesh}(a-b), respectively. The smooth topological interface (i.e., zero-level set) is remarkably grasped. 
Since Formulation 2 showed its proficiency versus other methods, for the remaining numerical examples, we investigate topology optimization due to brittle/ductile fracture by means of only Formulation 2.

\sectpb[Section52]{Example 2: L-shaped panel test under mixed-mode fracture}

The second example is concerned with a mixed-mode fracture problem, reported in \cite{jansari2019adaptive,mesgarnejad2015validation}. The main objective of this example is to evaluate the mesh sensitivity in the proposed numerical method. To this end, an L-shaped specimen with \textit{three different mesh sizes} under mixed-mode fracture is considered. By the L-shaped specimen, we aim to demonstrate the effects of stress concentration on fracture resistance in the optimal topology.

\begin{figure}[!t]
	\centering
	{\includegraphics[clip,trim=1cm 22cm 0cm 1cm, width=16cm]{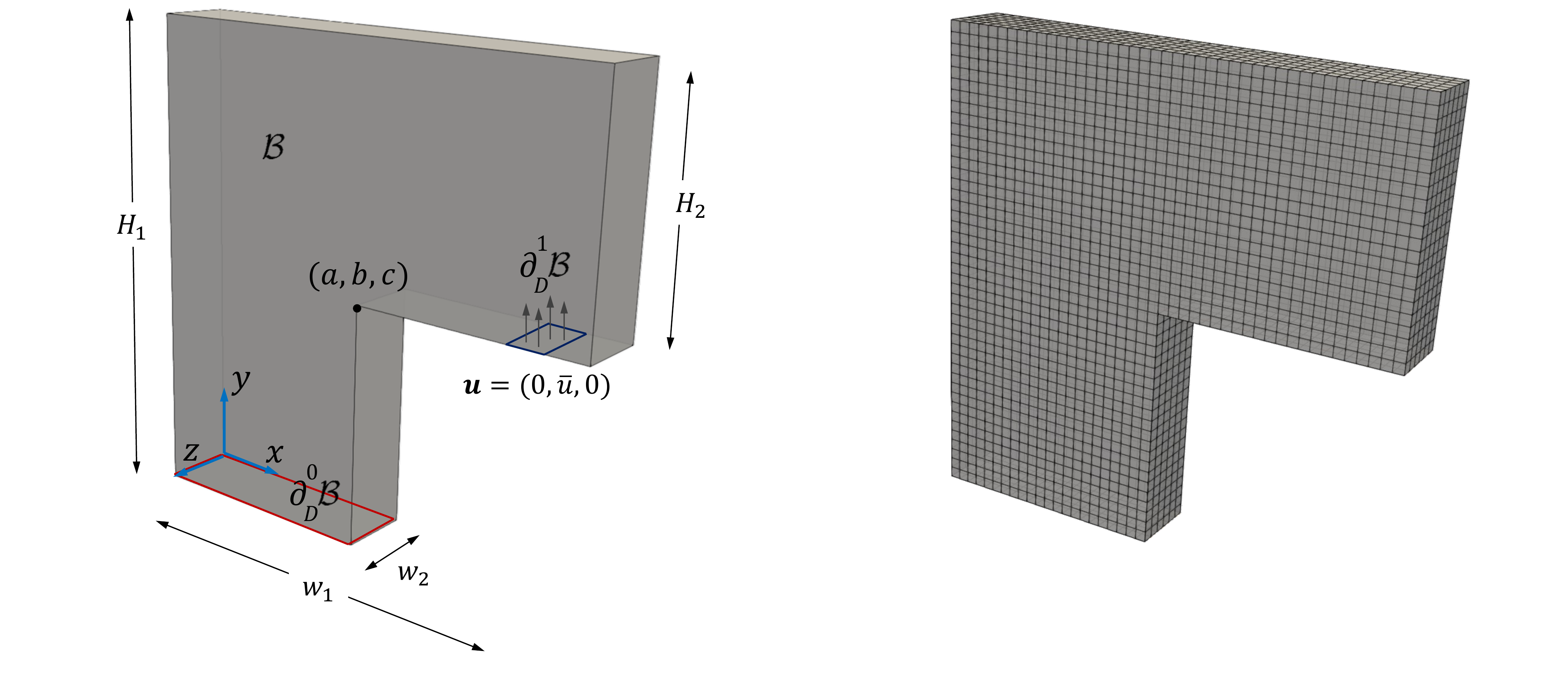}}  
	\vspace*{-0.7cm}
	\caption*{\hspace*{4.3cm}(a)\hspace*{8cm}(b)\hspace*{2cm}}
	\caption{Example 2. The representation of the (a) geometry and boundary conditions, and (b) finite element discretization.}
	\label{Exm2_bvp}
\end{figure}

These mesh sizes are denoted as coarse, medium, and fine discretization mesh. A boundary value problem is shown in Figure \ref{Exm2_bvp}(b), while the bottom surface is constrained in all directions $x-y-z$ directions. The geometrical dimensions are set as $w_1=500\;mm$, $w_2=75\;mm$,  $H_1=500\;mm$, and $H_2=250\;mm$. The reference point \grm{is} located in $(a,b,c)=(250,250,75)\;mm$. The numerical example is performed by applying a monotonic displacement increment ${\Delta \bar{u}}_x=1\times10^{-2}\;mm$ in the $area=(461,250,75)\times (481,250,0)$ of the specimen for 140-time steps. Thus, the final prescribed displacement load for this minimization problem is set as ${\bar{u}}_y=1.4\;mm$. The minimum finite element size in the design domains for the coarse mesh is $h=1.8$ which contains 3840 elements, for the medium mesh is $h=1.4$ which contains 7200 elements, and lastly for the fine mesh is $h=1.2\;mm$ which contains 8712 trilinear hexahedral elements. \noii{We note that for linear elasticity, we use the same discretization space as the fine mesh.}

\begin{figure}[t!]
 	\caption*{\underline{(a) Coarse mesh}}
 	\caption*{\hspace*{1cm}\underline{$\chi_v=0.52$}\hspace*{4.5cm}\underline{$\chi_v=0.42$}\hspace*{3.5cm}\underline{$\chi_v=0.4$}\hspace*{1.5cm}}
 	\vspace{-0.1cm}
 	\subfloat{\includegraphics[clip,trim=9.5cm 7.6cm 8cm 5.6cm, width=5.5cm]{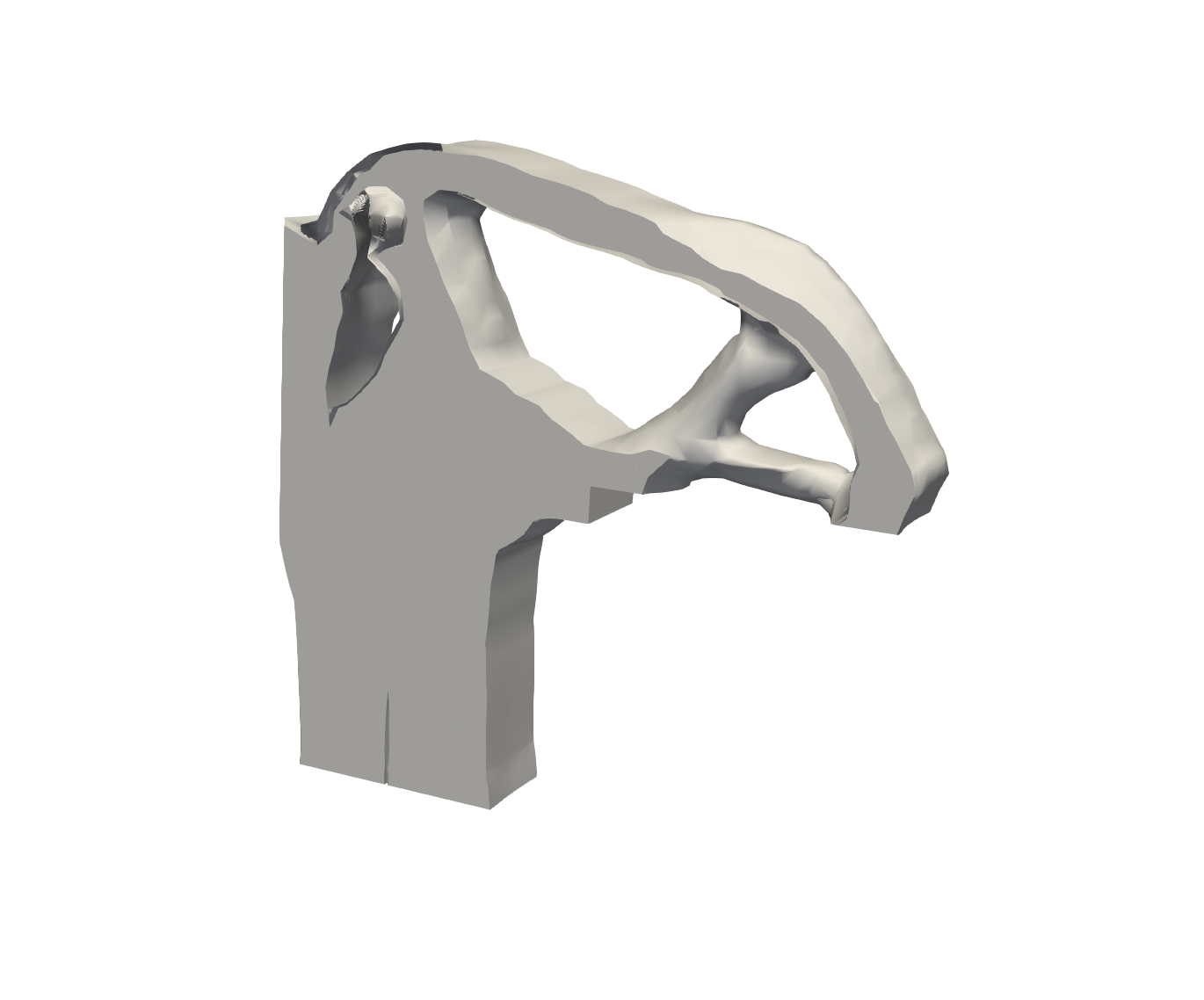}}   		\subfloat{\includegraphics[clip,trim=9.5cm 7.6cm 8cm 5.6cm, width=5.5cm]{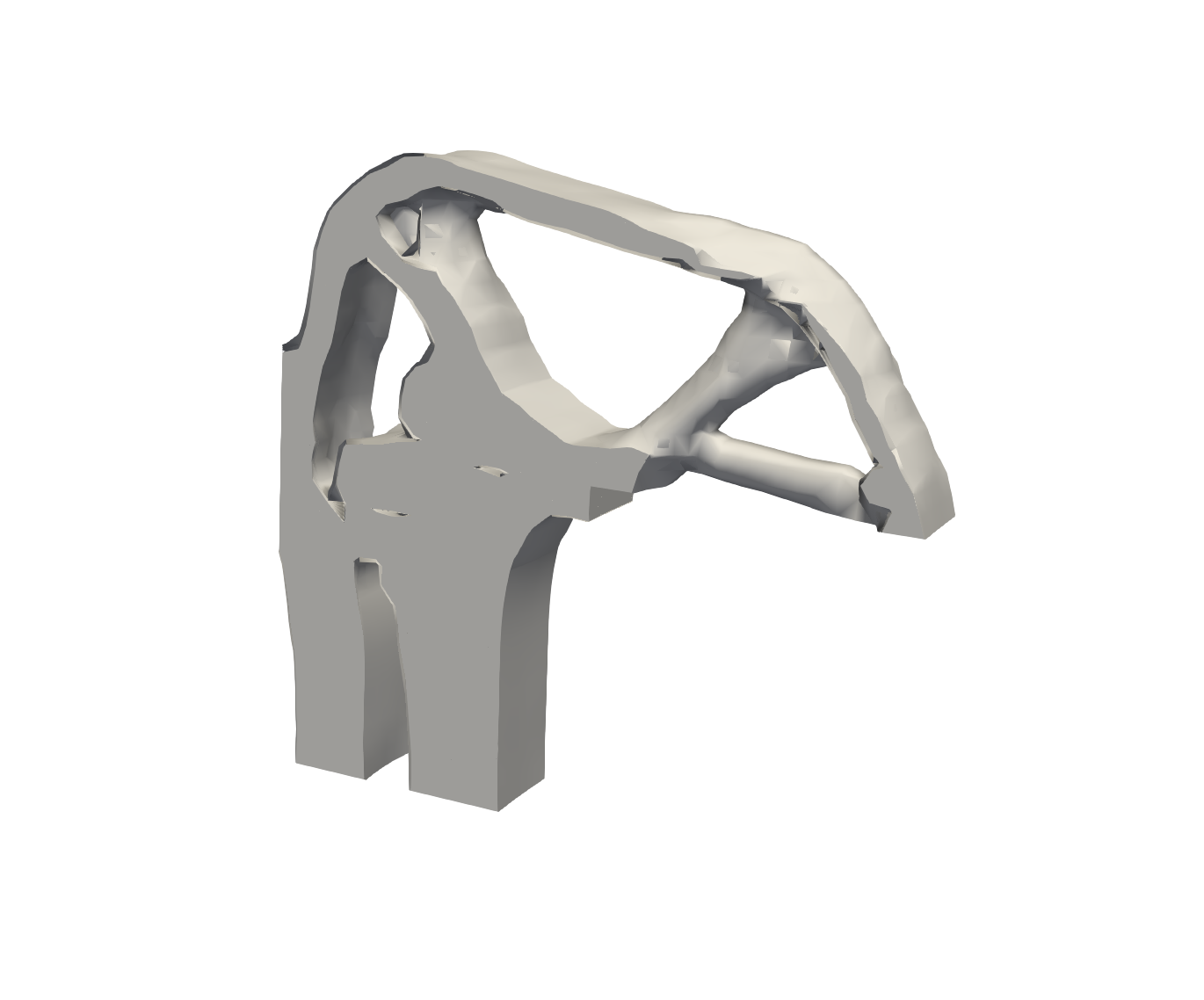}}   			\subfloat{\includegraphics[clip,trim=9.5cm 7.6cm 8cm 5.6cm, width=5.5cm]{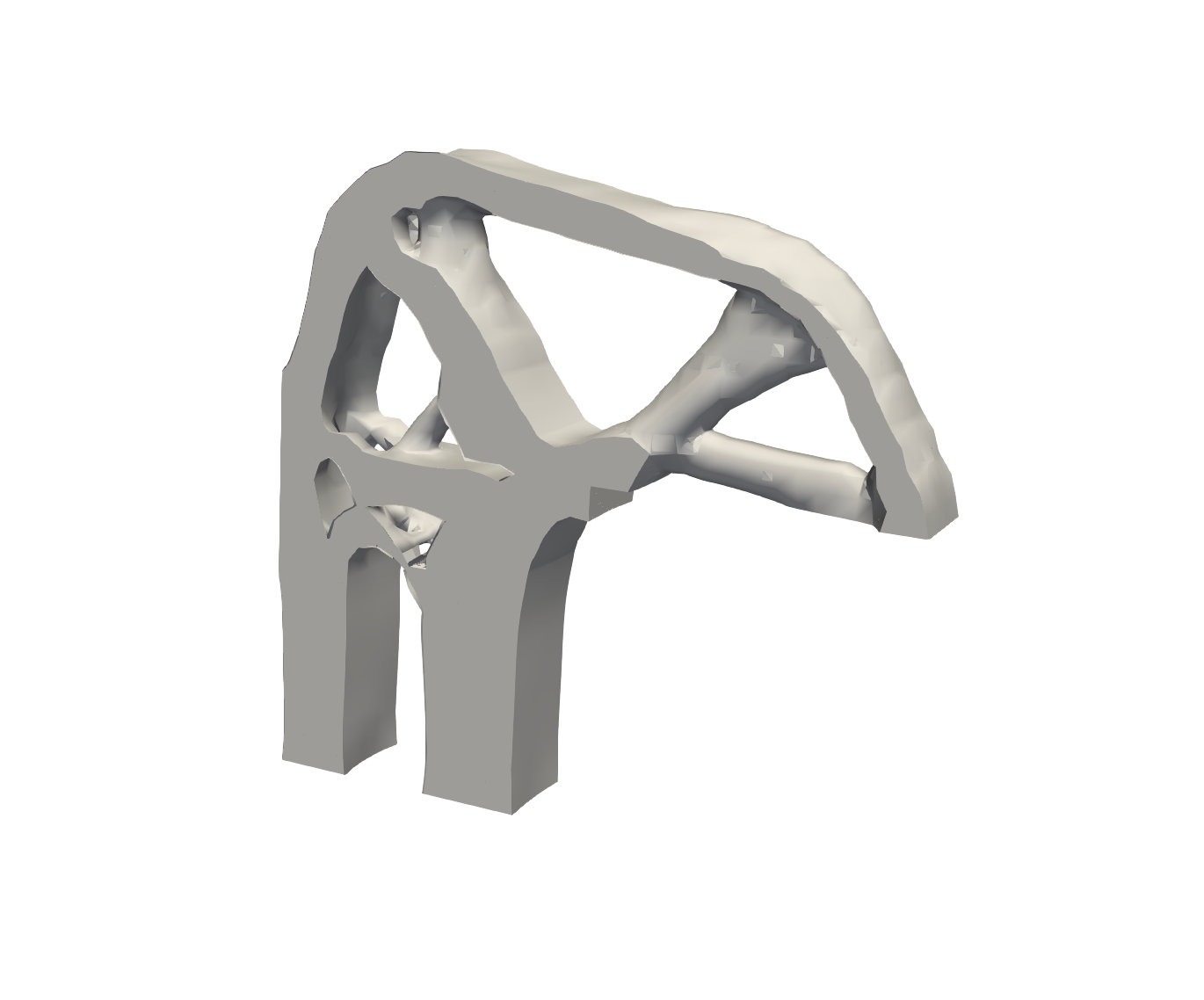}}
 	\caption*{\underline{(b) Medium mesh}}
 	\caption*{\hspace*{1cm}\underline{$\chi_v=0.52$}\hspace*{4.5cm}\underline{$\chi_v=0.42$}\hspace*{3.5cm}\underline{$\chi_v=0.4$}\hspace*{1.5cm}}
 	\vspace{-0.1cm}
 	\subfloat{\includegraphics[clip,trim=9.5cm 7.6cm 8cm 5.6cm, width=5.5cm]{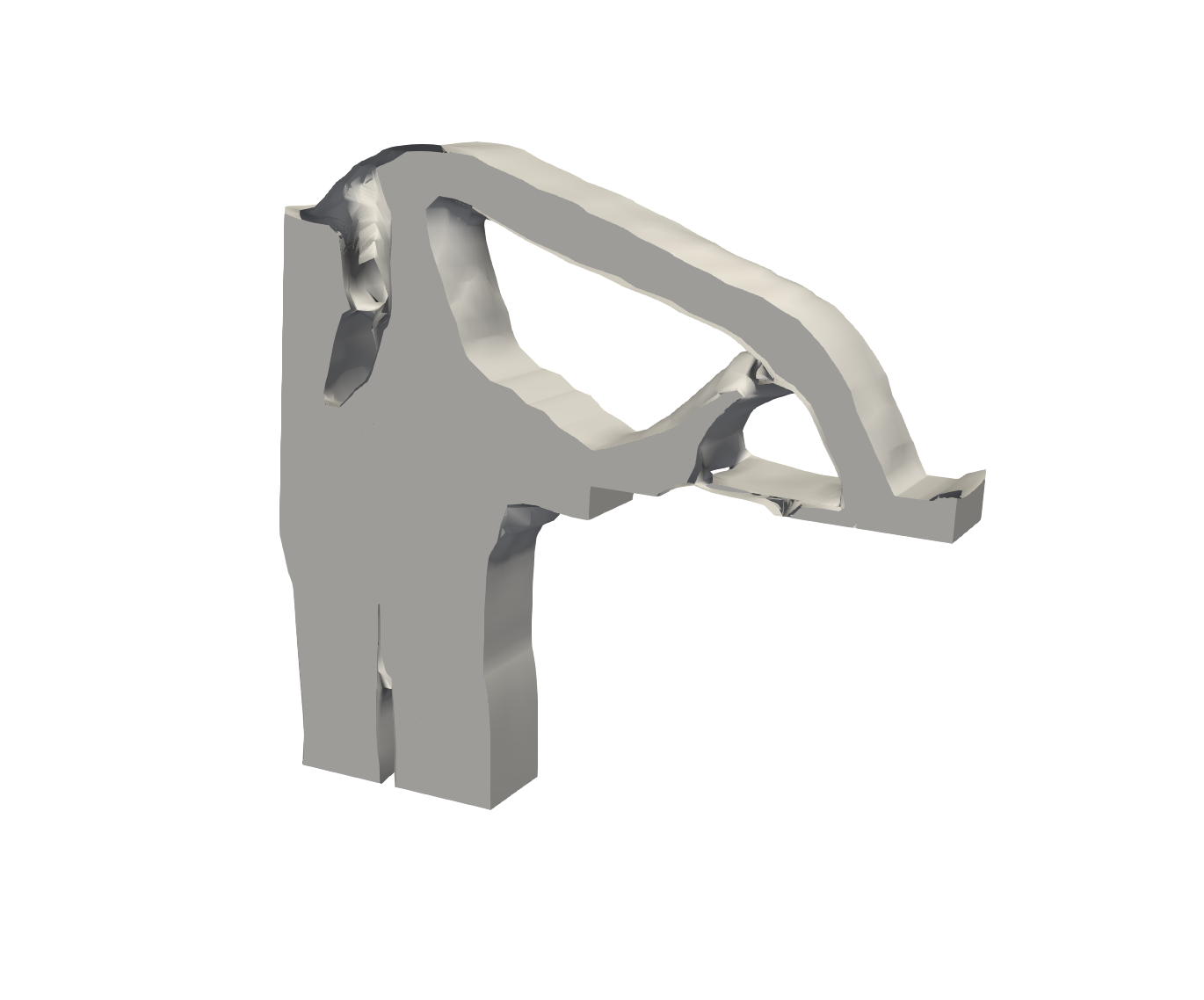}}   		\subfloat{\includegraphics[clip,trim=9.5cm 7.6cm 8cm 5.6cm, width=5.5cm]{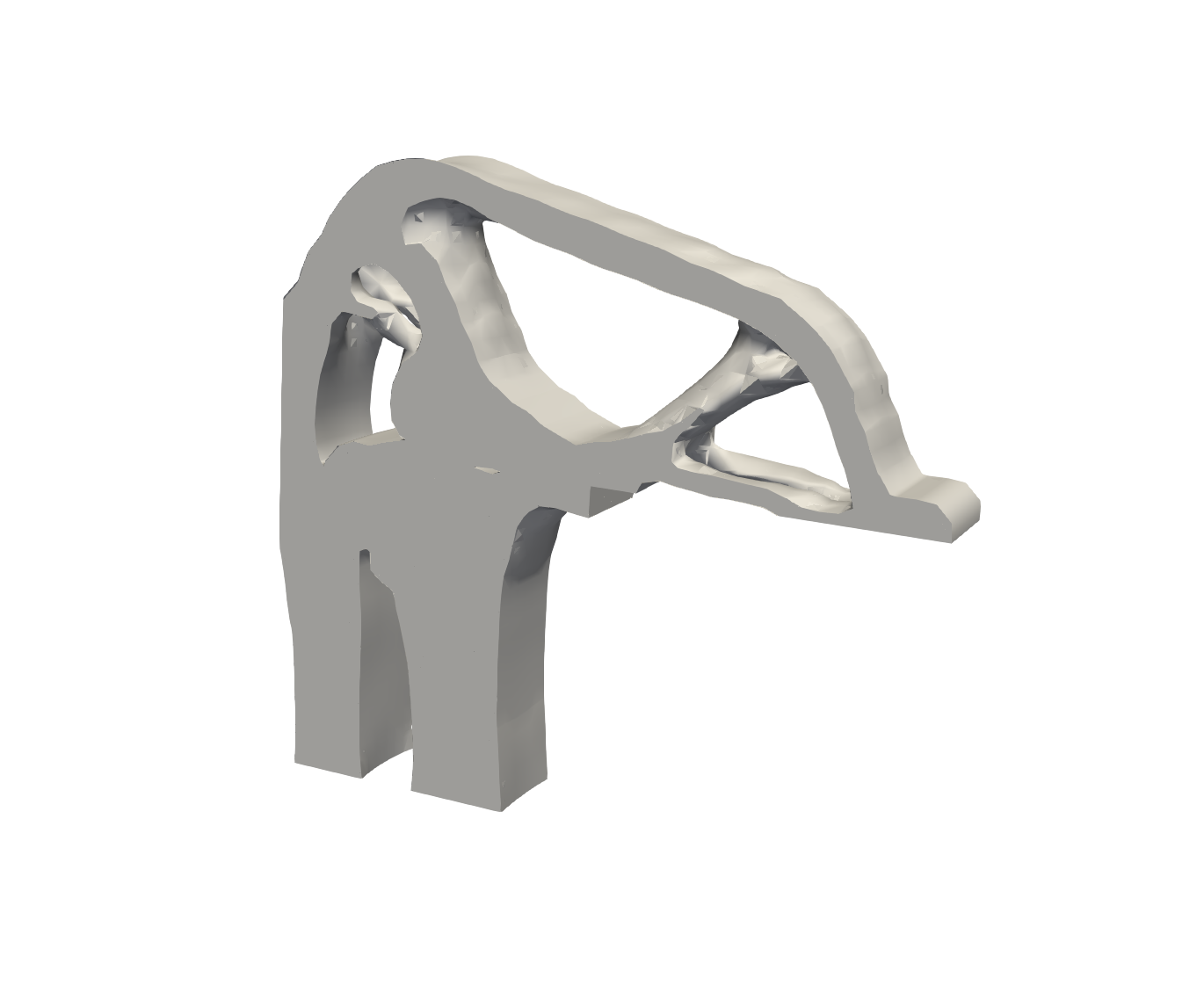}}   			\subfloat{\includegraphics[clip,trim=9.5cm 7.6cm 8cm 5.6cm, width=5.5cm]{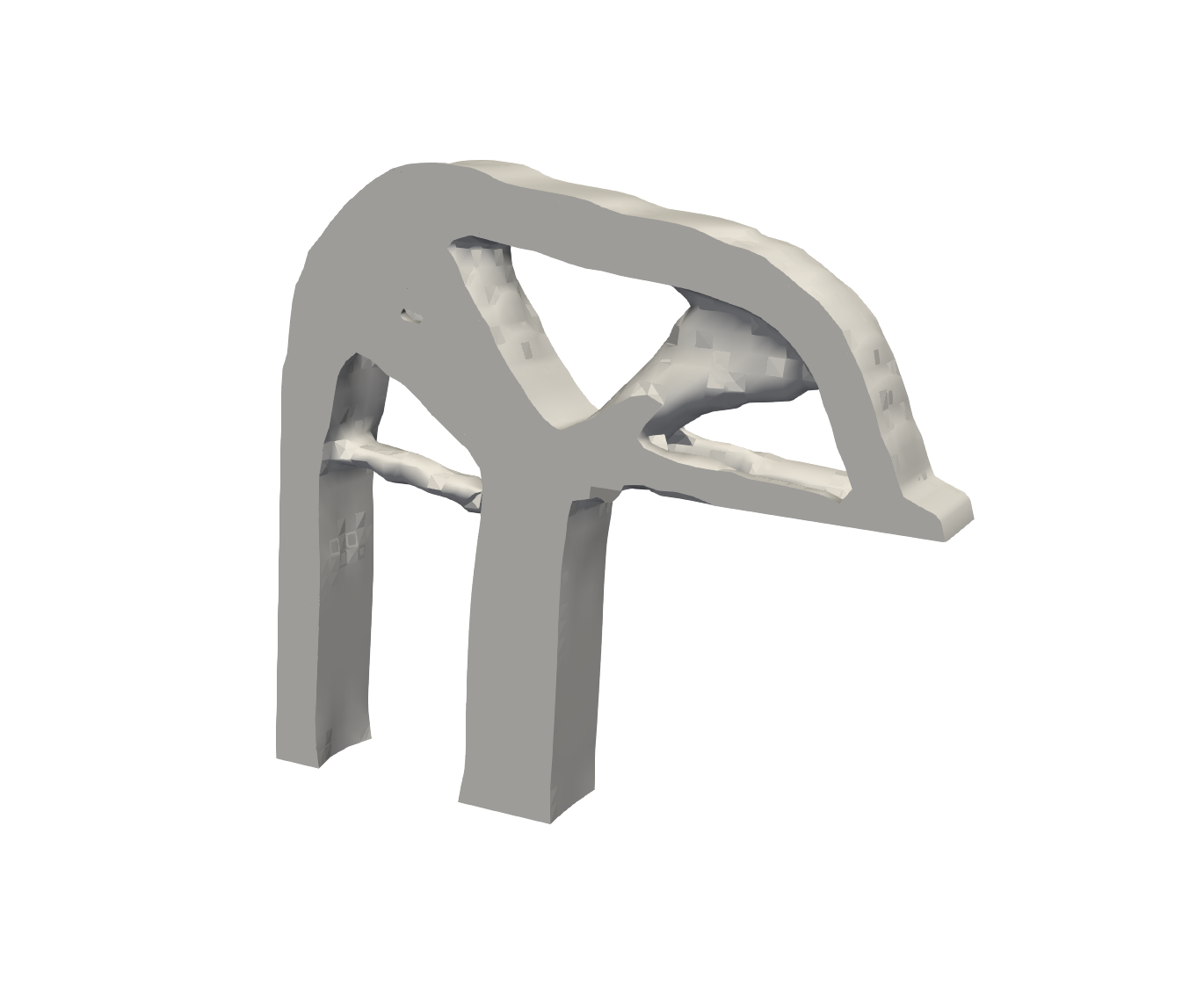}}
 	\caption*{\underline{(c) Fine mesh}}
 	\caption*{\hspace*{1cm}\underline{$\chi_v=0.52$}\hspace*{4.5cm}\underline{$\chi_v=0.42$}\hspace*{3.5cm}\underline{$\chi_v=0.4$}\hspace*{1.5cm}}
 	\vspace{-0.1cm}
 	\subfloat{\includegraphics[clip,trim=9.5cm 7.6cm 8cm 5.6cm, width=5.5cm]{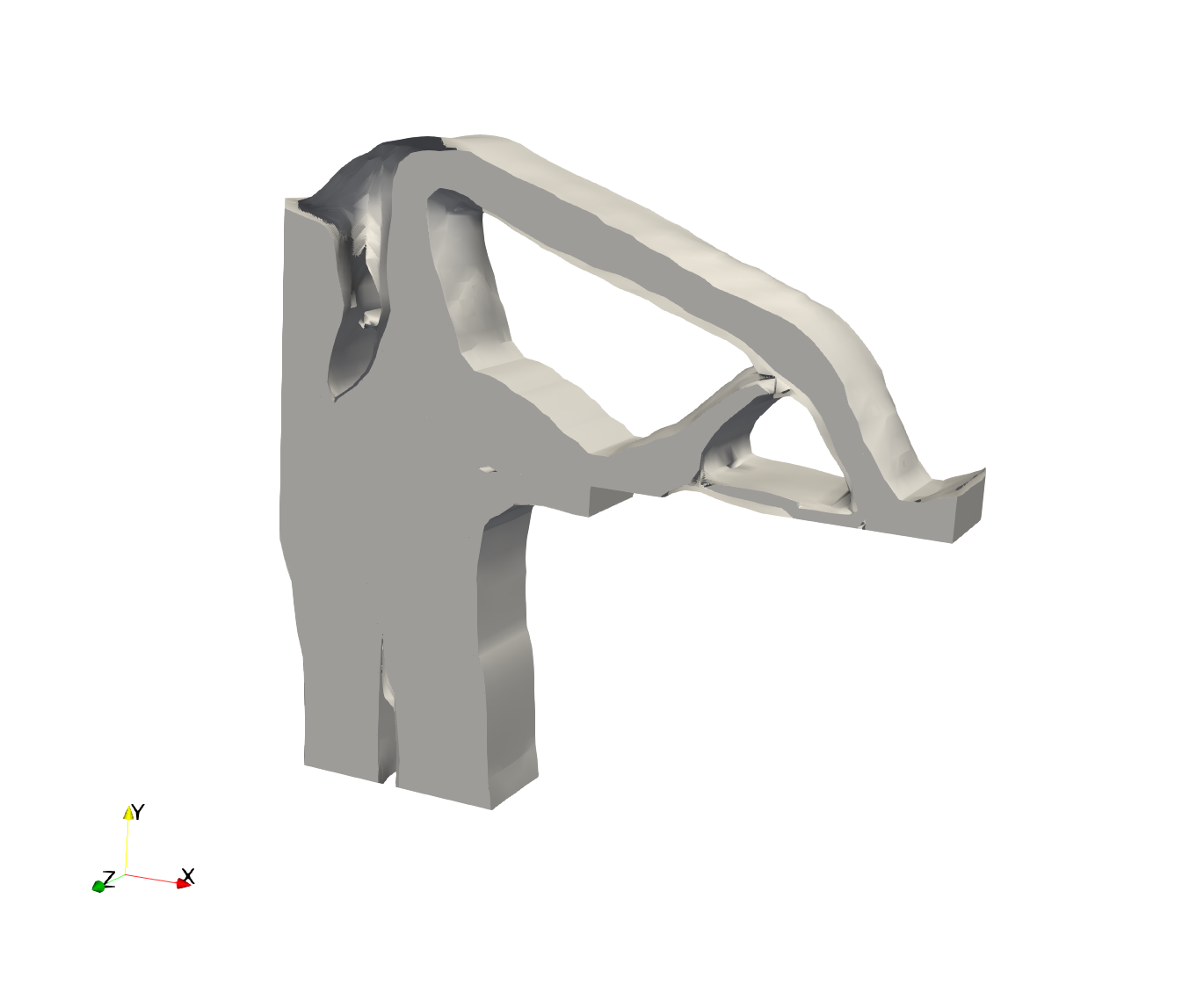}}   		\subfloat{\includegraphics[clip,trim=9.5cm 7.6cm 8cm 5.6cm, width=5.5cm]{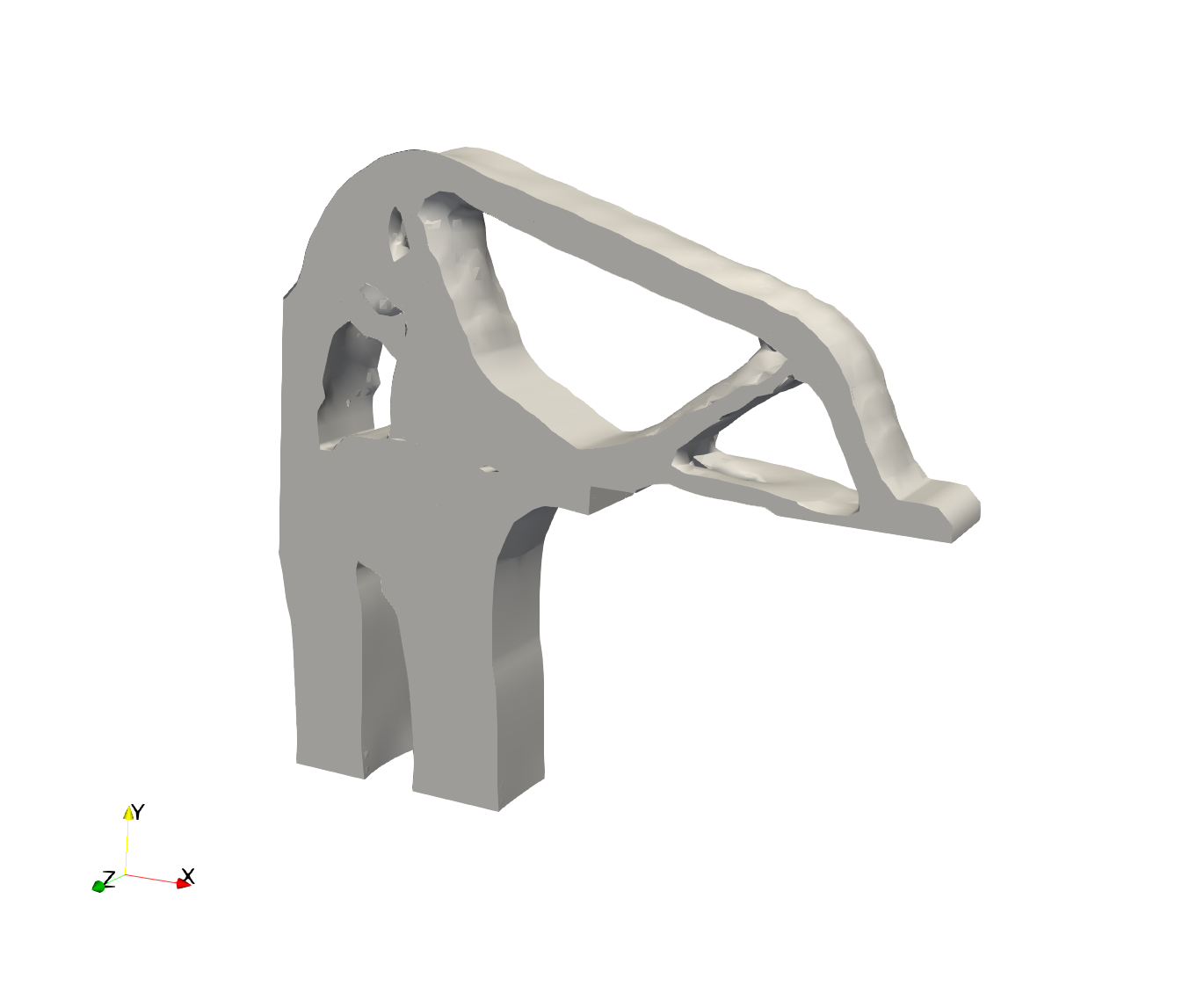}}   			\subfloat{\includegraphics[clip,trim=9.5cm 7.6cm 8cm 5.6cm, width=5.5cm]{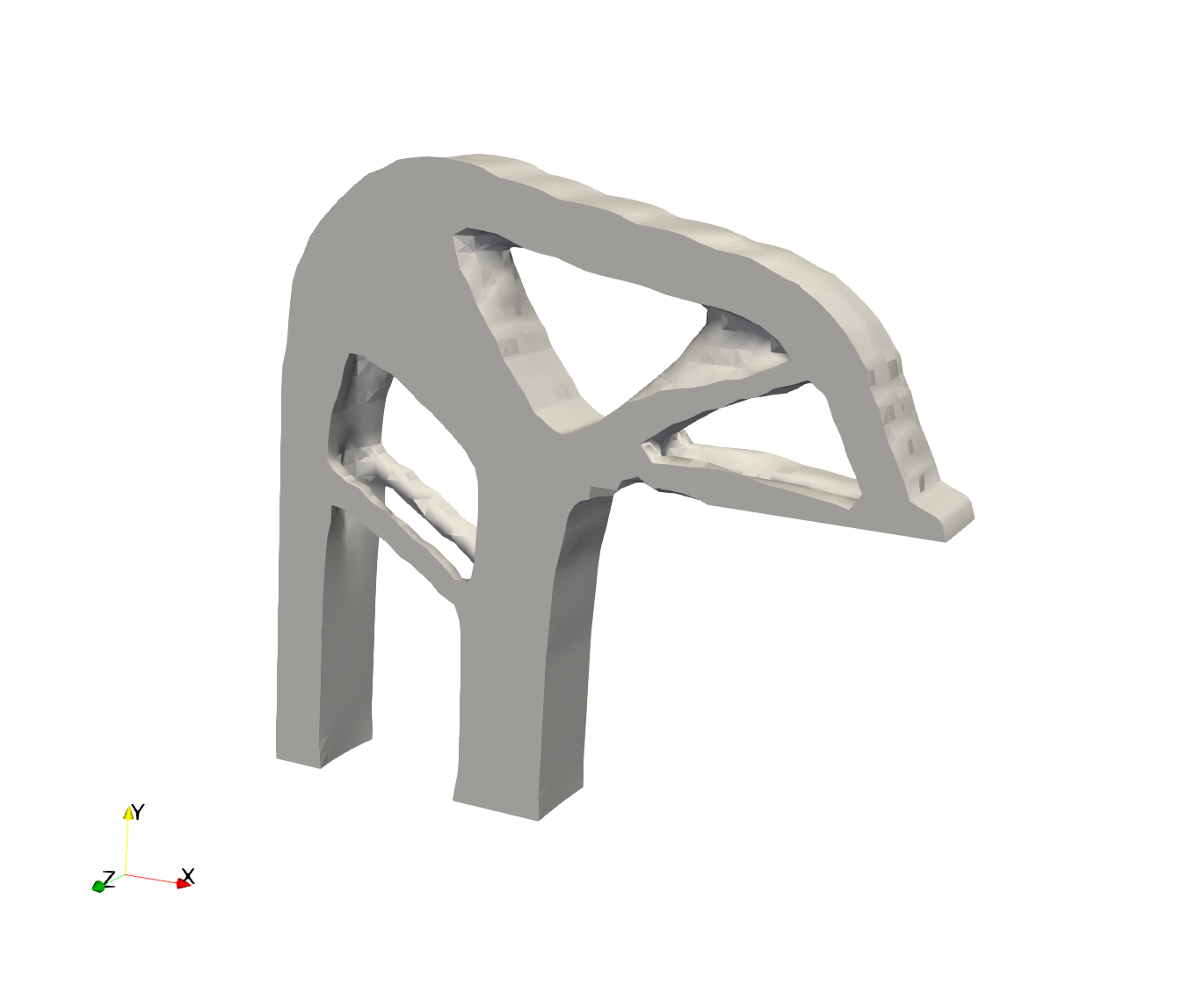}}
 	\caption{Example 2. Evolution history of the optimal layouts for different volume ratio of the L-shaped panel test for different discretization size.}
 	\label{Exm2_phi}
 \end{figure}

\begin{figure}[t!]
 	\caption*{\hspace*{1cm}\underline{$\chi_v=0.69$}\hspace*{4.5cm}\underline{$\chi_v=0.55$}\hspace*{3.5cm}\underline{$\chi_v=0.40$}\hspace*{1.5cm}}
 	\vspace{-0.1cm}
 	\subfloat{\includegraphics[clip,trim=9.5cm 7.6cm 7cm 4cm, width=5.5cm]{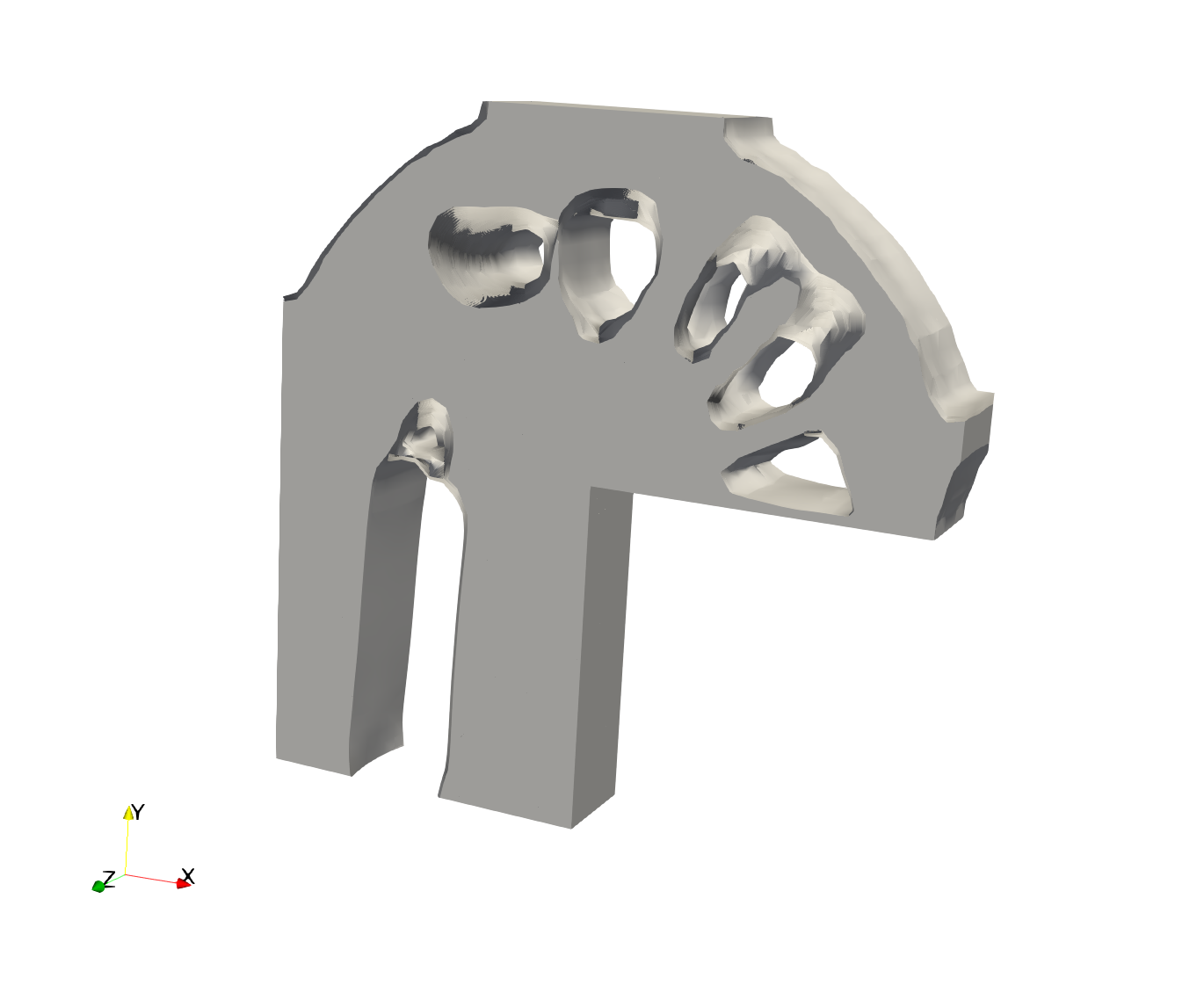}}   		\subfloat{\includegraphics[clip,trim=9.5cm 7.6cm 7cm 4cm, width=5.5cm]{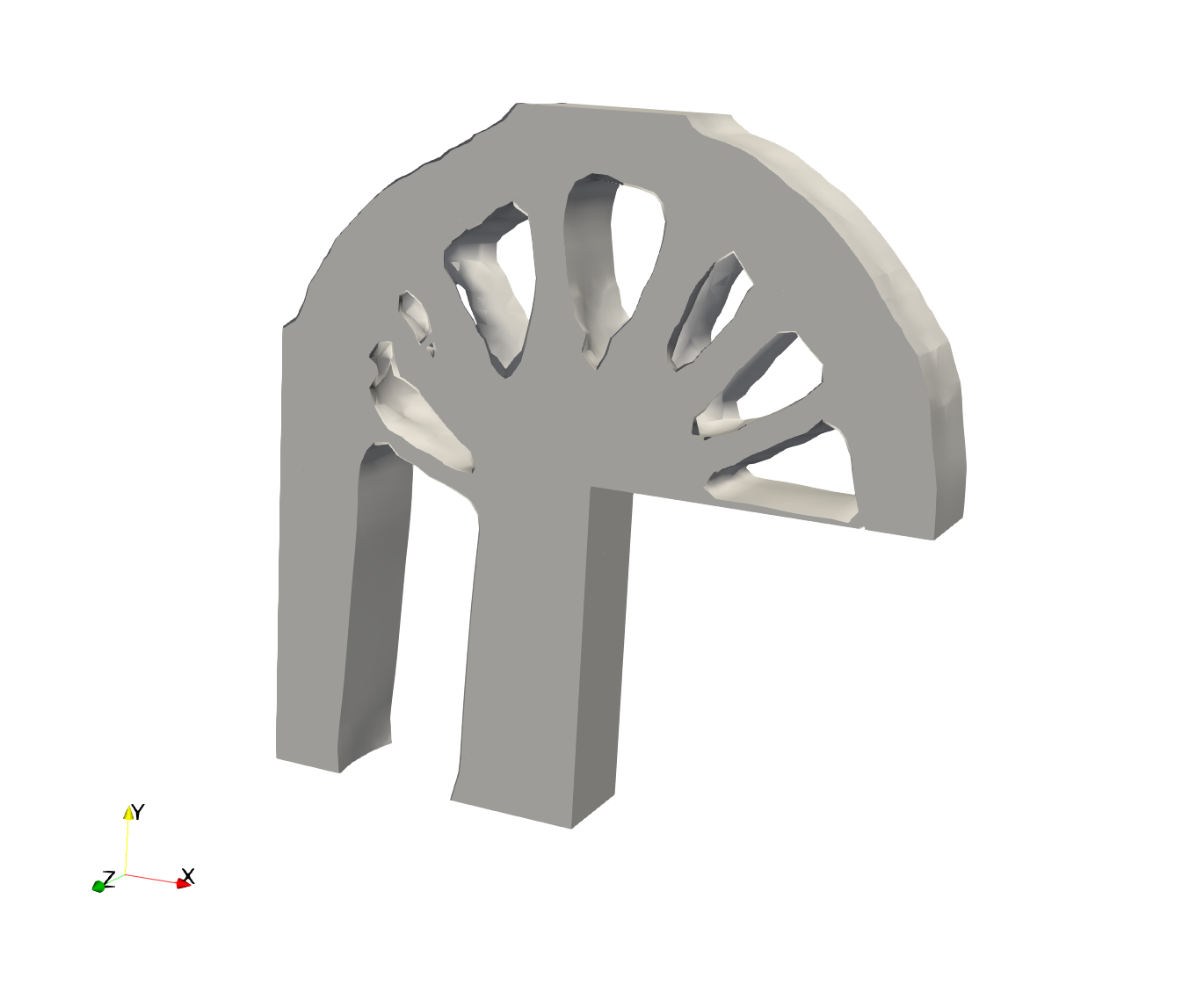}}   			\subfloat{\includegraphics[clip,trim=9.5cm 7.6cm 7cm 4cm, width=5.5cm]{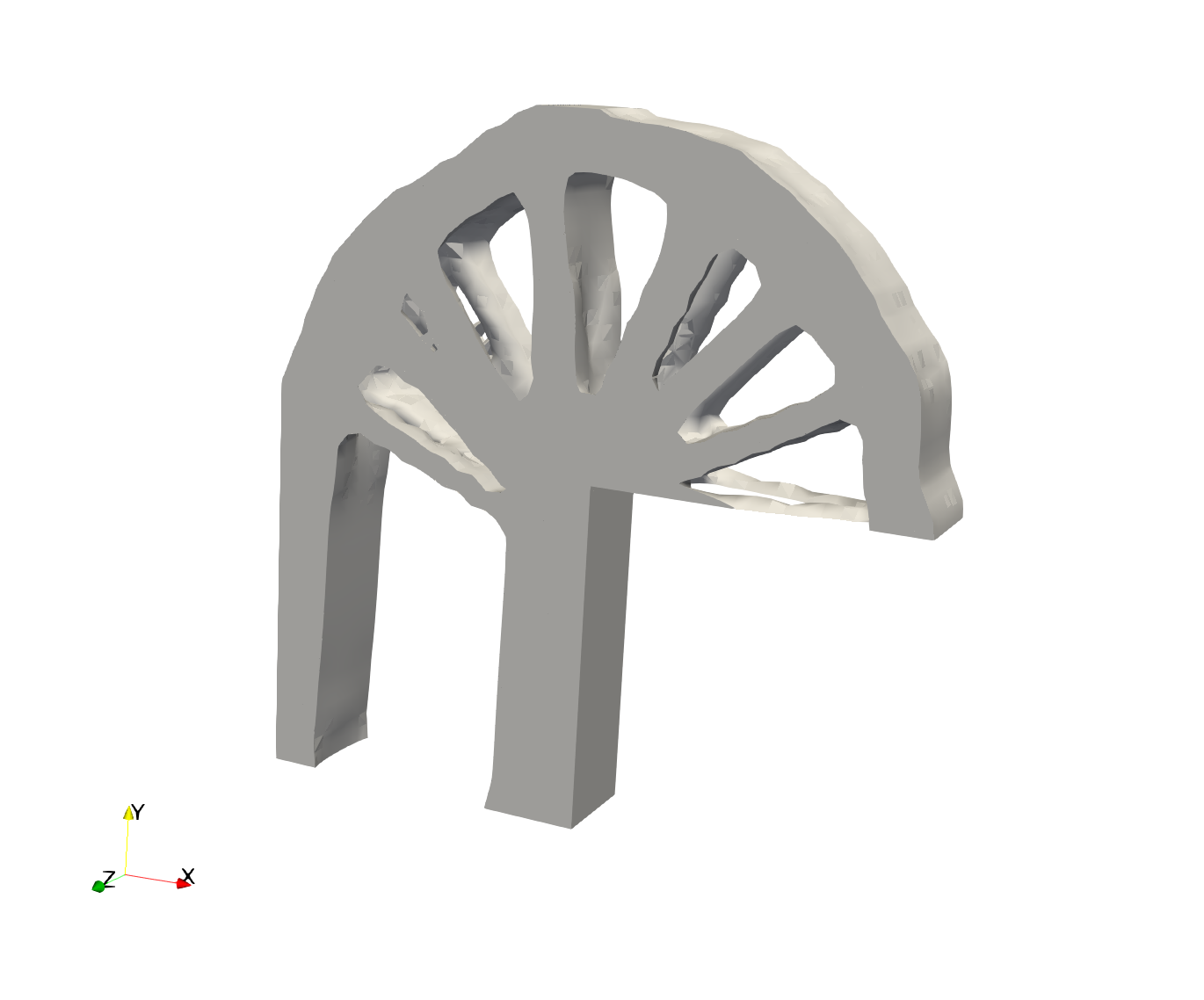}}	
 	\caption{Example 2. Evolution history of the optimal layouts for different volume ratio of the L-shaped panel test based on linear elasticity results.}
 	\label{Exm2_E}
 \end{figure}

We start with the presentation of the evolutionary history of the optimal layouts for different mesh discretization. A qualitative representation of the topological field is shown in Figure \ref{Exm2_phi}. Accordingly, for the sake of illustration, an optimum topology for different volume ratios of the elasticity results is depicted in Figure \req{Exm2_E}. It is clear that for the coarse mesh the obtained results show \grm{an} obvious difference between the final layouts, which indicates the mesh dependency in the method. However, the results obtained for medium and fine mesh are almost the same. It is worth noting that as the finite element discretization space is reduced, the final result improves significantly, to deal with crack propagation.
Another impacting factor that should be noted, in line with stress-constrained topology optimization problems, is that the optimum material distribution reduces stress singularity developed around the re-entrant corner, due to the fracture-resistant formulation, as can be seen in Figure \ref{Exm2_phi}. Indeed, this highlights the role of imposing \grm{an} additional constraint (i.e., phase-field equilibrium equation) on the topology optimization problem, while this can not be seen in the reference point (nearby the re-entrant corner point) of the L-shaped for the elasticity results depicted in Figure \ref{Exm2_E}. Another impacting factor that should be noted is that the sensitivity analysis which considers the phase-field fracture seeks to find the best pattern of material distribution in the design domain while reducing the fracture effects, simultaneously. This leads to an obvious difference in the normal velocity distribution obtained for the boundaries evolution compared to problems that ignore the effects of the phase-field fracture, especially in the re-entrant corner.

\begin{figure}[t!]
 	\vspace{-0.1cm}
 	\subfloat{\includegraphics[clip,trim=9cm 6.5cm 6.5cm 2.8cm, width=5.5cm]{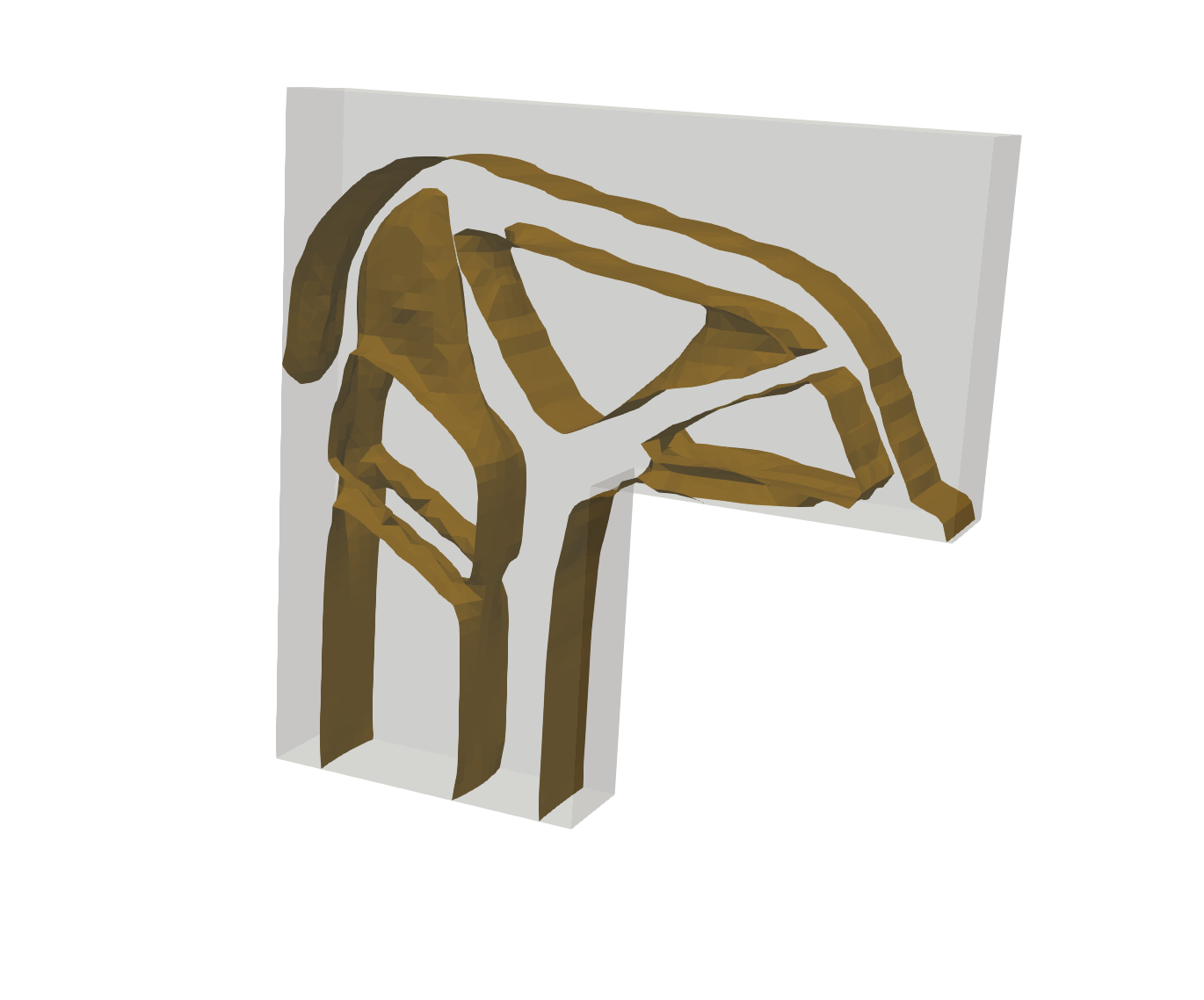}}   		\subfloat{\includegraphics[clip,trim=9cm 6.5cm 6.5cm 2.8cm, width=5.5cm]{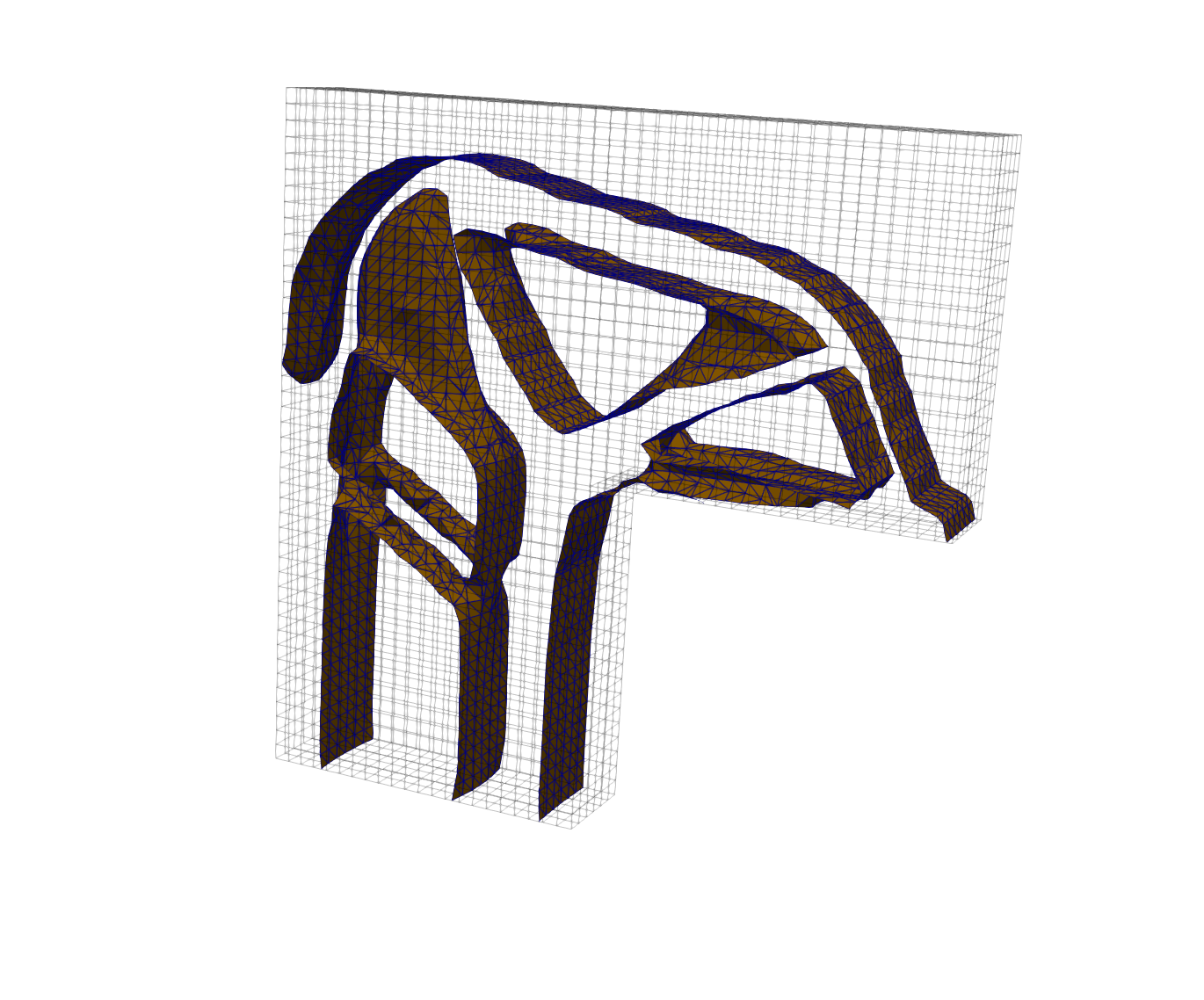}}   			\subfloat{\includegraphics[clip,trim=9cm 6.5cm 6.5cm 2.8cm, width=5.5cm]{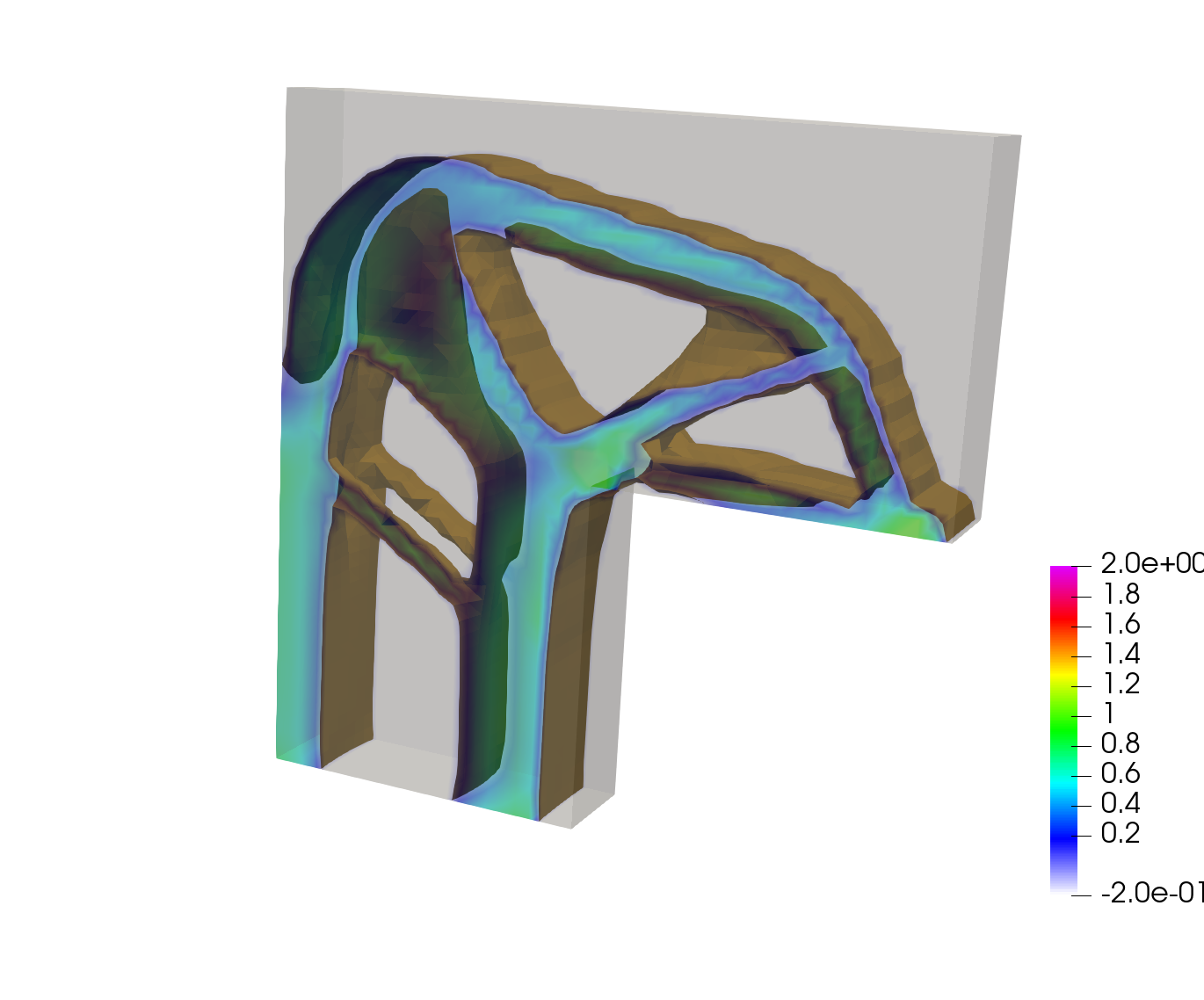}}
 	\vspace{-0.15cm}
 	\caption*{\hspace*{1.3cm}(a)\hspace*{5cm}{(b)}\hspace*{5cm}{(c)}}
 	\caption{Example 2. Level-Set representation for L-shaped panel test. Zero contour of level set surface $\Gamma_\Phi$ in (a) the continuum space, (b) the discretized space, and (c) computed topological field $\Phi$. We note that the gray area represents a non-material domain where $\Phi<0$.}
 	\label{Exm2_LSM}
 \end{figure}

The computed regularized topological field $\Phi(\Bx,t)$ by means of the reaction-diffusion evolution equation for fine mesh is outlined in Figure \ref{Exm2_LSM}. This includes zero contours of level set surface $\Gamma_{\Phi}$ in the continuous space, and in discretized space along with the evolved level set surface. We note that the gray color implies a non-material point, while the colorful area represents a material point. The load-displacement curve as a global indicator for the non-optimized, elasticity, and different mesh sizes are outlined in Figure \ref{Exm1_LD}(b). Evidently, it can be grasped that the maximum load-carrying capacity of the final layout is significantly increased compared to the non-optimized one, while this did not happen in elasticity when employing a coarse mesh. Thus, choosing an appropriate mesh discretization is essential. 
\begin{figure}[!t]
	\centering
	{\includegraphics[clip,trim=0cm 28cm 0cm 0cm, width=17cm]{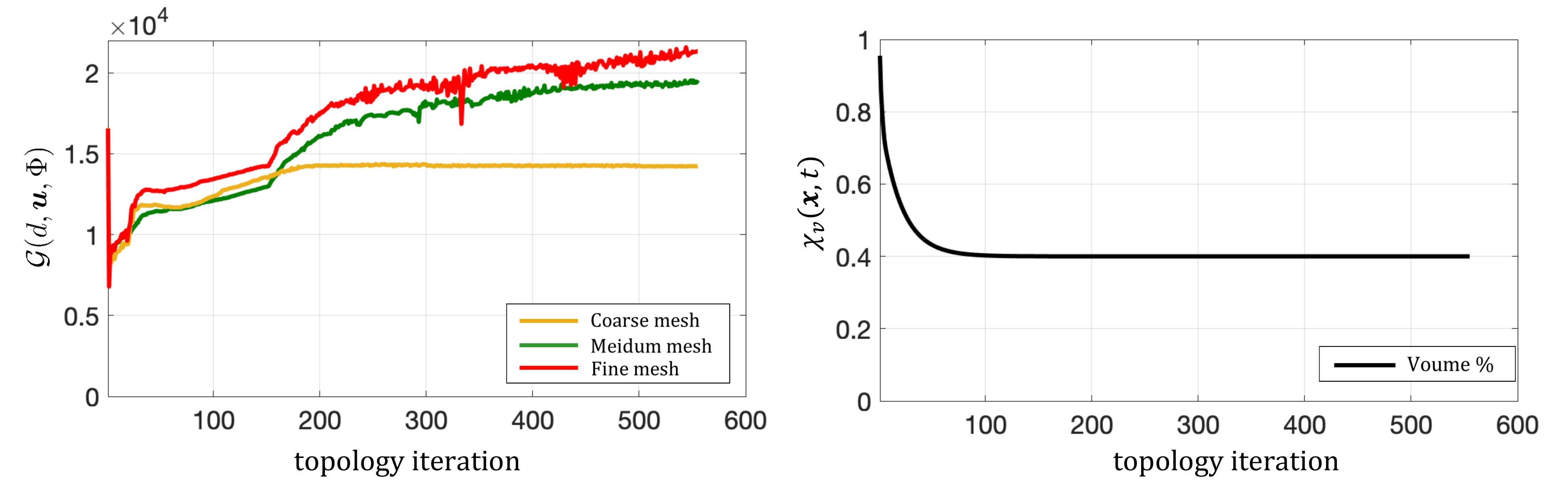}}  
	\vspace*{-0.7cm}
	\caption*{\hspace*{4.3cm}(a)\hspace*{8cm}(b)\hspace*{2cm}}
	\caption{Example 2. Convergence history for L-shaped panel test through (a) objective function, and (b) volume constraint function.}
	\label{Exm2_conv}
\end{figure}
 \begin{figure}[t!]
 	\centering
 	\vspace{-0.1cm}
 	\includegraphics[clip,trim=0cm 25cm 0cm 0cm, width=17cm]{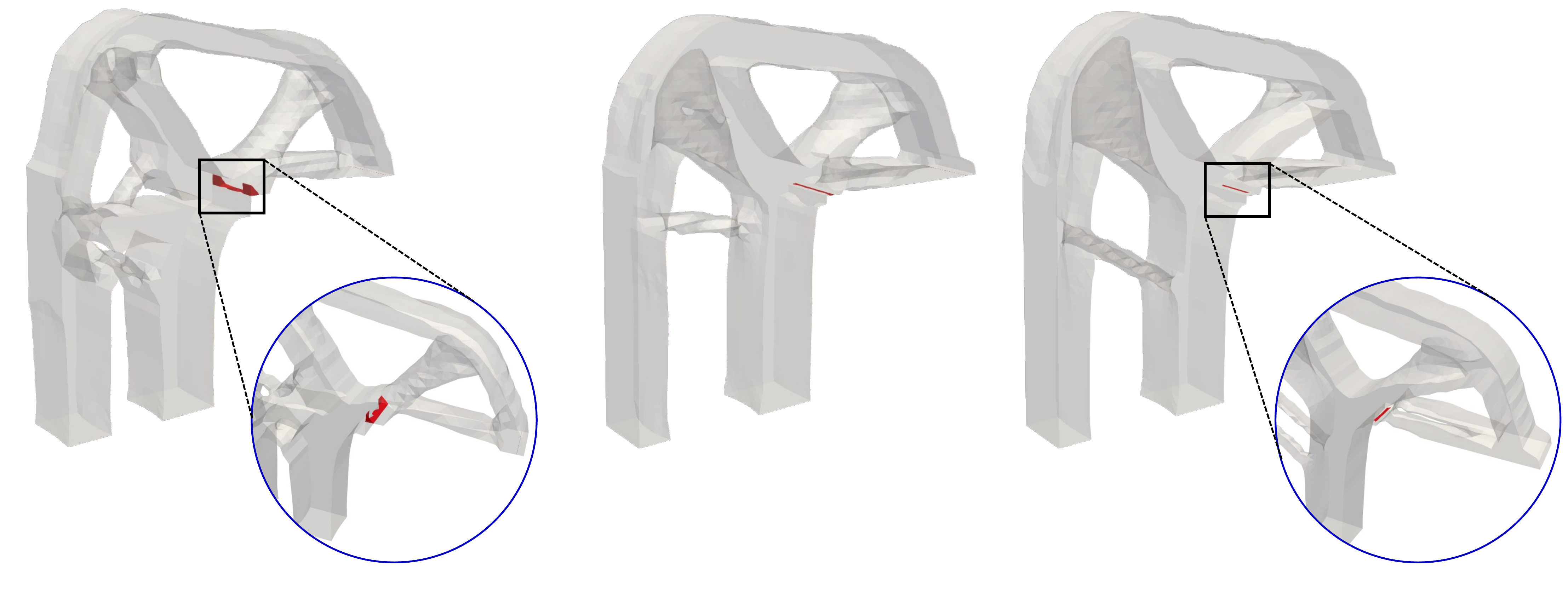}  		
 	\caption*{\hspace*{1.3cm}(a)\hspace*{5cm}{(b)}\hspace*{5cm}{(c)}}
 	\caption{Example 2. The crack phase-field profile for different discretization size in (a) coarse mesh, (b) medium mesh, and (c) fine mesh based on Formulation 2.}
 	\label{Exm2_phasefield_1}
 \end{figure}  

Additionally, it can be observed that the fracture limit point regardless of the mesh size, is $50.2\%$ more than the non-optimized case, while the elasticity result requires almost the same displacement load to observe the onset of fracture. This shows the great importance of considering fracture constraint in our topology minimization approach, and thus, it has successfully delayed fracture initiation.

\begin{figure}[t!]
 	\centering
 	\vspace{-0.1cm}
 	 \hspace*{1cm} 	
 	\subfloat{\includegraphics[clip,trim=10cm 3cm 7.7cm 5cm, width=5cm]{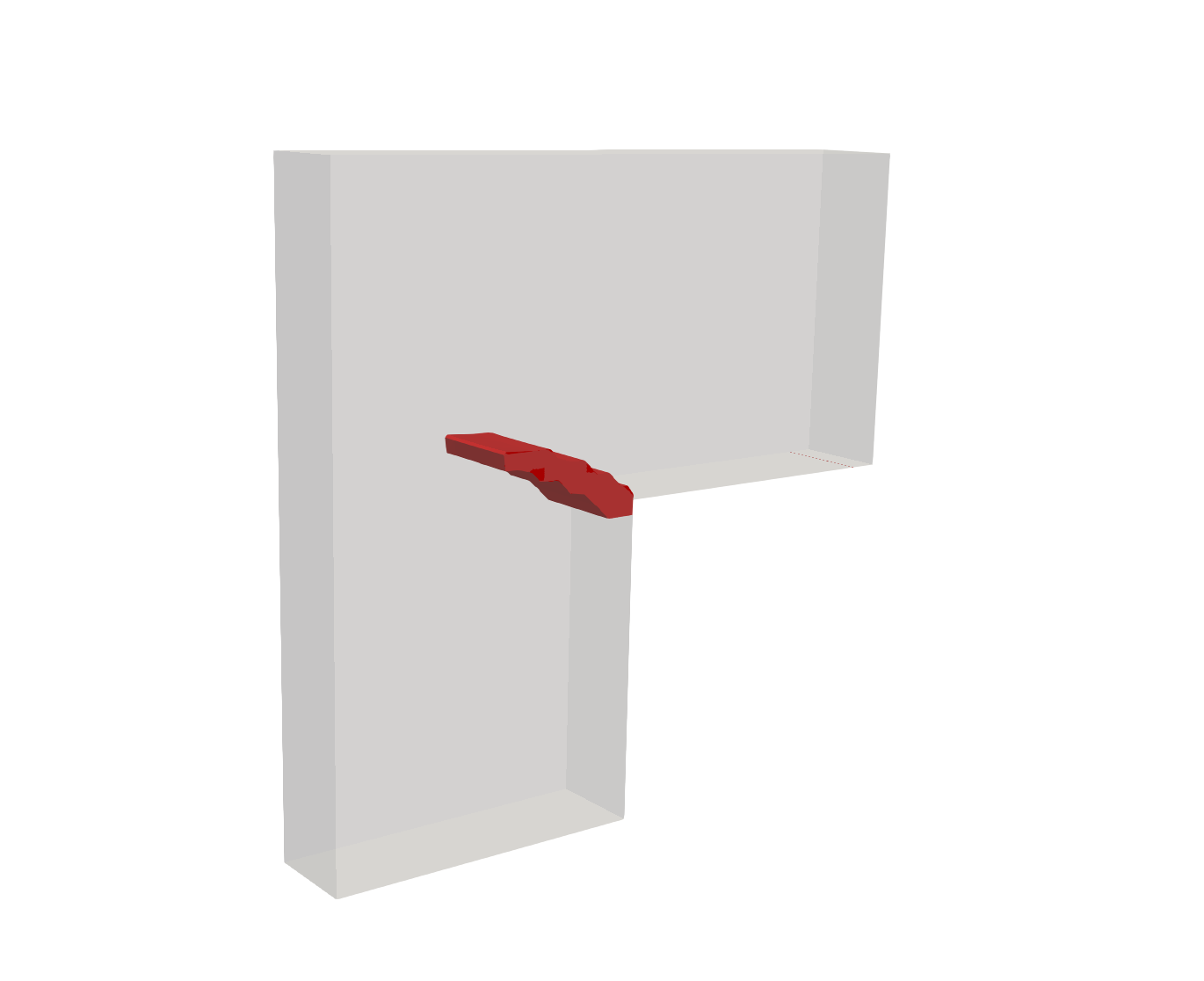}}	
 	 \hspace*{2cm} 
 	\subfloat{\includegraphics[clip,trim=10cm 3cm 7.7cm 5cm, width=5cm]{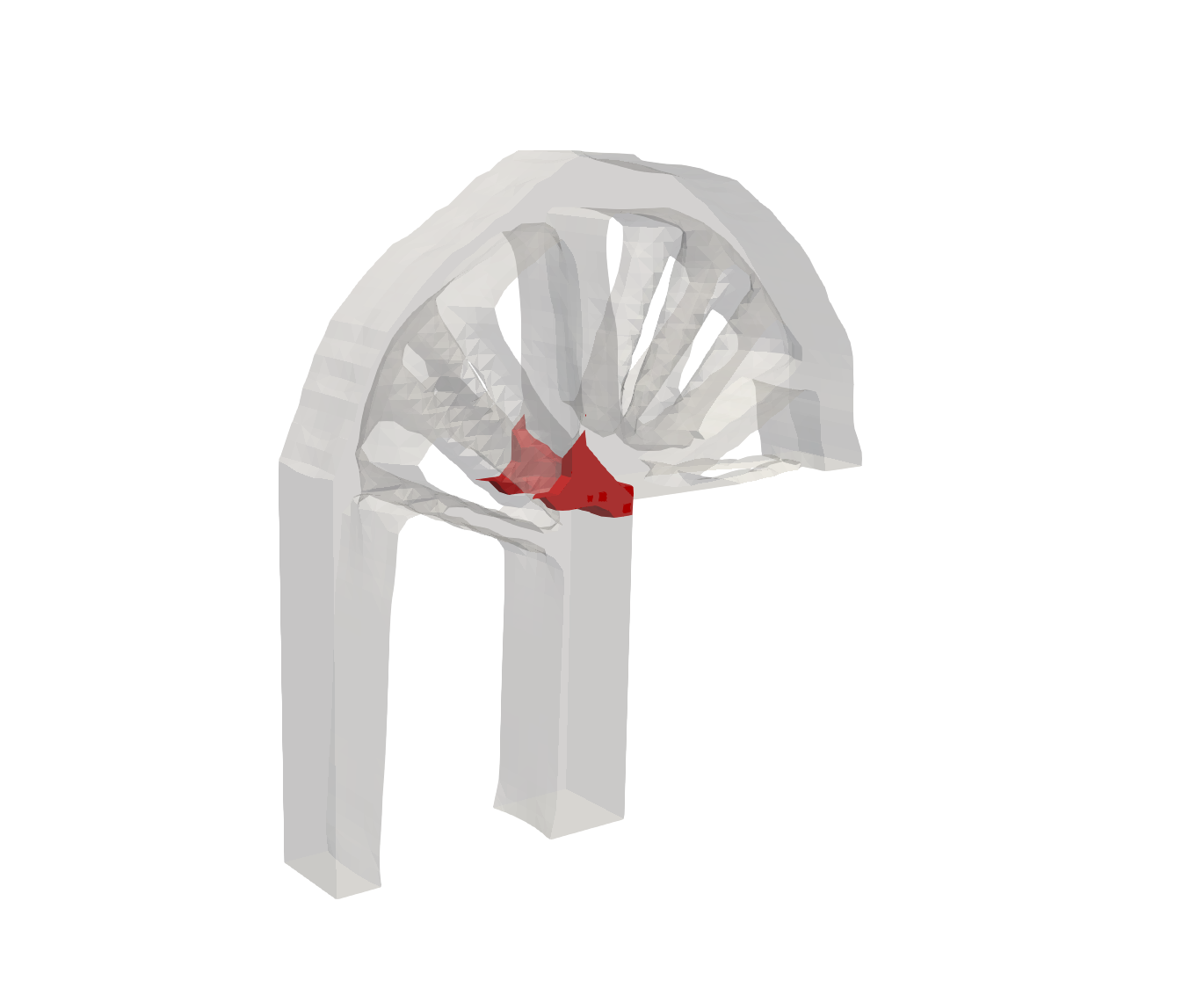}} 	
 	\caption*{\hspace*{1cm}(a)\hspace*{7cm}(b)\hspace*{1cm}}
 	\caption{Example 2. The crack phase-field profile based on (a) non-optimized, and (b) linear elasticity results undergoes brittle fracture.}
 	\label{Exm2_phi_full_E}
 \end{figure}

Another observation is the convergence history of the objective and the volume constraint functions which are depicted in Figures \ref{Exm2_conv}(a-b), respectively. It can be seen that with a finer mesh, the optimal value of the objective function is increased, thus finer mesh exhibits a stiffer response in its final optimum layout (since the critical stress state of the onset of fracture has an inverse relation to the discretization size, see \cite{borden2012phase}). Subsequently, the crack phase-field profiles for different discretization spaces are illustrated in Figure \ref{Exm2_phasefield_1}, while non-optimized and elasticity results are shown in Figure \ref{Exm2_phi_full_E}. The first important observation is the inadequacy of this topology obtained from both non-optimized and elasticity results against the fractured state. Unlikely, topology optimization by considering the fracture constraint regardless of the mesh size exhibits very limited cracked region, thus superior efficiency is observed.

Finally, the discretized final optimum layouts are depicted in Figure \ref{Exm23_mesh}(a). It is trivial that the smooth topological interface can be remarkably depicted.

\sectpb[Section53]{Example 3: Portal frame structure under compression loading}
\begin{figure}[!t]
	\centering
	{\includegraphics[clip,trim=1cm 24cm 0cm 0cm, width=16cm]{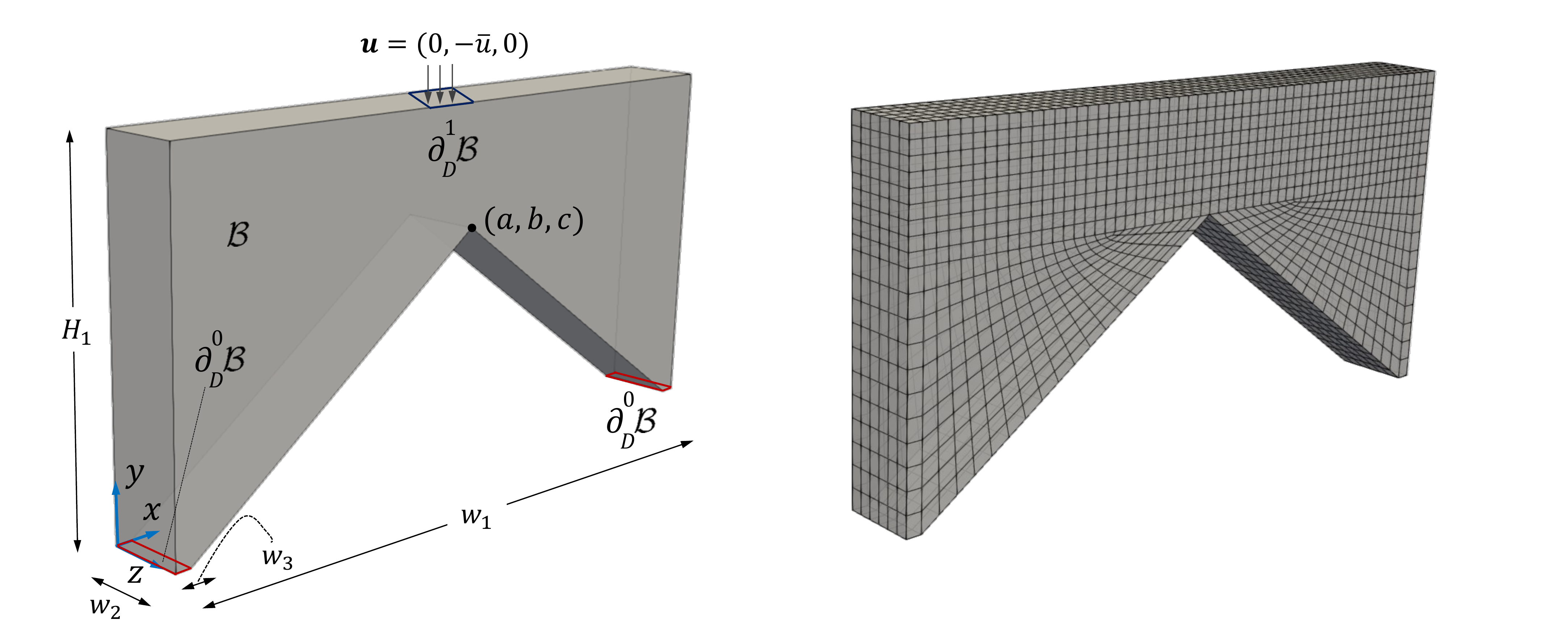}}  
	\vspace*{-0.7cm}
	\caption*{\hspace*{4.3cm}(a)\hspace*{8cm}(b)\hspace*{2cm}}
	\caption{Example 3. The representation of the (a) geometry and boundary conditions, and (b) finite element discretization.}
	\label{Exm3_bvp}
\end{figure}

\begin{figure}[t!]
	\caption*{\hspace*{1cm}\underline{$\chi_v=0.88$}\hspace*{4.5cm}\underline{$\chi_v=0.63$}\hspace*{3.5cm}\underline{$\chi_v=0.40$}\hspace*{1.5cm}}
	\vspace{-0.1cm}
	\subfloat{\includegraphics[clip,trim=6cm 6cm 6cm 10cm, width=5.6cm]{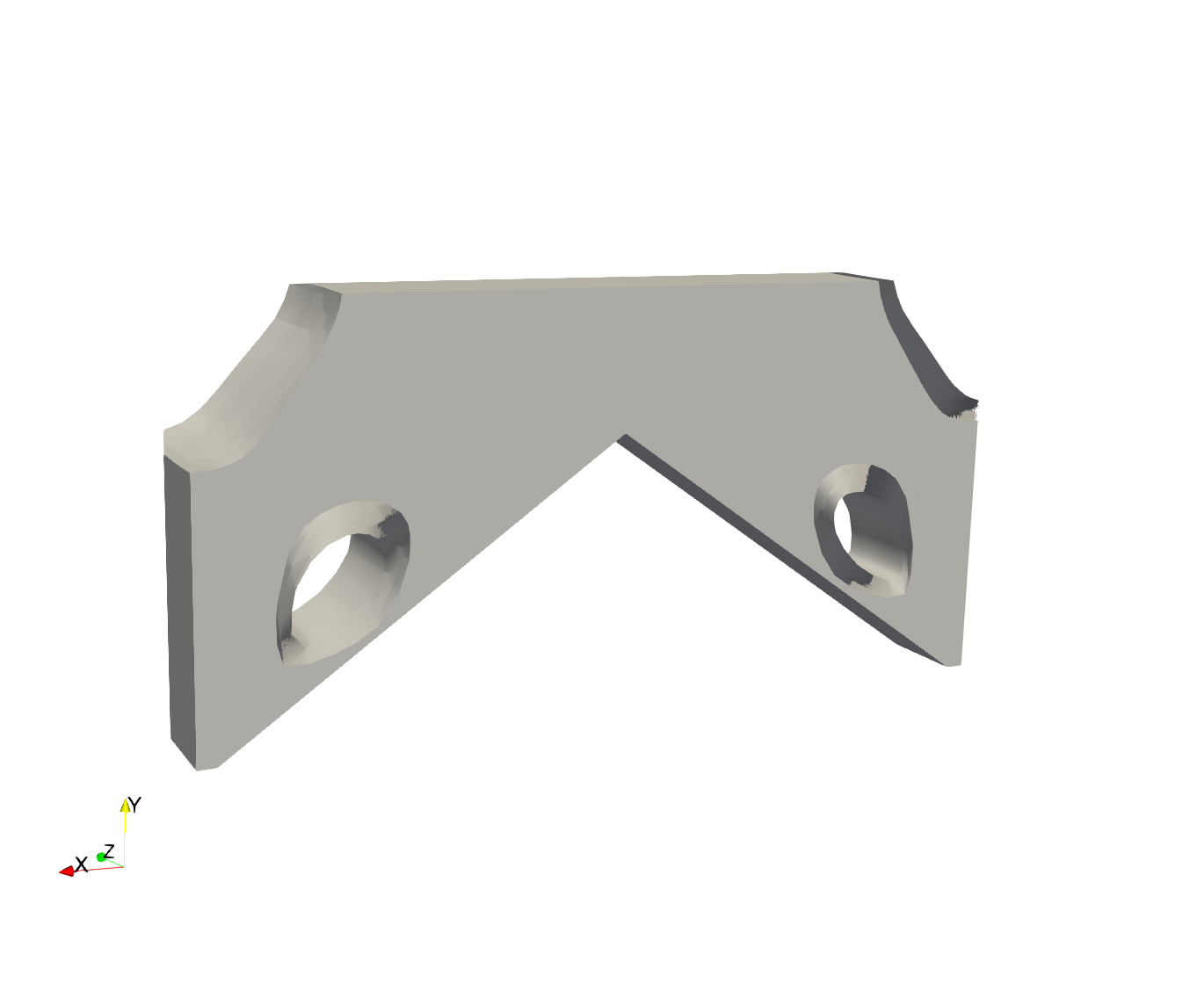}}   		\subfloat{\includegraphics[clip,trim=6cm 6cm 6cm 10cm, width=5.6cm]{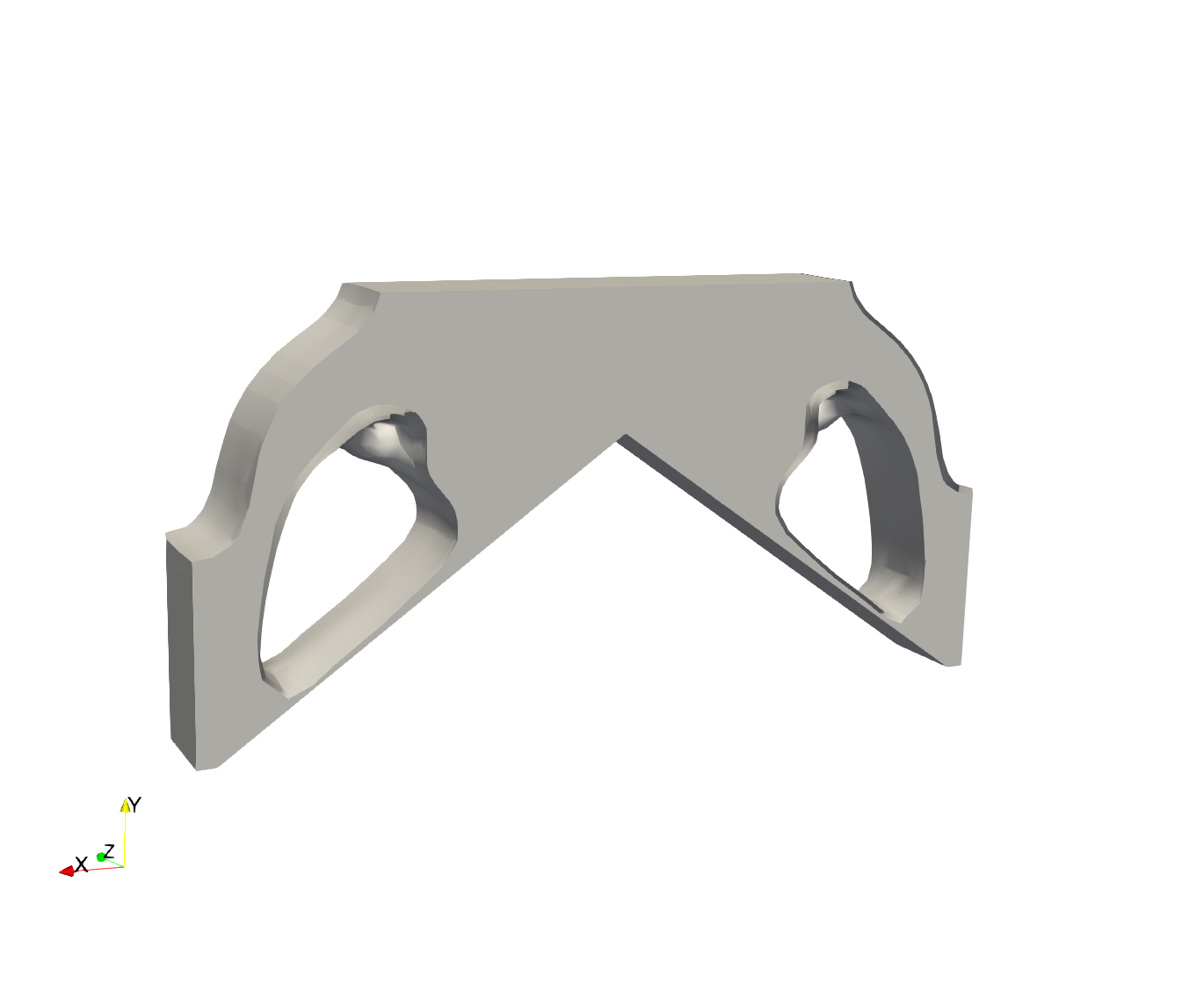}}   			\subfloat{\includegraphics[clip,trim=6cm 6cm 6cm 10cm, width=5.6cm]{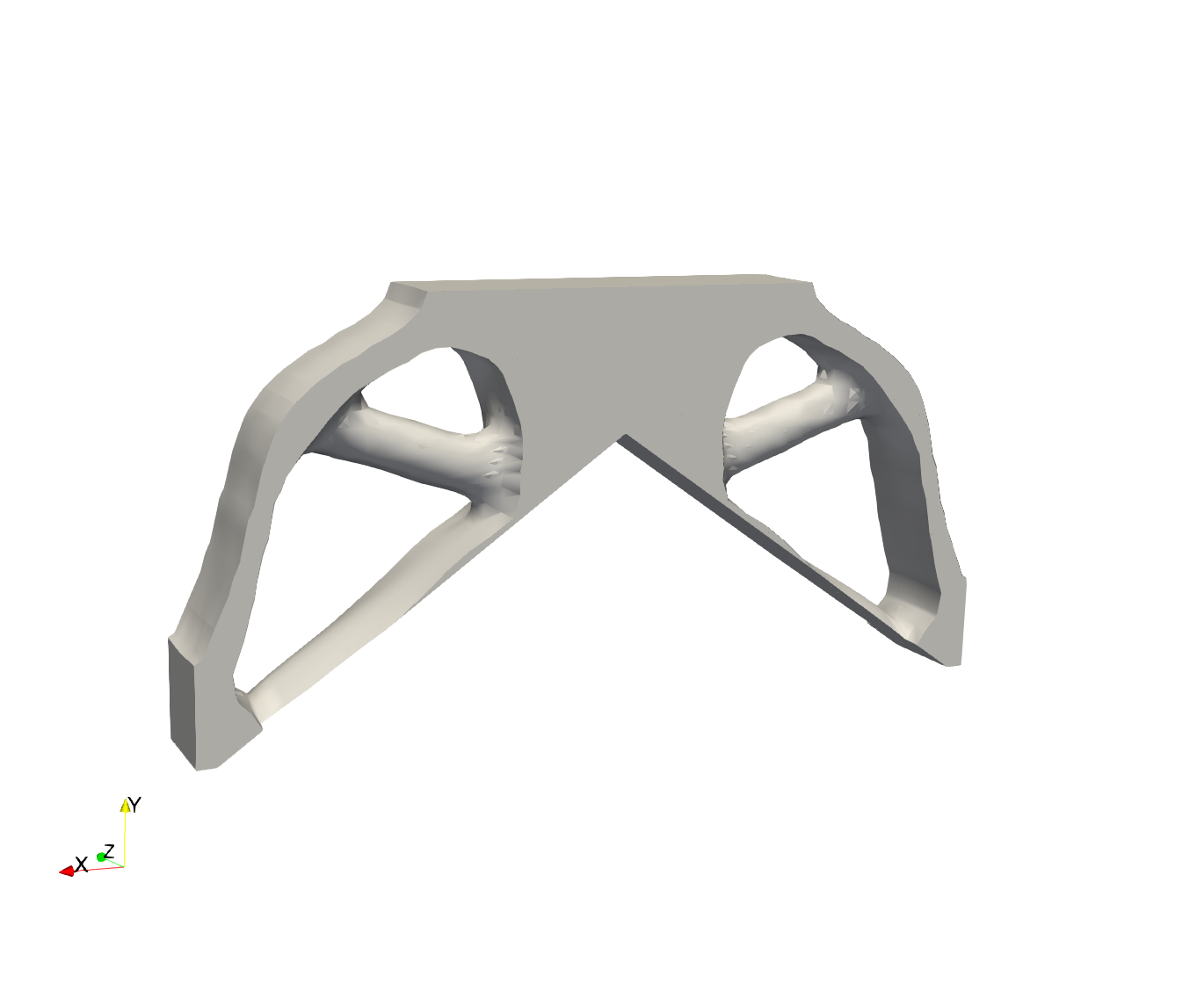}}	
	\caption{Example 3. Evolution history of the optimal layouts for different volume ratio of the Portal frame structure.}
	\label{Exm3_phi}
\end{figure}  
The last two examples concern the performance of the proposed topology optimization in the case of ductile failure. In this example, a three-dimensional portal frame structure shown in Figure \ref{Exm3_bvp}(a) is investigated, which is subject to a prescribed downward displacement distributed over a narrow area at the top of the frame. This benchmark example has been reported in many contributions \cite{Allaire1,Russ}. The geometrical configuration is shown in Figure \ref{Exm3_bvp}(a). The edge on its leftmost support is constrained in all directions while the edge on its rightmost is fixed for displacement in $y$ direction. The geometrical dimensions shown in Figure \ref{Exm3_bvp}(a) are set as $w_1=4.5\;mm$, $w_2=0.5\;mm$, $w_3=0.1\;mm$, and $H_1=2.25\;mm$. The reference point \grm{is} located in $(a,b,c)=(2.25,1.5,0.5)\;mm$. 

A monotonic displacement increment  ${\Delta \bar{u}}_y=-2.5\times^{-4}\;mm$ is prescribed downward in a vertical direction in a part of the top page $area=((2.15,2.25,0.5)\times (2.35,2.25,0)$ of the specimen for 150-time steps. 
Thus, the final prescribed displacement load for this optimization problem is set as ${\bar{u}}_y=-0.037\;mm$. The minimum finite element size in the solid domains is selected to be $h_{min}=0.09\;mm$ which implies 5520 trilinear hexahedral elements.

\begin{figure}[t!]
	\caption*{\underline{$\chi_v=1$}\hspace*{6cm}\underline{$\chi_v=0.75$}\hspace*{4cm}\underline{$\chi_v=0.62$}\hspace*{1cm}}
	\vspace{-0.1cm}
	\subfloat{\includegraphics[clip,trim=6cm 6cm 6cm 10cm, width=5.6cm]{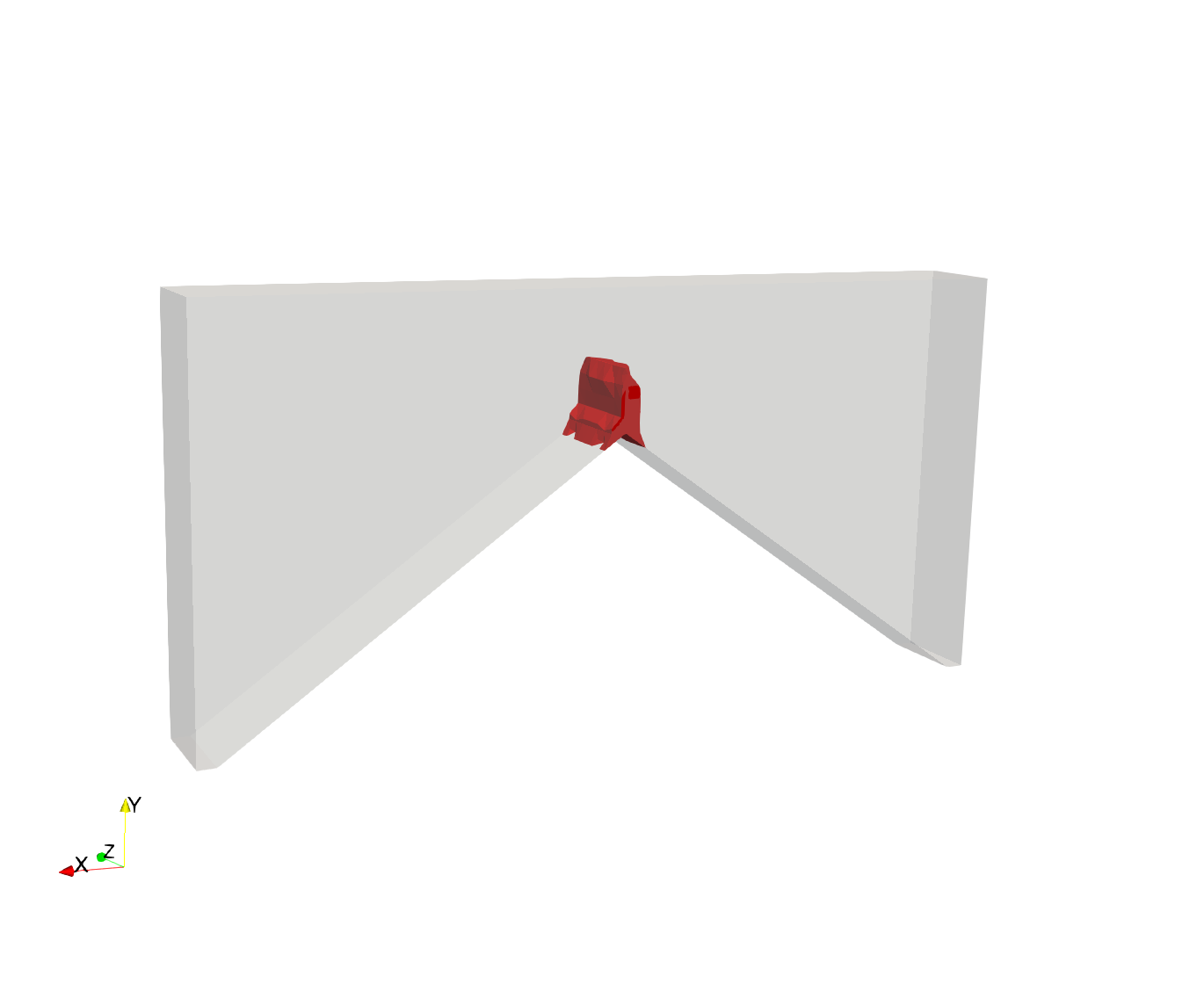}}   		\subfloat{\includegraphics[clip,trim=6cm 6cm 6cm 10cm, width=5.6cm]{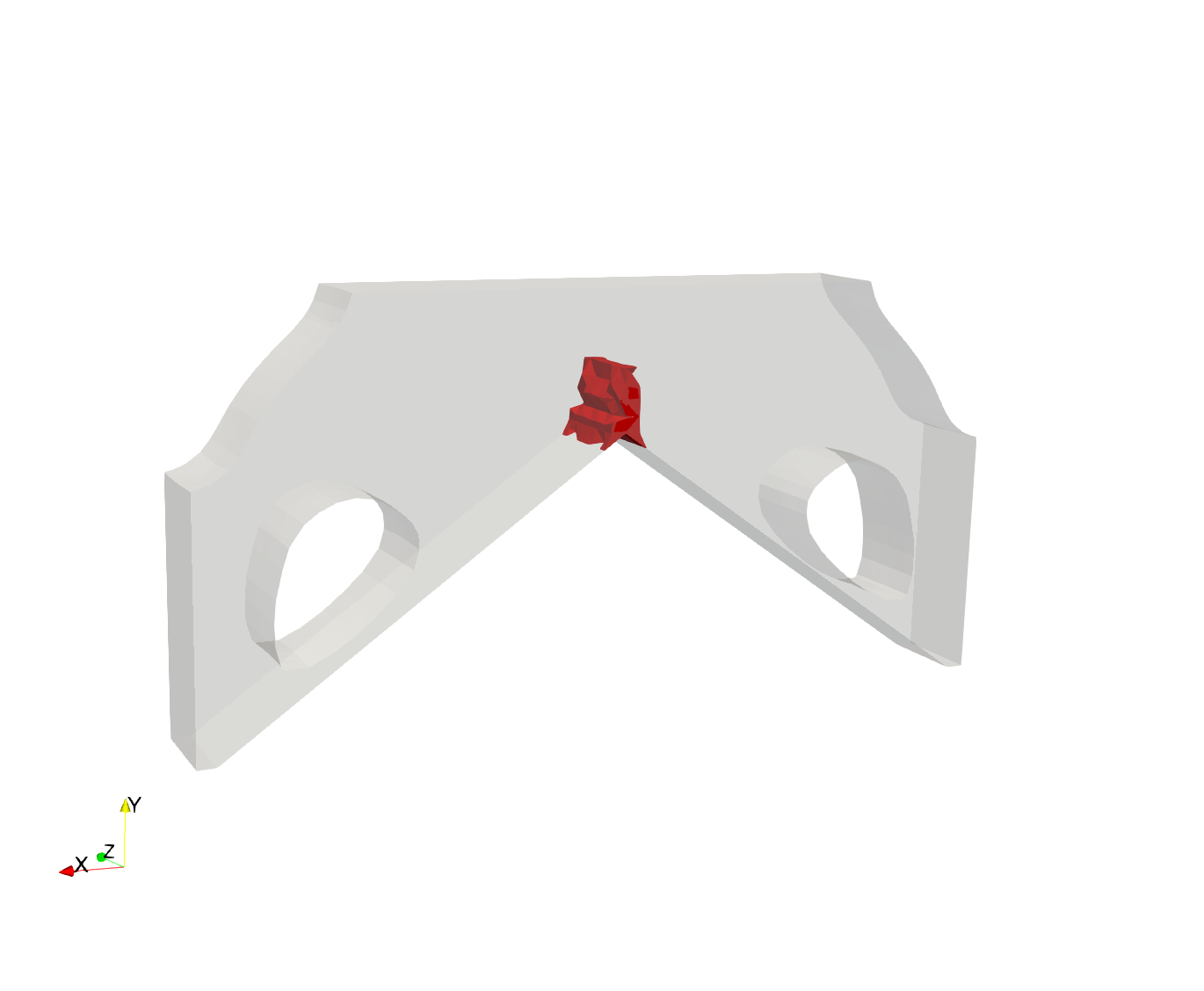}}   			\subfloat{\includegraphics[clip,trim=6cm 6cm 6cm 10cm, width=5.6cm]{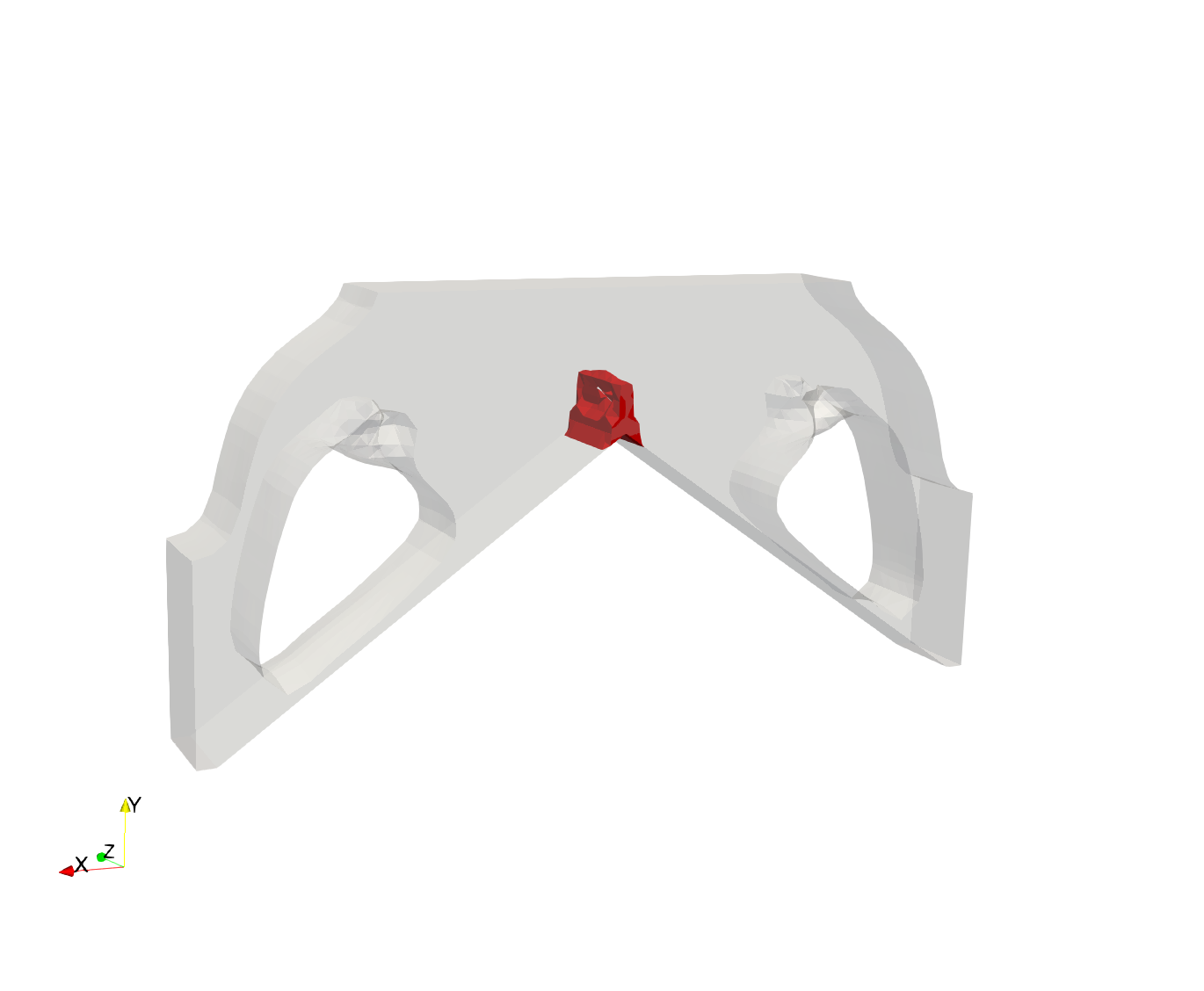}}	
	\caption*{\underline{$\chi_v=0.59$}\hspace*{6cm}\underline{$\chi_v=0.56$}\hspace*{4cm}\underline{$\chi_v=0.55$}\hspace*{1cm}}
	\vspace{-0.1cm}
	\subfloat{\includegraphics[clip,trim=6cm 6cm 6cm 10cm, width=5.6cm]{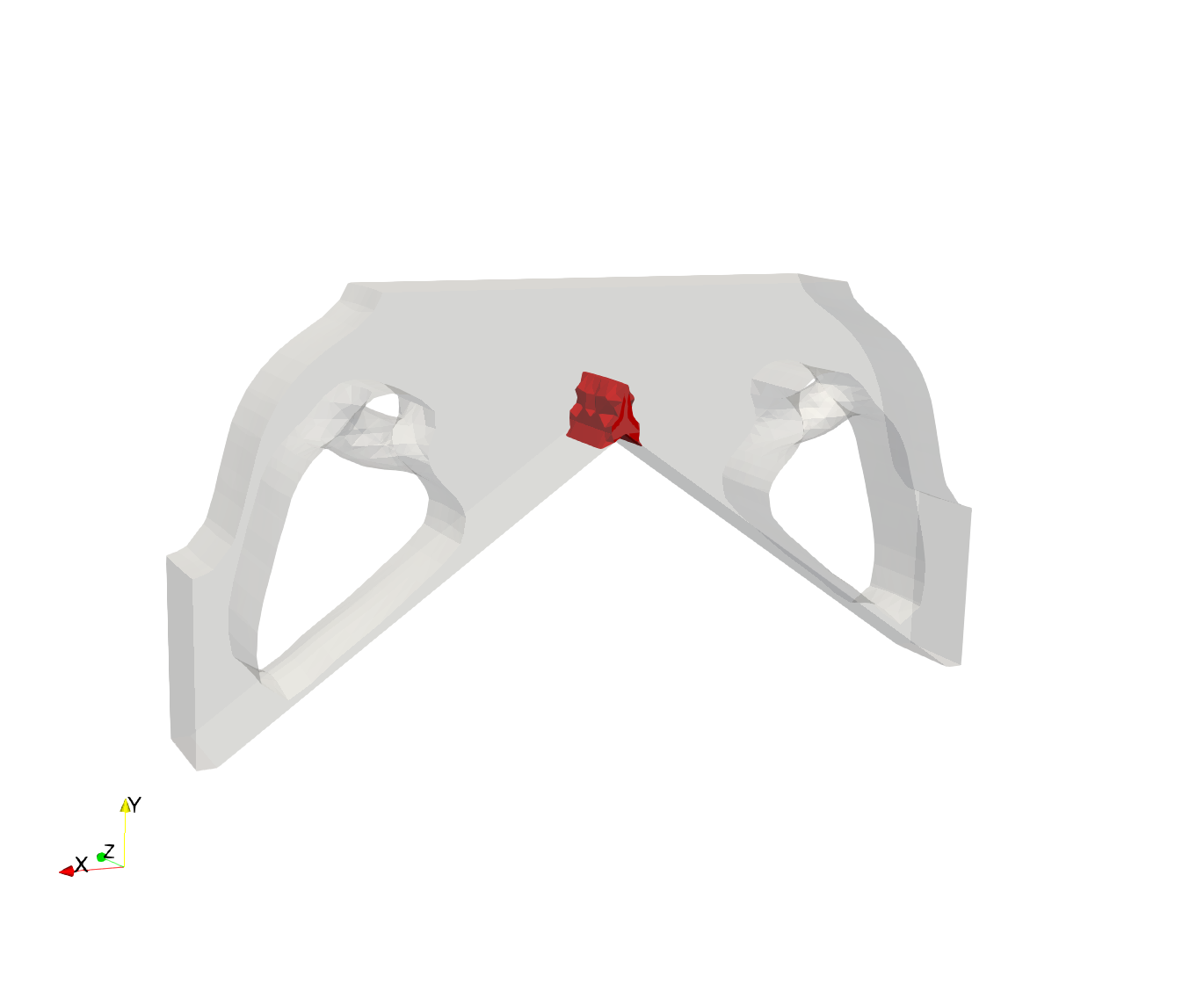}}   		\subfloat{\includegraphics[clip,trim=6cm 6cm 6cm 10cm, width=5.6cm]{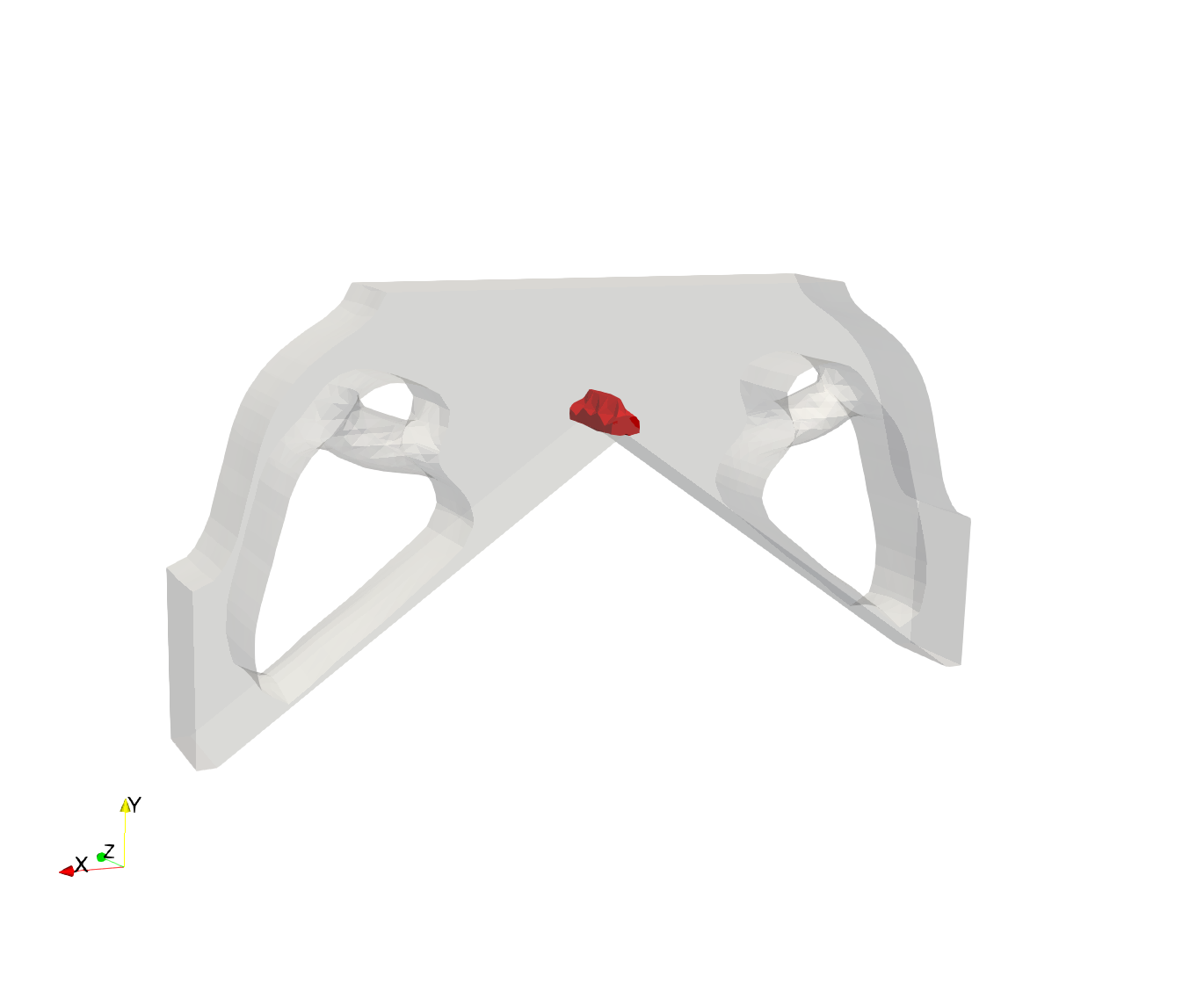}}   			\subfloat{\includegraphics[clip,trim=6cm 6cm 6cm 10cm, width=5.6cm]{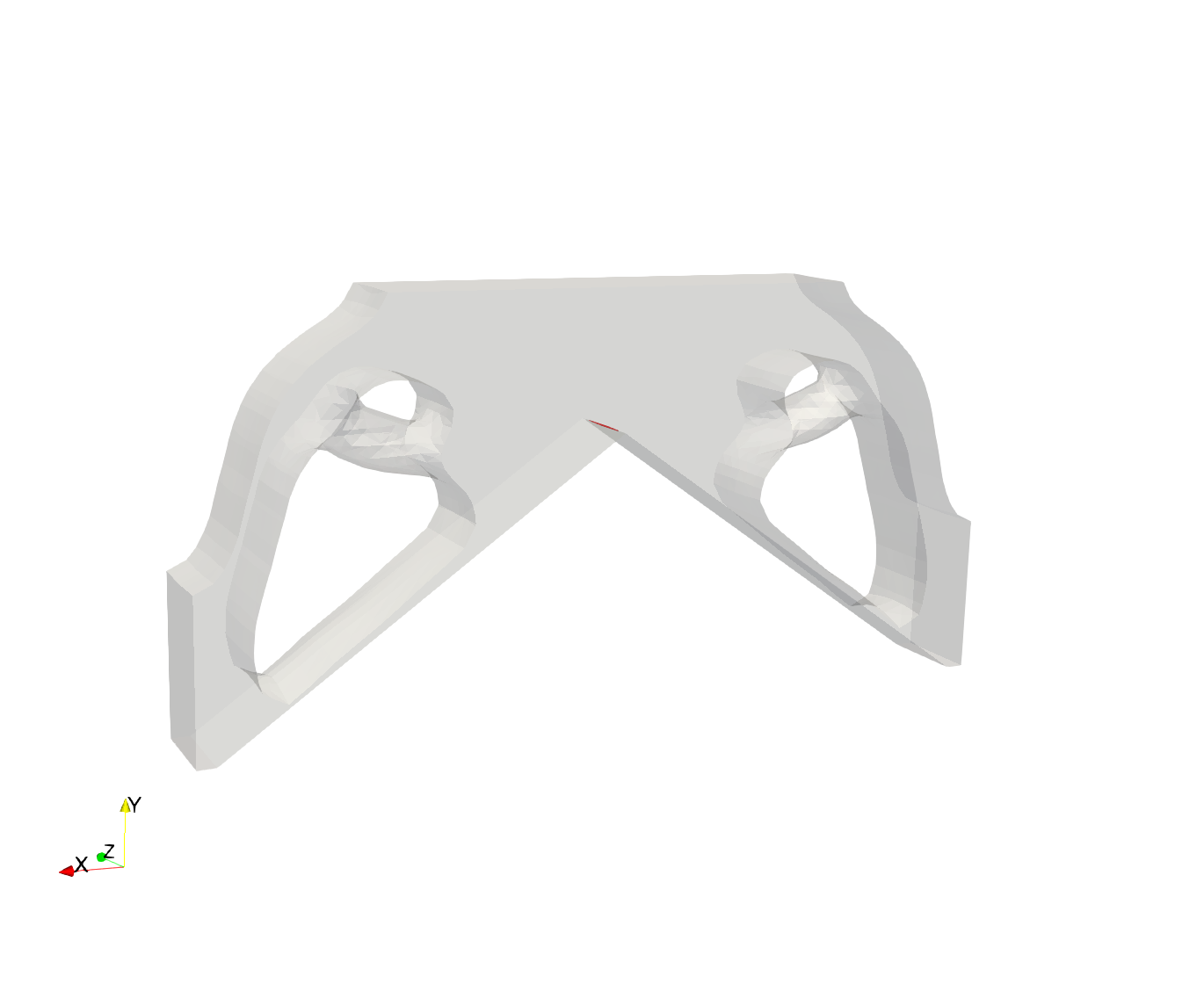}}	
	\caption{Example 3. The crack phase-field profile for different optimal topology layout of volume ratio}
	\label{Exm3_d}
\end{figure}

\begin{figure}[t!]
	\caption*{\underline{$\chi_v=1$}\hspace*{6cm}\underline{$\chi_v=0.57$}\hspace*{4cm}\underline{$\chi_v=0.40$}\hspace*{1cm}}
	\vspace{-0.1cm}
	\subfloat{\includegraphics[clip,trim=6cm 6cm 6cm 10cm, width=5.6cm]{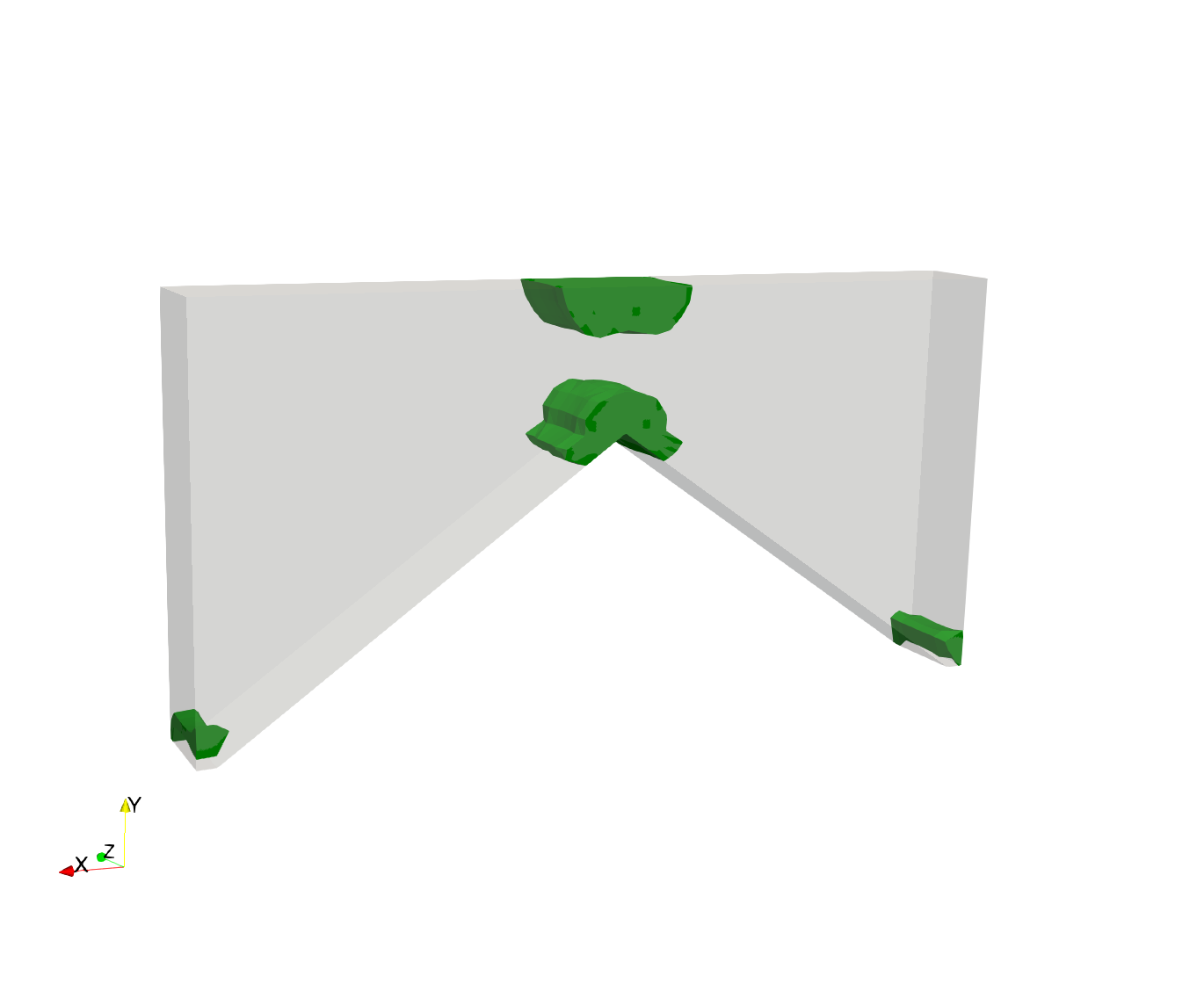}}   		\subfloat{\includegraphics[clip,trim=6cm 6cm 6cm 10cm, width=5.6cm]{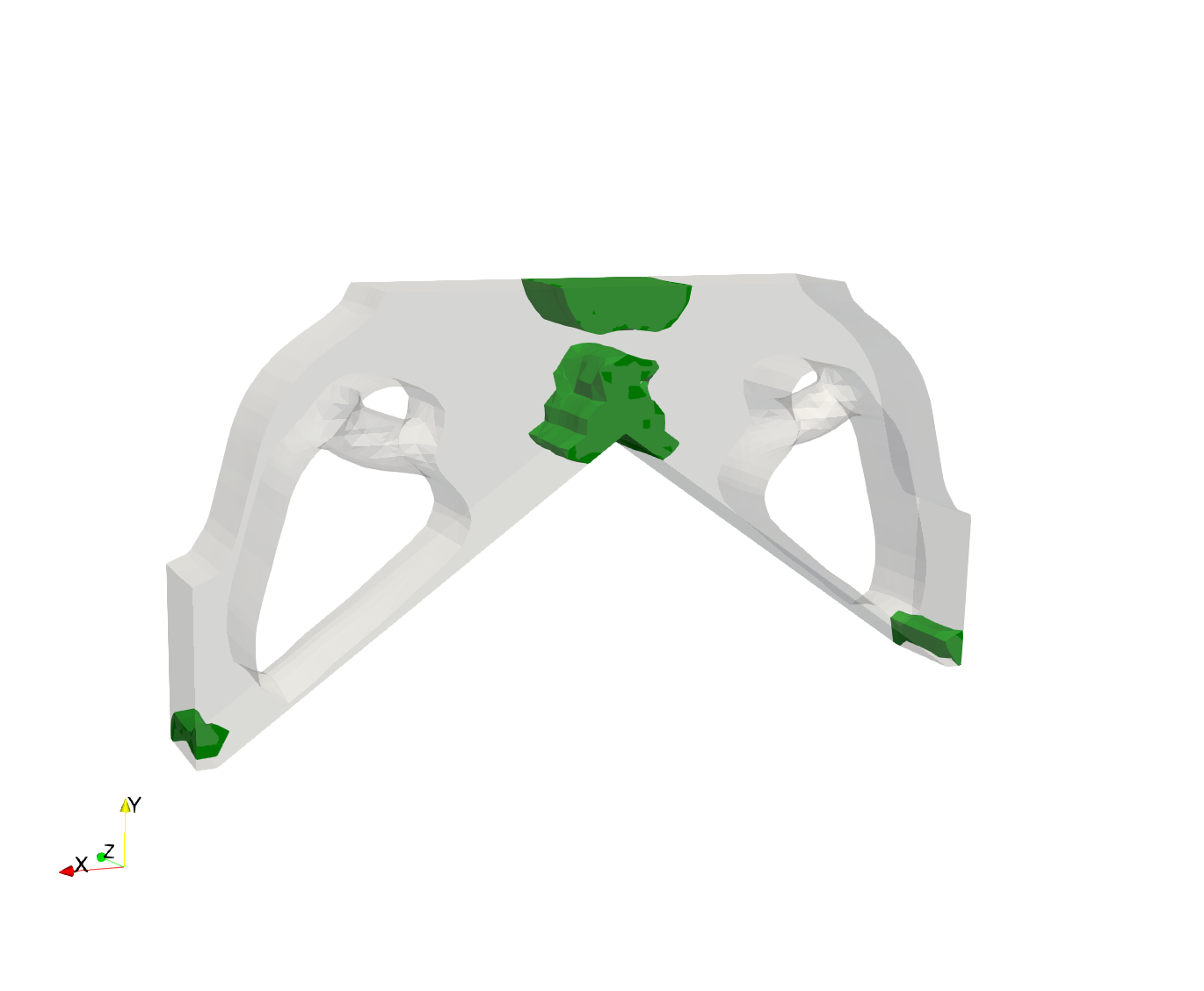}}   			\subfloat{\includegraphics[clip,trim=6cm 6cm 6cm 10cm, width=5.6cm]{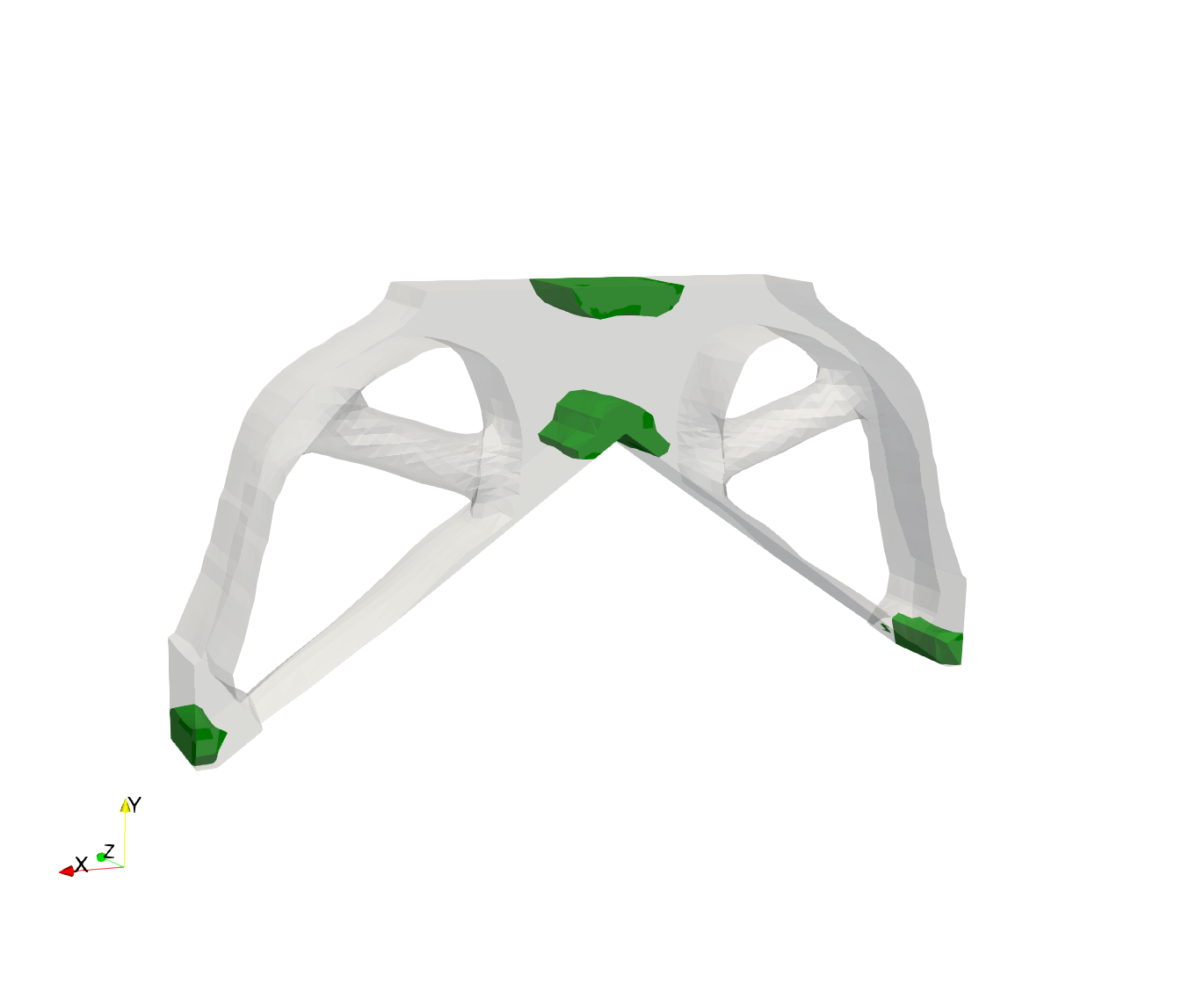}}	
	\caption{Example 3.  Computed hardening value (plastic zone) for evolution history of the different optimal layouts of volume ratio for ductile fracture model.}
	\label{Exm3_alpha}
\end{figure}

The evolution history of the optimal layouts for the portal frame structure is presented in Figure \ref{Exm3_phi}. Additionally, the crack phase-field profiles due to optimal layout results are depicted in Figure \ref{Exm3_d}. An important observation is that the evolutionary history of the optimal layout by approaching the volume fraction of $\chi_v=0.40$ leads to a reduction of the fracture zone. Thus, the fracture-resistance topology optimization formulation has prevented crack propagation in the material domain. So, the final optimum layout presents the stiffest response due to the fractured state, although it has less volume ratio. Besides the fracture patterns, Figure \ref{Exm3_alpha} illustrates the equivalent plastic strain $\alpha(\Bx,t)$ of the evolutionary history of the optimal layout. Accordingly, the plastic zone is shown by the green color, in which maximum equivalent plastic strain occurs \grm{nearby the} prescribed load, supports, and corner surface in the middle. It is worth noting that as long as the optimum layout approaches $\chi_v=0.40$, the plastic area will be reduced. This shows another interesting advantage of using the proposed method.

\begin{figure}[!t]
	\centering
	\subfloat{{\includegraphics[clip,trim=0cm 28cm 0cm 0cm, width=17cm]{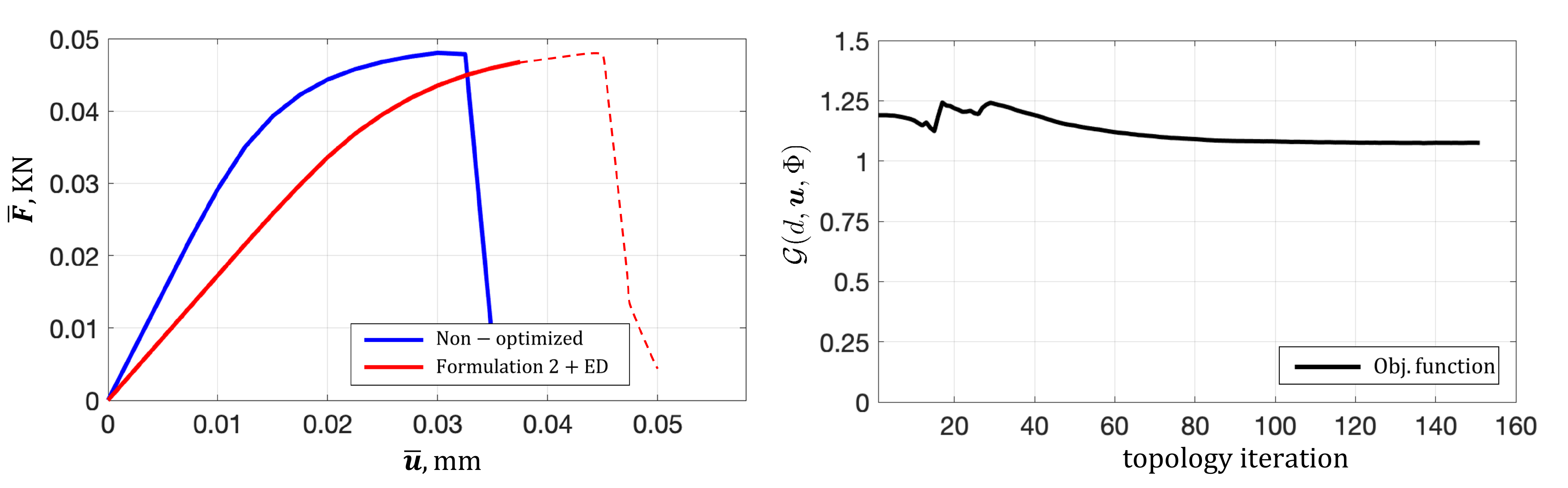}} } 
	\vspace{-0.15cm}
	\caption*{\hspace*{4.3cm}(a)\hspace*{8cm}(b)\hspace*{2cm}}
	\subfloat{{\includegraphics[clip,trim=0cm 28cm 0cm 0cm, width=17cm]{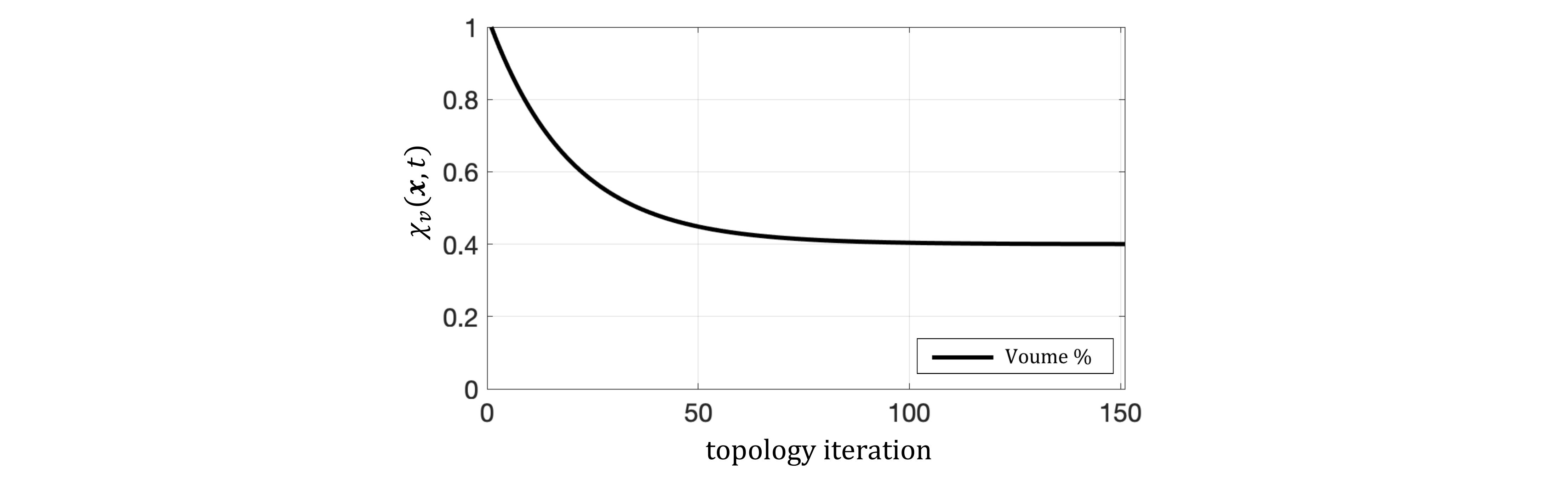}}} 
	\vspace{-0.25cm}
    \caption*{\hspace*{2cm}(c)}
	\caption{Example 3.  Computed response for the  portal frame structure under compression loading. (a) Load-displacement curves, (b) convergence history of the objective function, and (c) volume constraint function.}
	\label{Exm3_LD}
\end{figure}

Another impacting factor that should be noted is the load-displacement curve for non-optimized and optimized layouts, illustrated in Figure \ref{Exm3_LD}(a). \noii{The failure softening region is further depicted as a dashed line, to observe the delay in the fracture state.} Indeed, using topology optimization toward fracture avoids crack initiation, thus there will be no softening due to fracture as expected from the proposed model. The convergence performance for the objective and volume constraint functions are depicted in Figure \ref{Exm3_LD}(b-c), respectively. Evidently, by means of objective behavior in Figure \ref{Exm3_LD}(b), after 100 optimization iterations, a stable behavior is observed. It can be grasped that the volume constraint function asymptotically approaches its desired value, i.e., $\chi_v=0.40$.

\begin{figure}[!b]
	\centering
	{\includegraphics[clip,trim=1.3cm 14.5cm 0cm 1cm, width=17cm]{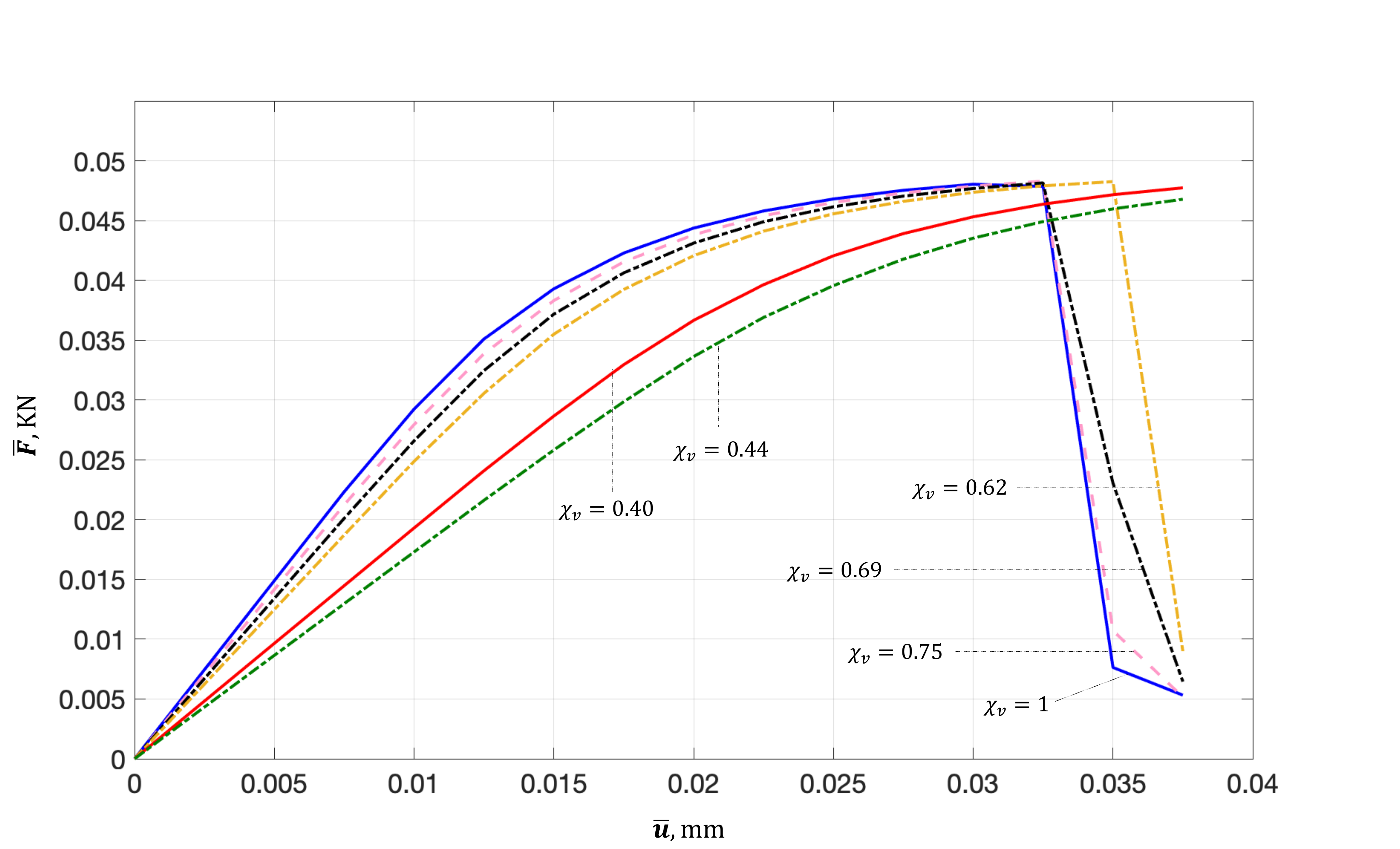}}  
	\caption{Example 3. Comparison load-displacement curves for different volume ratio
		of topology optimization iteration up to final optimum layout.}
	\label{Exm3_inc}
\end{figure}

\begin{figure}[t!]
	\centering
	\vspace{-0.1cm}
	\subfloat{\includegraphics[clip,trim=10cm 7.1cm 7.7cm 5.6cm, width=7cm]{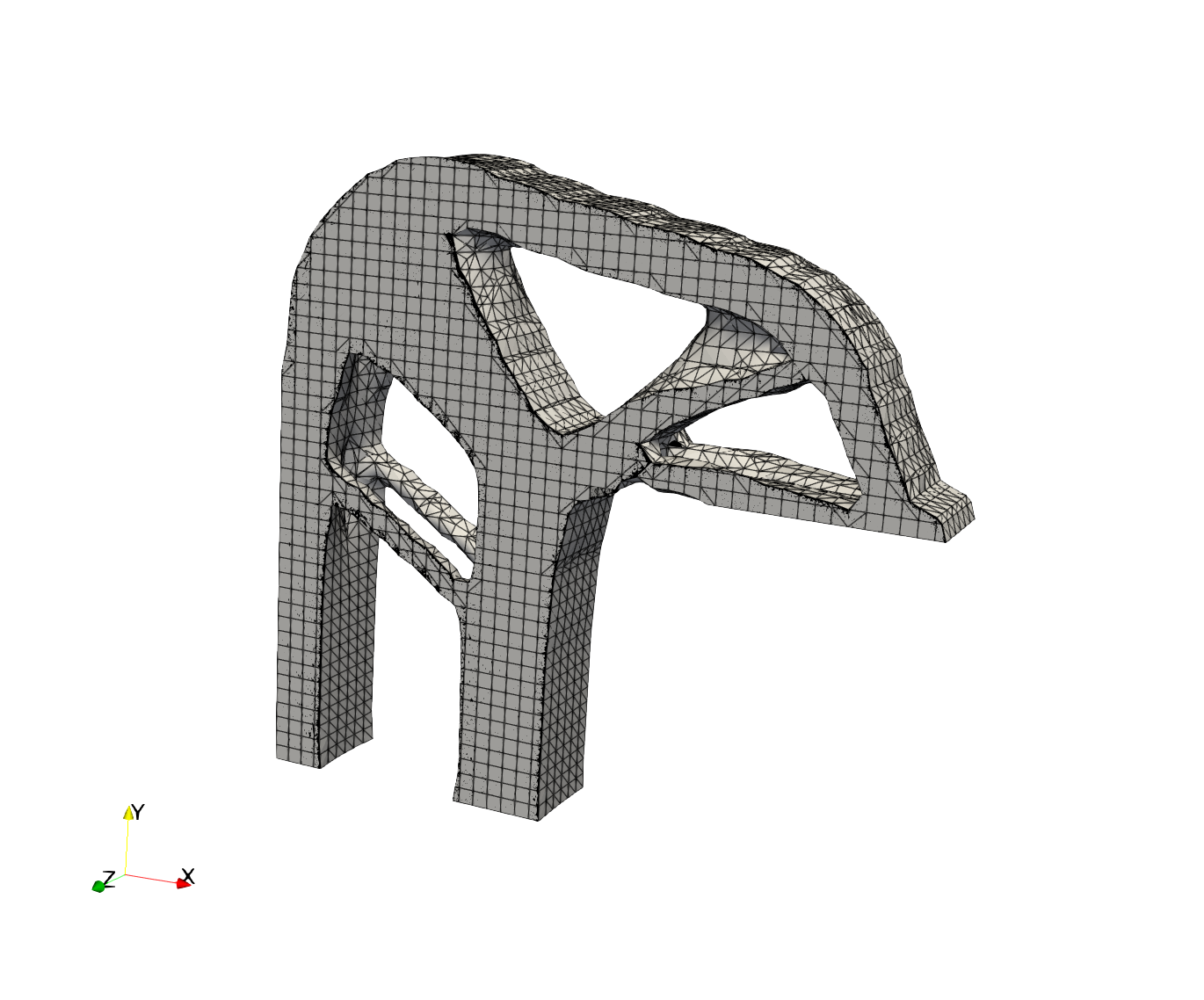}}   	\subfloat{\includegraphics[clip,trim=6cm 6cm 7cm 10cm, width=7cm]{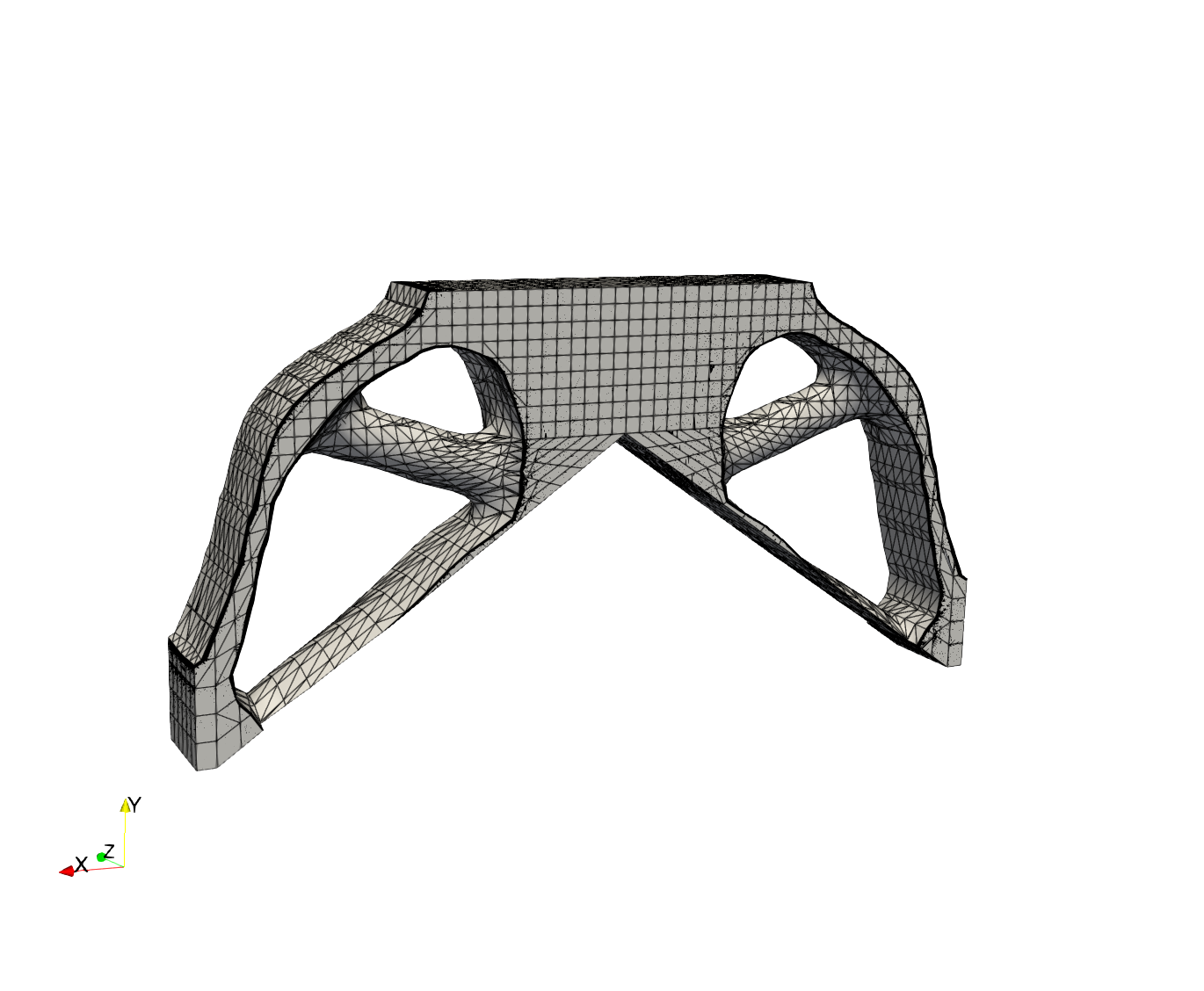}}		
	\caption*{\hspace*{1.8cm}(a)\hspace*{8cm}(b)}
	\caption{Finite element discretization of the optimal layouts when $\chi_v=0.4$ for (a) Example 2, and (b) Example 3.}
	\label{Exm23_mesh}
\end{figure}

Next, we investigate the evolution of load-displacement curves for different volume ratios up to the convergence of the optimization to the final topology. This is outlined in Figure \ref{Exm3_inc}. Indeed, this evolution highlights the effects of the damage response within topology optimization, which is reduced for every new optimum layout, thus softening region will gradually vanish. As a result, the proposed model shows its proficiency approximately after $\chi_v\ge0.44$. Lastly, the smooth discretized final optimum layouts for this example are depicted in Figure \ref{Exm23_mesh}(b). 

\sectpb[Section53]{Example 4: Short cantilever beam under shear loading}

Finally, we investigate the optimum layout of a three-dimensional cantilever beam under shear loading. This is a common topology optimization problem which has been reported in \cite{desai_Allaire,Russ,marino2021mixed,dede2012isogeometric}. Here, we determine an optimum topology configuration due to both brittle and ductile phase-field fracture in order to compare the two. We will discuss this in detail.

\begin{figure}[!t]
	\centering
	{\includegraphics[clip,trim=1cm 24cm 0cm 0cm, width=16cm]{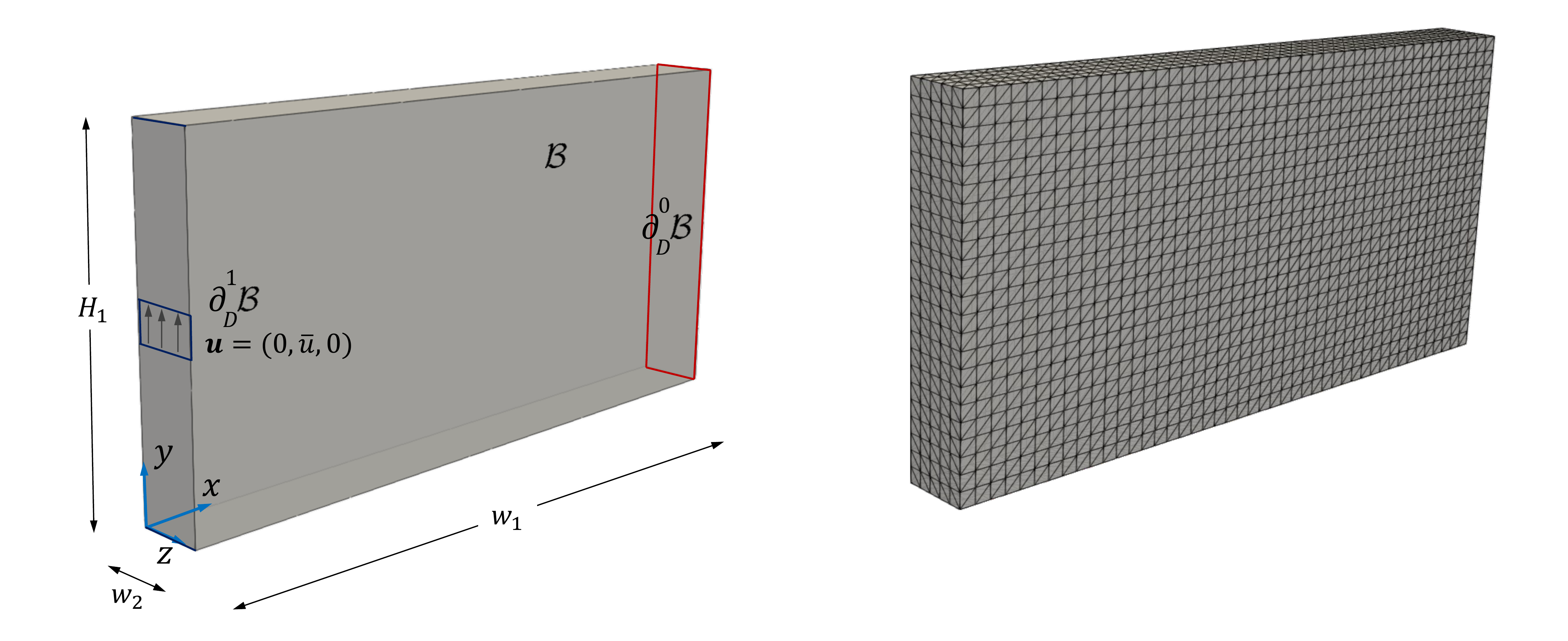}}  
	\vspace*{-0.7cm}
	\caption*{\hspace*{4.3cm}(a)\hspace*{8cm}(b)\hspace*{2cm}}
	\caption{Example 4. The representation of the (a) geometry and boundary conditions, and (b) finite element discretization.}
	\label{Exm4_bvp}
\end{figure}

A boundary value problem is depicted in Figure \ref{Exm4_bvp}(b). The surface on its rightmost support is constrained in all directions. The geometrical dimensions for Figure \ref{Exm4_bvp}(b) are set as $w_1=1\;mm$, $w_2=0.1\;mm$, and $H_1=0.5\;mm$. The reference point \grm{is} located in $(a,b,c)=(2.25,1.5,0.5)\;mm$. The numerical example is performed by applying a monotonic displacement increment  ${\Delta \bar{u}}_y=0.5\times^{-4}\;mm$ which is prescribed \noii{upward} in a vertical direction in a part of the top page $area=(0,0.23,0)\times (0,0.27,0.1)$ of the specimen for 146-time steps. Thus, the final prescribed displacement load for this minimization problem is ${\bar{u}}_y=7.3\times^{-3}\;mm$. The minimum finite element size in the solid domains is $h_{min}=0.02\;mm$ which contains 32161 linear tetrahedral elements.

We start with the presentation of the quantitative and qualitative topology optimization constrained by failure response for cantilever beam considering \textit{two different} scenarios.

\begin{figure}[t!]
	\caption*{\hspace*{1cm}\underline{$\chi_v=0.56$}\hspace*{4.5cm}\underline{$\chi_v=0.50$}\hspace*{3.5cm}\underline{$\chi_v=0.45$}\hspace*{1.5cm}}
	\vspace{-0.1cm}
	\subfloat{\includegraphics[clip,trim=5cm 10cm 6.5cm 11cm, width=5.5cm]{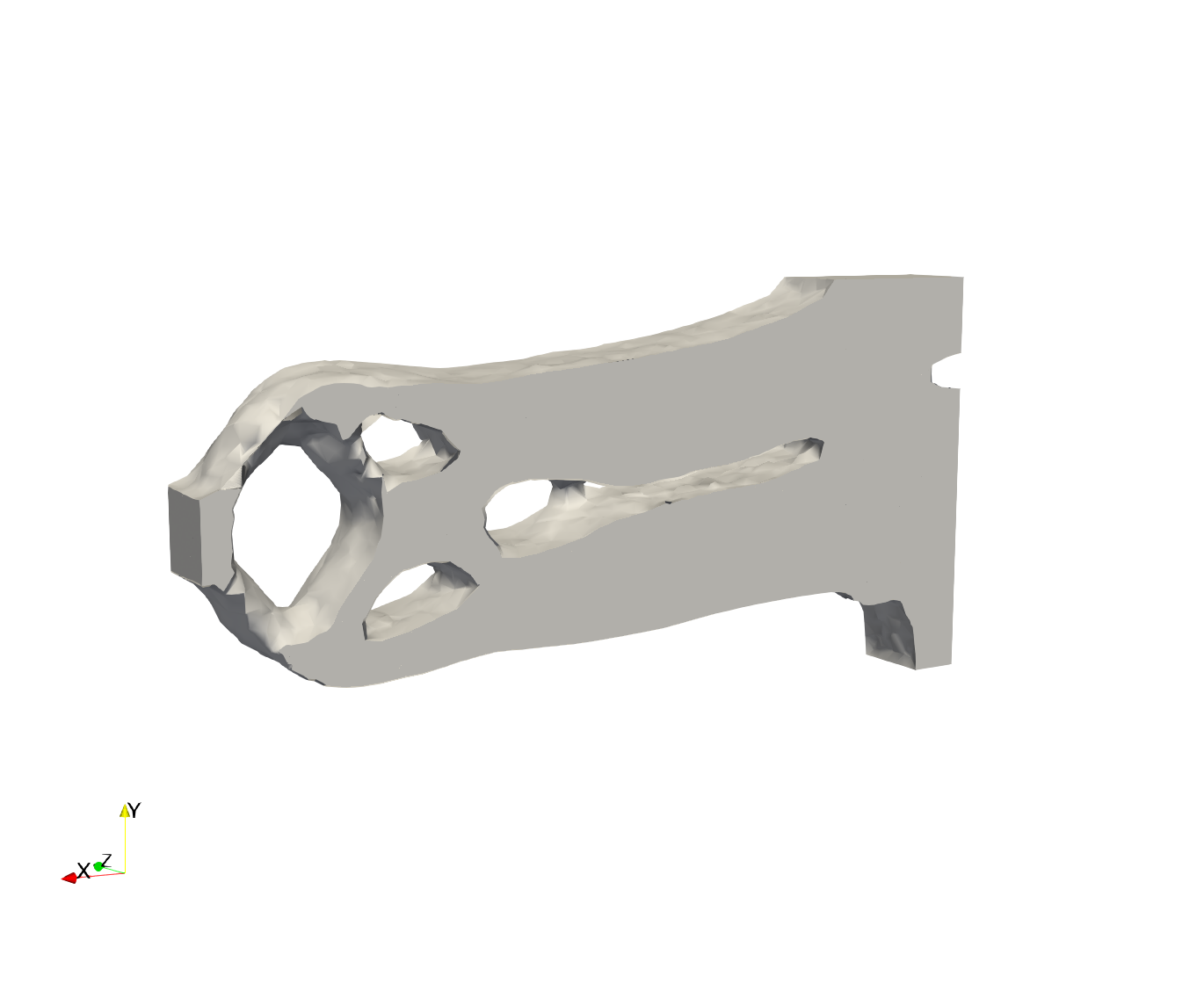}}   		\subfloat{\includegraphics[clip,trim=5cm 10cm 6.5cm 11cm, width=5.5cm]{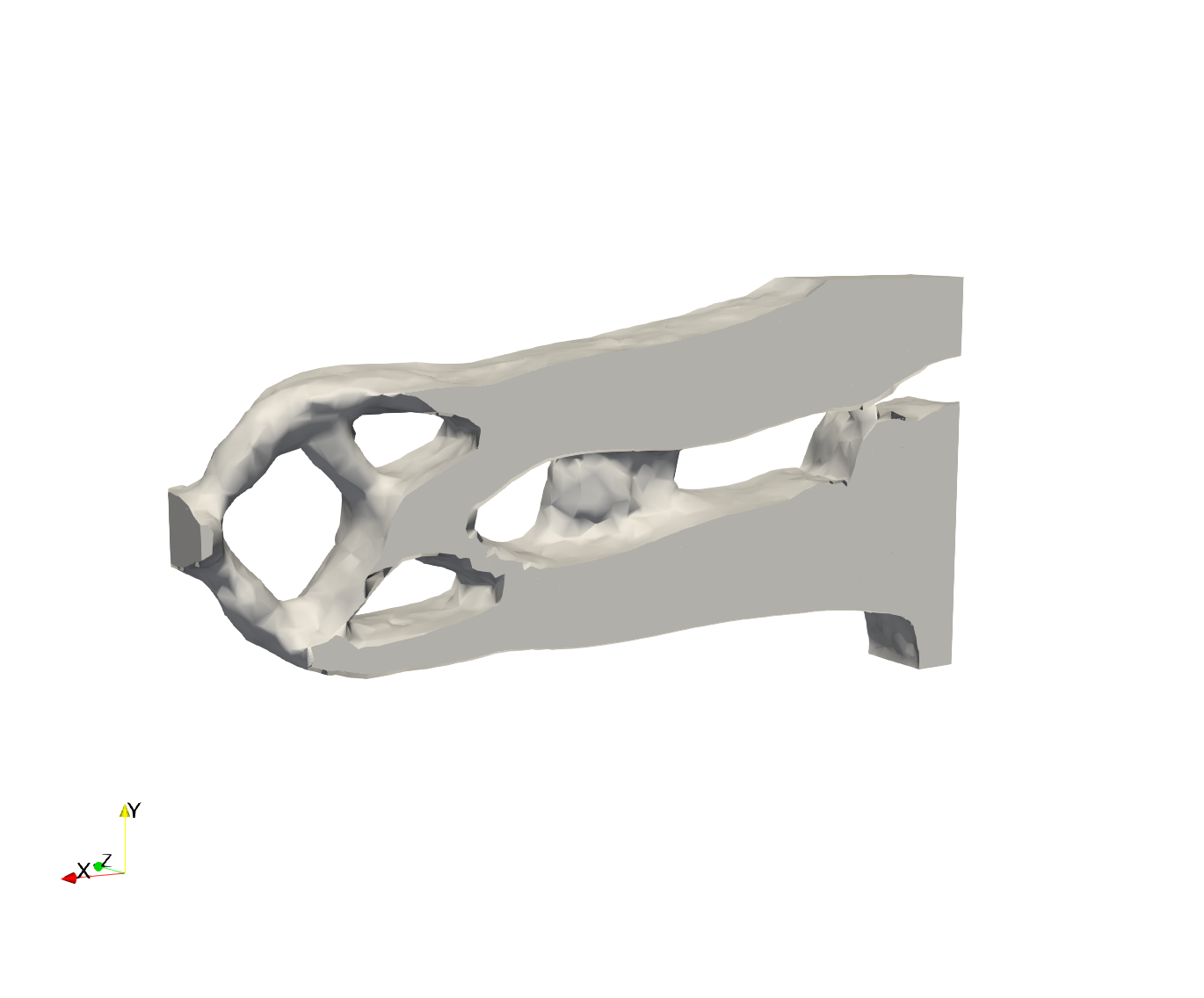}}   			\subfloat{\includegraphics[clip,trim=5cm 10cm 6.5cm 11cm, width=5.5cm]{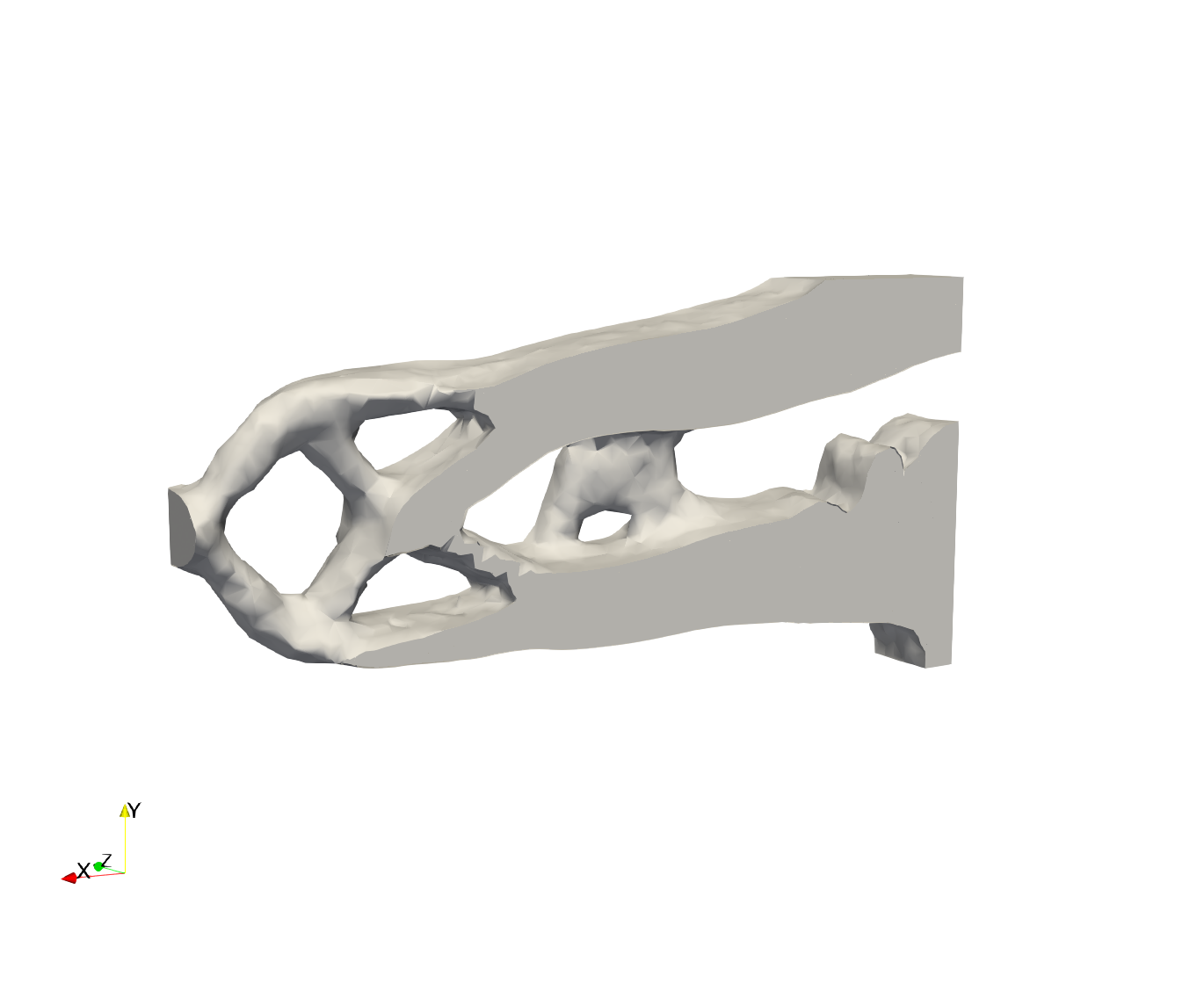}}
	\caption*{\hspace*{1cm}\underline{$\chi_v=0.41$}\hspace*{4.5cm}\underline{$\chi_v=0.405$}\hspace*{3.5cm}\underline{$\chi_v=0.40$}\hspace*{1.5cm}}
	\vspace{-0.1cm}
	\subfloat{\includegraphics[clip,trim=5cm 10cm 6.5cm 11cm, width=5.5cm]{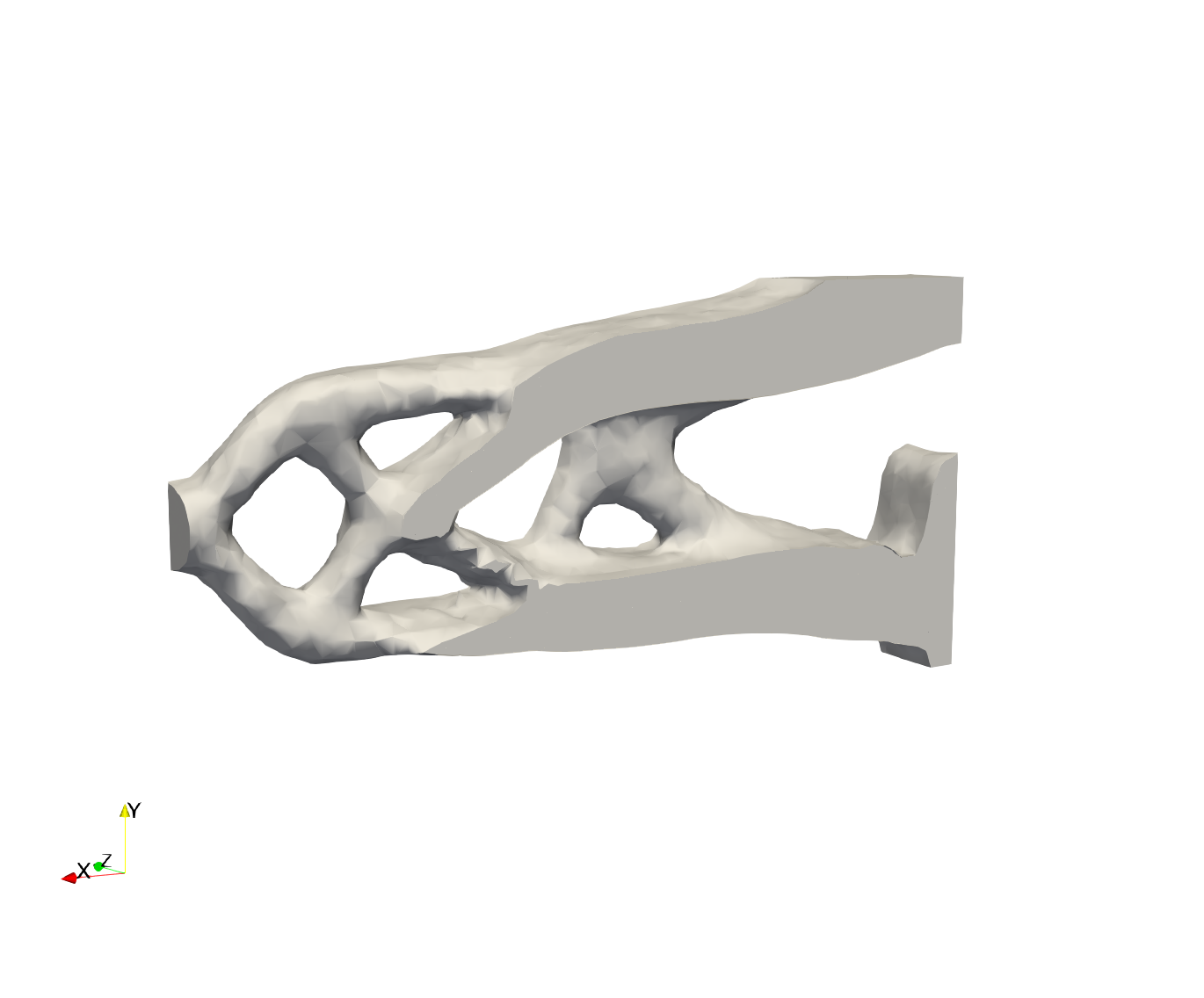}}   		\subfloat{\includegraphics[clip,trim=5cm 10cm 6.5cm 11cm, width=5.5cm]{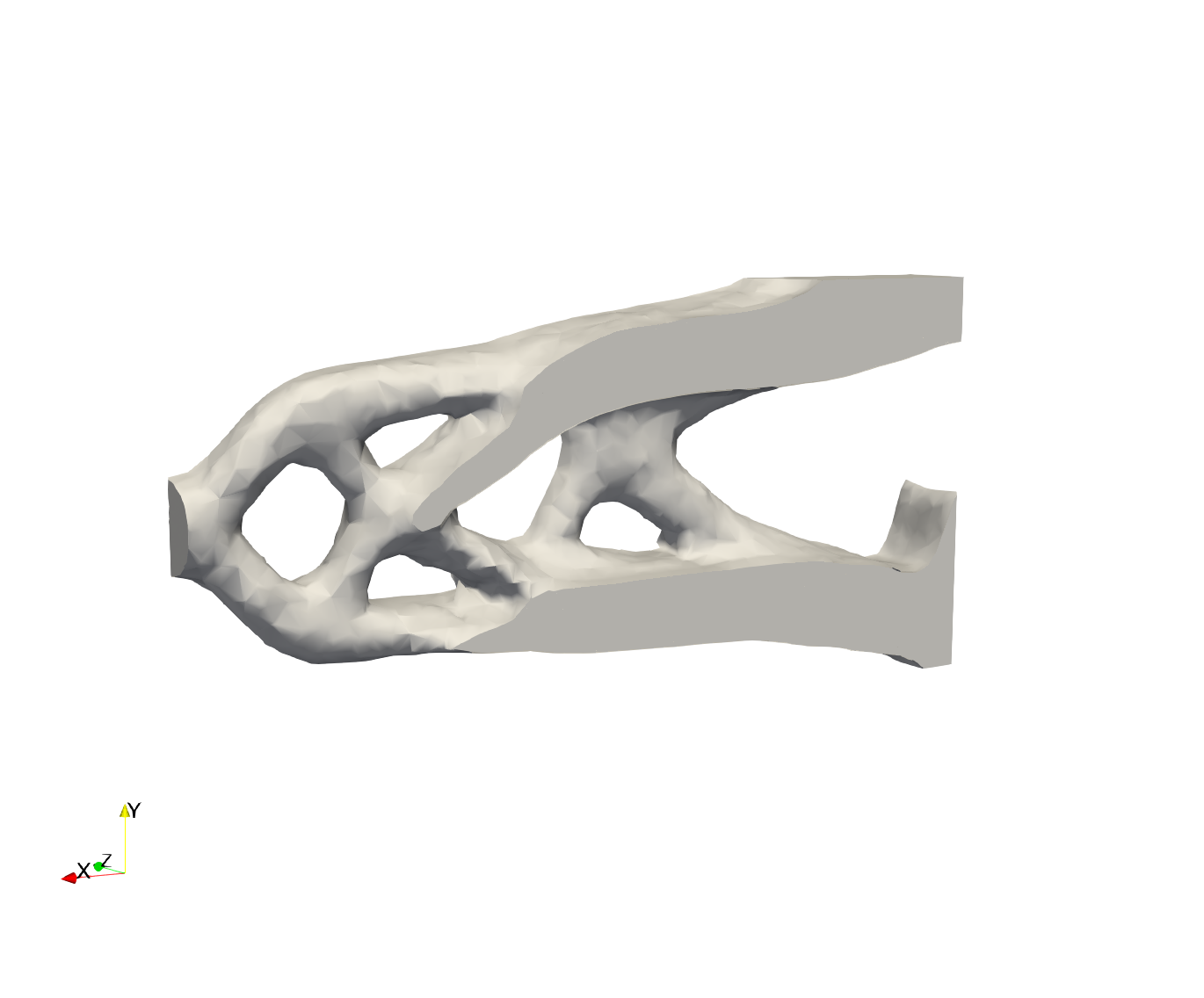}}   			\subfloat{\includegraphics[clip,trim=5cm 10cm 6.5cm 11cm, width=5.5cm]{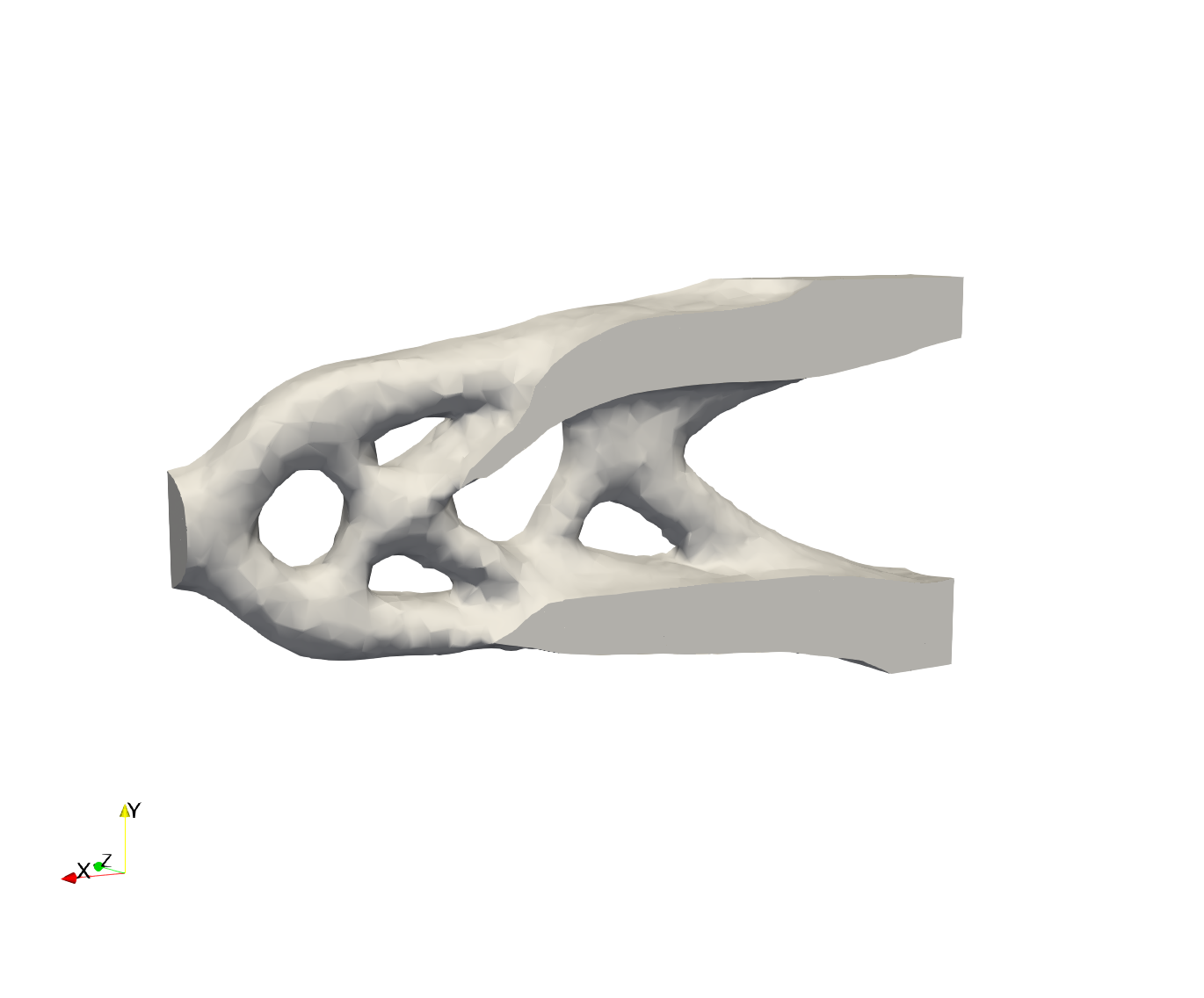}}
	\caption{Example 4 (Case a). Evolution history of the optimal layouts for different volume ratio of the  short cantilever beam for brittle fracture.}
	\label{Exm4_ED_phi}
\end{figure}

\begin{figure}[t!]
	\caption*{\hspace*{1cm}\underline{$\chi_v=0.81$}\hspace*{4.5cm}\underline{$\chi_v=60$}\hspace*{3.5cm}\underline{$\chi_v=0.40$}\hspace*{1.5cm}}
	\vspace{-0.1cm}
	\subfloat{\includegraphics[clip,trim=5cm 10cm 6.5cm 11cm, width=5.5cm]{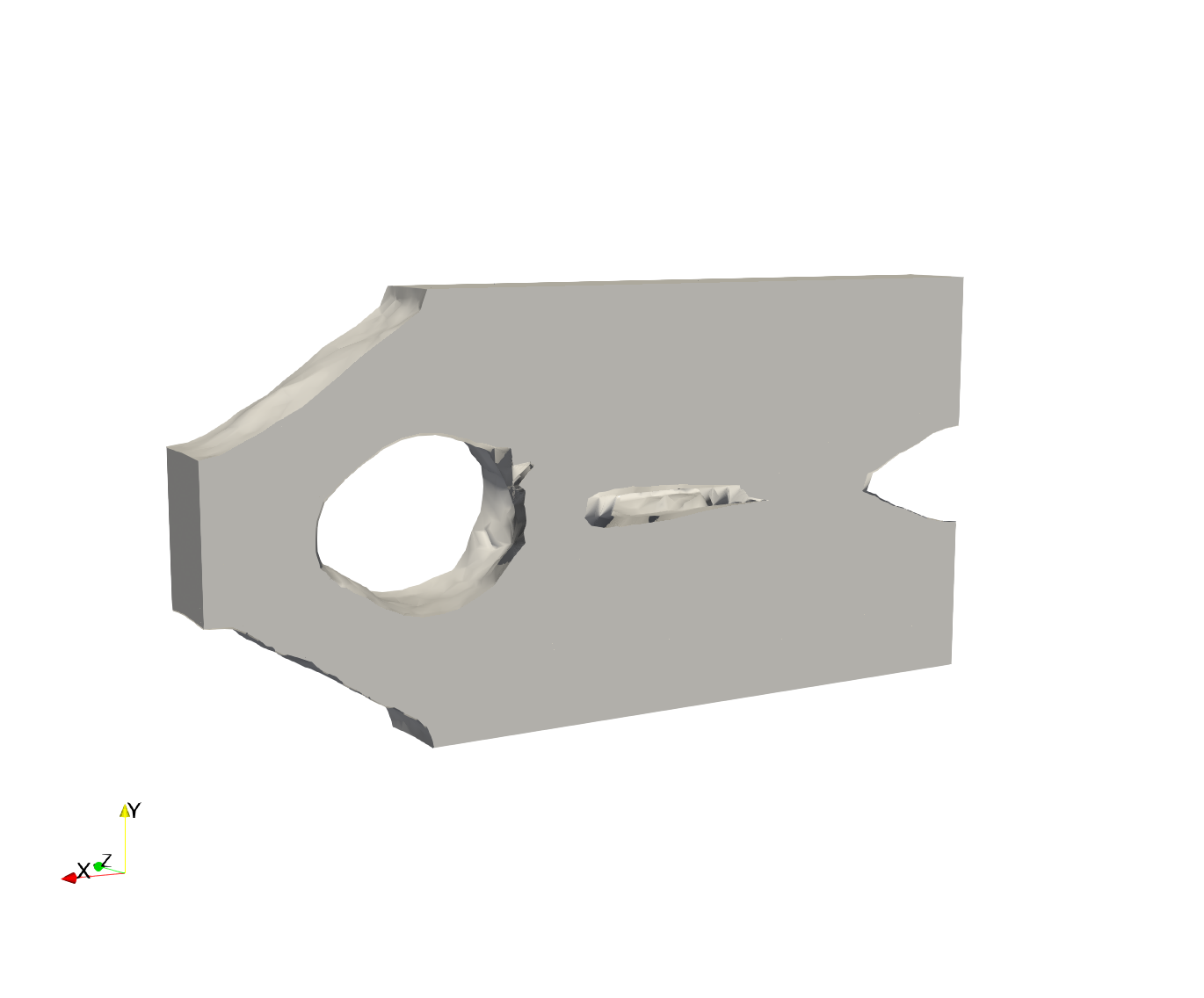}}   		\subfloat{\includegraphics[clip,trim=5cm 10cm 6.5cm 11cm, width=5.5cm]{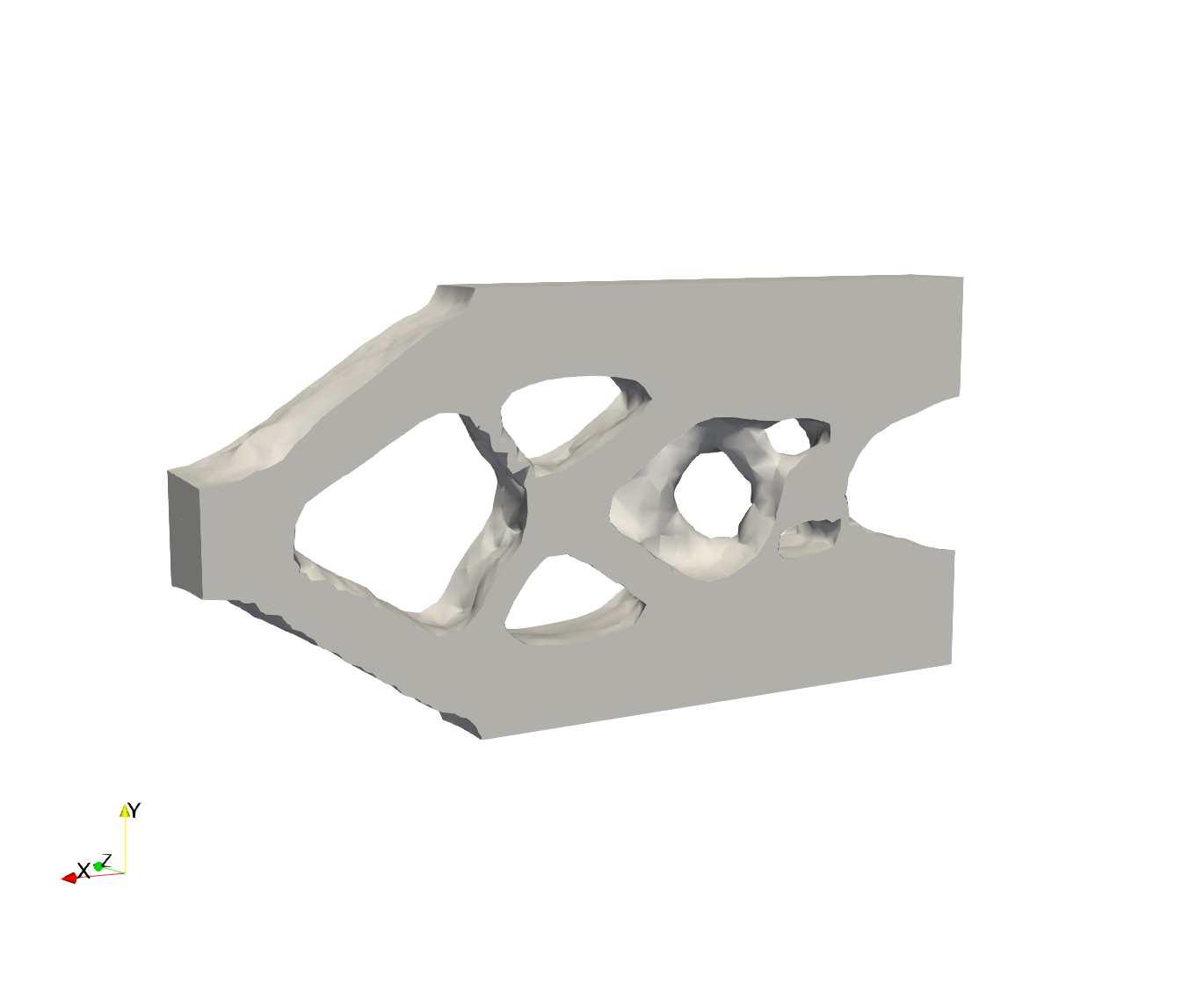}}   			\subfloat{\includegraphics[clip,trim=5cm 10cm 6.5cm 11cm, width=5.5cm]{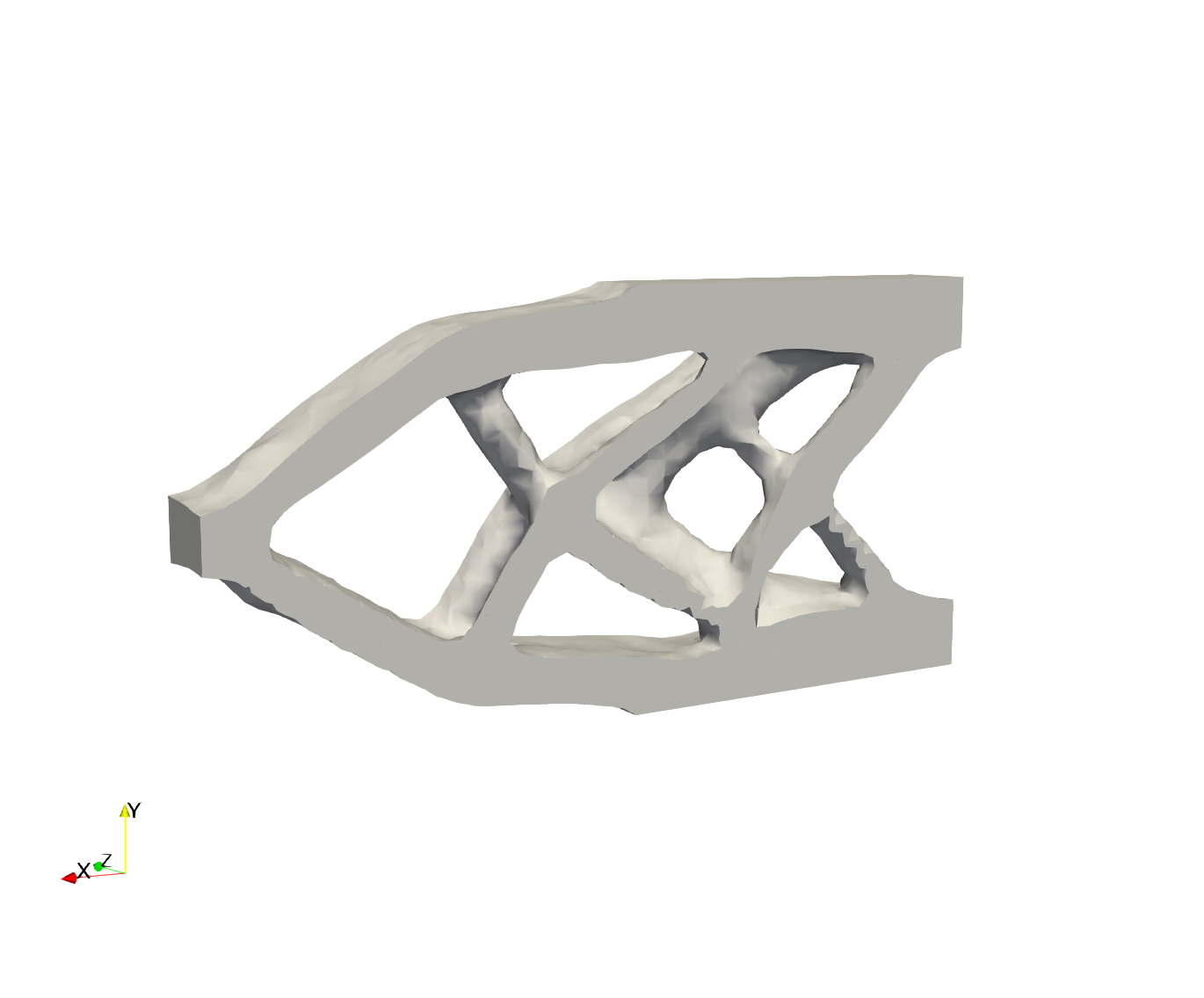}}
	\caption{Example 4. Evolution history of the optimal layouts for different volume ratio of the  short cantilever beam for linear elasticity results.}
	\label{Exm4_E_phi}
\end{figure}

\begin{figure}[t!]
	\caption*{\hspace*{1cm}\underline{$\chi_v=0.81$}\hspace*{4.5cm}\underline{$\chi_v=0.60$}\hspace*{3.5cm}\underline{$\chi_v=0.40$}\hspace*{1.5cm}}
	\vspace{-0.1cm}
	\subfloat{\includegraphics[clip,trim=5cm 10cm 6.5cm 11cm, width=5.5cm]{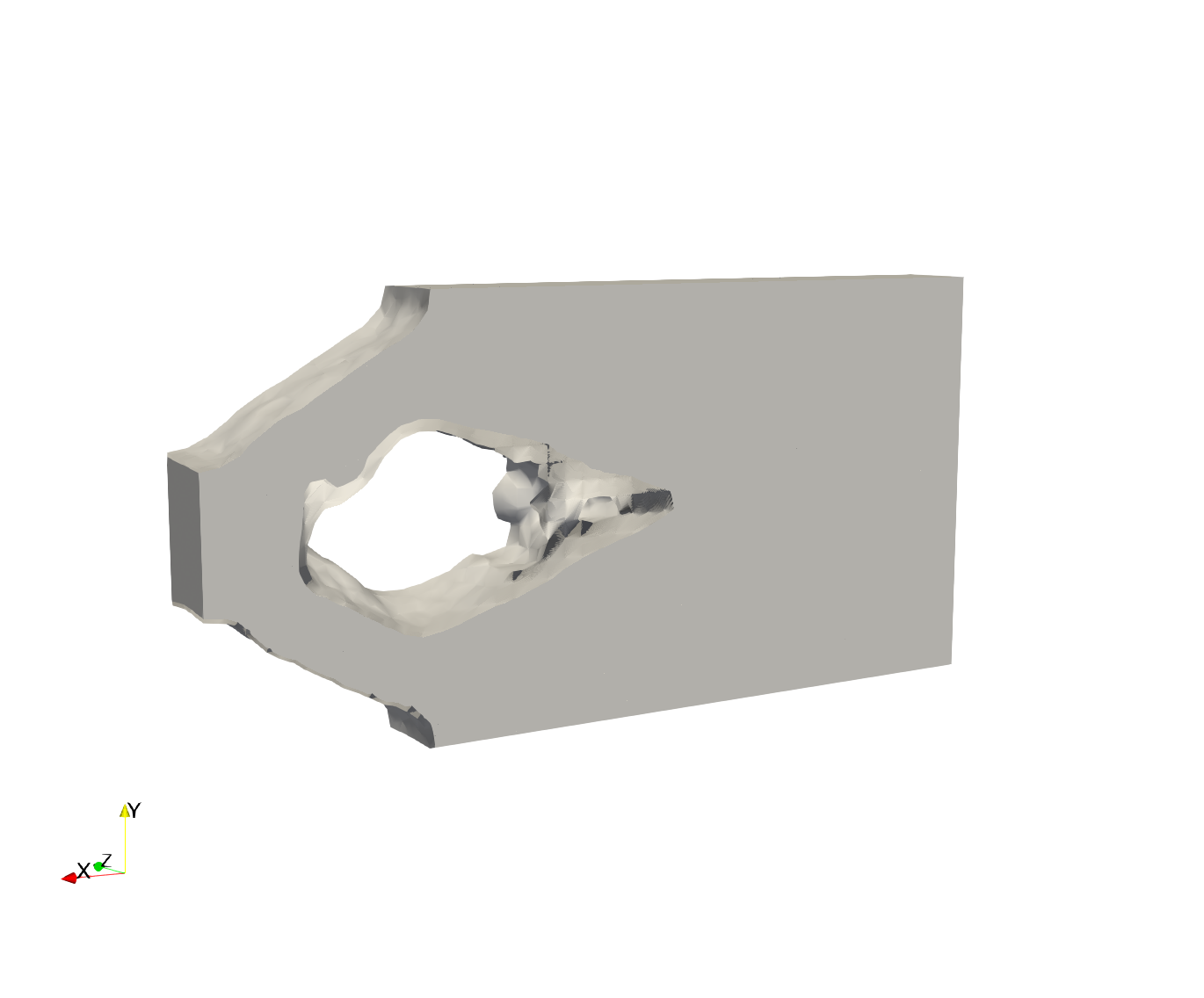}}   		\subfloat{\includegraphics[clip,trim=5cm 10cm 6.5cm 11cm, width=5.5cm]{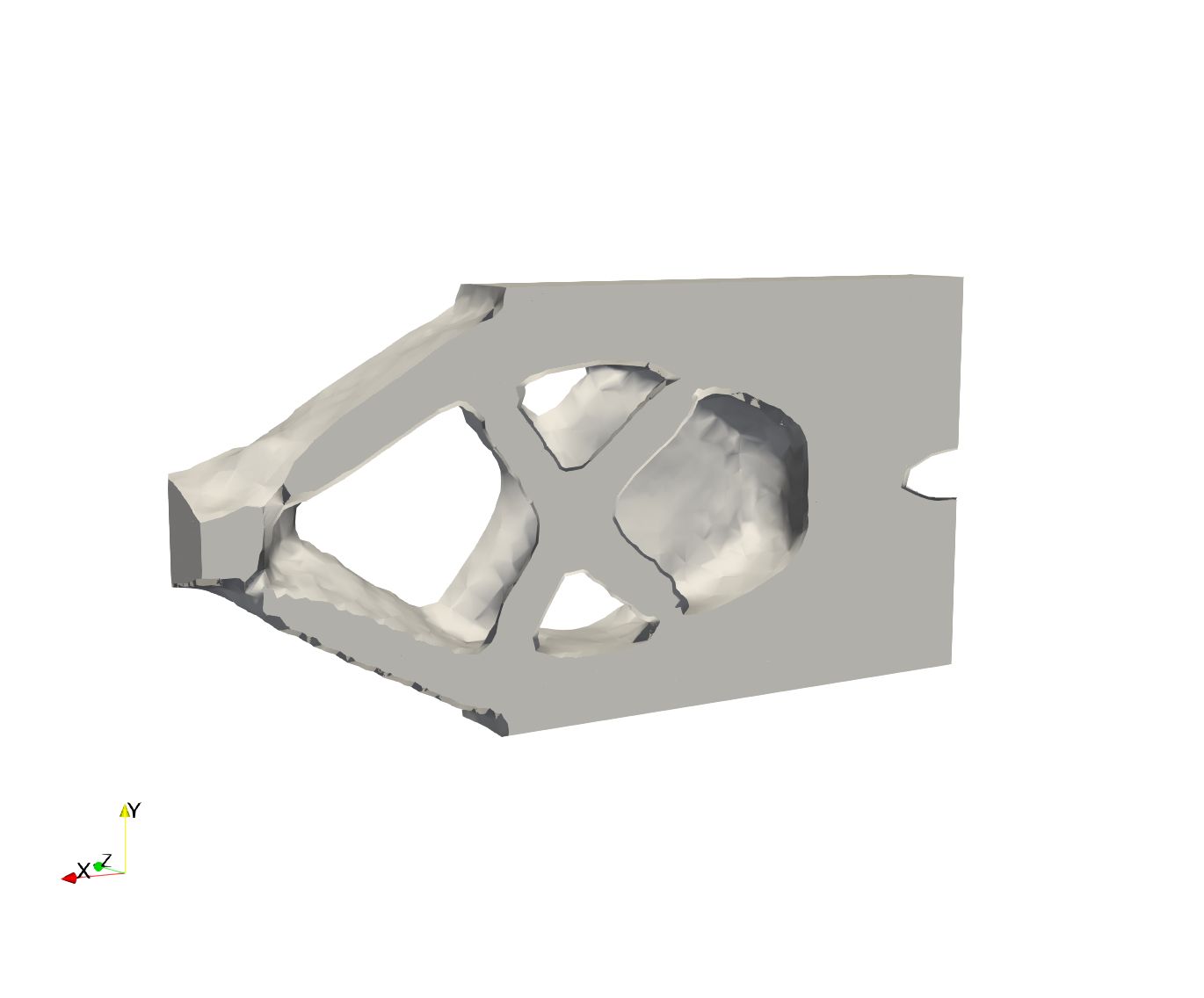}}   			\subfloat{\includegraphics[clip,trim=5cm 10cm 6.5cm 11cm, width=5.5cm]{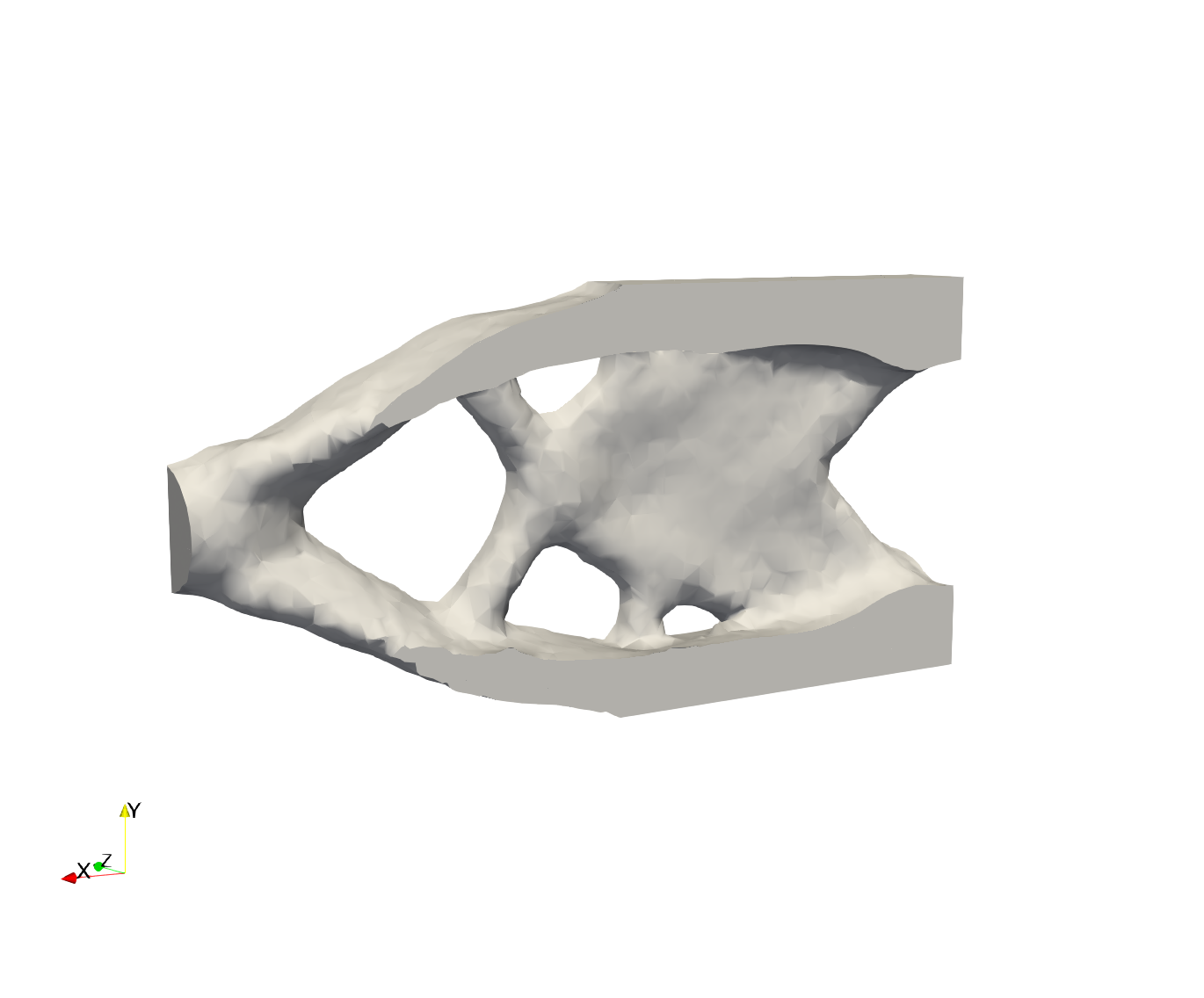}}
	\caption{Example 4 (Case b). Evolution history of the optimal layouts for different volume ratio of the  short cantilever beam for ductile fracture.}
	\label{Exm4_EPD_phi}
\end{figure}

\begin{figure}[!t]
	\centering
	{\includegraphics[clip,trim=0cm 28cm 0cm 0cm, width=17cm]{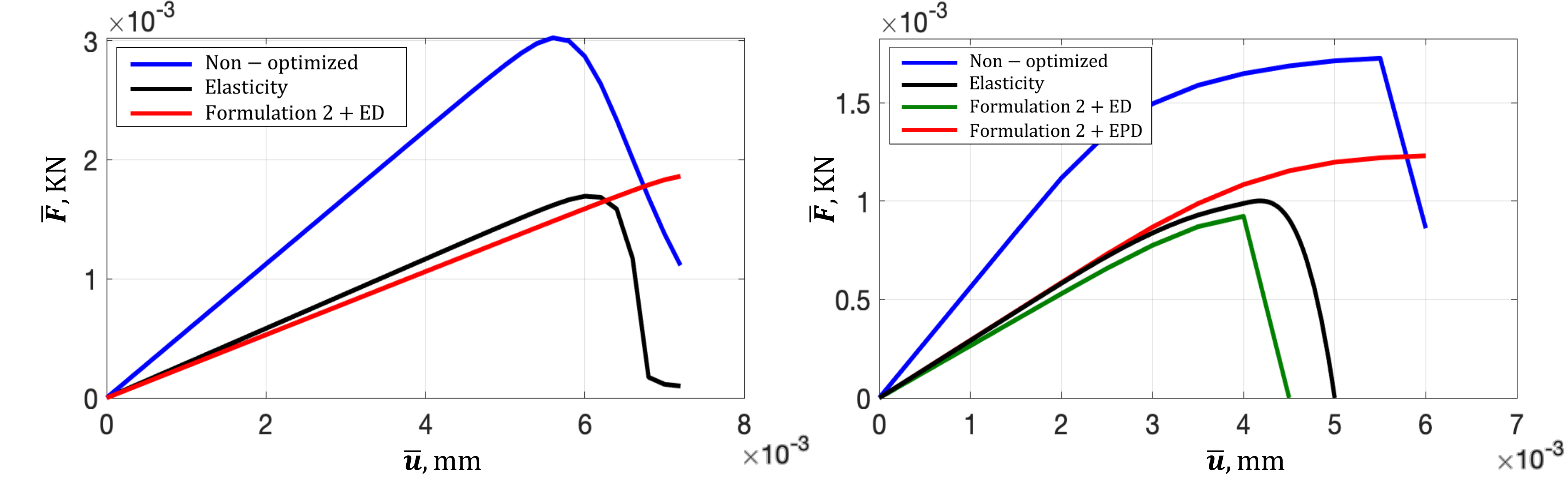}}  
	\vspace*{-0.7cm}
	\caption*{\hspace*{4.3cm}(a)\hspace*{8cm}(b)\hspace*{2cm}}
	\caption{Example 4. Comparison load-displacement curves for the short cantilever beam of the, (a) brittle fracture model described for Case a, (b) ductile fracture model illustrated in Case b.}
	\label{Exm4_LD}
\end{figure}

\begin{figure}[!t]
	\centering
	{\includegraphics[clip,trim=0cm 28cm 0cm 0cm, width=17cm]{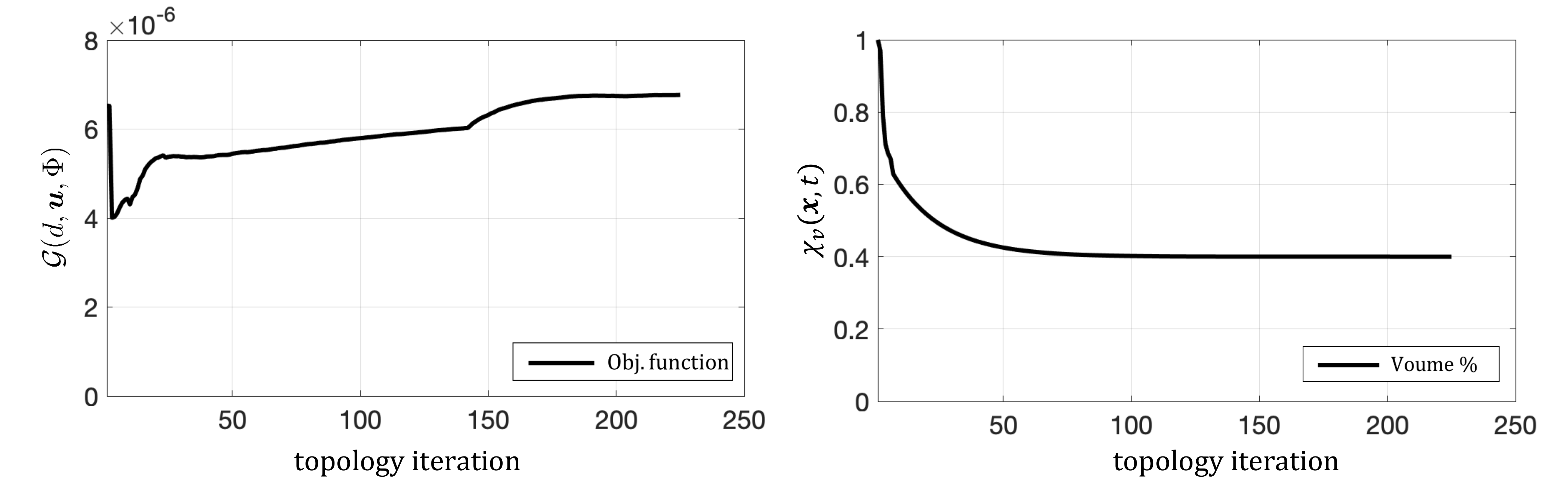}}  
	\vspace*{-0.7cm}
	\caption*{\hspace*{4.3cm}(a)\hspace*{8cm}(b)\hspace*{2cm}}
	\caption{Example 4 (Case a). Convergence history of the brittle fracture for the (a) objective
		function, and (b) volume constraint function.}
	\label{Exm4_conv1}
\end{figure}

\begin{figure}[!t]
	\centering
	{\includegraphics[clip,trim=0cm 28cm 0cm 0cm, width=17cm]{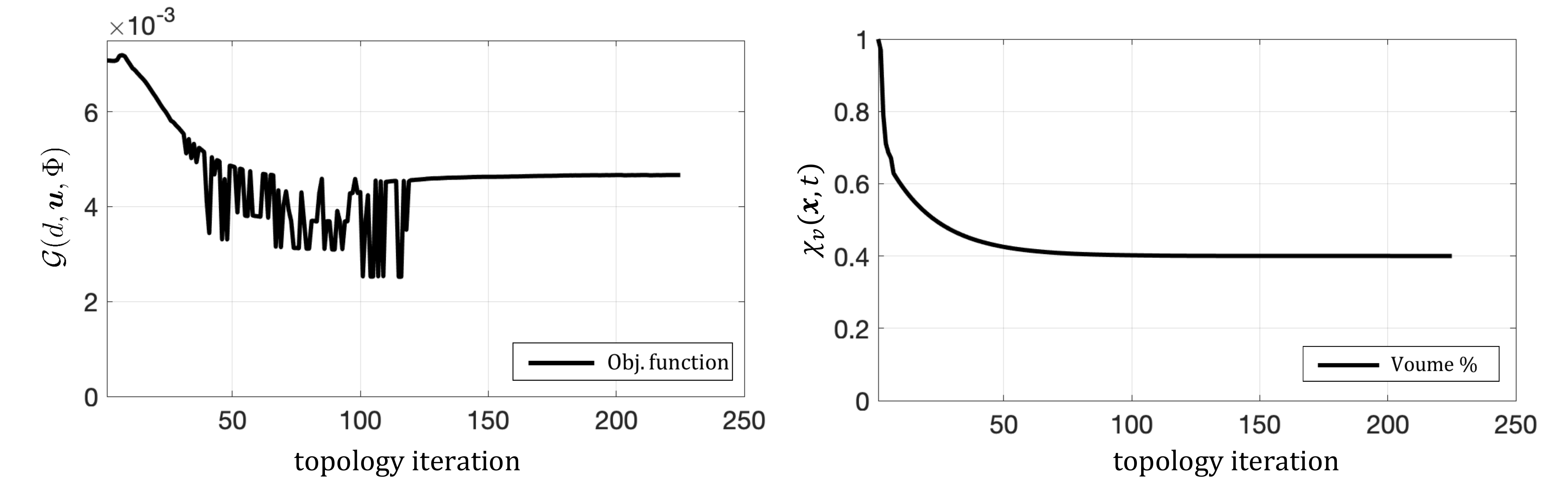}}  
	\vspace*{-0.7cm}
	\caption*{\hspace*{4.3cm}(a)\hspace*{8cm}(b)\hspace*{2cm}}
	\caption{Example 4 (Case b). Convergence history of the ductile fracture for the (a) objective
		function, and (b) volume constraint function.}
	\label{Exm4_conv2}
\end{figure}

\begin{figure}[!t]
	\subfloat{\includegraphics[clip,trim=5cm 8cm 6.5cm 11cm, width=5.5cm]{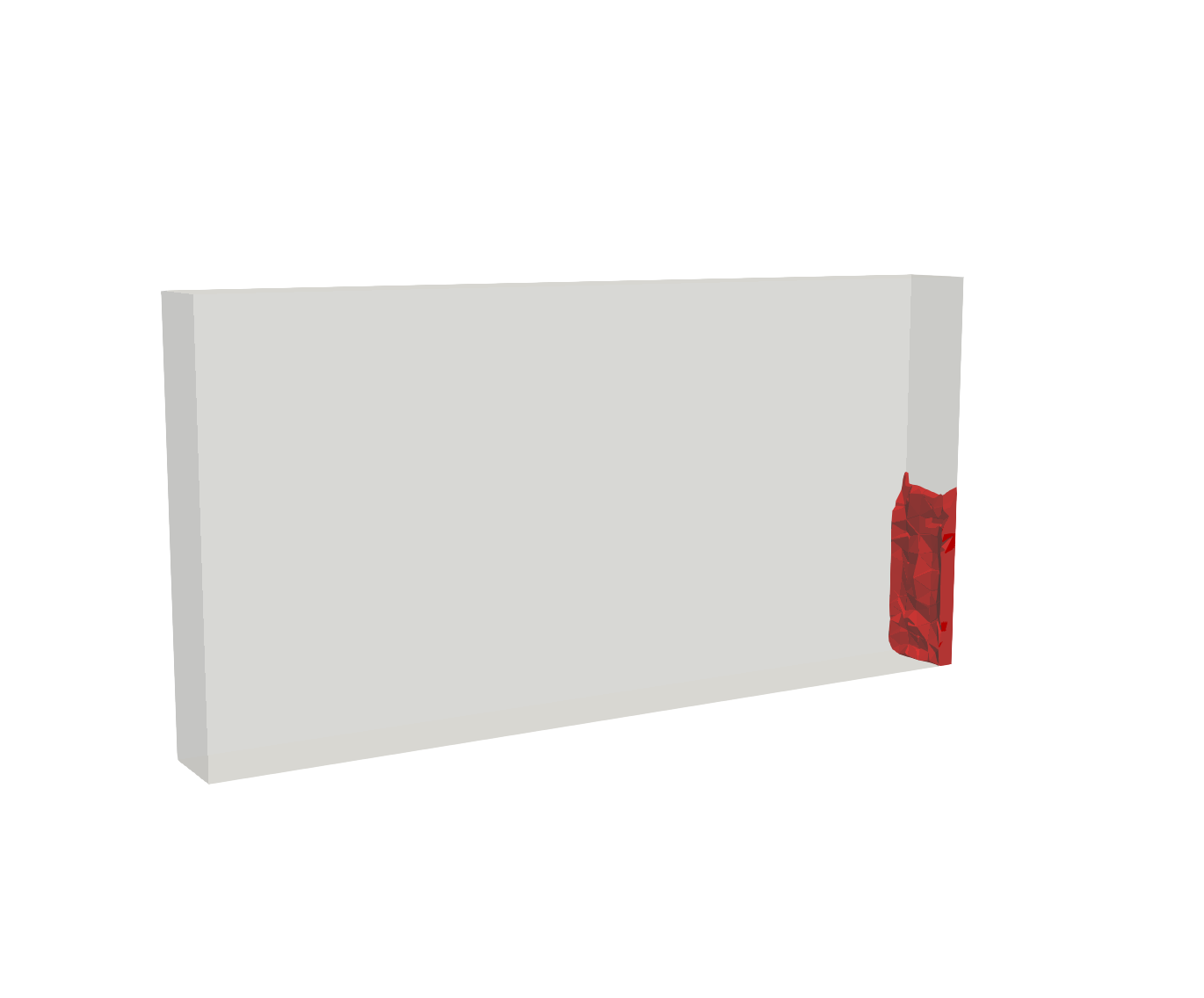}}   		\subfloat{\includegraphics[clip,trim=5cm 8cm 6.5cm 11cm, width=5.5cm]{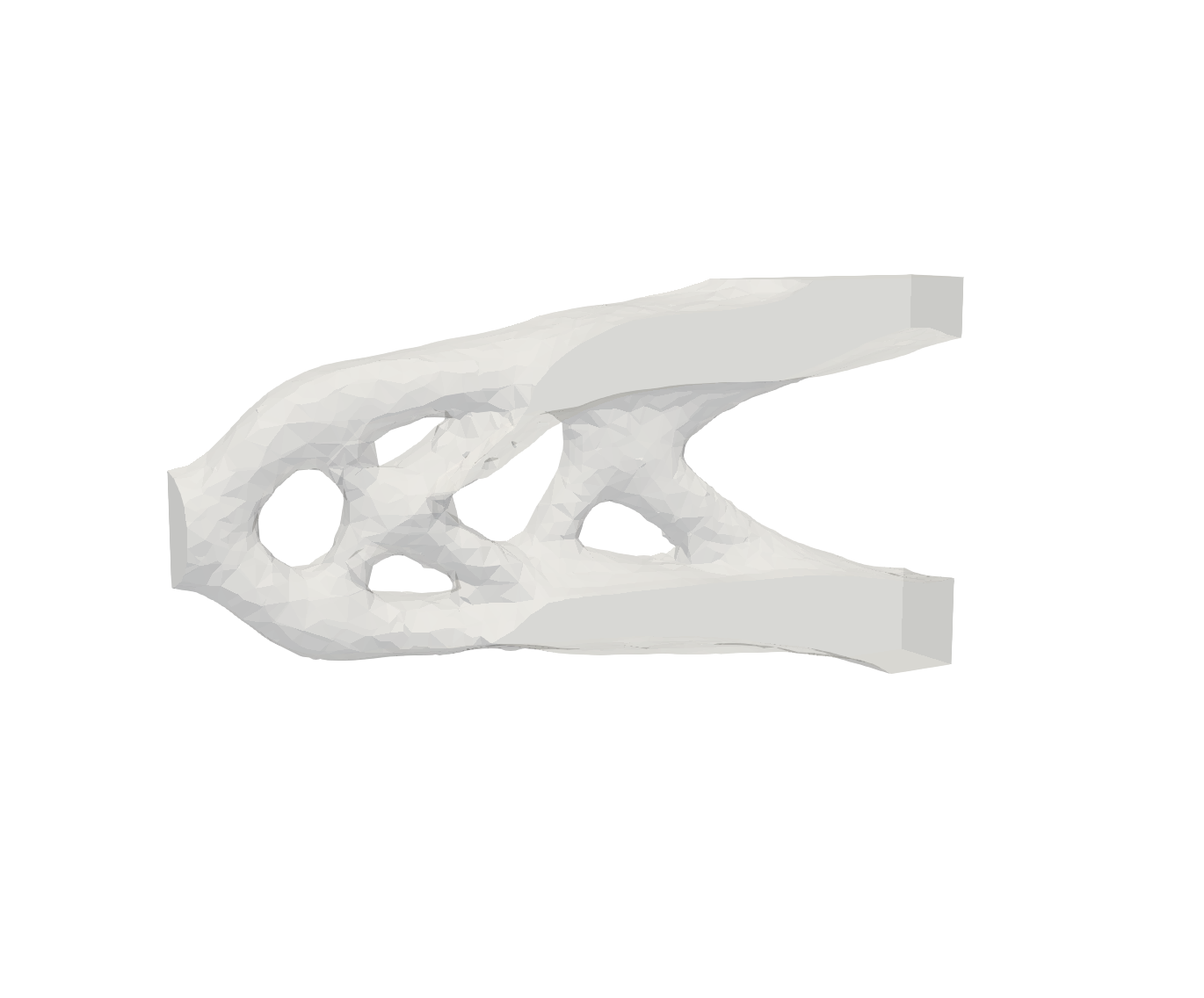}}   			\subfloat{\includegraphics[clip,trim=5cm 8cm 6.5cm 11cm, width=5.5cm]{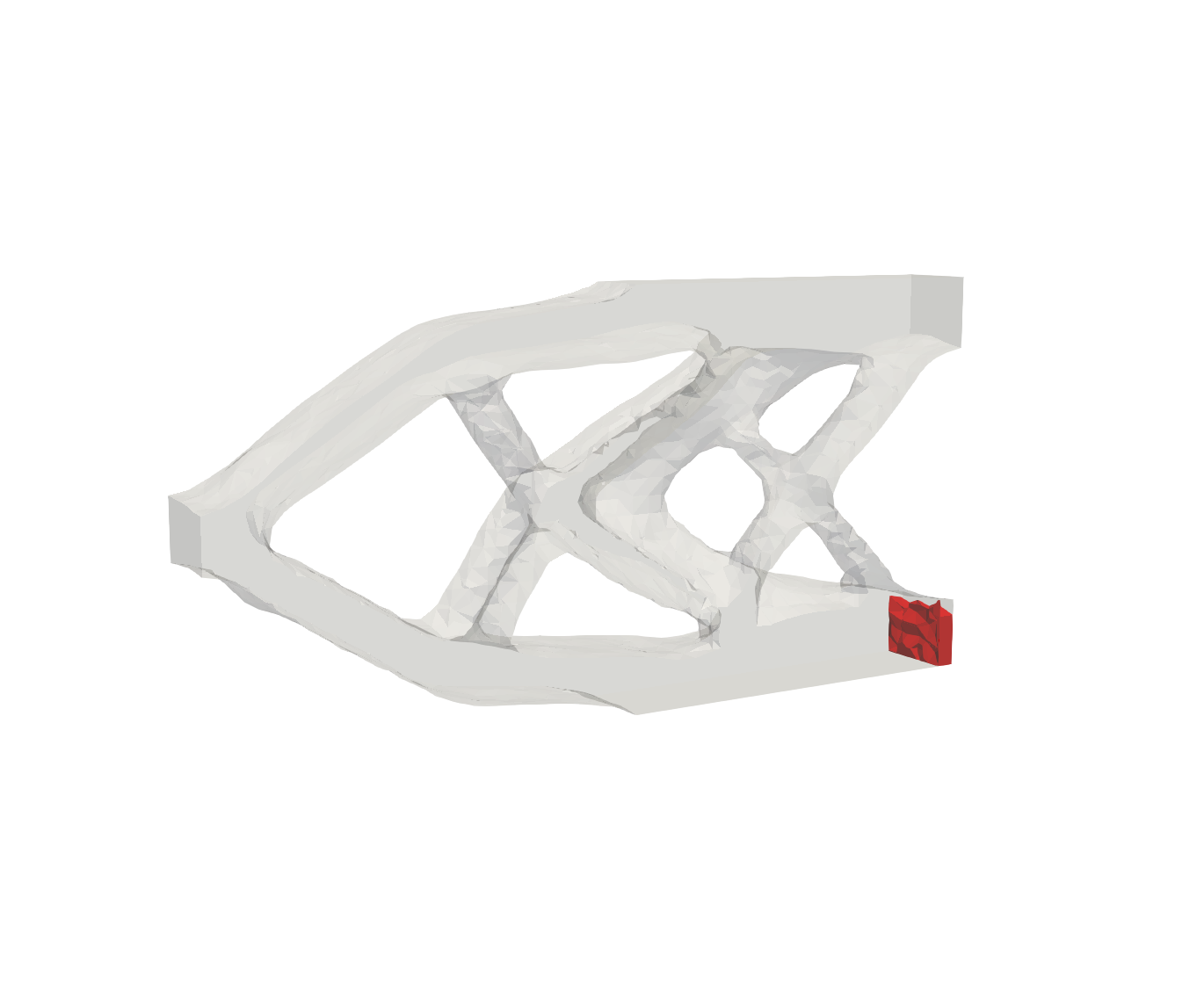}}
	\vspace{-0.3cm}
	\caption*{\hspace*{1.3cm}(a)\hspace*{5cm}{(b)}\hspace*{5cm}{(c)}}
	\caption{Example 4 (Case a). The crack phase-field profile for different optimal topology layout of volume ratio based on (a) non-optimized, (b) brittle fracture model, and (c) linear elasticity results undergoes brittle fracture.}
	\label{Exm4_ED_d}
\end{figure}

\begin{figure}[t!]
	\caption*{\underline{$\chi_v=1$}\hspace*{6cm}\underline{$\chi_v=0.94$}\hspace*{4cm}\underline{$\chi_v=0.91$}\hspace*{1cm}}
	\vspace{-0.1cm}
	\subfloat{\includegraphics[clip,trim=6cm 8cm 6cm 10cm, width=5.6cm]{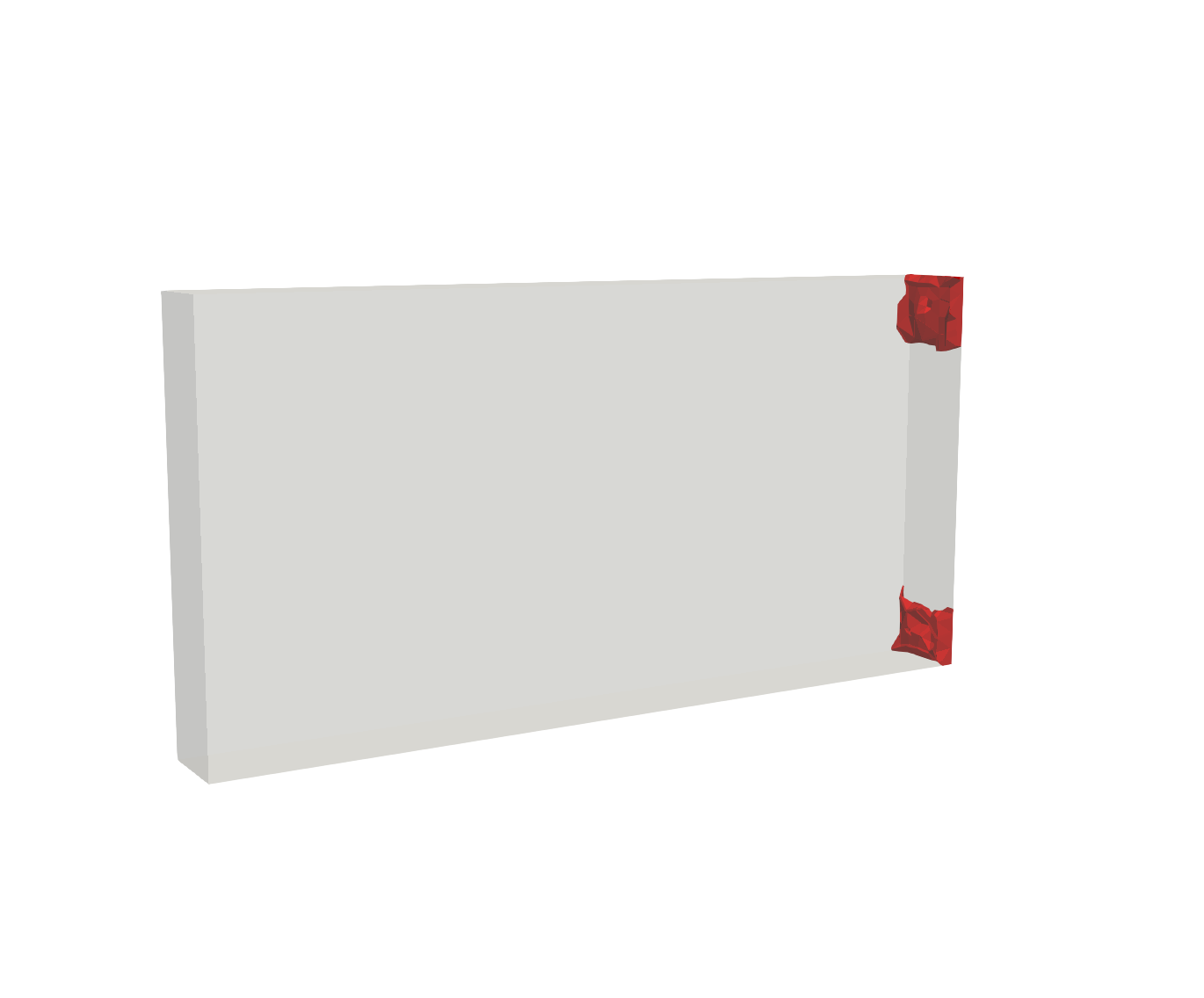}}   		\subfloat{\includegraphics[clip,trim=6cm 8cm 6cm 10cm, width=5.6cm]{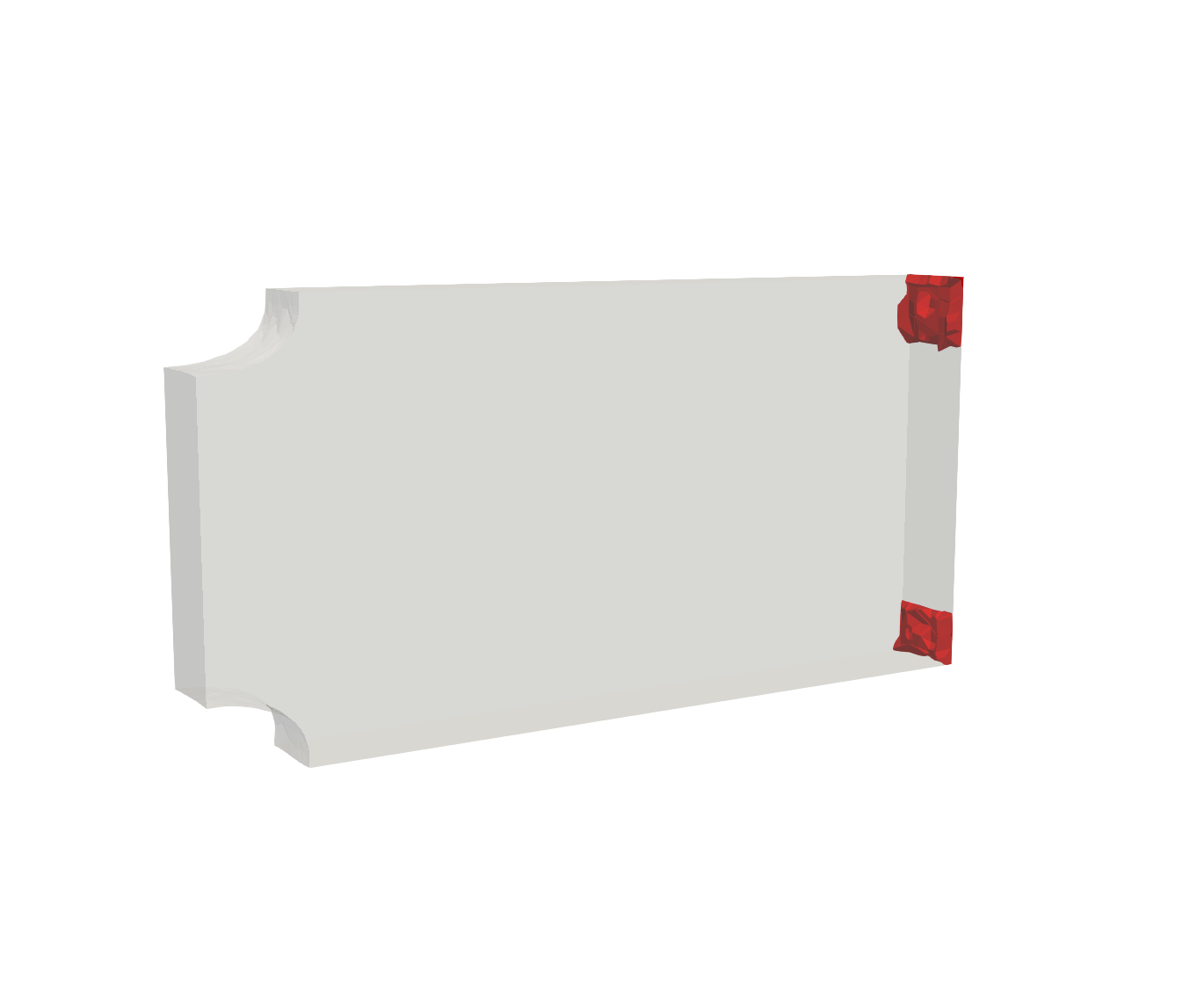}}   			\subfloat{\includegraphics[clip,trim=6cm 8cm 6cm 10cm, width=5.6cm]{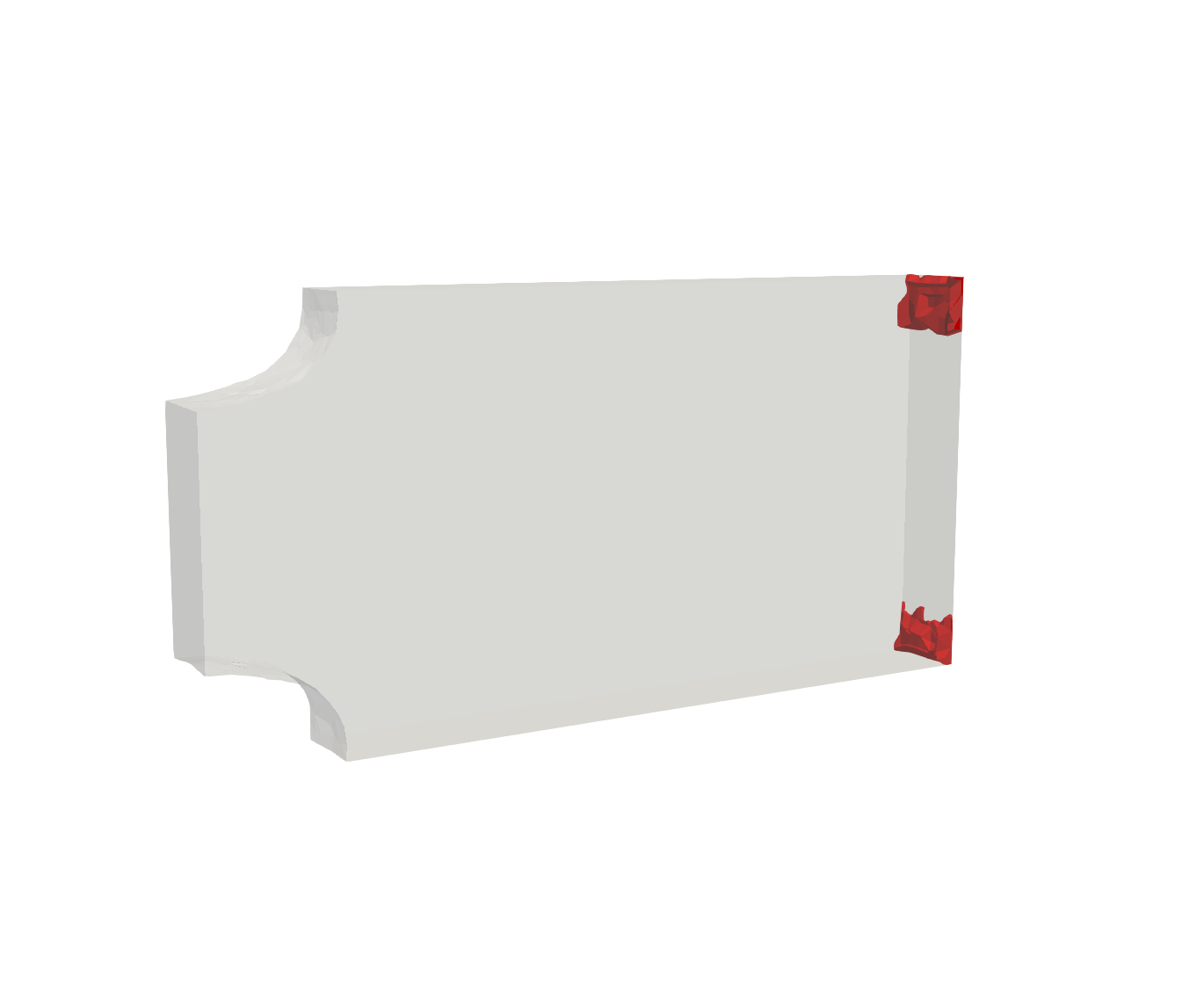}}
	\caption*{\underline{$\chi_v=0.88$}\hspace*{5cm}\underline{$\chi_v=0.40$}}
	\centering
	\vspace{-0.2cm}
	\subfloat{\includegraphics[clip,trim=6cm 8cm 6cm 10cm, width=5.6cm]{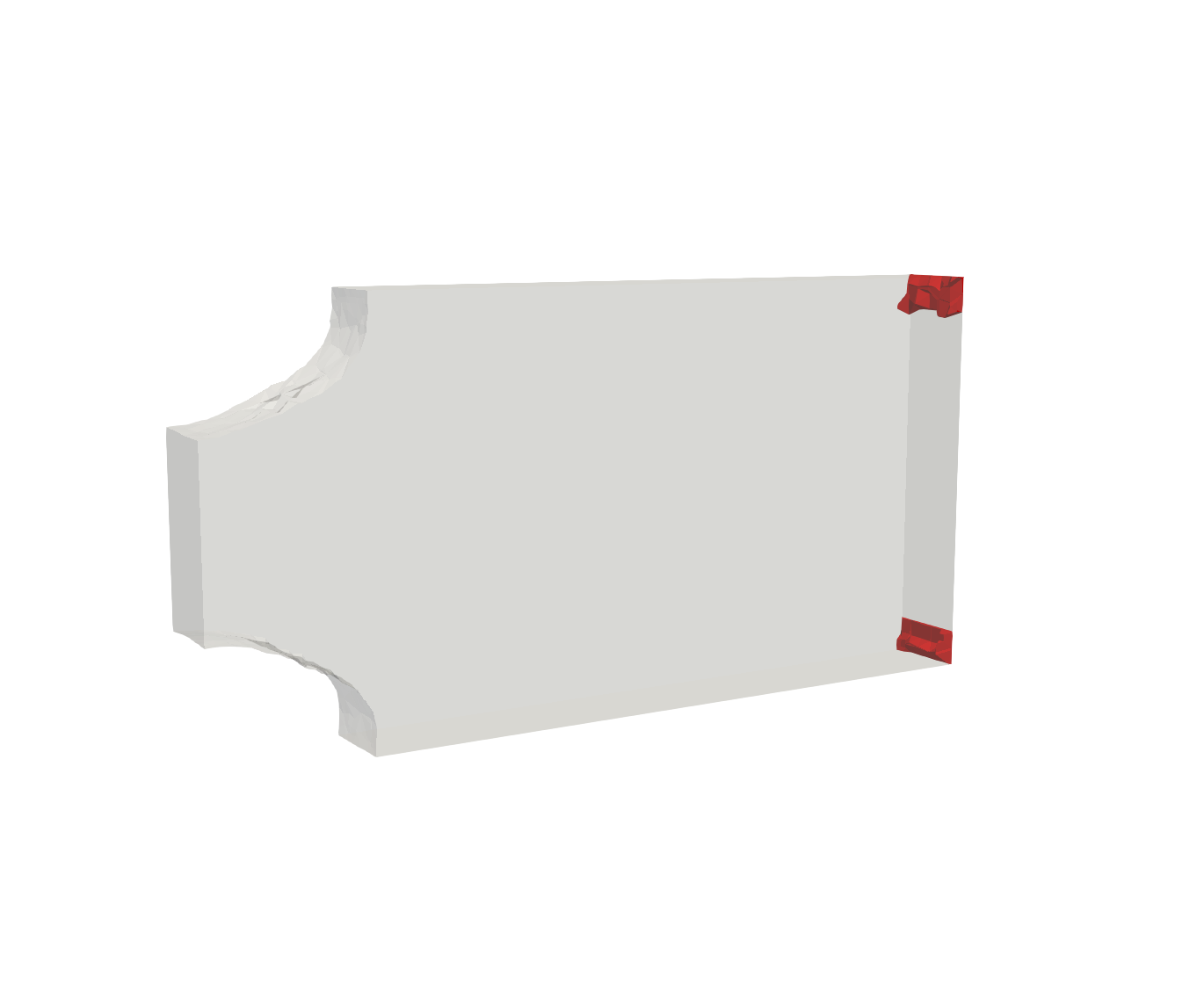}}   	  			\subfloat{\includegraphics[clip,trim=6cm 8cm 6cm 10cm, width=5.6cm]{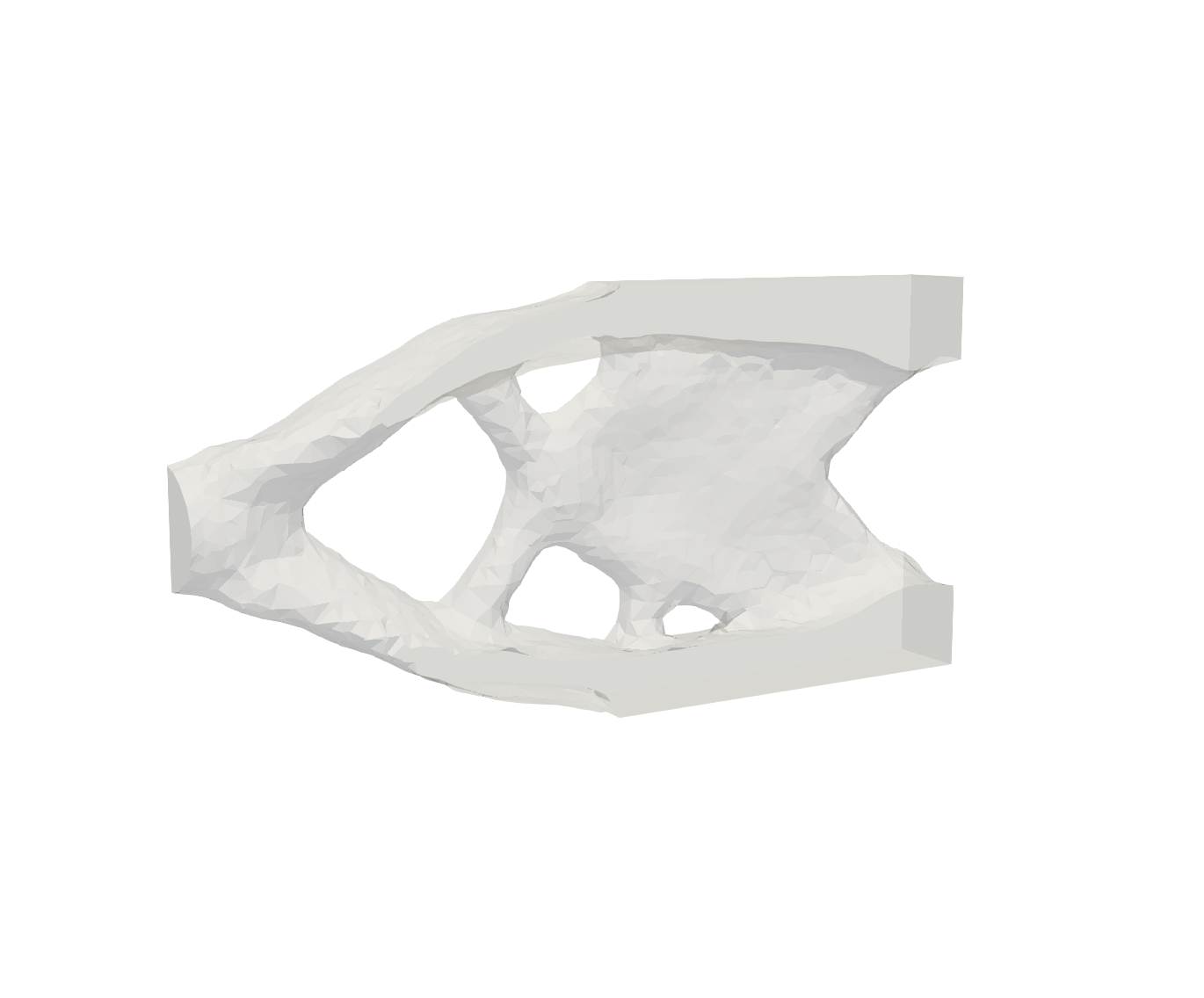}}	
	\caption{Example 4 (Case b). The crack propagation patterns for evolution history of different optimal layouts of volume ratio for ductile fracture model.}
	\label{Exm4_EPD_d}
\end{figure}

\begin{figure}[t!]
	\caption*{\underline{$\chi_v=1$}\hspace*{6cm}\underline{$\chi_v=0.94$}\hspace*{4cm}\underline{$\chi_v=0.91$}\hspace*{1cm}}
	\vspace{-0.2cm}
	\subfloat{\includegraphics[clip,trim=6cm 8cm 6cm 10cm, width=5.6cm]{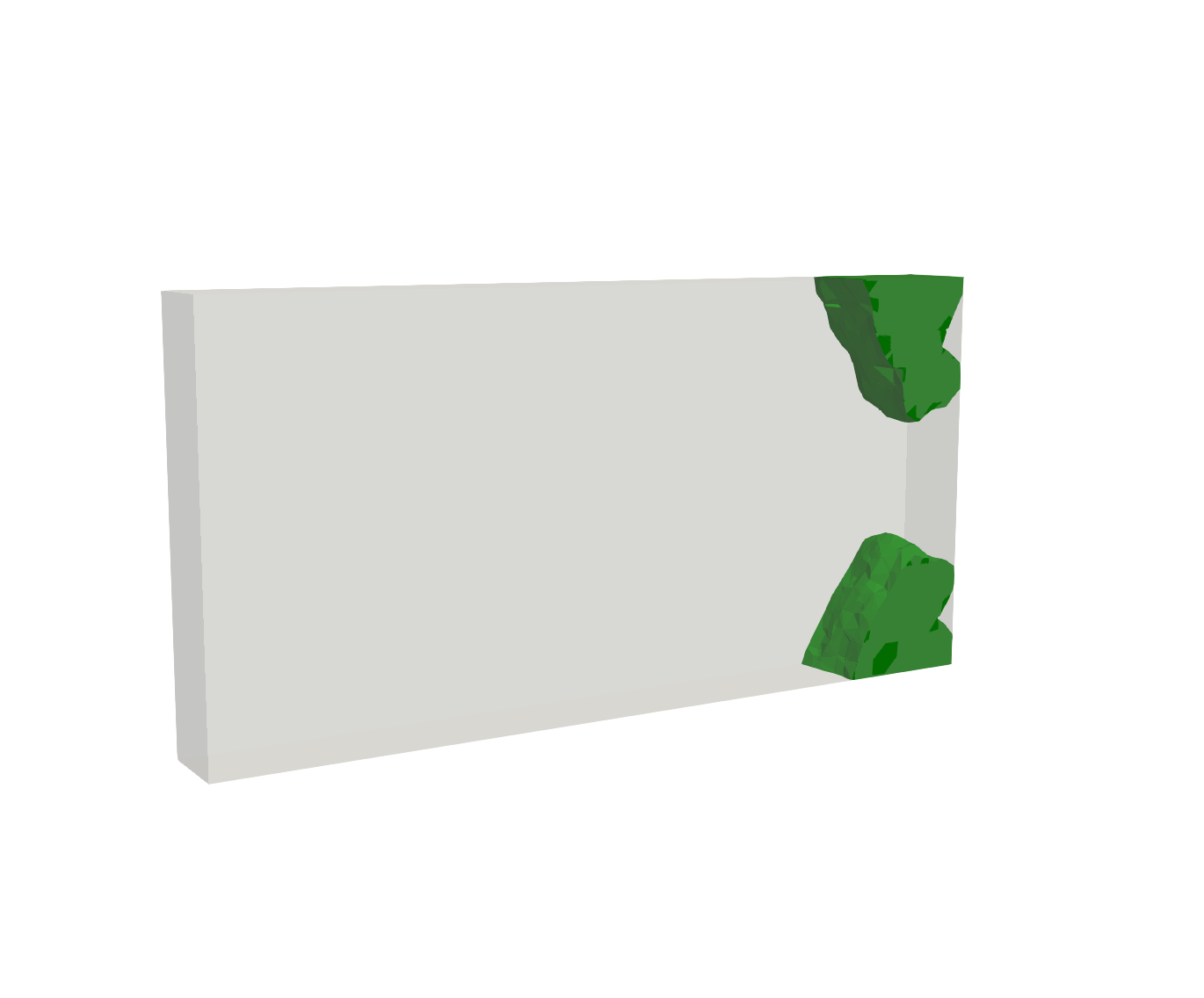}}   		\subfloat{\includegraphics[clip,trim=6cm 8cm 6cm 10cm, width=5.6cm]{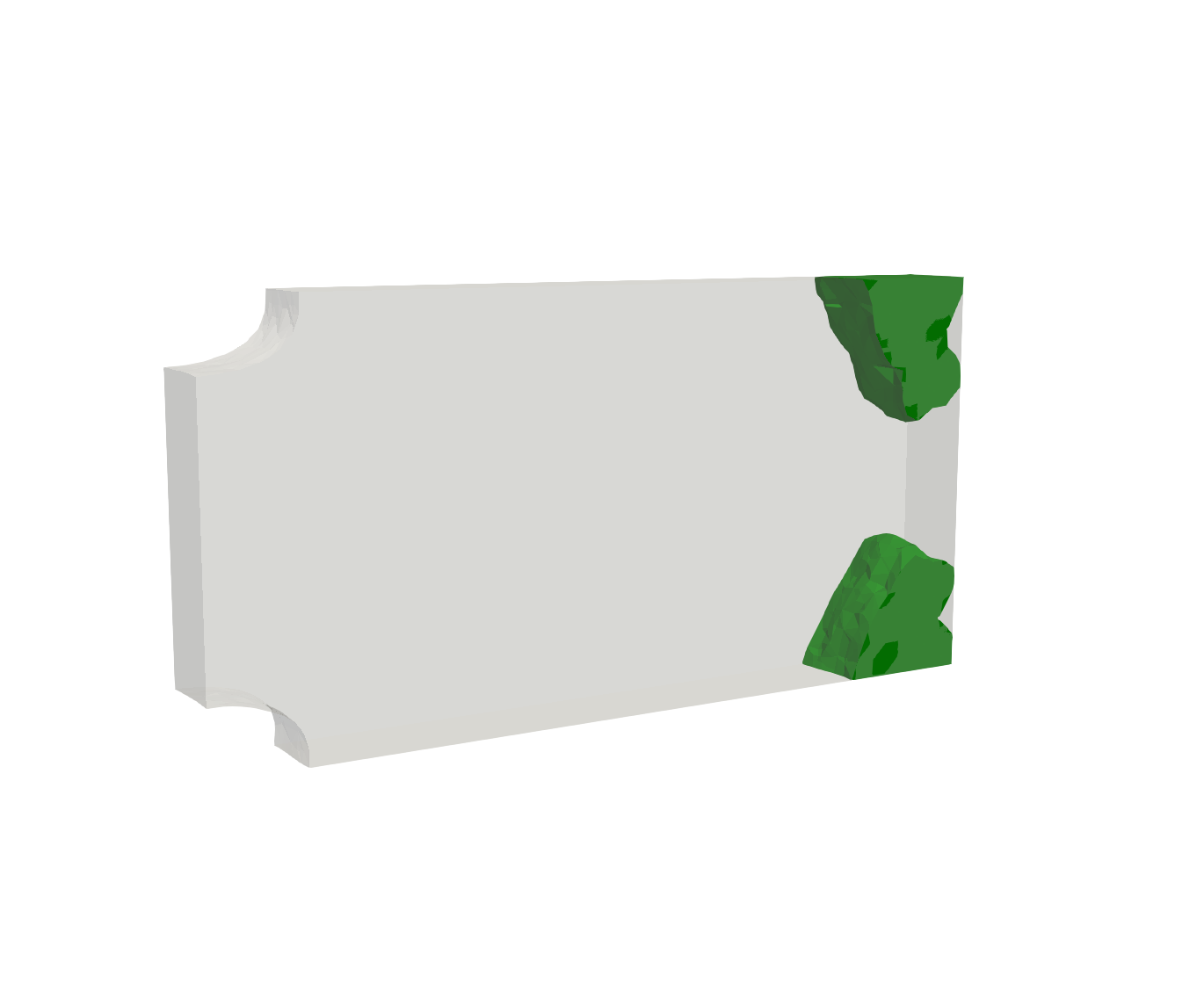}}   			\subfloat{\includegraphics[clip,trim=6cm 8cm 6cm 10cm, width=5.6cm]{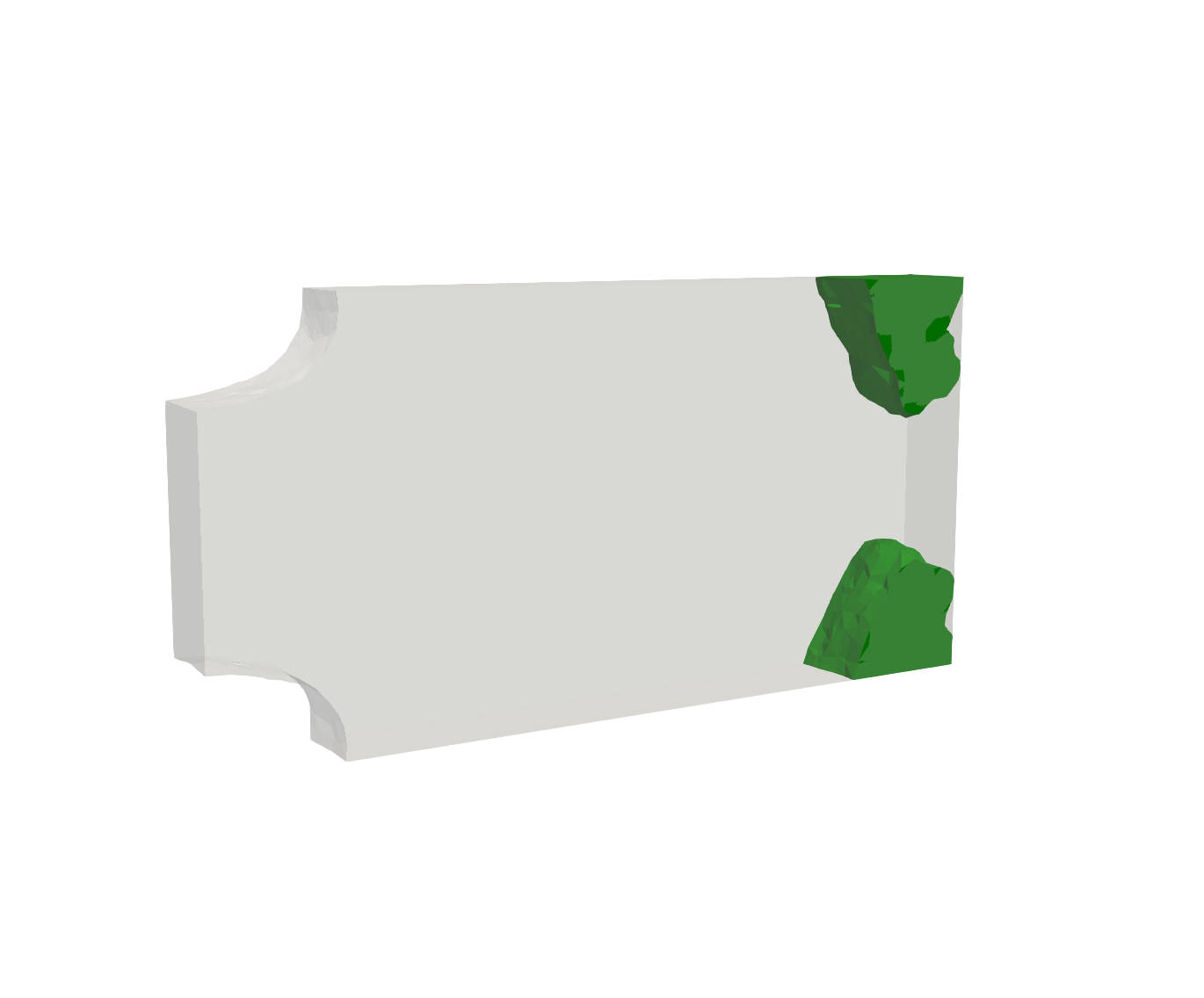}}	
	\caption*{\underline{$\chi_v=0.88$}\hspace*{6cm}\underline{$\chi_v=0.81$}\hspace*{4cm}\underline{$\chi_v=0.40$}\hspace*{1cm}}
	\vspace{-0.2cm}
	\subfloat{\includegraphics[clip,trim=6cm 8cm 6cm 10cm, width=5.6cm]{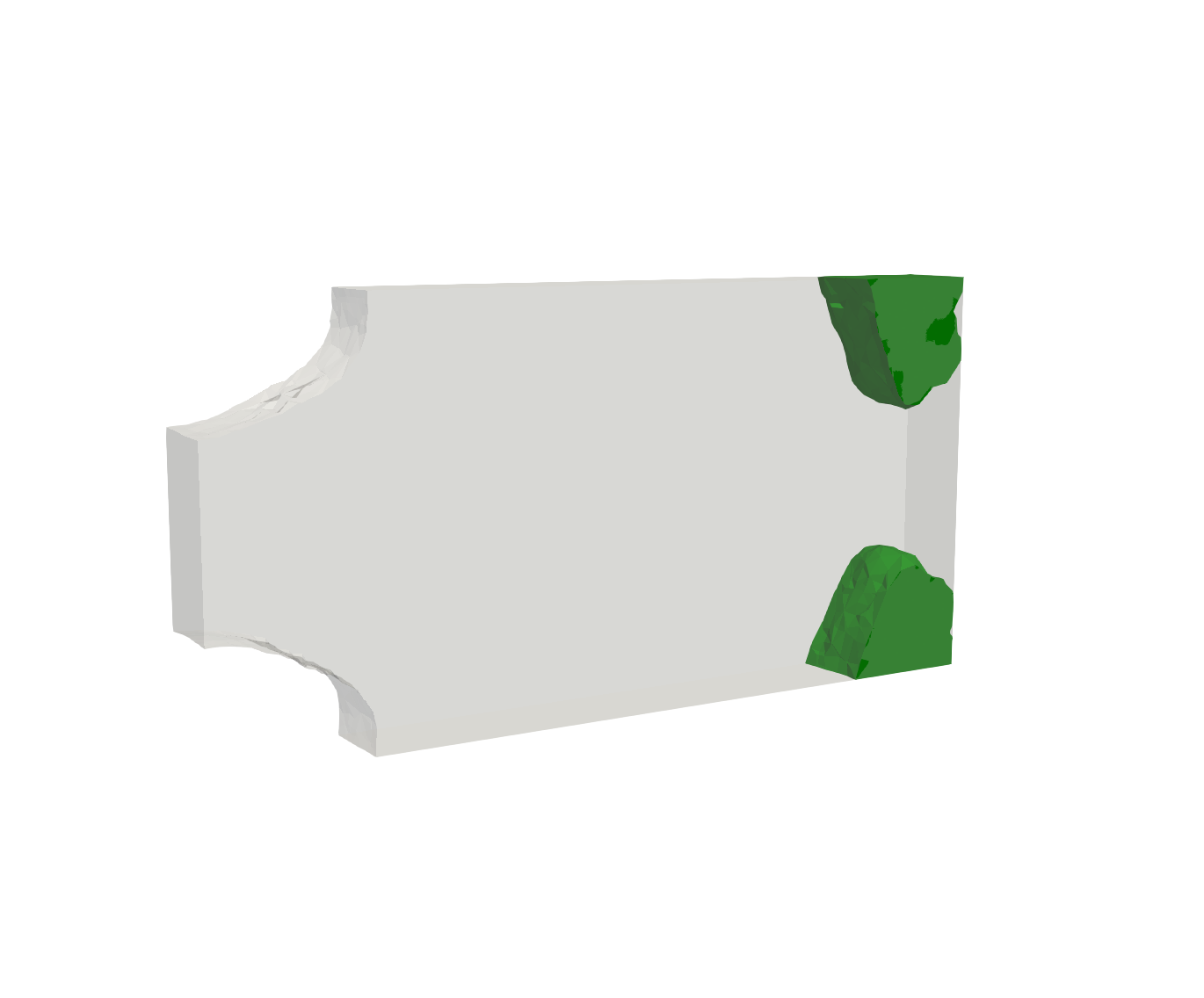}}   		\subfloat{\includegraphics[clip,trim=6cm 8cm 6cm 10cm, width=5.6cm]{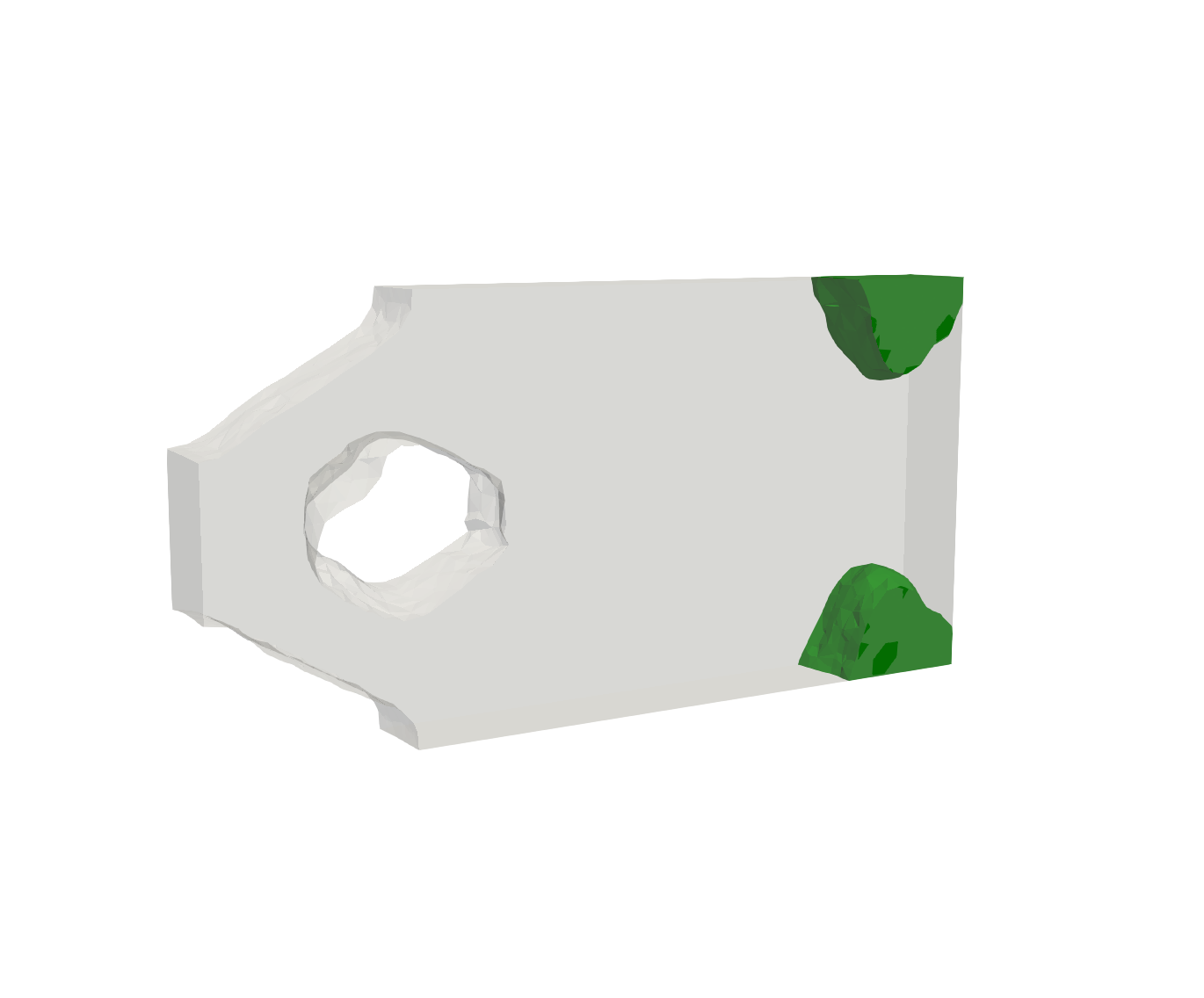}}   			\subfloat{\includegraphics[clip,trim=6cm 8cm 6cm 10cm, width=5.6cm]{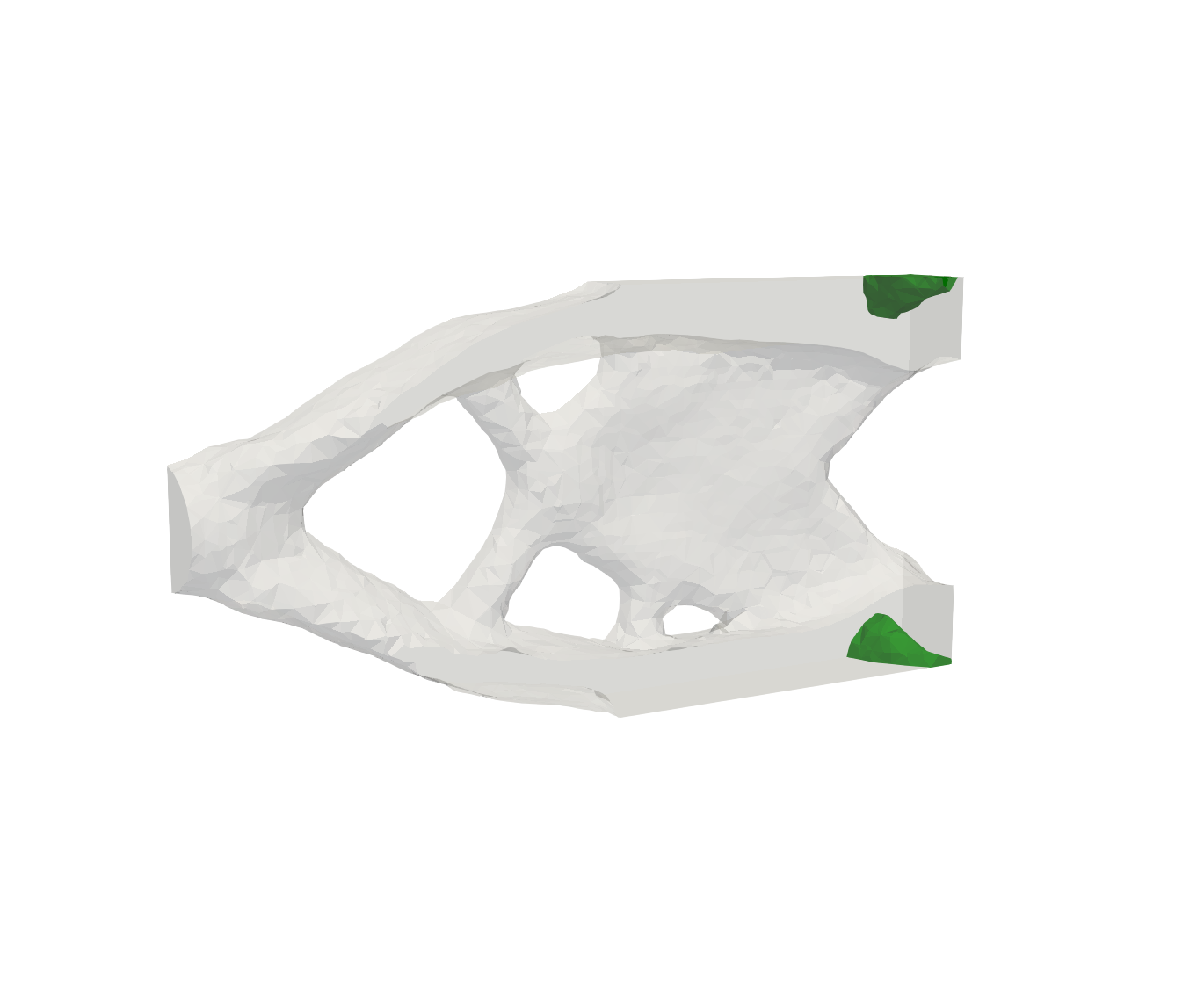}}	
	\caption{Example 4 (Case b).  Computed hardening value (plastic zone) for evolution history of the different optimal layouts of volume ratio for ductile fracture model.}
	\label{Exm4_EPD_alpha}
\end{figure}

\begin{figure}[t!]
	\centering
	\vspace{-0.1cm}
	\subfloat{\includegraphics[clip,trim=6cm 12cm 7.7cm 5.6cm, width=8.2cm]{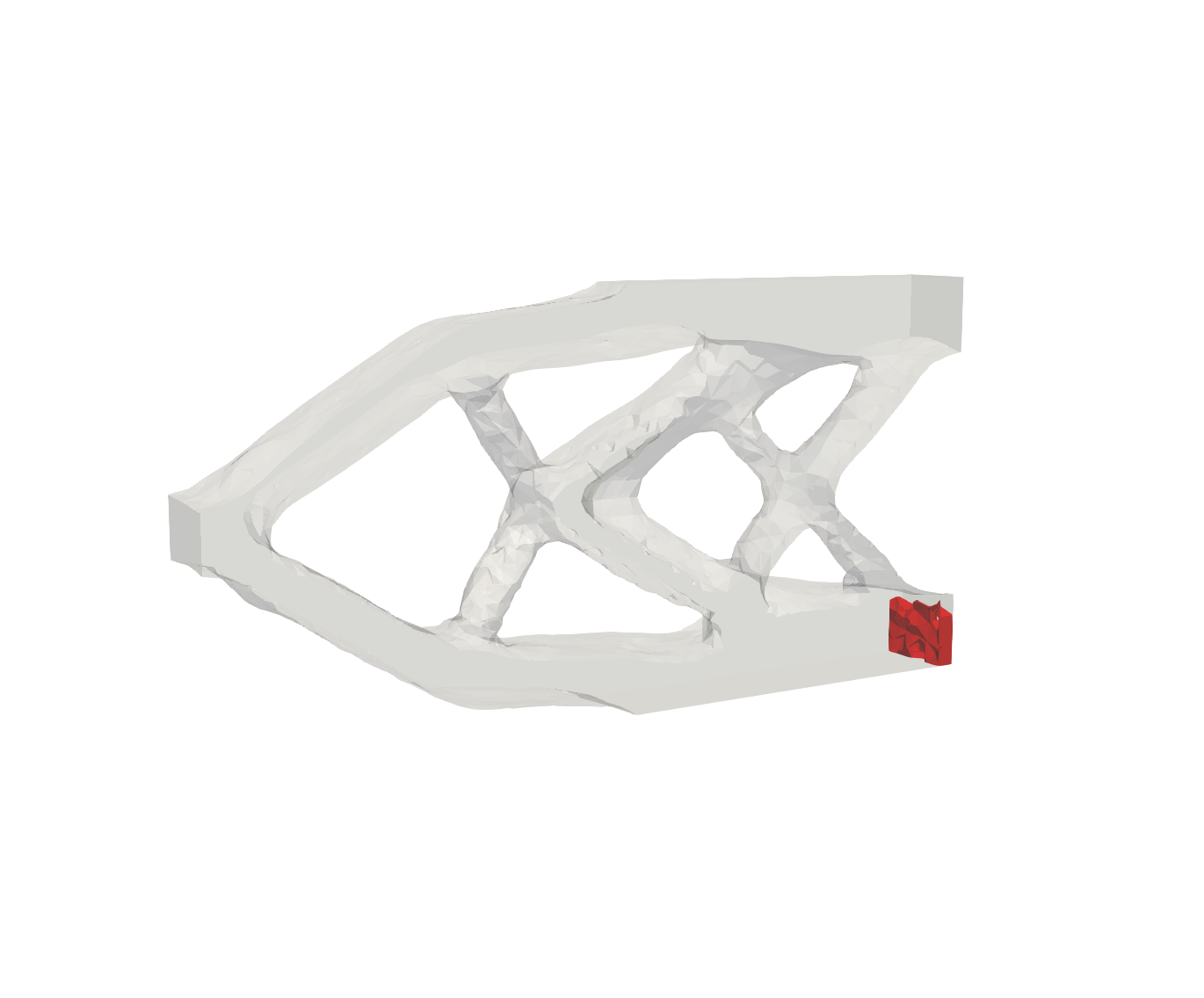}}   	\subfloat{\includegraphics[clip,trim=6cm 12cm 7cm 10cm, width=8.2cm]{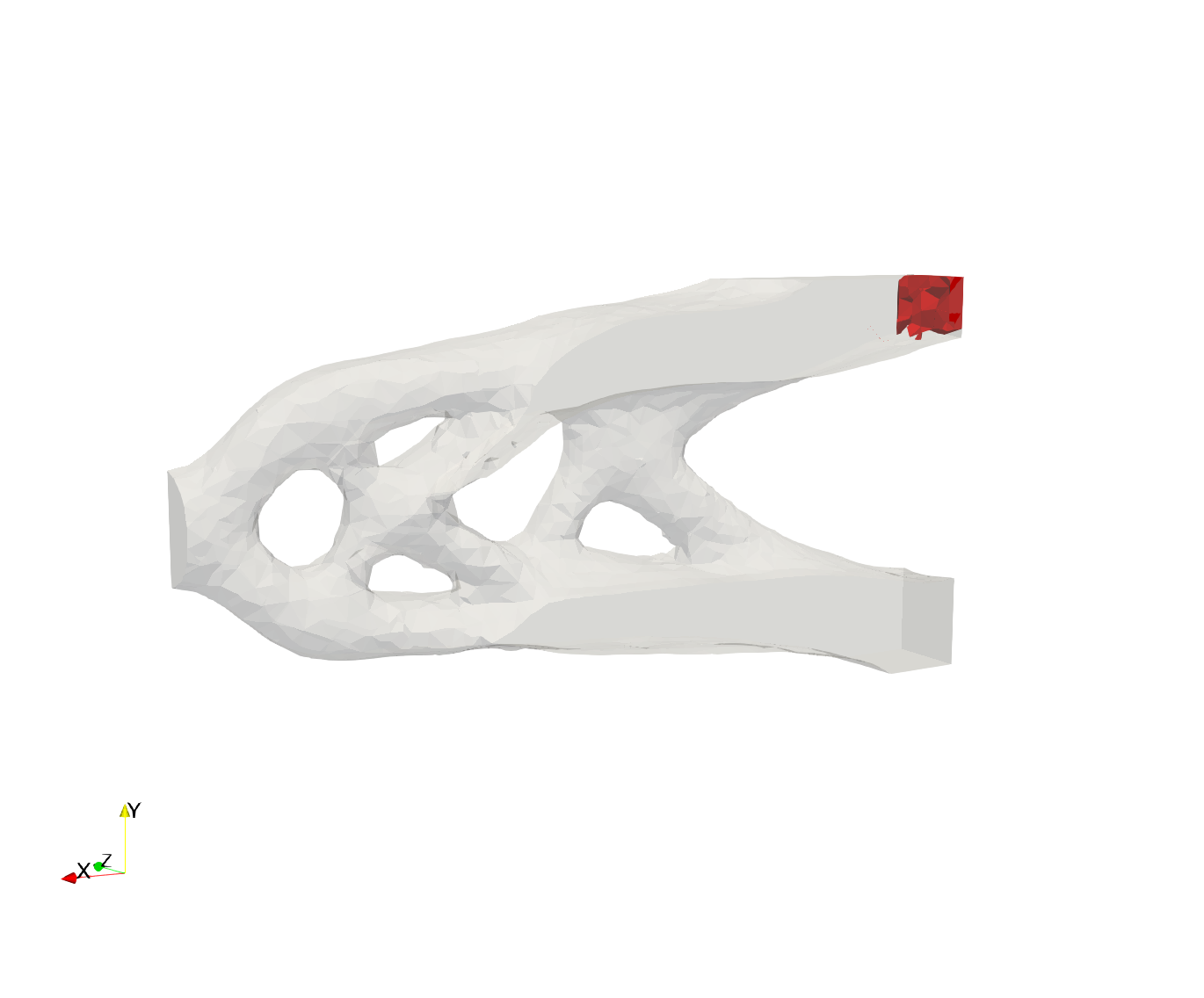}}		
	\caption*{\hspace*{1.8cm}(a)\hspace*{8cm}(b)}
	\caption{Example 4 (Case b). Crack propagation patterns for ductile fracture model in Example 4: (a) Pure elasticity, and (b) brittle fracture model.}
	\label{Exm4_d_E_ED}
\end{figure}

\begin{figure}[!t]
	\centering
	{\includegraphics[clip,trim=2cm 23cm 0cm 2cm, width=16cm]{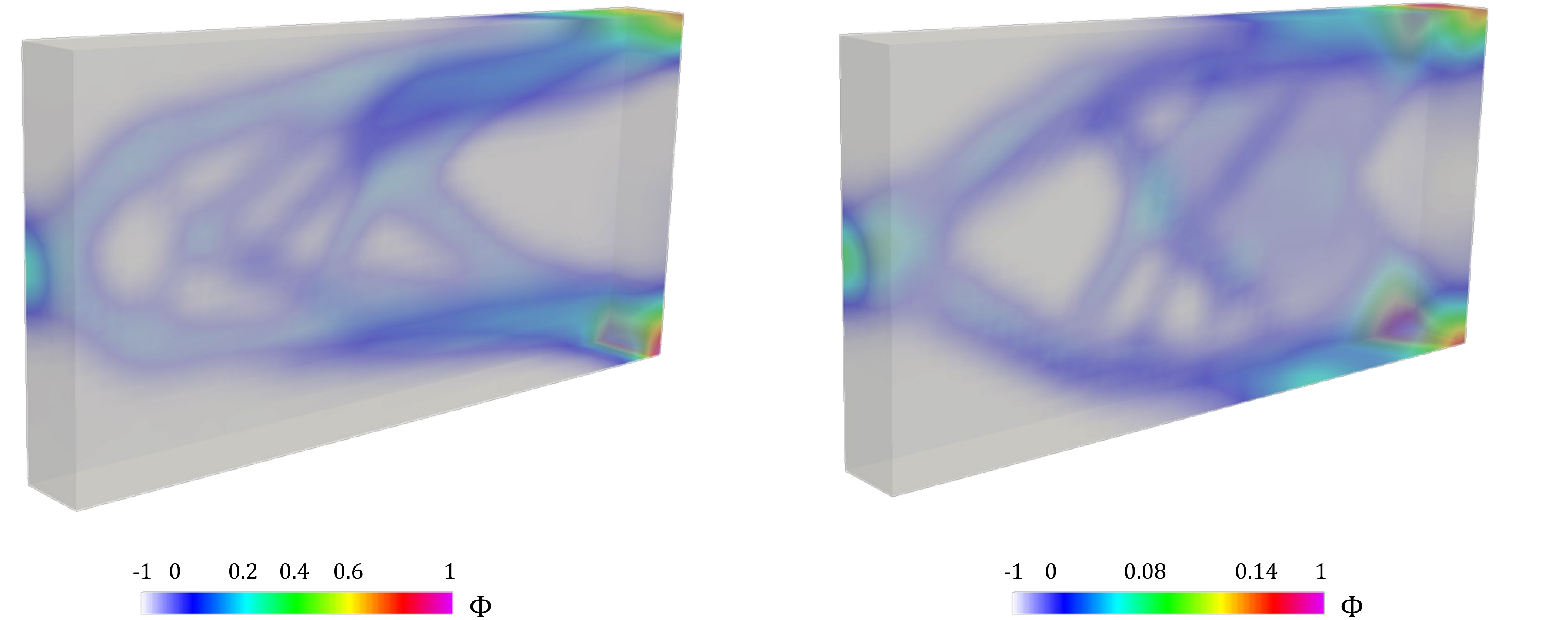}}  
	\vspace*{-0.7cm}
	\caption*{\hspace*{3cm}(a)\hspace*{9cm}(b)\hspace*{2cm}}
	\caption{Example 4. Computed topological field $\Phi$ at the final topology optimization iteration at $\chi_v=0.4$ for (a) Case a, and (b) Case b.}
	\label{Exm4_LSM}
\end{figure}

\begin{figure}[!t]
	\centering
	\vspace{-0.1cm}
	\subfloat{\includegraphics[clip,trim=6cm 12cm 7.7cm 5.6cm, width=8.2cm]{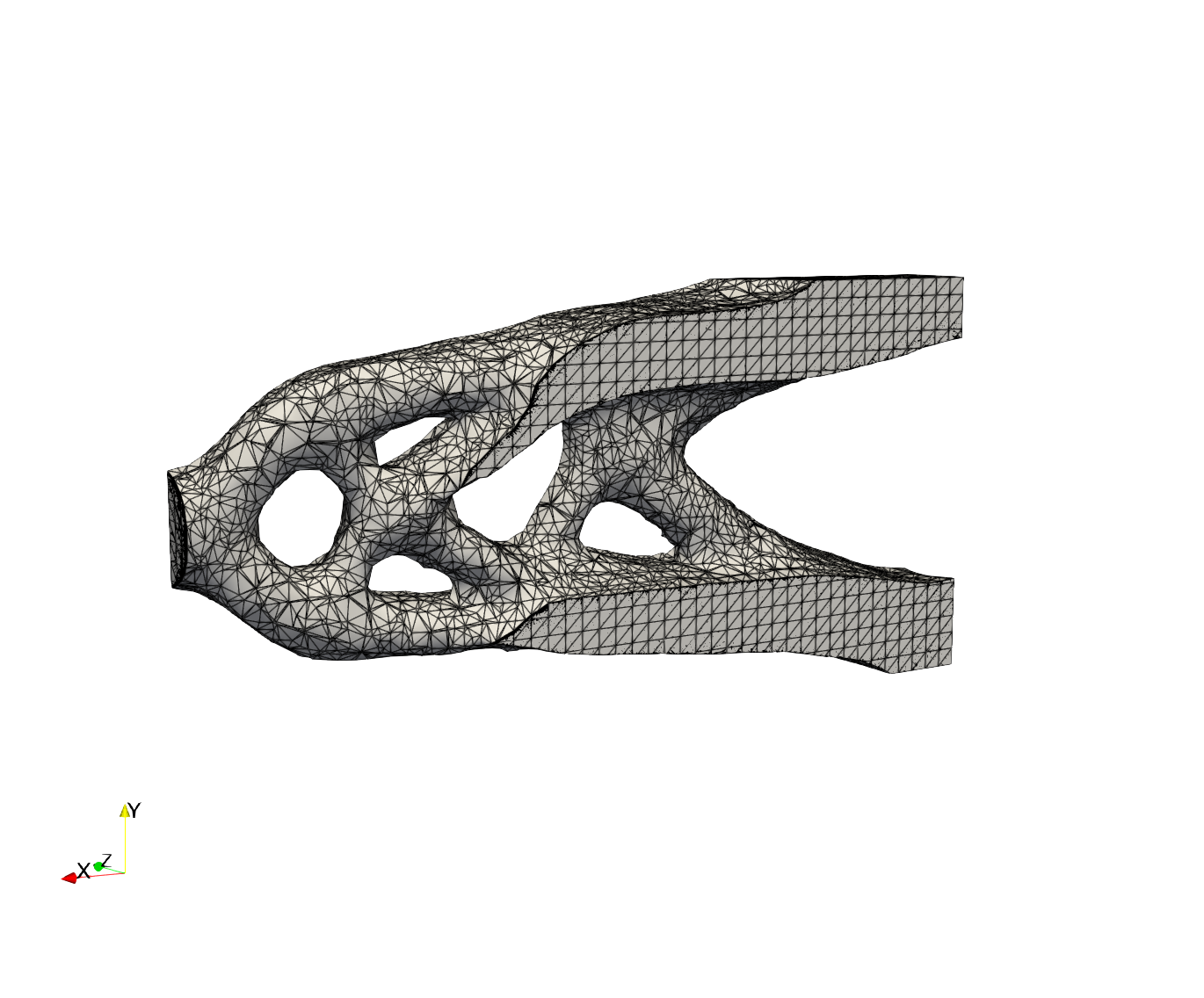}}   	\subfloat{\includegraphics[clip,trim=6cm 12cm 7cm 10cm, width=8.2cm]{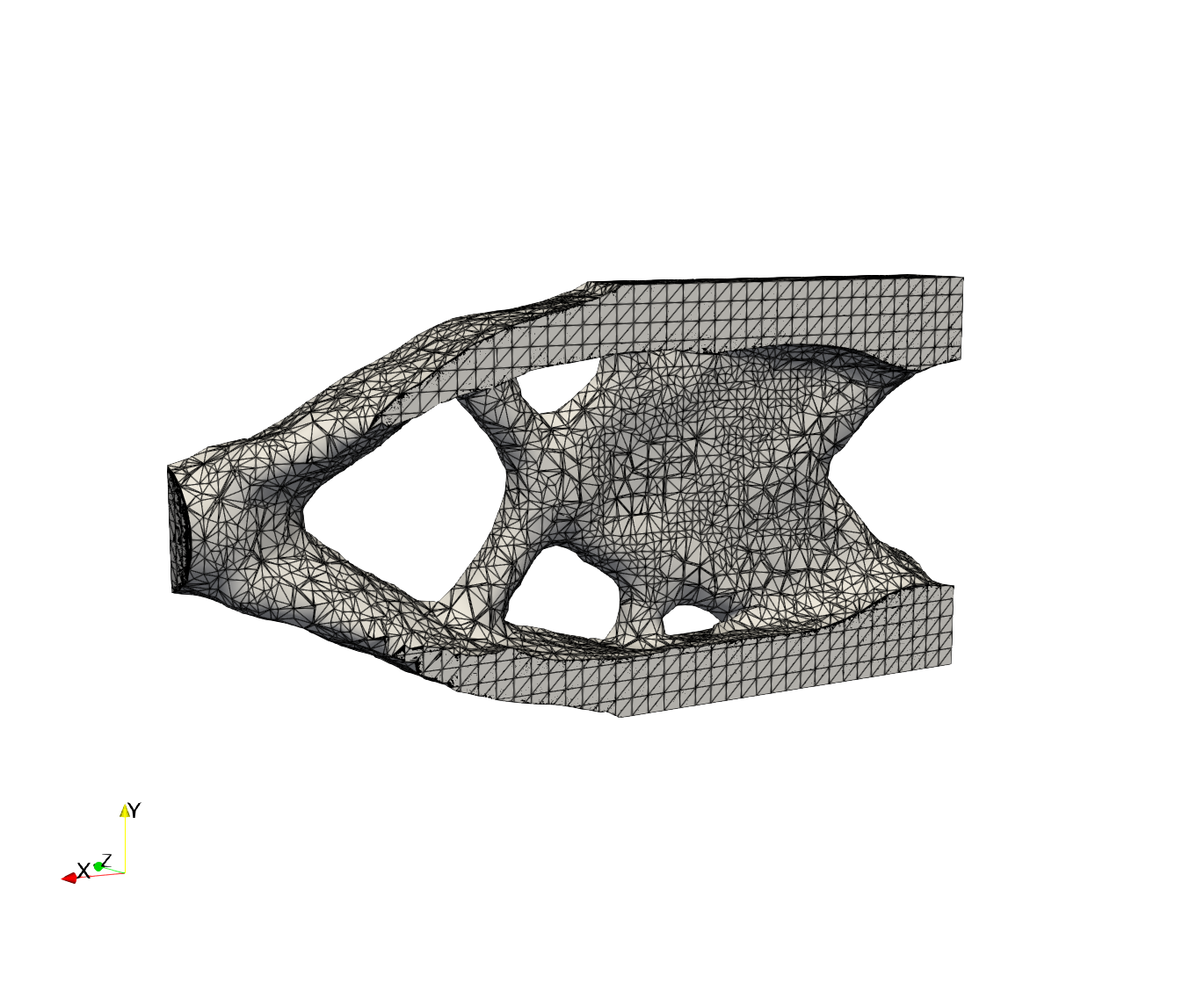}}		
	\caption*{\hspace*{1.8cm}(a)\hspace*{8cm}(b)}
	\caption{Example 4. Finite element discretization of the optimal layouts when $\chi_v=0.4$ (a) Case a, and (b) Case b.}
	\label{Exm4_mesh}
\end{figure}

\begin{itemize}
	\item \textbf{Case a.} First, topology optimization of brittle fracture for a cantilever beam under shear loading is investigated. The material properties in Table \ref{material-parameters} are used, however we set $\sigma_Y=10^{16}$. The evolution history of the optimal topology is shown in Figure \ref{Exm4_ED_phi} for different volume ratios up to $\chi_v=0.40$. For sake of comparison, topology optimization of linear elasticity results are also shown in Figure \ref{Exm4_E_phi}. Quantitative load-displacement curves for the non-optimized, elasticity and topology optimization due to brittle fracture are shown in Figure \ref{Exm4_LD}. At this point, it is necessary to remark that if we only obtain optimum topology due to linear elasticity, it does not help avoiding a fractured state in the structure, and thus we should consider fracture constraint in our optimization problem.  The convergence history of the objective and volume constraint functions are depicted in Figure \ref{Exm4_conv1}, which clearly reach a converged state. Finally, we investigate the crack phase-field profiles due to the optimal layout obtained through the brittle fracture, and elasticity results along with non-optimized results which are shown in Figure \ref{Exm4_ED_d}. \noii{We note that since, the loading is applied upward to the geometry, the bottom areas of the cantilever beam is stretched  due to tensile stresses which is below the neutral axis.}
	\item \textbf{Case b.} Second, topology optimization of ductile failure for the cantilever beam is examined.  Accordingly, the evolutionary history of the optimal topology is shown in Figure \ref{Exm4_EPD_phi}. For a better insight into the proposed model, quantitative load-displacement curves for the non-optimized, elasticity, and topology optimization due to brittle fracture, and ductile fracture are studied in Figure \ref{Exm4_LD}(b). Here, by topology optimization due to brittle fracture, we mean first to obtain the final optimum layout from the optimization process due to brittle fracture, then compute the load-displacement curve through the ductile phase-field. This curve is shown in black color in Figure \ref{Exm4_LD}(b). In fact, we claim that topology optimization obtained due to ductile fracture here has the best performance to avoid fracture state for elastic-plastic material, which is shown in Figure \ref{Exm4_LD}(b). Thus, topology optimization obtained due to ductile fracture (the red color in Figure \ref{Exm4_LD}(b)) will exhibit no softening due to fracture (no crack initiation), showing the great potential of the proposed model. 
	
	Another impacting factor that should be noted is the convergence history of objective and volume constraint functions, depicted in Figure \ref{Exm4_conv2}. The crack phase-field profiles due to the optimal layout of different volume ratios are shown in Figure \ref{Exm4_EPD_d}. Additionally, the equivalent plastic strain $\alpha(\Bx,t)$ of the evolutionary history of the optimal layout for different volume ratios \grm{is} depicted in Figure \ref{Exm4_EPD_alpha}. We note that the maximum equivalent plastic strain which is shown by the green color occurs near fixed supports on the left surface. It is worth noting that as long as the optimum layout approaches $\chi_v=0.40$, the fracture area and the plastic zone will be drastically reduced. Indeed, this demonstrates the great efficiency of the proposed method. 
	
	We further examine the optimum layout obtained due to ductile fracture. To this end, the final optimum layout obtained from the optimization process due to elasticity result, and brittle fracture are used, to compute the ductile failure response. Results \grm{are} depicted in Figure \ref{Exm4_d_E_ED}. \noii{We note that since, we have the plastic deformation, the deviatoric stress plays a significant rule, while the hydrostatic pressure has no contribution in the plastic stage, i.e., $\text{tr}(\Bve_p)=0$. In fact, in the right top corner, we have the maximum deviatoric stress (and also the largest hardening quantity) thus fracture initiate/propagate in this area.}
	Evidently, both optimum configurations yield crack propagation near supports, while this is not the case for topology optimization due to ductile fracture, see Figure \ref{Exm4_EPD_d} for $\chi_v=0.40$. Hence, the result demonstrates the necessity of using an appropriate topology optimization framework constrained with essential terms (in this case, fracture equilibrium force vector). Since only the ductile fracture-resistance optimal layout has been able to eliminate crack propagation in the material domain. 
\end{itemize}
Accordingly, it is worth showing the computed regularized topological field $\Phi(\Bx,t)$ through the reaction-diffusion evolution equation for both brittle, and ductile fracture which is outlined in Figures \ref{Exm4_LSM}. Here,  the gray implies a non-material point, while the colorful area \grm{represents the} material point, and the purple color depicts the material boundary (zero contour of level set function). Lastly, the discretized final optimum layouts due to brittle and ductile fracture are depicted in Figure \ref{Exm4_mesh}(a-b), respectively. Evidently, the smooth topological interface is grasped.

    \cleardoublepage
\sectpa[Section7]{Conclusion}

This study outlined a robust and efficient topology optimization framework to enhance the brittle and ductile phase-field fracture resistance of structures. Herein, we have used a level-set based topology optimization as a mathematical method that seeks to find the best material distribution that satisfies the equilibrium, objective, and constraint functions. Accordingly, two different formulations were proposed. One required only the residual force vector of the deformation field as a constraint, while the second formulation accounted for the residual force vector of the deformation and phase-field fracture simultaneously. The first type enables computations performed with less implementation effort, however, it doesn't eliminate the fracture. In contrast, the second formulation directly constraints the fracture state and clearly leads to more optimal designs. Here, the optimal material distribution is realized by the reaction-diffusion based level-set topology optimization method, along with the consistent adjoint sensitivity analysis which has been derived for both formulations. Here, since we are dealing with a rate-dependent nonlinear boundary value problem the sensitivities are computed at every time increment, results in path-dependent sensitivity analysis.

Four numerical examples were presented to substantiate our algorithmic developments, in which the first two are related to brittle fracture and the last two deal with ductile phase-field fracture.  In the first example, we have shown that Formulation 2 has superior efficiency in comparison to Formulation 1, with optimum layout obtained due to linear elasticity, and non-optimized domain (the original one) which undergoes brittle failure. Thus, Formulation 2 showed its proficiency compared to others. In the second example, we intended to evaluate the mesh sensitivity in the proposed numerical method. To that end, three different mesh sizes under mixed-mode fracture loading were considered. The first important observation was that finer meshes lead to a qualitative distribution of materials with higher numbers of holes, which is capable \grm{of} providing more optimal results than coarser meshes. Nevertheless, finer meshes may lead to more complex structures that are more difficult to manufacture, as well as they require higher computational costs for both the forward and optimization problems. Additionally, it can be observed that the fracture limit point regardless of the mesh size,  is about 50.2\% greater than non-optimized results, and also optimum elasticity results undergoing brittle fracture. This shows the great importance of considering fracture constraint in our topology minimization approach, thus superior accuracy is observed. 
Findings showed that similar to stress-constrained topology optimization problems, it can be grasped from the results that optimal material distribution is such that stress singularity around the re-entrant is avoided as the optimizer removed the material around the re-entrant.

In the third and fourth examples, It was concluded that the final optimum layout due to topology optimization presents the stiffest response due to the ductile fractured state, although it has less volume ratio. Additionally, the equivalent plastic strain $\alpha(\Bx,t)$, and plastic zone were reduced. Thus, ductile fracture-resistant topology optimization prevented crack propagation in the material domain. In particular, in the last example, we examined the optimum layout obtained through topology optimization for both linear elasticity and brittle fracture undergoing ductile failure. It was shown that only a formulation that accounts for ductile fracture can avoid crack propagation (thus no degradation of the material). 

Several topics for further research emerging from the present study. Firstly, an efficient level-set based topology optimization due to ductile phase-field fracture could be extended towards large deformations. Secondly, by considering additional physics and in particular thermal effects involved in the model which can play as an additional source for the objective and constraint functions need to be taken into consideration. Lastly, since for the topology optimization of the fracturing material, the finite element treatment of the phase-field formulation is computationally demanding, thus, an idea of an adaptive multi-scale (spatial reduced-order modeling) formulation is particularly appealing. We would like to tackle these issues in our future work. 

\subsection*{Acknowledgment}
N. Noii  acknowledges the  Deutsche Forschungsgemeinschaft which was founded by the Priority Program \texttt{DFG-SPP 2020} within its second funding phase. 

\begin{Appendix}
	\setcounter{equation}{0}
	\renewcommand{\theequation}{A.\arabic{equation}}
	\setcounter{figure}{0}
	\renewcommand{\thefigure}{A.\arabic{figure}}
\subsection*{\noii{Appendix A. Sensitivity verification using central difference method }}
In order to verify the sensitivity analysis outlined in Section 3.4, we employ the portal frame structure under compression loading given in Example 3, with boundary value problems depicted in Figure 23. Here, the sensitivities are verified via a typical finite difference method through central difference approximation for the Lagrangian functional of the optimization problem in (90). More specifically, we compare the sensitivity results obtained in (106) with central difference method applied in (90).
\begin{equation}
\texttt{FEM:} \quad
{\widehat{v}_{\Phi}}\left| _{\Phi = 0} \right.=-\mathcal{G}=-D_{\Phi} \BfrakL
=-\sum_{\text{n}=1}^{N} \Big(
D_{\Phi}  \BfrakL^i\Big)
	\label{app1_eq:fem1}
\end{equation}
where as the sensitivity analysis is obtained by (106).  Thereafter, with central difference method, we have following approximation for each element $e\in(1,n_e)$ in our portal domain through:
\begin{equation}
\texttt{FDM:} \quad
		{\widehat{v}_{\Phi}}\left| _{\Phi = 0} \right.
		\approx
        -\sum_{\text{n}=1}^{N}\frac{\BfrakL^i (\Phi_1,\cdot\cdot\cdot,\Phi_e+\Delta\Phi_e,\cdot\cdot\cdot,\Phi_{n_e})-\BfrakL^i (\Phi_1,\cdot\cdot\cdot,\Phi_e-\Delta\Phi_e,\cdot\cdot\cdot,\Phi_{n_e})}{2\Delta \Phi_e}
	\label{app1_eq:fdm1}
\end{equation}
such that we set $\Delta \Phi_e=10^{-4}$. We now recapitulate the comparison of the finite element method with finite difference method. To do so, Figure \ref{Exm1_sens1} shows the sensitivity obtained by the analytical \req{app1_eq:fem1}, and numerical schemes \req{app1_eq:fdm1}. Due to the sensitivity analysis of the Lagrangian function, the values with negative signs have been appeared for some nodal points in domain and it can be seen in Figure  \ref{Exm1_sens1}. In fact, this indicates that void nucleation regions in the topology optimization process.
For a more detailed comparison of the results, the relative error is shown in Figure  \ref{Exm1_sens2}. The results indicate that there is a good agreement between analytical and numerical sensitivities, while analytical approach is computationally much faster.

It is worth noting that the source of the residual error in in Figure  \ref{Exm1_sens2} can be illustrated because of (i) choice of pseudo-density perturbation value $\Delta \Phi_e$, (ii) the use of the exact Heaviside step function in \req{eq:heavydide_def} versus its regularized version, and (iii) the use of the regularized Dirac delta function in \req{eq:dirac_approx} (see Figure \ref{Figure_dirac_filter}a). It is worth noting that due to the fact that the exact Heaviside step function does not have an analytical derivative, therefore, the regularized Heaviside function, which is the integral of the regularized Dirac delta function, has been used in the finite difference method. While in the analytical method, the exact Heaviside step function has been used.\\

For further comparison, the sensitivity of the Lagrangian function for a constant value of the volume constraint Lagrange multiplier which enters in ${\widehat{v}_{\Phi}}\left| _{\Phi = 0} \right.$ is depicted in Figure  \ref{Exm1_sens3}. This is the cross-sectional view in $x-y$ direction, in which the maximum quantity for the  ${\widehat{v}_{\Phi}}\left| _{\Phi = 0} \right.$ appeared in middle singular point (where the crack initiation is observed, see Figure \ref{Exm3_d}) and the load applied. Evidently, the reasonable agreement between FEM  and FDM has been observed, which demonstrates the accuracy of the sensitivity analysis has been done in Section 3.4.
\begin{figure}[!b]
		\centering
	\vspace{-0.1cm}
	\includegraphics[clip,trim=0cm 21.5cm 3cm 0cm, width=14cm]{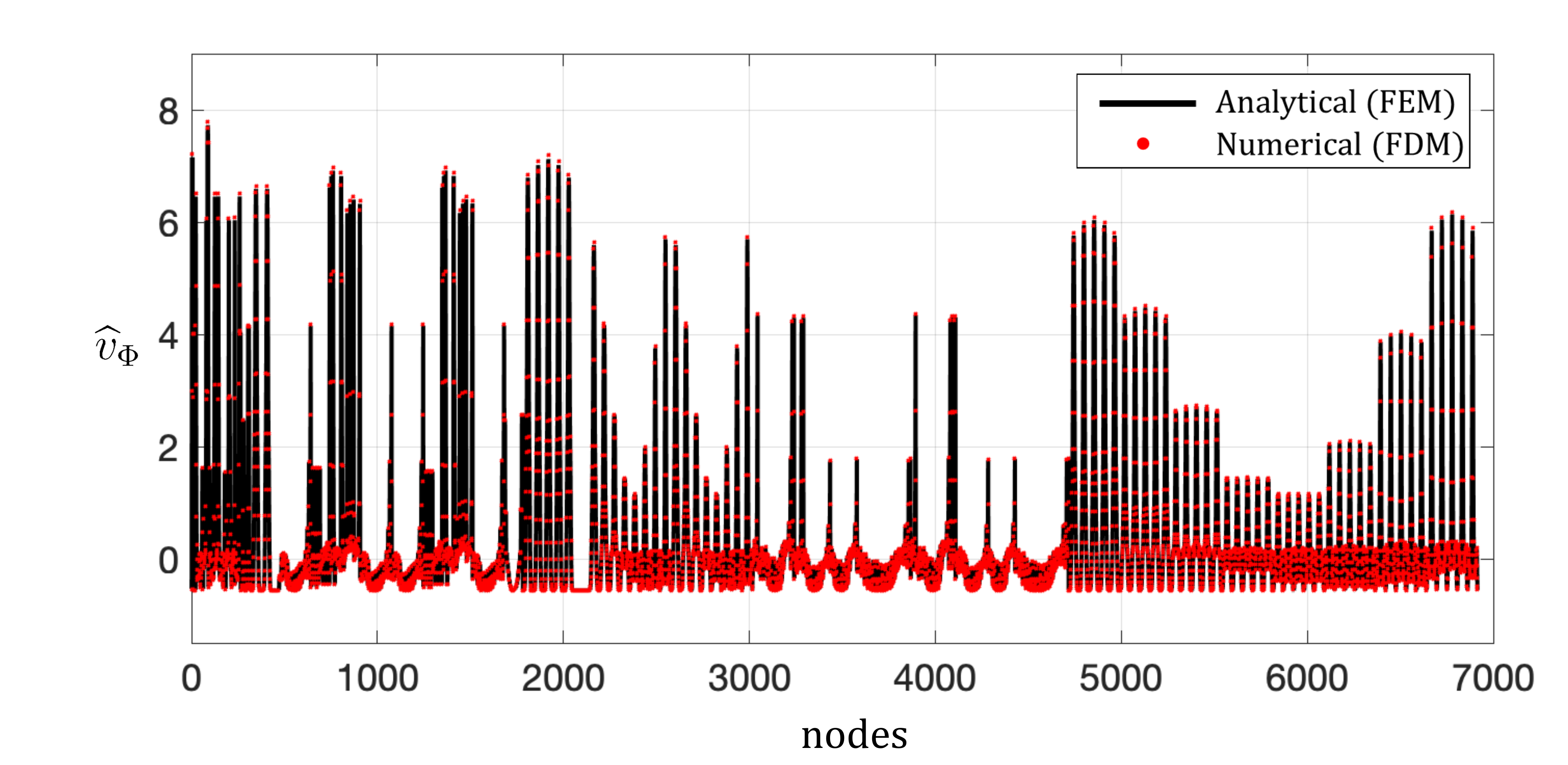}		
	\caption{Example 3. Sensitivity analysis performed through analytical approach by FEM given in (106) versus the numerical sensitivity by FDM given in \req{app1_eq:fdm1}.}
	\label{Exm1_sens1}
\end{figure}

\begin{figure}[!ht]
		\centering
	\vspace{-0.1cm}
	\includegraphics[clip,trim=0cm 21.5cm 3cm 0cm, width=14cm]{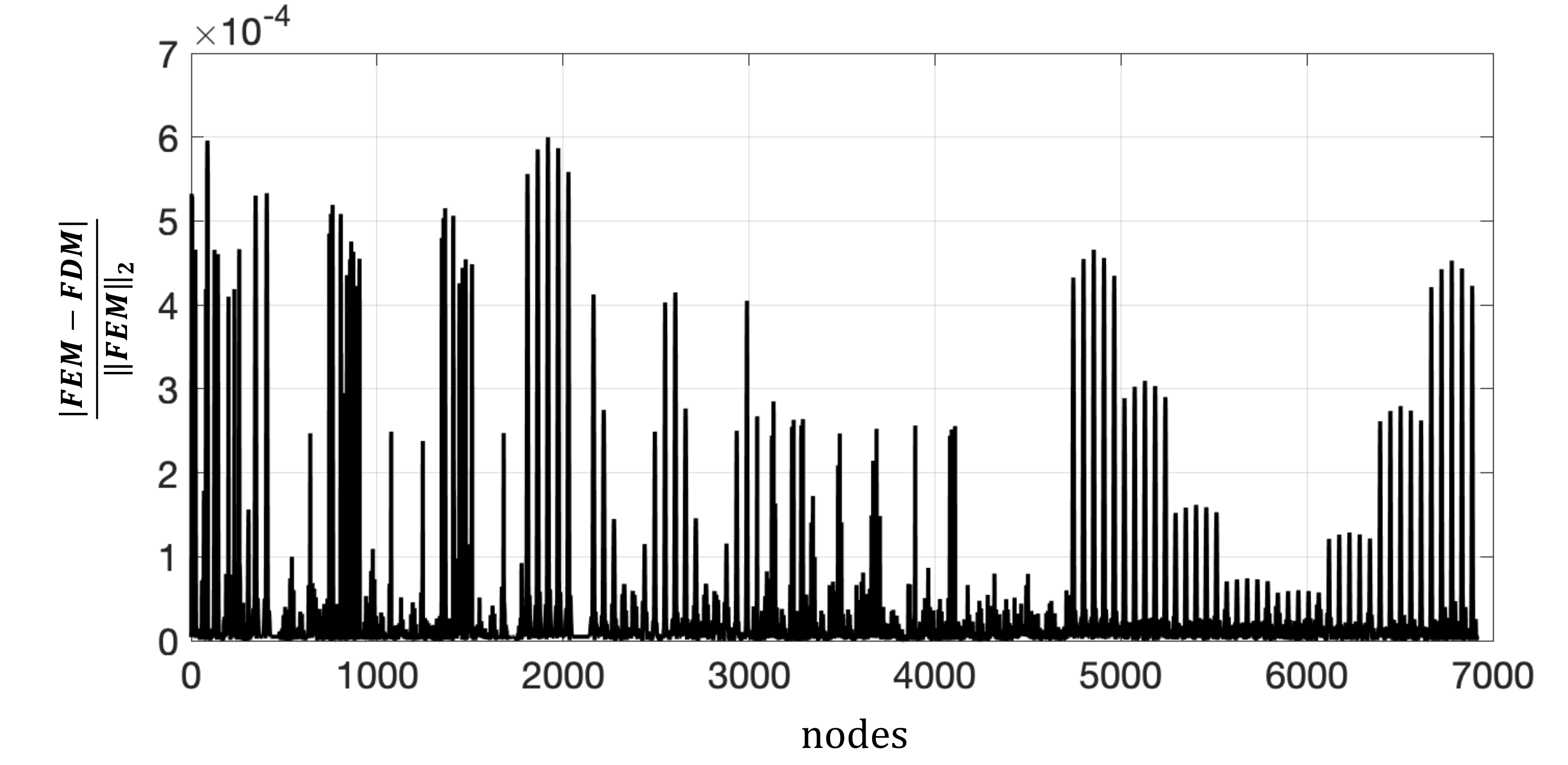}		
	\caption{Example 3. The relative sensitivity errors between analytical and numerical schemes.}
	\label{Exm1_sens2}
\end{figure}  

\begin{figure}[!b]
	\centering
	\vspace{-0.1cm}
	\includegraphics[clip,trim=1cm 29cm 1.5cm 1cm, width=16cm]{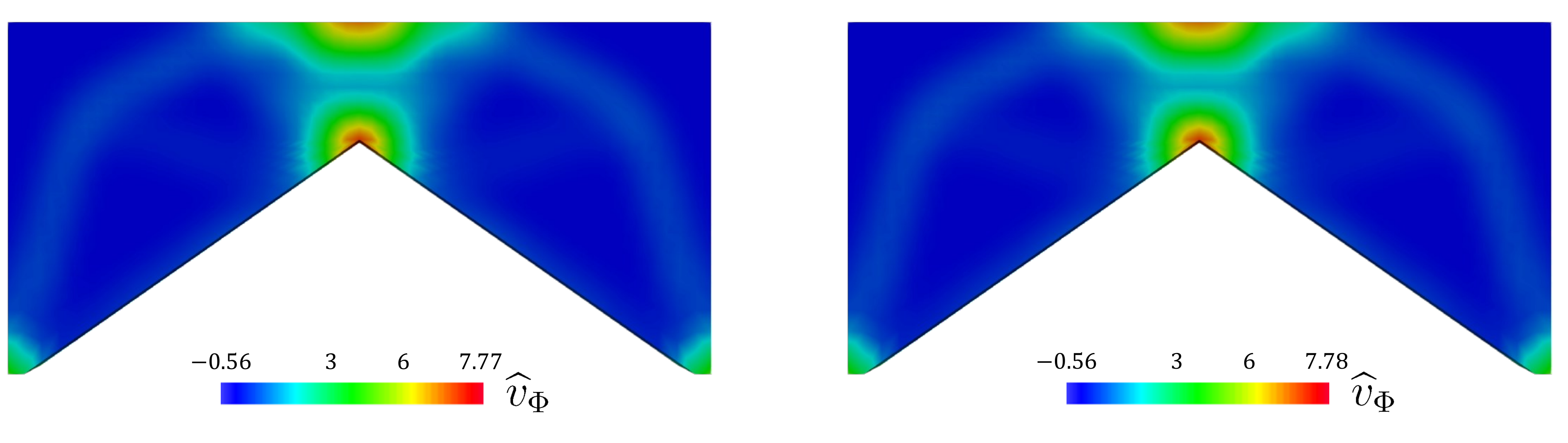}		
	\vspace*{-0.7cm}
	\caption*{\hspace*{2cm}(a)\hspace*{8cm}(b)\hspace*{2cm}}
	\caption{Example 3. The sensitivity distribution enters in ${\widehat{v}_{\Phi}}\left| _{\Phi = 0} \right.$ for the cross-sectional view in $x-y$ direction of the portal frame (a) FEM, and (b) FDM.}
	\label{Exm1_sens3}
\end{figure}  
\end{Appendix}
\cleardoublepage
    {\normalsize
	\begin{spacing}{0.8}
		\bibliographystyle{ieeetr}
		\bibliography{./lit}
\end{spacing}}
\end{document}